\documentclass[10pt,leqno]{amsart}
\usepackage{amsmath,amsthm,amsfonts}
\usepackage {latexsym}
\usepackage{amssymb}

\usepackage[dvips]{epsfig}

\newtheorem{theorem}{Theorem}[section]
\newtheorem{lemma}[theorem]{Lemma}
\newtheorem{prop}[theorem]{Proposition}
\newtheorem{defi}[theorem]{Definition}
\newtheorem{corollary}[theorem]{Corollary}
\newtheorem{cor}[theorem]{Corollary}
\newtheorem{remark}[theorem]{Remark}

\renewcommand{\theequation}{\thesection .\arabic{equation}}

\let\sect\section
\renewcommand\section{\setcounter{equation}{0}
\gdef\theequation{\thesection .\arabic{equation}}\sect}

\newcommand{\cA}{{\mathcal{A}}}
\newcommand{\cB}{{\mathcal{B}}}
\newcommand{\cD}{{\mathcal{D}}}
\newcommand{\cE}{{\mathcal{E}}}
\newcommand{\cG}{{\mathcal{G}}}
\newcommand{\cH}{{\mathcal{H}}}
\newcommand{\cJ}{{\mathcal{J}}}
\newcommand{\cN}{{\mathcal{N}}}
\newcommand{\cR}{{\mathcal{R}}}
\newcommand{\cU}{{\mathcal{U}}}

\newcommand{\cM}{{\mathcal{M}}}
\newcommand{\cT}{{\mathcal{T}}}
\newcommand{\cS}{{\mathcal{S}}}
\newcommand{\cF}{{\mathcal{F}}}

\newcommand{\cP}{{\mathcal{P}}}

\newcommand{\IC}{{\mathbb{C}}}

\newcommand{\IR}{{\mathbb{R}}}
\newcommand{\TT}{{\mathbb{T}}}

\newcommand{\tor}{\TT}
\newcommand{\ZZ}{{\mathbb{Z}}}
\newcommand{\IZ}{{\mathbb{Z}}}

\newcommand{\uE}{{\underline{E}}}
\newcommand{\oE}{{\overline{E}}}
\newcommand{\oN}{{\overline{N}}}
\newcommand{\uN}{{\underline{N}}}

\newcommand{\uu}{{\underline{u}}}
\newcommand{\uv}{{\underline{v}}}
\newcommand{\ux}{{\underline{x}}}
\newcommand{\ox}{{\bar{x}}}
\newcommand{\uw}{{\underline{w}}}
\newcommand{\uz}{{\underline{z}}}

\newcommand{\be}{\begin{eqnarray}}
\newcommand{\ee}{\end{eqnarray}}

\newcommand{\degg}{\mathop{\rm{deg}}}
\newcommand{\dist}{\mathop{\rm{dist}}}
\newcommand{\disc}{\mathop{\rm{disc}}}

\newcommand{\mes}{\mathop{\rm{mes}\, }}
\newcommand{\compl}{\mathop{\rm{compl}}}
\renewcommand{\mod}{{\rm{mod}\, }}
\newcommand{\rtr}{\mathop{\rm{tr}\, }}

\newcommand{\Res}{\mathop{\rm{Res}}}

\newcommand{\spec}{\mathop{\rm{spec}}}
\newcommand{\const}{\mathop{\rm{const}}}

\newcommand{\tilw}{\tilde{w}}

\newcommand{\ul}{\underline{\ell}}

\newcommand{\xo}{(x,\omega)}
\newcommand{\hnxo}{H_N\xo}

\newcommand{\zoe}{(z, \omega, E)}
\newcommand{\capo}{\cA_{\rho_0}}

\newcommand{\vep}{{\varepsilon}}
\newcommand{\ve}{{\varepsilon}}

\newcommand{\nn}{\nonumber}

\newcommand{\notint}{{-\!\!\!\!\!\!\int}}
\newcommand{\nint}{\mathop{\notint}}
\newcommand{\la}{\langle}
\newcommand{\ra}{\rangle}
\def\beeq{\begin{equation}}
\def\eneq{\end{equation}}
\def\eps{\varepsilon}
\def\les{\lesssim}

\def\cZ{{\mathcal Z}}
\def\cS{{\mathcal S}}
\def\bm{\begin{matrix}}
\def\endm{\end{matrix}}

\def\Im{{\rm Im}}
\def\cC{{\mathcal C}}
\def\wt{\widetilde}

\def\Dioph{{\rm Dioph}}

\def\un{\underline{n}}

\newcommand{\car}{\mathop{\rm{Car}}\nolimits}
\newcommand{\mybigcup}{\mathop{\textstyle{\bigcup}}\limits}
\begin{document}

\title[Resonances and the
formation of gaps]{On resonances and the formation of gaps in the spectrum\\
of quasi-periodic Schr\"odinger equations}
\author{Michael Goldstein
\and Wilhelm Schlag}

\address{Department of Mathematics, University of Toronto, Toronto, Ontario, Canada M5S 1A1}

\email{gold@math.toronto.edu}

\address{The University of Chicago, Department of Mathematics, 5734 South University Avenue, Chicago, IL 60637, U.S.A.}

\email{schlag@math.uchicago.edu}

\thanks{The first author was partially supported by a  Guggenheim Fellowship and by an NSERC grant. The second author was partially
supported by the NSF, DMS-0617854. The first author is very grateful
to Galina Perelman for several illuminating discussions during her visit at
the University of Toronto in July of 2005.}

\maketitle

\begin{abstract}
We consider one-dimensional difference Schr\"odinger equations
$$
\bigl[H(x,\omega)\varphi\bigr](n)\equiv -\varphi(n-1)-\varphi(n+1) +
V(x + n\omega)\varphi(n) = E\varphi(n)\ ,
$$
$n \in \IZ$, $x,\omega \in [0, 1]$ with real-analytic potential
function $V(x)$.  If $L(E,\omega_0)>0$ for all $E\in(E', E'')$ and
some Diophantine $\omega_0$, then the integrated density of states
is absolutely continuous for almost every $\omega$ close to
$\omega_0$, see~\cite{Gol Sch2}. In this work we apply the methods
and results of \cite{Gol Sch2} to establish the formation of a dense
set of gaps in $\spec( H(x,\omega))\cap (E',E'')$. Our
approach is based on an induction on scales argument, and is therefore both
constructive as well as quantitative. Resonances between
eigenfunctions of one scale lead to "pre-gaps" at a larger scale.
To pass to actual gaps in the spectrum, we show that  these pre-gaps cannot be filled more than a finite
(and uniformly bounded) number of times. To accomplish this, one
relates a pre-gap to pairs of complex zeros of the Dirichlet
determinants off the unit circle using the techniques of~\cite{Gol
Sch2}. Amongst other things, we establish in this work a non-perturbative version
of the co-variant parametrization of the eigenvalues and eigenfunctions via the
phases in the spirit of
Sinai's (perturbative) description of the spectrum~\cite{Sin1} via his function~$\Lambda$.
This allows us to relate
the gaps in the spectrum with the graphs of the eigenvalues
parametrized by the phase. Our infinite volume theorems hold for all Diophantine frequencies~$\omega$ 
up to a set of Hausdorff dimension zero. 
\end{abstract}

\section{Introduction and statement of the main results}\label{sec:intro}

The main goal of this work is to establish a multiscale
description of the structure of the spectrum of
quasi-periodic Schr\"odinger equations
\begin{equation}\label{eq:1.Sch}
\bigl[H(x,\omega)\varphi\bigr](n) \equiv - \varphi(n-1) -
\varphi(n+1) +  V(x + n\omega)\varphi(n) = E\varphi(n)
\end{equation}
in the regime of exponentially localized eigenfunctions.  We assume
that $V(x)$ is a $1$-periodic, real-analytic function, and that
$\omega \in [0, 1]$.  Let $H_N(x,\omega)$ be the restriction of
$H(x,\omega)$ to the finite interval $[1, N]$ with zero boundary
conditions.  Consider the union $\cS_N = \mybigcup_x\spec( \hnxo)$,
where $\spec (\hnxo)$ stands for the spectrum of $\hnxo$.  The set
$\cS_N$ is closed, so \[\cS_N = \bigl[\underline{E}(N), \overline
E(N)\bigr] \setminus\mybigcup_k \bigl(\uE(N, k), \oE(N, k) \bigr),\quad \uE(N) = \min\limits_{\cS_N}\, E,\; \oE(N) =
\max\limits_{\cS_N}\, E\]
where $\bigl(\uE(N, k), \oE(N, k) \bigr)$
are the maximal intervals of $\bigl[\uE(N), \oE(N)\bigr] \setminus
\cS_N$.  More specifically, the goals of this work are as follows:

\begin{enumerate}
\item[(a)] To relate the intervals $\bigl(\uE(N, k), \oE(N, k) \bigr)$ and
$\bigl(\uE(N', k'), \oE(N', k')\bigr)$ for ``consecutive scales'' $N
\gg N'$.

\item[(b)] To ``label'' the interval $\bigl(\uE(N, k),$ $\oE(N, k)\bigr)$ relative  to the intervals $\bigl(\uE(m, \ell), \oE(m, \ell)
\bigr)$ of the previous scales.

\item[(c)] To describe the mechanism responsible for the formation of
the intervals $\bigl(\uE(N, k), \oE(N, k)\bigr)$ inside  the set
$\cS_{N'}$, $N' \ll N$, independently of any $\bigl(\uE(N', k'),
\oE(N', k')\bigr)$.
\end{enumerate}

Our interest in these properties is largely motivated by possible
applications to inverse spectral problems for the quasi-periodic
Schr\"odinger equation and the Toda lattice with quasi-periodic
initial data~\cite{To}. This paper relies heavily on the methods developed
in~\cite{Gol Sch2}. For the
convenience of the reader, we recall -- and expand upon -- some of
the material of that paper in Sections~\ref{sec:basictools}-\ref{sec:separation}.
The eigenvalues \[E_1^{(N)}(x,\omega) < E_2^{(N)}(x,\omega) < \dots
< E_N^{(N)}(x,\omega)\] of $\hnxo$ are real analytic functions of $x
\in [0, 1]$.  Although the graphs of the functions
$E_j^{(N)}(x,\omega)$ can be very complicated, the following was
proved in~\cite{Gol Sch2} for Diophantine~$\omega$,
see~\eqref{eq:1.Diaph}, and positive Lyapunov exponents: there exist
intervals $\bigl(E'_{N,k}, E''_{N, k}\bigr)$, $k = 1, 2,\dots,
k_N$, with
\[ \max\limits_k\bigl(E''_{N, k} - E'_{N, k}\bigr) \le
\exp\bigl(-\bigl(\log N\bigr)^A\bigr), \qquad k_N \le \exp
\bigl(\bigl(\log N\bigr)^B\bigr),\] with constants $1 \ll B \ll A$
depending on~$\omega$, such that if \[E_j^{(N)}(x,\omega) \notin
\cE_N = \mybigcup_k \bigl(E'_{N, k}, E''_{N, k}\bigr),\] for some
$j$ and $x$, then \[\bigl|\partial_x E_j^{(N)} (x,\omega) \bigr| >
\exp\bigl(-N^\delta\bigr)\] Here $0 < \delta \ll 1$ is an arbitrary
but fixed small parameter.  In other words,  the
graphs of $E_j(x,\omega)$ have controlled slopes off a small set
$\cE_N$.

\begin{figure}[ht]
\centerline{\hbox{\vbox{ \epsfxsize= 13.0 truecm \epsfysize=6.0
truecm \epsfbox{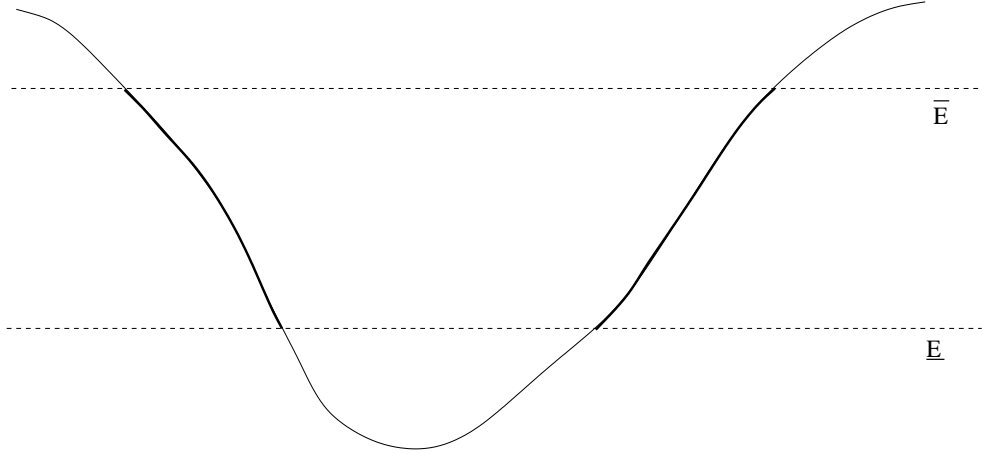}}}}
\label{fig:1} \caption{$I$--segments}
\end{figure}

The segments of the graph where $E_j^{(N)}(x,\omega) \in I$ and $I =
(\uE, \oE)$ is an interval disjoint from $\cE_N$, are called
$I$--{\em segments}. They are denoted by
$\bigl\{E_j^{(N)}(x,\omega), \ux,\ox\bigr\}$, where
\[E_j^{(N)}(\ux,\omega) = \uE, \quad E_j^{(N)}(\ox,\omega) = \oE\]  The
$I$--segments are important for our purposes, because they allow us
to locate the resonances and to describe the graphs of the functions
$E_j^{(\oN)}(x,\omega)$ for $\oN \gg N$ in the region where the
resonance occurs. A possible definition of a resonance is as
follows: With a constant $A\gg1$ depending on~$\omega$,
\begin{equation}\label{eq:1.res}
\tau = \bigl|E_{j_1}^{(N)} \xo - E_{j_2}^{(N)} (x + m\omega, \omega)
\bigr| < m^{-A}
\end{equation}
for some $x \in \tor$, $1 \le j_1$, $j_2 \le N$ and $m > N$. In
fact, there is some stability in the constant~$A$ with regard to
small perturbations of~$\omega$.

The significance of such resonances was explained in the work by
Sinai~\cite{Sin1} on quasi-periodic Anderson localization for
potentials $V(x)=\lambda \cos(2\pi x)$ in the regime of large $|\lambda|$,
see~\eqref{eq:1.Sch}. Sinai developed a  KAM-type scheme to analyze
the functions $E_j^{(N)}(x,\omega)$ and the corresponding
eigenvectors. The critical points of $E_j^{(\oN)}(x,\omega)$ with
$\oN \gg N$ were proved to be closely related to resonances as
in~\eqref{eq:1.res}. It is very important for the analysis of the
resonances~\eqref{eq:1.res} in~\cite{Sin1} that given $x \in \tor$
and $j_1$ there exist at most one $j_2$ and $m \le \oN$ so that
\eqref{eq:1.res} holds. {\sl For that reason the function $V(x)$ in
\cite{Sin1} is assumed to have two monotonicity intervals with
non-degenerate critical points.} That allows one to reduce the
analysis of $E_j^{(\oN)}(x,\omega)$ to an eigenvalue problem for a
$2\times 2$ matrix function of the form
\begin{equation}\label{eq:1.twobytwo}
A(x) = \begin{bmatrix}
E_1(x-x_0) & \vep(x)\\
\vep(x) & E_2(x - x_0)\end{bmatrix}\ ,
\end{equation}
where $E_1(0) = E_2(0)$, $\partial_x E_1 < 0$, $\partial_x E_2 > 0$
locally around zero, and $\vep(x)$ is small together with its
derivatives. It is easy to check that the eigenvalues $E^+(x)$,
$E^-(x)$ of $A(x)$ plotted against $x$ are as in Figure~2, at least
locally around~$x_0$.

\begin{figure}[ht]\label{fig:2}
\centerline{\hbox{\vbox{ \epsfxsize= 13.0 truecm \epsfysize= 5.0
truecm \epsfbox{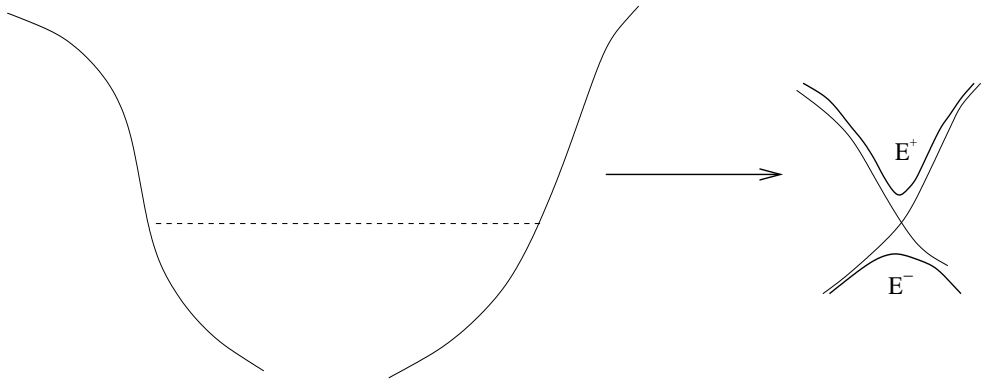}}}} \caption{Classical formation
of the resonant eigenvalues}
\end{figure}


We would like to emphasize that some of the conclusions which we
reach in this paper are similar in spirit to those of
Sinai~\cite{Sin1}. This is particularly true in regards to the main
result involving gaps and the aforementioned pictures describing the
splitting of eigenvalues. At the same time, we stress that we use
entirely nonperturbative methods (i.e., we are only assuming
positive Lyapunov exponent rather than large $|\lambda|$) and we
work with more general potentials than cosine. In this respect we
would like to mention the recent breakthrough by Puig~\cite{Puig},
who established the Cantor structure of the spectrum for the almost
Mathieu case (cosine potential) and Diophantine~$\omega$. Earlier,
Choi, Elliott, and Yui~\cite{CEY} had obtained gaps for the case of
Liouville rotation numbers~$\omega$. The remaining cases of
irrational rotation numbers (i.e., those with behavior intermediate
to Diophantine and Liouville) was settled by Avila and
Jitomirskaya~\cite{AJ} (but this again only applies to the cosine).

The major objective in this work is to locate those segments of the
graphs of some $E_{k_1}^{(\oN)}(x,\omega)$,
$E_{k_2}^{(\oN)}(x,\omega)$ which look like $E^+(x,\omega)$ and
$E^-(x,\omega)$ in Figure~2. Ultimately, such regions give rise to
gaps in the spectrum. Before we state the main result of this work
let us recall the central notions involved in it.

It is convenient to replace $V(x)$ in \eqref{eq:1.Sch} by
$V\bigl(e(x)\bigr)$ (with $e(x)=e^{2\pi ix}$), where $V(z)$ is an
analytic function in the annulus $\cA_{\rho_0} = \bigl\{z \in \IC: 1
- \rho_0 < |z| < 1 + \rho_0\bigr\}$ which assumes only real values
for $|z| = 1$. The monodromy matrices are as follows
\begin{equation}\label{eq:1.monodr}
\begin{aligned}
M_{[a, b]}(z, \omega, E) & = \prod^a_{k = b}\, A\bigl(ze(k\omega), \omega, E\bigr)\\[5pt]
A(z,\omega, E) & = \begin{bmatrix} V(z) -E & -1\\ 1 & 0\end{bmatrix}
\end{aligned}
\end{equation}
$a, b \in \IZ$, $a < b$, $E \in \IC$.  For $M_{[1, N]}(z, \omega,
E)$ we reserve the notation $M_N(z, \omega, E)$.  For almost all $z
= e(x + iy) \in \cA_{\rho_0}$ the limit
\begin{equation}\label{eq:1.Lyapdef}
\lim_{N\to \infty}\ N^{-1} \log \big \|M_N(z, \omega, E)\big \|
\end{equation}
exists; if $\omega$ is irrational, then the limit does not depend on
$x$ a.s.~and it is denoted by $L(y, \omega, E)$.  The most important
case is $y=0$, and  we reserve the notation $L(\omega, E)$ for the
Lyapunov exponents $L(0, \omega, E)$.  We always assume that the
frequency $\omega$ satisfies the same Diophantine condition as in
\cite{Gol Sch2}, namely
\begin{equation}\label{eq:1.Diaph}
\|n\omega \| \ge \tfrac{c}{n(\log n)^a}\qquad\text{for all $n \ge
1$}
\end{equation}
and some $a > 1$.  We denote the class of $\omega$ satisfying
\eqref{eq:1.Diaph} by $\tor_{c, a}$ and further define 
\[
 \Dioph := \bigcup_{a>1,c>0} \tor_{c,a}
\]
Let $\omega\in \Dioph$.
By a theorem of Avron and Simon~\cite{AvrSim} the spectrum $\spec (H(x, \omega))$ does not depend on~$x$
and we denote it by~$\Sigma_\omega$.

\begin{theorem}\label{th:1.mainth}
Assume that $L( \omega, E)  > 0$
for any $\omega \in (\omega',\omega'')$ and any  $E\in(E', E'')$.
There exists
 a set $\Omega$ of Hausdorff dimension zero such that for any
$\omega \in (\omega', \omega'') \cap \Dioph\backslash \Omega$, the intersection
$\Sigma_\omega\bigcap (E', E'')$ is a Cantor set.
\end{theorem}

We remark that it follows from this theorem that  if $L( \omega_0, E) \ge \gamma
> 0$ for some $\omega_0 \in \tor_{c,a}$ and $E_0\in \IR$
then there exists $\rho^{(0)} = \rho^{(0)}(V, c, a, \gamma) > 0$,
and a set $\Omega$ of Hausdorff dimension zero such that for any
$\omega \in (\omega_0 - \rho^{(0)}, \omega_0 + \rho^{(0)}) \cap
\Dioph \backslash \Omega$ the spectrum  $\spec (H(x, \omega))\bigcap (E_0-
\rho^{(0)}, E_0 + \rho^{(0)})$ is a Cantor set.  This is due to the fact that $L(\omega,E)>\gamma/2$
for all $|\omega-\omega_0|<\rho^{(0)}$ (see \cite{Gol Sch1} or~\cite{BourJit}).

Concerning the statement of Theorem~\ref{th:1.mainth}, note that the removal
of a set of Hausdorff dimension zero cannot be achieved by a ``Fubini''-type argument;
rather, it requires some information on the  complexity of a suitable
cover of sets of bad frequencies~$\omega$ (relative to finite volume).
Throughout this paper we rely heavily
on the notion of ``complexity'' of a set real or complex numbers: if
$\cS\subset\IR$, then
\[
 \mes(\cS)<\eps, \quad \compl(\cS)<K
\]
mean that for some intervals $I_k$,
\[
 \cS\subset \bigcup_{k=1}^K I_k, \quad \sum_{k=1}^K |I_k|<\eps
\]
In all cases considered here, $K\eps\ll 1$ and we will often replace the
latter condition by the stronger $\max_k |I_k|<\eps^2$.  In the complex
case, replace `interval' by `disk'.
We derive Theorem~\ref{th:1.mainth} as a simple corollary of our
analysis of the gap development in finite volume.
 The following theorem should be thought of as an (important) representative
of the finite volume analysis -- the reader will find a more
detailed  description in later sections. As usual, $[y]$ stands for
the entire part of $y$ and $\cH^s_\alpha$ is the $s$-dimensional Hausdorff outer
measure of scale $\alpha$, see e.g.~Falconer~\cite{Falc}. Finally, we introduce the notation
 $H^{(P)}_{N}(x,\omega)$  for the
Schr\"odinger operator on $[-N+1,N]$ with periodic boundary
conditions and \[ \cS^{(P)}_{N, \omega}:=\bigcup_{x\in \TT} \spec(
H^{(P)}_N(x,\omega))\]

\begin{theorem}\label{th:1.2}
Assume that $L(\omega, E) \geq \gamma > 0$ for any $\omega \in
(\omega', \omega'')$ and any $E\in(E', E'')$. Given $c>0$, and
$a>1$, $0<s<1$ there exist positive integers $N_0 = N_0(V, c, a,
\gamma,s)$ and $T_0 = T_0(V, c, a, \gamma), A=A(V, c, a, \gamma)$ such that for any
$N_1 \geq N_0$ there exists a subset $\Omega_{N_1,s} \subset \tor$
\[
\cH^s_{\alpha(N_1)}(\Omega_{N_1,s})\le 1,\quad \alpha(N_1)=\exp(- (\log \log N_1)^A)
\]
such that for all $\omega \in \tor_{c, a} \cap (\omega',
\omega'') \backslash \Omega_{N_1,s}$ the following statement holds:
Set \[ \oN(N_1,1):=N_1,\quad \oN(N_1,t+1):=[\exp((\oN(N_1,t))^{\delta})] \quad \forall\; 
t\ge1 \] with some small $0 < \delta=\delta(V, c,
a, \gamma) \ll 1$. Then there exists $N \le \oN(N_1,T_0)$  depending on
$\omega$, such that for any interval $I = (\uE, \oE)$, $I \subset (E',
E'')$ with $|I|
> \exp(-(\log N_1)^C)$  there exists a subinterval $I^{(1)} =
(\uE^{(1)}, \oE^{(1)}) \subset I$ such that $|I^{(1)}|>
\exp(-\oN(N_1,T_0))$ and $ \cS^{(P)}_{N, \omega} \cap I^{(1)} =
\emptyset $. 
\end{theorem}

In the previous theorem, we use the Hausdorff outer measure just for simplicity.
In fact, one has the following bounds on the measure and complexity of $\Omega_{N_1,s}$:
\[
\Omega_{N_1,s}=\bigcup_{1\le t \le T_0+1} \cB_{N_1,t}, \quad  \mes (
\cB_{N_1,t}) \le \mu(t), \quad \compl ( \cB_{N_1,t}) \le C(t)
\]
where $\mu(t)=\exp(-(\log \log \oN(N_1,t))^A)$, $C(t)=\mu(t)^{-s}$. 

\smallskip It is not clear how to pass from Theorem~\ref{th:1.2}
 to Theorem~\ref{th:1.mainth} via the basic definition of the
 spectrum in $\ell^2$ alone. The well-known description
 of the spectrum via polynomially bounded solutions (via the Schnol--Simon theorem)
 appears not to be too helpful in this context, either, since it is an existence
theorem and thus non-effective.
 Let us recall in passing that the  ``proper'' spectrum~$\Sigma_\omega$, which is the
 closure of the set $\{ E_j(x,\omega)\}_j$ of the eigenvalues of $H(x,\omega)$, does {\em not}
 depend on $x$, whereas the set $\{ E_j(x,\omega)\}_j$ of eigenvalues itself
 does.
 The methods developed in \cite{Gol Sch2}, and which are expanded upon here,
 are centered around the parametrization of the eigenvalues by the phase.
 Amongst the properties of these parametrizations we single out the crucial {\em separation }
as the most important; it says that the eigenvalues $E_j^{(N)}(x,\omega)$ of $H_{[-N,N]}(x,\omega)$
satisfy
\[
 |E_j^{(N)}(x,\omega) - E_k^{(N)}(x,\omega)|> e^{-N^\delta}
\]
for any $j\ne k$ provided $E_j^{(N)}(x,\omega)\not \in \cE_{N,\omega}$ where the latter set
is small both with regard to measure and complexity, see Section~\ref{sec:separation} for more
details.  In his seminal
paper~\cite{Sin1} on cosine-like potentials, Sinai introduced a (multi-valued) function
$\Lambda$, defined almost everywhere
on~$\tor$, which allows for
the parametrization of the eigenvalues in the {\em infinite volume} in a co-variant fashion. This means
that for almost every $x\in\tor$, the set $\{\Lambda(x+j\omega)\}_{j=-\infty}^\infty$ is the complete set of
eigenvalues of $H(x,\omega)$ and the associated orthonormal basis of eigenfunctions
$\{\psi_j(\cdot,x,\omega)\}_{j=-\infty}^\infty$ satisfies
\[
 \psi_j(\cdot,x+\omega,\omega) = \psi_j(\cdot+1,x,\omega)
\]
The following theorem on infinite volume Anderson localization arises as part
of our construction of gaps. It is proved via an induction on scales argument with a suitable
finite volume localization statement at its core. Amongst other things, it shows that the
set of exceptional frequencies~$\omega$ which need to be removed from all Diophantine~$\omega$ in
order to ensure Anderson localization in 
the work of Bourgain and the first author, see~\cite{BouGol}, 
is of Hausdorff dimension zero. The positive measure statement in the theorem improves on~\cite{Boupos} for the same reason. 
However, our proof of localization is very different technically
speaking from the one in~\cite{BouGol} and property (4) in the following theorem is new. In essence, this property 
controls the number of monotonicity intervals of Sinai's function. 

\begin{theorem}\label{th:1.sinai}
 Assume that $L(\omega_0, E) \geq \gamma > 0$ for some $\omega_0 \in \tor_{c, a}$
and any $E\in \IR$. Then there exist $\rho^{(0)} = \rho^{(0)}(V, c,
a, \gamma) >0$, and a set $\Omega\in \tor$ of Hausdorff dimension zero
such that for any $\omega \in (\omega_0-\rho^{(0)}, \omega_0+\rho^{(0)})
\cap \Dioph\backslash \Omega$ the spectrum satisfies 
$
\mes (\Sigma_{\omega}) > 0
$. 
Furthermore, there exists a set $\cB_{\omega}\subset \tor$ of
Hausdorff dimension zero,
 such that for any $x\in \tor \setminus \cB_{\omega}$ the following conditions hold:
\begin{itemize}
\item[(1)] There exists an orthonormal basis
$\{\psi_j(x,\omega,\cdot)\}_{j\ge1}$ of eigenfunctions of
$H(x,\omega)$ in $\ell^2(\ZZ)$,
\[H(x,\omega)\psi_j(x,\omega,\cdot)=E_j(x,\omega)\psi_j(x,\omega,\cdot)
\]
Moreover, each function $\psi_j(x,\omega,\cdot)$ is exponentially
localized and 
\[
\lim_{|N|\rightarrow \infty} (2|N|)^{-1}\log \big( |\psi_j(x,\omega,N)|^2+|\psi_j(x,\omega,N-1)|^2\big)  =
-L(\omega,E_j(x,\omega))
\]

\item[(2)] The eigenvalues $E_j(x,\omega)$ are simple
\item[(3)] The set $\cB_{\omega}$ is invariant under
the shifts $x\mapsto x+m\omega\; (\mod 1)$, $m\in \ZZ$. For any $x\in
\tor \setminus \cB_{\omega}$, $j\ge1$, and $m\in \ZZ$ ,
$\psi_j(x,\omega,\cdot+m)$ is the eigenfunction of
$H(x+m\omega,\omega)$ with the eigenvalue $E_j(x,\omega)$

\item[(4)] For each $E\in \IR$ the set
\[
\cT(E):=\{x\in \tor \setminus \cB_{\omega} \;:\; \exists j \text{\ \ so that\ \ }E=E_j(x,\omega) \}
\]
is either empty or consists of a union of trajectories
$\Gamma(x(E,k))$, $x(E,k)\in \tor \setminus \cB_{\omega}$,
$1\le k\le k(E)$ where $k(E)\le C(V)<\infty$ and 
\[
\Gamma(x):=x+\omega\ZZ\quad (\mod 1)
\]
If $V$ is a trigonometric polynomial of degree\footnote{This means that $V(x)=\sum_{k=-k_0}^{k_0} a_k e(kx)$ with $a_{-k}=\overline{a_k}$.} $k_0$, then $C(V)\le
2k_0$.
\end{itemize}

\end{theorem}

We feel that the methods of this paper, combined with some ``soft'' measure
theoretic consideration, should allow for the construction of a true Sinai's function. I.e., 
we claim that there exists a function~$\Lambda:\tor\to\IR$, defined up to a set
of Hausdorff dimension zero and with at most $C(V)$ monotonicity intervals where $C(V)$
is as in part (4) above such that 
\[
 E_j(x,\omega) = \Lambda(x+j\omega)  \qquad\forall j\in \ZZ,\;\forall\; x\in\tor\setminus\Omega
\]
where $\Omega$ is of Hausdorff dimension zero. However, we have chosen not to pursue this issue here. 

\noindent 
In Section~\ref{sec:sinai} we derive a detailed finite volume 
version of  Theorem~\ref{th:1.sinai}.   Amongst
other things, this derivation gives an effective quantitative
description of the spectrum of the problem \eqref{eq:1.Sch} on the
whole lattice $\IZ$ in terms of the spectrum on finite volume
and also allows for a simple transition from Theorem~\ref{th:1.2} to
Theorem~\ref{th:1.mainth}. The co-variant parametrization of the
eigenvalues and eigenfunctions via the phases is based on the
description of the exponentially localized eigenfunctions on the
interval $[-\oN, \oN]$ by means the eigenfunction on the interval
$[-N, N]$ with $N \leq \oN$, combined with the aforementioned
quantitative repulsion property of the Dirichlet eigenvalues on a
finite volume. We discuss these results in
Sections~\ref{sec:resultant}--\ref{sec:separation}.
As already mentioned
before, we produce gaps (on finite volume) from resonances of the
previous scale. This requires restricting the graphs of the
eigenfunctions to segments which have a controlled slope (in the
sense of a favorable lower bound). Thus, in Section~\ref{sec:parameter} we
introduce $I$--segments $\bigl\{E_j^{(N)}(x,\omega),
\ux,\ox\bigr\}$ of the graphs of the eigenfunctions on a finite
volume $[-N,N]$ which have slopes bounded below (in absolute value)
by $e^{-N^\delta}$. Amongst those we single out {\em regular}
$I$-segments which have the property that the eigenfunctions with
eigenvalues $E_j^{(N)}(x,\omega)$ are supported away from the
boundary of $[-N,N]$ for all $\ux<x<\ox$.
 This is needed in order to assure that crossing $I$-segments do indeed
form a resonance at a larger scale as in the figures above. For that
we use the spectrum and the eigenfunctions of the Schr\"odinger
operator on a finite interval with periodic (respectively,
antiperiodic) boundary conditions. The key analytical tool for the
study of periodic boundary conditions consists of a large deviation
estimate for the trace of the propagator matrix $M_N$ which we
derive in Section~\ref{sec:trace}. We now describe the strategy
behind the proof of Theorem~\ref{th:1.2} in more detail. One first
shows that segments $E_1(x), E_2(x)$ as in the
matrix~\eqref{eq:1.twobytwo} exist. Invoking the estimates for the
separation of the Dirichlet eigenvalues and the zeros of the
Dirichlet determinants established in~\cite{Gol Sch2}, one next
shows that the resonance defined by $E_1, E_2$ leads to two new
eigenvalues $E^+(x), E^-(x)$ of the ``next scale'', see Figure~2. We
call the interval
\[\bigl(\max_{x\in J} E^-(x), \min_{x\in J} E^+(x)\bigr)
\]
a {\it pre-gap} at scale~$\oN$. The interval $J$ here is the common domain
of $E_1$ and $E_2$.  At this point one faces the obstruction of  a so-called
{\it triple resonance}. Recall that the resonance defined by~\eqref{eq:1.res} is
called a {\em double resonance} if, with $B\gg A$
from~\eqref{eq:1.res},
\begin{equation}\label{eq:1.doubleres}
\Big |E_{j_1}^{(N)} \xo - E_{j_3}^{(N)} (x + m'\omega, \omega) \Big
| > ({m'})^{-B}
\end{equation}
for any pair $(j_3, m') \ne (j_2, m)$, with $ m\ll m' \le \oN$,
where $\oN \asymp \exp\bigl(N^{\delta}\bigr)$ is the ``next scale''.
Otherwise it is called a triple (or higher order) resonance.

\noindent As  mentioned above the triple resonance obstruction already appears in
Sinai's perturbative method~\cite{Sin1}, see also Bourgain's paper on almost Mathieu~\cite{BouKAM}.
In fact, by the choice of a cosine-like
potential and for large $|\lambda|$ this type of resonance is
excluded in~\cite{Sin1} and~\cite{BouKAM}.
For general potentials,  it was shown by J.~Chan~\cite{Chan} that  if \[\big |\partial_x E_j^{(N)}(x,\omega)
\big | + \big |\partial_{xx} E_j^{(N)}(x,\omega)\big
|>\lambda(N,\omega)
>0\] with a suitable function~$\lambda(N,\omega)$, then triple resonances do
not occur for most $\omega$. In Section~\ref{sec:tripleelim} we follow a similar approach in  the case
where the graphs $E_j^{(N)}(x,\omega)$
have controlled slopes. Moreover, for the case of analytic potentials,
one can show that the triple (or higher) resonance obstruction can occur only
for a set of frequencies of Hausdorff dimension zero.

In order to run an ``induction on scales'' argument that shows how pre-gaps
in finite volume eventually lead to gaps in infinite volume, we invoke
the mechanism of counting complex zeros of the characteristic determinants of~$H_N\xo$
as developed in~\cite{Gol Sch2}. Using complexified
notations, the characteristic  determinants are as follows:

\begin{equation}\label{eq:1.Dirdet}
\begin{aligned}
f_N(z,\omega, E)  &= \det\bigl(H_N(z, \omega)-E\bigr)\\
& =
\begin{vmatrix}
 V\bigl(ze(\omega)\bigr) - E & -1 & 0  &\cdots &\cdots & 0\\[5pt]
-1 & V\bigl(ze(2\omega)\bigr) - E & -1 & 0 & \cdots & 0\\[5pt]
\vdots & \vdots & \vdots & \vdots & \vdots & \vdots\\
&&&&&-1 \\[5pt]
0 & \dotfill & 0 & -1 && V\bigl(ze(N\omega)\bigr) - E
\end{vmatrix}
\end{aligned}
\end{equation}
where
\begin{equation}\label{eq:1.Dirdet'}
f_{[a, b]}(z, \omega, E) = f_{b-a+1} \bigl(ze(a\omega), \omega,
E\bigr)
\end{equation}
For future reference, we remark that for any interval
$\Lambda\subset \ZZ$, $H_\Lambda(z,\omega)$ denotes the matrix
obtained from $H(z,\omega)$ by restriction to~$\Lambda$ with
Dirichlet boundary conditions; in addition, we let $f_\Lambda:=\det
(H_\Lambda-E)$. As already remarked above, we write $H_N$ and $f_N$
for $H_{[1,N]}$ and $f_{[1,N]}$, respectively (although occasionally,
the same notation will also be used relative to the interval $[-N,N]$).
It is well--know that
these functions are closely related to the monodromy (or propagator) matrices~\eqref{eq:1.monodr}.
In fact,
\begin{equation}\label{eq:1.mondet}
M_N(z,\omega, E) = \begin{bmatrix}
f_N(z,\omega, E) & - f_{N-1}\bigl(ze(\omega), \omega, E\bigr)\\
f_{N-1}(z,\omega, E) & - f_{N-2}\bigl(ze(\omega),\omega, E\bigr)
\end{bmatrix}
\end{equation}
By means of this relation, large deviation estimates and an
avalanche principle expansion  for the function $\log \big
|f_N(z,\omega, E)\big |$ were developed in \cite{Gol Sch2}.  In
Section~\ref{sec:basictools} we recall the statements of these results and
prove some corollaries. These corollaries, combined with a suitable
version of the Jensen formula (see (e) in Section~\ref{sec:basictools})
enable one to locate and count the zeros of $f_N(\cdot, \omega, E)$
in the annulus $\capo$ and its subdomains. In particular, this
technique allows one to claim that if \[E \in \bigl(\max_x E^-(x) +
\exp\bigl(-\oN^{1/2}\bigr), \min_x E^+(x) -
\exp\bigl(-\oN^{1/2}\bigr)\bigr),\] where $\bigl(\max E^-(x), \min
E^+(x)\bigr)$ is a pre-gap at scale $\oN$, then $f_{\oN}(\cdot,
\omega, E)$ has two complex zeros $\zeta_\ell = e(x_\ell +
iy_\ell)$, with $\exp\bigl(-N^\delta\bigr) > |y_\ell| >
\exp\bigl(-\oN^\delta\bigr)$, $\ell = 1, 2$.  This is due to the
absence of triple resonances and the stability of the number of
zeros of $f_{\oN}(\cdot, \omega, E)$ under small perturbations of
$E$. The most effective form of the last property consists of the
Weierstrass preparation theorem for $f_N(\cdot, \omega, E)$, which
is described in part~(g) of Section~\ref{sec:basictools}.  To complete the
description of the formation of a gap from a pre-gap we use the
translations of the segments $\bigl\{E_{j_1}^{(N)}(x), \ux,
\ox\bigr\}$ under the shifts $x \mapsto x + k\omega$. Using the
localization property of eigenfunctions on a finite interval (see
Section~\ref{sec:anderson}), we show that if a double resonance
\eqref{eq:1.res} occurs then the same is true for a sequence of
segments which are ``almost'' identical with the shifts $E_{j_1}(x+
k\omega), E_{j_2}(x+ k\omega)$, $1 \le k \le \oN\bigl(1-
0(1)\bigr)$. This method is explained in great detail in
Sections~\ref{sec:segments}-\ref{sec:formation}. In particular, the possible
locations of the "center of localization" of an eigenfunction plays
an important role with the ``bad case'' being when this center is
too close to the boundary of the finite volume interval. Due to this
method we obtain a whole sequence of complex zeros $\zeta_{k,\ell}
\cong e\bigl(x_\ell + k\omega + iy_\ell\bigr)$ of $f_N(\cdot,
\omega, E)$. So, the numbers
$$
\cM_N(E) = N^{-1} \#\bigl\{z: 1-\rho_N < |z| < 1 + \rho_N,
f_N(z,\omega, E) = 0\bigr\}
$$
$\rho_N = \exp\bigl(-N^\delta\bigr)$ decrease at least by $2-o(1)$
by going from scale $N$ to scale $\oN$, provided $E$ is in the
pre-gap.  After a finite number of inductive steps one can locate a gap
and complete the proof of Theorem~\ref{th:1.2}. Needless to say, the
zero counting mechanism for the Dirichlet determinants from~\cite{Gol Sch2}
 plays a crucial role here.
The relevant sections in this regard are~\ref{sec:jensen} and~\ref{sec:zerocount},
where the applications of the Jensen formula and the avalanche principle expansions,
respectively, are presented. For a summary of~\cite{Gol Sch2} see~\cite{Gol Sch4}, and
for a heuristic discussion of gaps in the context of this paper see~\cite{Gol Sch3}.

\section{A review of the basic tools}\label{sec:basictools}

In this section we give a sketch of the main ingredients of the
method developed in~\cite{Gol Sch2}. We of course do not reproduce
all the material from that paper in full detail, and refer the
reader for most proofs to~\cite{Gol Sch2}.
We start our discussion with the classical Cartan estimate for
analytic functions.

\bigskip
\noindent {\em (a)~Cartan Estimate}

\medskip
\begin{defi}\label{def:2.1}
Let $H \gg 1$.  For an arbitrary subset $\cB \subset \cD(z_0,
1)\subset \IC$ we say that $\cB \in \car_1(H, K)$ if $\cB\subset
\bigcup\limits^{j_0}_{j=1} \cD(z_j, r_j)$ with $j_0 \le K$, and
\begin{equation}\label{eq:2.1}
\sum_j\, r_j < e^{-H}\ .
\end{equation}
If $d$ is a positive integer greater than one and $\cB \subset
\prod\limits_{i=1}^d \cD(z_{i,0}, 1)\subset \IC^d$ then we define
inductively that $\cB\in \car_d(H, K)$ if for any $1 \le j \le d$ there
exists $\cB_j \subset \cD(z_{j,0}, 1)\subset \IC, \cB_j \in \car_1(H,
K)$ so that $\cB_z^{(j)} \in \car_{d-1}(H, K)$ for any $z \in \IC
\setminus \cB_j$,  here $\cB_z^{(j)} = \bigl\{(z_1, \dots, z_d) \in \cB:
z_j = z\bigr\}$.
\end{defi}

 \begin{remark}\label{rem:2.1} (a) This definition is consistent with the notation of
Theorem~4 in Levin's book~\cite{levin}, p.~79. \\
(b) It is important in the definition of $\car_d(H,K)$ for $d>1$
that we control both the measure and the complexity $K$ of each
slice $\cB_z^{(j)}$, $1\le j\le d$.
\end{remark}

The following lemma is a straightforward consequence of this definition.

\begin{lemma}
\label{lem:2.cart_12}\qquad

\begin{enumerate}
\item[{\rm{(1)}}] Let $\cB_j \in \car_d(H, K)$, $\cB_j \subset
\prod\limits^d_{j=1} \cD(z_{j,0}, 1)$, $j = 1, 2, \dots, T$.  Then $\cB
= \bigcup\limits_j\, \cB_j \in \car_d\bigl(H - \log T, TK\bigr)$.

\item[{\rm{(2)}}] Let \[ \cB \in \car_d(H, K),\quad \cB \subset
\prod\limits^d_{j=1} \cD\bigl(z_{j,0}, 1\bigr)\] Then there exists
\[\cB' \in \car_{d-1}(H, K),\quad \cB' \subset \prod\limits^d_{j=2}
\cD\bigl(z_{j, 0}, 1\bigr)\]such that $\cB_{(w_2, \dots, w_d)} \in
\car_1(H, K)$, for any $(w_2, \dots, w_d) \in \cB'$.
\end{enumerate}
\end{lemma}

\noindent Next, we generalize the usual Cartan estimate to several
variables.

\begin{lemma}
\label{lem:2.high_cart}
 Let $\varphi(z_1, \dots, z_d)$ be an analytic function defined
in a polydisk $\cP = \prod\limits^d_{j=1} \cD(z_{j,0}, 1)$, $z_{j,0} \in
\IC$.  Let $M \ge \sup\limits_{\uz\in\cP} \log |\varphi(\uz)|$,  $m \le \log
\bigl |\varphi(\uz_0)\bigr |$, $\uz_0 = (z_{1,0},\dots, z_{d,0})$.  Given $H
\gg 1$ there exists a set $\cB \subset \cP$,  $\cB \in
\car_d\bigl(H^{1/d}, K\bigr)$, $K = C_d H(M - m)$,  such that
\beeq
\label{eq:2.cart_bd}
\log \bigl | \varphi(z)\bigr | > M-C_d H(M-m)
\eneq
for any $z \in \prod^d_{j=1} \cD(z_{j,0}, 1/6)\setminus \cB$.
\end{lemma}
\begin{proof}
The proof goes by induction over $d$. For $d = 1$ the
assertion is Cartan's estimate for analytic functions. Indeed, Theorem~4 on page~79 in~\cite{levin}
applied to $f(z)=e^{-m}\varphi(z)$ yields that
\[ \log \bigl | \varphi(z)\bigr | > m-C H(M-m)=M -(CH+1)(M-m) \]
holds outside of a collection of disks $\{\cD(a_k,r_k)\}_{k=1}^K$
with $\sum_{k=1}^K r_k\lesssim \exp(-H)$. Increasing the constant
$C$ leads to~\eqref{eq:2.cart_bd}. Moreover, $K/5$ cannot exceed the
number of zeros of the function $\varphi(z)$ in the disk
$\cD(z_{1,0},1)$ counted with multiplicity, which is in turn
estimated by Jensen's formula, as $\lesssim M-m$, see the following
section. Although this bound on $K$ is not explicitly stated in
Theorem~4 in~\cite{levin}, it can be deduced from the proofs of
Theorems~3 and~4 in \cite{levin}. Indeed, one can assume that each
of the disks $\cD(a_k,r_k)$ contains a zero of $\varphi$, and it is
shown in the proof of Theorem~3 in \cite{levin} that no point is
contained in more than five of these disks. Hence we have proved the
$d=1$ case with a bad set $\cB\in \car_1(H,C(M-m))$, which is
slightly better than stated above (the $H$ dependence of $K$ appears
if $d>1$ and we will ignore some slight improvements that are
possible to the statement of the lemma due to this issue).

In the general case take $1 \le j \le d$ and consider $\psi(z) =
\varphi\bigl(z_{1,0}, \dots, z_{j-1,0}, z, z_{j+1,0}, \dots,
z_{d,0}\bigr)$.  Due to the $d=1$ case there exists $\cB^{(j)} \in
\car_1\bigl(H^{{1/d}}, C_1(M-m)\bigr)$,  such that
$$
\log \bigl |\psi(z)\bigr | > M-C_1 H^{1/d}(M-m)
$$
for any $z \in \cD\bigl(z_{j,0}, 1/6\bigr) \setminus \cB^{(j)}$.  Take
arbitrary $z_{j,1} \in \cD\bigl(z_{j,0}, 1/6\bigr) \setminus \cB^{(j)}$
and consider the function
\[\chi\bigl(z_1, z_2, \dots, z_{j-1}, z_{j+1}, \dots, z_d\bigr) =
\varphi\bigl(z_1, \dots, z_{j-1}, z_{j,1}, z_{j+1}, \dots,
z_d\bigr)\]
 in the
polydisk $\cP':=\prod\limits_{i\ne j} \cD\bigl(z_{i, 0}, 1\bigr)$.
Then
\begin{align*}
\sup_{\cP'} \log
\bigl | \chi(z_1, \dots, z_{j-1}, z_{j+1}, \dots, z_d)\bigr | &\le M, \\
\log \bigl |\chi(z_{1,0}, \dots, z_{j-1,0}, z_{j+1,0}, \dots,
z_{d,0}) \bigr | &>  M-CH^{1/d}(M-m).
\end{align*}
Thus $\chi$ satisfies the conditions
of the lemma with the same $M$ and with $m$ replaced with
\[M-CH^{1/d}(M-m) .\] We now apply the inductive assumption for $d-1$ and with
$H$ replaced with $H^{\frac{d-1}{d}}$ to finish the proof.
\end{proof}

Later we will need the following general assertion which is a combination of the Cartan-type
estimate of the previous lemma and Jensen's formula on the zeros of analytic functions,
see (e) of the present section.

\begin{lemma}
\label{lem:2.cart_zero} Fix some $\uw_0=(w_{1,0}, w_{2,0}, \dots,
w_{d,0})\in{\mathbb C}^d$ and suppose that $f(\uw)$ is an analytic
function in $\cP = \prod\limits^d_{j=1} \cD(w_{j,0},1)$.  Assume
that
$ M \ge \sup_{\uw\in \cP} \log |f(\uw)|$, 
and let
$m \le  \log |f(\uw_1)|$ for some $ \uw_1 = (w_{1,1}, w_{2,1}, \dots, w_{d,1}) \in \prod\limits^d_{j=1} \cD(w_{j,0}, 1/2)$.
Given $H \gg 1$ there exists $\cB'_H \subset \cP' = \prod\limits^d_{j=2} \cD(w_{j,0}, 3/4)$, $\cB'_H \in
\car_{d-1} \bigl(H^{1/d}, K\bigr)$, $K = CH(M - m)$ such that for any $\uw' = (w_2, \dots, w_d) \in \cP' \setminus \cB'_H$ the
following holds: if
\[ \log|f(\tilw_1, \uw') | < M-C_dH(M-m) \text{\ \ for some\ \ }\tilde{w}_1\in\cD(w_{1,0},1/2),\]
then there exists $\hat w_1$ with $|\hat w_1 - \tilw_1| \lesssim e^{-H^{\frac1d}}$ such that $f(\hat w_1, \uw') = 0$.
\end{lemma}
\begin{proof} Due to Lemma~\ref{lem:2.high_cart},
there exists $\cB_H \subset \cP$, $\cB_H \in
\car_d\bigl(H^{1/d}, K\bigr)$, $K = C_dH(M-m)$ such that for any $\uw \in \prod\limits^d_{j=1} \cD(w_{j, 0}, 3/4)\setminus \cB_H$
one has
\begin{equation}
\label{eq:2.fuw}
\log \big | f(\uw) \big | > M - C_d H(M - m)\ .
\end{equation}
By Lemma~\ref{lem:2.cart_12}, part (2), there exists $\cB'_H \subset
\prod\limits^d_{j=2} \cD\bigl(w_{j, 0}, 1\bigr)$, $\cB'_H \in
\car_{d-1} (H^{\frac{1}{d}},K)$ such that $\bigl(\cB_H\bigr)_{\uw'} \in
\car_1(H^{\frac{1}{d}},K)$ for any $\uw'=(w_2, \dots, w_d) \in \cB'_H$. Here $(\cB)_{\uw'}$
stands for the $\uw'$--section of $\cB$. Assume
$$
\log \big |f(\tilde w_1, \uw') \big | < M-C_dH(M-m)
$$
for some $\tilde w_1 \in \cD(w_{1,0},1/2)$, and $\uw'\in
\cP' \setminus \cB'_H$. Since $\bigl(\cB_H\bigr)_{\uw'} \in
\car_1(H^{\frac1d},K)$ there exists $r \lesssim \exp\bigl(-H^{1/d}\bigr)$ such that
\begin{equation}\nn
\bigl\{z: |z - \tilde w_1| = r\bigr\} \cap \bigl(\cB_H\bigr)_{\uw'} = \emptyset\ .
\end{equation}
Then in view of \eqref{eq:2.fuw},
$$
\log \big | f(z,\uw')|  > M - C_d H(M - m)
$$
for any $|z - \tilde w_1| = r$.  It follows from the maximum
principle  that $f(\cdot, \uw')$ has at least one zero in the disk
$\cD(\tilde w_1, r)$, as claimed.
\end{proof}

\medskip \noindent {\em (b)~Large deviation theorem for the
monodromies and their entries}

\medskip
Let $M_n(z,\omega, E)$ be the monodromies defined as in
\eqref{eq:1.monodr}. The entries of $M_n\zoe$ are the determinants
$f_{[1+a, N-b]} \zoe$, $a, b \in \{0,1\}$, see \eqref{eq:1.Dirdet},
\eqref{eq:1.Dirdet'}.  Let \[L(y,\omega, E)= \lim_{N\to\infty}
N^{-1}\int \log \|M_N(e(x+iy),\omega,E)\|\, dx
\] be the Lyapunov exponent. We shall assume throughout this paper
that the Lyapunov exponents are bounded away from zero; the positive
lower bounded on the Lyapunov exponent will typically be denoted
by~$\gamma$. We shall also adhere to the following convention
regarding constants for the remainder of the paper:

\begin{defi}
  \label{def:constants} Constants appearing in the paper will be denoted by $A,B,C$ as well as  $A_j, B_j, C_j$, $j\ge0$.
As a rule, they will be allowed to depend on $\omega,\gamma, V, E$. The dependence on~$V$ will only
be exclusively through $\rho_0>0$ and $\|V\|_{L^\infty(\cA_{\rho_0})}$ where $V$ is
analytic on the annulus~$\cA_{\rho_0}$.  Moreover, the dependence on $\omega$ will be only through
$a,c$ where $\omega\in\tor_{c,a}$. Finally, constants depending
on~$E$ will be uniform for $E$ ranging over bounded sets.  For any positive
numbers $a,b$ we let $a\les b$ denote $a\le Cb$ and $a\ll b$ denote
$a\le C^{-1}b$. Finally, $a\asymp b$ stands for $a\les b$ and $b\les
a$.
\end{defi}

We now state the large deviation estimates (LDEs) which are
fundamental to the arguments of this paper. It will be assumed
tacitly that $V$ is analytic on~$\cA_{\rho_0}$ for some $\rho_0>0$ and
Definition~\ref{def:constants} will be in force. However, we shall
for now {\em not assume} that $V$ is real-valued.

\begin{prop}\label{prop:2.ldtst}  Assume that $L(y,\omega, E)\ge\gamma > 0$
for some $y \in (-\rho_0/10, \rho_0/10)$, $\omega \in \tor_{c, a}$,
$E \in \IC$. Then  for any $N \ge 2$,
\begin{align}
 \mes\bigl\{x\in\tor\::\: \bigl|\log \big \|M_N\bigl(e(x+iy),\omega,
E\bigr)\big\| - NL(y,\omega, E)\bigr| > H\bigr\} &\le
\label{eq:2.ldtm} C\exp\bigl(-H\big /\bigl(\log N\bigr)^{C_0}\bigr) \\
 \mes\bigl\{x\in\tor\::\: \bigl|\log \big |f_N\bigl(e(x+iy),\omega, E\bigr)\big |-
NL(y,\omega, E)\bigr|> H\bigr\} &\le \label{eq:2.ldtd}
 C\exp\bigl(-H\big /\bigl(\log N\bigr)^{C_0}\bigr)
\end{align}
for all $H> \bigl(\log N\bigr)^{C_0}$.
\end{prop}

We remark that it makes no difference here whether we write $NL$ or
$NL_N$. This is due to the estimate
\[
0\le L_N(\omega,E) - L(\omega,E) \le \frac{C}{N} \quad \forall\;
N\ge1
\]
from~\cite{Gol Sch1}.
 The bound \eqref{eq:2.ldtm} for the monodromies (in
an even sharper form) is in~\cite{Gol Sch1}. In~\cite{Gol Sch2} it
is shown how to pass from~\eqref{eq:2.ldtm} to~\eqref{eq:2.ldtd},
see Sections~2, 3
 of that paper.

\bigskip
\noindent
{\em (c)~The avalanche principle expansion for the Dirichlet determinants}

\medskip
Another basic tool in this paper is the following avalanche
principle, see \cite{Gol Sch1} and~\cite{Gol Sch2}.

\begin{prop}
\label{prop:AP} Let $A_1,\ldots,A_n$ be a sequence of  $2\times
2$--matrices whose determinants satisfy
\begin{equation}
\label{eq:detsmall} \max\limits_{1\le j\le n}|\det A_j|\le 1.
\end{equation}
Suppose that \be
&&\min_{1\le j\le n}\|A_j\|\ge\mu>n\mbox{\ \ \ and}\label{large}\\
   &&\max_{1\le j<n}[\log\|A_{j+1}\|+\log\|A_j\|-\log\|A_{j+1}A_{j}\|]<\frac12\log\mu\label{diff}.
\ee Then
\begin{equation}
\bigl|\log\|A_n\cdot\ldots\cdot A_1\|+\sum_{j=2}^{n-1}
\log\|A_j\|-\sum_{j=1}^{n-1}\log\|A_{j+1}A_{j}\|\bigr| <
C\frac{n}{\mu} \label{eq:AP}
\end{equation}
with some absolute constant $C$.
\end{prop}

Combining this with the large deviation theorems from above yields
the following expansion for the determinants. Let $f_N\zoe$ be the
determinants defined as in (\ref{eq:1.Dirdet}), and let $L(\omega,
E)$ be the Lyapunov exponent as above but with $y=0$. As before, $V$
is analytic on~$\cA_{\rho_0}$. For convenience we will assume that $V$
is real-valued. As mentioned above, any constant depending on~$V$
depends only on~$\rho_0$ and $\|V\|_{L^\infty(\cA_{\rho_0})}$.

\begin{cor}\label{cor:apldt} Assume that $L(\omega, E)\ge \gamma > 0$ for
some $\omega \in \tor_{c,a}$, $E \in \IC$.  There exists $N_0 =
N_0(V, \omega, \gamma,E)$, $\rho^{(0)} = \rho^{(0)}(V, \omega,
\gamma,E) > 0$ such that for any $N \ge N_0(V, \omega, \gamma,E)$
and any integers $\ell_1, \dots, \ell_n$, $\bigl(\log N\bigr)^{C_0}
< \ell_j < cN$ (where $C_0=C_0(a)$ is a large constant),
$\sum\limits_j \ell_j = N$ the following expansion is valid:
\begin{equation}\label{eq:2.apexp}
 \log \bigl|f_N\bigl(e(x+iy),\omega, E\bigr) \bigr| = \sum^{n-1}_{j=1}\,
\log \big \|A_{j+1}(z)\, A_j(z)\big \| -
\sum^{n-1}_{j=2}\, \log \big \|A_j(z)\big \| +
O\bigl(\exp\bigl(-\ul^{1/2}\bigr)\bigr)
\end{equation}
for any $z = e(x + iy) \in \capo\setminus \cB_{N,\omega, E}$, where
\begin{align*}
\cB_{N,\omega, E} &= \mybigcup^{k_0}_{k=1} \cD\bigl(\zeta_k,
\exp\bigl(-\ul^{1/2}\bigr)\bigr)\ ,\quad \ul = \min_j\,
\ell_j,\quad  k_0 \lesssim N, \\
   A_m(z)
  &= M_{\ell_m}\bigl(ze(s_m\omega), \omega, E\bigr), \qquad m = 2,\dots,
  n -1 \\
  A_1(z) &= M_{\ell_1}(z, \omega, E) \begin{bmatrix} 1 & 0\\ 0
  & 0\end{bmatrix} \\
  A_n(z) &= \begin{bmatrix} 1 & 0\\ 0 &
  0\end{bmatrix}M_{\ell_n}\bigl(ze(s_n\omega), \omega, E\bigr)
  \end{align*}
and with $s_m = \sum\limits_{j< m} \ell_j$.
\end{cor}

A detailed derivation of this theorem can be found in Sections~2
and~3 of \cite{Gol Sch2}.

\bigskip
\noindent {\em (e)~Uniform upper estimates on the norms of monodromy
matrices} \label{sec:2.upper}

\medskip\noindent The proof of the uniform upper estimate is based on
an application of the avalanche principle expansion in combination
with the following useful general property of averages of
subharmonic functions.

\begin{lemma}
\label{lem:2.lipsub}
Let $1>\rho>0$ and suppose $u$ is subharmonic on~$\cA_\rho$ such that
$\sup_{z\in \cA_\rho} u(z)\le 1$ and $\int_{\tor} u(e(x))\,dx\ge0$. Then
for any $r_1,r_2$ so that $1-\frac{\rho}{2} < r_1,r_2 < 1+\frac{\rho}{2}$ one has
\[ |\langle u(r_1 e(\cdot)) \rangle - \langle u(r_2e(\cdot)) \rangle| \le C_\rho\,|r_1-r_2|,\]
here $\la v(\cdot)\ra = \int^1_0 v(\xi)\, d\xi$.
\end{lemma}

For the proof see Lemma~4.1 in~\cite{Gol Sch2}.
 This assertion immediately implies the following corollary regarding the continuity of
$L_N$ in~$y$.

\begin{corollary}
\label{cor:2.liplyap}
Let $L_N(y,\omega,E)$ and $L(y,\omega,E)$ be defined as above. Then
with some constant $\rho>0$ that is determined by the potential,
\[
|L_N(y_1,\omega,E) - L_N(y_2,\omega,E)| \le C|y_1-y_2|
 \text{\ \ \ for all\ \ \ }|y_1|,|y_2| < \rho
\]
uniformly in $N$. In particular, the same bound holds for $L$ instead of~$L_N$ so that
\[ \inf_{E} L(\omega,E)>\gamma>0 \]
implies that
\[ \inf_{E,|y|\ll \gamma} L(y,\omega,E) > \frac{\gamma}{2}.\]
\end{corollary}

The following result improves on the uniform upper bound on the
monodromy matrices from~\cite{BouGol} and~\cite{Gol Sch1}. The
$(\log N)^A$ error here (rather than $N^\sigma$, say, as
in~\cite{BouGol} and~\cite{Gol Sch1}) is crucial for the study of
the distribution of the zeros of the determinants and eigenvalues,
see Proposition~4.3 in~\cite{Gol Sch2}. We remind the reader of our
convention regarding constants, see Definition~\ref{def:constants}.

\begin{prop}
\label{prop:2.unifb} Assume $L(\omega, E)\ge\gamma > 0$, $\omega \in
\tor_{c,a}$. Then
\begin{equation}\nn
\sup\limits_{x\in\tor} \log \| M_N (x,\omega,E) \| \le
NL_N(\omega,E)+ C(\log N)^{C_0} \ ,
\end{equation}
for all $N\ge2$.
\end{prop}

We now list some applications of this upper bound. See Section~4
of~\cite{Gol Sch2}.

\begin{corollary}
\label{cor:2.uniflocb} Fix $\omega_1 \in \tor_{c,a}$ and $E_1 \in
\IC$, $|y| < \rho_0$. Assume that $L(y,\omega_1,E_1)\ge\gamma > 0$.
Then
\begin{align*}
 & \sup \bigl\{ \log \big \| M_N \bigl(e(x+iy), \omega,E\bigr) \big \| : |E - E_1| + |\omega - \omega_1|
< N^{-C},\, x \in \tor\bigr\}\\
&\qquad \le NL_N (y,\omega_1,E_1) + C(\log N)^{C_0}
\end{align*}
for all $N\ge2$.
\end{corollary}

The importance here lies with the large size of the perturbations: a
crude argument would only allow for perturbations of size~$e^{-CN}$.
To achieve the much larger size~$N^{-C}$ one needs to invoke the
avalanche principle with smaller factors of size $\ell\asymp \log N$
which is allowed by the sharp LDE (on scale~$\ell$) from~\cite{Gol
Sch1}.

\begin{corollary}
\label{cor:2.derlocb} Fix $\omega_1\in \tor_{c,a}$ and $E_1 \in
\IC$, $|y| < \rho_0$. Assume that $L(y,\omega_1,E_1)\ge\gamma > 0$.
Let $\partial$ denote any of the partial derivatives $\partial_x,
\partial_y,
\partial_E$ or $\partial_\omega$.  Then
\begin{align*}
 & \sup \bigl\{\log\big \| \partial M_N\bigl(e(x+ iy), \omega,E\bigr) \big \|: |E - E_1| + |\omega -
\omega_1| < N^{-C}, x \in \tor\bigr\}\\
&\qquad \le N L_N(y,\omega_1,E_1) + C(\log N)^{C_0}
\end{align*}
for all $N\ge2$. Here $C_1=C_1(a)$ and
$C=C(V,\rho_0,a,c,\gamma,E_1)$.
\end{corollary}
\begin{proof} Clearly, for all $x, y, \omega,E$,
\begin{align*}
&\partial M_N \bigl(e(x+ iy), \omega,E\bigr) \\& = \sum^N_{n=1}
M_{N-n} \bigl(e(x + n\omega + iy), \omega,E\bigr) \,\partial
\begin{bmatrix}
 V(e(x + n\omega + iy)) - E  & -1\\ 1 & 0\end{bmatrix} M_{n-1} \bigl(e(x + iy),
\omega,E\bigr)
\end{align*}
Since $|E - E_1| + |\omega - \omega_1| < N^{-C}$, the statement now
follows from Corollary~\ref{cor:2.uniflocb},
Corollary~\ref{cor:2.liplyap}, as well as the rate of convergence
estimate
\[
0\le L_N(\omega,E)- L(\omega,E) \le \frac{C}{N},\quad \forall \,
N\ge2
\]
from \cite{Gol Sch1}.
\end{proof}

The previous bound on the derivatives implies the following bound on
differences of propagator matrices.

\begin{corollary}
\label{cor:2.lipnorm}
Under the assumptions of the previous corollary,
\begin{align*}
& \bigl\| M_N \bigl(e(x + iy), \omega,E\bigr) - M_N \bigl(e(x_1 + iy_1),\omega_1,E_1\bigr) \bigr\|\\
\le & \bigl(|E - E_1| + |\omega - \omega_1| + |x - x_1| + |y -
y_1|\bigr) \cdot
 \exp \bigl(NL_N \bigl(y_1, \omega_1, E_1\bigr) + C(\log N)^{C_0}\bigr)
\end{align*}
provided $|E - E_1| + |\omega - \omega_1| + |x - x_1| < N^{-C}$,
$|y_1| < \rho_0/2$, $|y - y_1| < N^{-1}$. In particular,
\begin{equation}
\label{eq:2.lipnormt}
\begin{split}
\bigl| \log {\bigl | f_N\bigl(e(x+iy), \omega, E\bigr)\bigr|\over
\bigl| f_N\bigl(e(x_1 + iy_1), \omega_1, E_1\bigr)\bigr|} \bigr| &
\le
\bigl(|E -E_1| + |\omega - \omega_1| + |x - x_1|+ |y - y_1|\bigr)\\
&\qquad {\exp\bigl(NL(y_1,\omega_1, E_1\bigr) + C(\log
N)^{C_0}\bigr)\over \bigl| f_N\bigl(e(x_1 + iy_1), \omega_1,
E_1\bigr)\bigr|}\ ,
\end{split}
\end{equation} for all $N\ge2$
provided the right-hand side of~\eqref{eq:2.lipnormt} is less than $1/2$.
\end{corollary}
\begin{proof}
 For \eqref{eq:2.lipnormt} estimate
 \begin{equation}
 \label{eq:f_Ndiff2}
  \begin{split}
&| f_N\bigl(e(x+iy), \omega, E\bigr) - f_N\bigl(e(x_1 + iy_1),
\omega_1, E_1\bigr) |\\
&\les (|E -E_1| + |\omega - \omega_1| + |x - x_1|+ |y - y_1|) \sup
|d f_N(e(x'+iy'), \omega',E')|
  \end{split}
 \end{equation}
 where the supremum is taken over all $x',y',\omega',E'$ on the line
 joining $(x,y,\omega,E)$ to $(x_1,y_1,\omega_1,E_1)$ and $d$ stands for the derivative in all variables. By
 Corollary~\ref{cor:2.derlocb} we can bound
 \[
\sup |d f_N(e(x'+iy'), \omega',E')| \les \exp(NL(y_1,\omega_1, E_1)+
C(\log N)^{C_0})
 \]
Dividing \eqref{eq:f_Ndiff2} by $f_N\bigl(e(x_1 + iy_1), \omega_1,
E_1\bigr)$ therefore  yields
\begin{align*}
&  \Bigl|  {\bigl| f_N\bigl(e(x+iy), \omega,
E\bigr)\bigr|\over\bigl
| f_N\bigl(e(x_1 + iy_1), \omega_1, E_1\bigr)\bigr|} -1 \Bigr | \\
&\les \bigl(|E -E_1| + |\omega - \omega_1| + |x - x_1|+ |y -
y_1|\bigr) {\exp\bigl(NL(y_1,\omega_1, E_1\bigr) + C(\log
N)^{C_0}\bigr)\over \bigl| f_N\bigl(e(x_1 + iy_1), \omega_1,
E_1\bigr)\bigr|}
\end{align*}
By assumption, the right-hand side here is $<\frac12$. Hence,
\eqref{eq:2.lipnormt} follows by taking logarithms.
\end{proof}

A particular instance of this bound is the following one.

\begin{corollary}\label{cor:2.lognormcor} Using the notation of the previous corollary one has
\begin{equation}
\label{eq:2.lognormrat} \Bigl | \log {\big \|M_N\bigl(e(x + iy),
\omega, E\bigr) \big \|\over \big \|M_N \bigl(e(x_1 + iy_1),
\omega_1, E_1 \bigr)\big \|} \Bigr|  < C\exp \bigl(-(\log
N)^{C_0}\bigr)
\end{equation}
\begin{equation}\label{eq:2.logdetrat}
\Bigl| \log {\big | f_N\bigl(e(x+ iy), \omega, E\bigr) \big |\over
\big | f_N\bigl(e(x_1 + iy_1),\omega_1, E_1 \bigr)\big |} \Bigr |
< C\exp\bigl(-(\log N)^{C_0}\bigr)
\end{equation}
for any $|E - E_1| + |\omega - \omega_1| + |x - x_1| + |y - y_1| <
\exp \bigl(-(\log N)^{4C_0}\bigr)$, $x_1\in \cA_{\rho_0/2}
\setminus \cB_{y_1,\omega_1,E_1}$, where $\mes(
\cB_{y_1,\omega_1,E_1}) < \exp \bigl(-(\log N)^{C_0}\bigr)$,
$\compl(\cB_{y_1,\omega_1,E_1}) \le CN$. In
particular, 
\begin{equation} \label{eq:2.lexpregul} |L(y,\omega,
E) -L(y_1,\omega_1, E_1)| \le C\exp(-(\log N)^{C_0})
\end{equation}
provided $|E - E_1| + |\omega - \omega_1| +  |y - y_1| < \exp
\bigl(-(\log N)^{4C_0}\bigr)$. 
\end{corollary}

An important application of the uniform upper bounds is the following
analogue of Wegner's estimate from the random case. We provide the
proof here just to demonstrate how the previous corollaries can be
applied.

\begin{lemma}
\label{lem:wegner} Let $V$ be analytic and real-valued on $\tor$
as in the previous result. Suppose $\omega\in \tor_{c,a}$. Then for
any $E \in \IR$, $H \ge (\log N)^{C_0}$ one has
\begin{equation}
\label{eq:2.wegn} \mes \bigl\{x \in \tor\::\: \dist\bigl(\spec(
H_N(x,\omega)), E\bigr) < \exp(-H)\bigr\} \les \exp \bigl(-H/(\log
N)^{C_0}\bigr)
\end{equation}
for all $N\ge2$.  Moreover, the set on the left-hand side is
contained in the union of $\les N$ intervals each of which does not
exceed the bound stated in~\eqref{eq:2.wegn} in measure.
\end{lemma}
\begin{proof} By Cramer's rule
\begin{equation}
  \label{eq:cramerG} \bigl|\bigl(H_N(x, \omega) - E\bigr)^{-1} (k, m)\bigr| = {\big
|f_{[1, k]}\bigl(e(x), \omega, E\bigr)\big |\, \bigl| f_{[m+1, N]}
\bigl(e(x), \omega, E \bigr)| \over \big | f_N\bigl(e(x), \omega,
E\bigr)\big |}
\end{equation}
By Proposition~\ref{prop:2.unifb}
\begin{equation}
\nn \log \big | f_{[1, k]}\bigl(e(x), \omega, E\bigr)\big | + \log
\big |f_{[m+1, N]}\bigl(e(x), \omega, E\bigr)\big | \le NL(\omega,E)
+ C(\log N)^{C_0}
\end{equation}
for any $x \in \tor$.  Therefore,
\begin{equation}\label{eq:distHN}
  \big\| \bigl(H_N(x, \omega) - E\bigr)^{-1} \big\| \le N^2\
{\exp\bigl(NL(\omega,E) + C(\log N)^{C_0}\bigr)\over \big
|f_N\bigl(e(x), \omega, E\bigr)\big|}
\end{equation}
for any $x \in \tor$.  Since
$$\dist\bigl(\spec\bigl(H_N(x,\omega),
E\bigr)\bigr) = \big\| \bigl(H_N(x,\omega) -
E\bigr)^{-1}\big\|^{-1}\ ,
$$
the lemma follows from Proposition~\ref{prop:2.ldtst}.
\end{proof}

Next, we derive an important application of
Lemma~\ref{lem:2.cart_zero} and Proposition~\ref{prop:2.unifb} to
the Dirichlet determinants $f_N$. The constants $C_0,C_1,C_2$ depend
on~$\omega$ as explained above, see Definition~\ref{def:constants}.

\begin{corollary}
\label{cor:2.lexcepzero} Suppose $\omega\in \tor_{c,a}$. Given
$E_0\in\IC$ and $H > (\log N)^{C_2}$, $N\ge2$,  there exists
\[ \cB_{N,E_0, \omega}(H) \subset \IC, \qquad\cB_{N,E_0,\omega}(H)\in \car_1(\sqrt{H},HN^2)\]
such that for any $z \in \IC \setminus \cB_{N,E_0, \omega}(H)$ with
 $|\Im z|< N^{-1}$, and large $N$ the
following holds: If
\[ \log \big | f_N \bigl(e(z), \omega, E_1\bigr) \big | < NL(\omega,E_1) - H(\log N)^{C_2},\quad |E_0-E_1|<\exp(-(\log N)^{C_2}), \]
then $f_N\bigl(e(z), \omega,E\bigr) = 0$ for some $|E - E_1|
\lesssim   \exp(-\sqrt{H})$. Similarly, given $x_0\in\tor$ and
$|y_0|<N^{-1}$, let $z_0=e(x_0+iy_0)$. Then for any $H \gg 1$, the
following  holds: if
\[ \log \big |f_N\bigl(z_0,\omega,E\bigr)\big | < NL(\omega,E) - H(\log N)^{C_2}, \]
then $f_N\bigl(z, \omega,E\bigr) = 0$ for some $|z-z_0|\lesssim
\exp(-{H})$.
\end{corollary}
\begin{proof} Set $r_0=\exp(-(\log N)^{C_0})$ with some (large) constant $C_0=C_0(a)$ as above. Fix any $z_0$ with $|z_0|=1$
and consider the analytic function
\[ f(z,E)=f_N(z_0+(z-z_0)N^{-1},E_0+(E-E_0)r_0,\omega)\]
on the polydisk $\cP=\cD(z_0,1)\times\cD(E_0,1).$
Then, by Proposition~\ref{prop:2.unifb},
\[ \sup_{\cP}\log|f(z,E)| \le NL(E_0,\omega) + C(\log N)^{C_0} =M\]
and by the large deviation theorem,
\[ \log|f(z_1,E_0)|> NL(E_0,\omega) - (\log N)^{C_0}=m\]
for some $|z_0-z_1|<1/100$, say. By Lemma~\ref{lem:2.cart_zero} there exists
\[ \cB_{z_0,E_0,\omega}(H)\subset\IC, \qquad \cB_{z_0,E_0,\omega}(H)\in \car_1(\sqrt{H}, H(\log N)^{C_1}) \]
so that for any $z\in \cD(z_0,1/2)\setminus\cB_{z_0,E_0,\omega}(H)$ the following holds: If
\[ \log|f(z,E_1)|< NL(E_0,\omega)-H(\log N)^{C_1} \]
for some $|E_1-E_0|<1/2$, then there is $E$ with $|E_1-E|\lesssim
\exp(-\sqrt{H})$ such that $f(z,E)=0$. Now let $z_0$ run over a
$N^{-\frac32}$-net on $|z|=1$ and define $\cB_{N,E_0,\omega}(H)$ to
be the union of the sets $z_0+N^{-1}\cB_{z_0,E_0,\omega}(H)$. The
first half of the lemma now follows by taking $C_2$ sufficiently
large and by absorbing some powers of $\log N$ into~$H$ if needed.

The second half of the lemma dealing with zeros in the $z$ variable
can be shown without appealing to Lemma~\ref{lem:2.cart_zero}.
Indeed, we apply Cartan's estimate in $d=1$ directly to $u(\cdot)=
\log |f_N(\cdot,\omega,E)|$ on the disk $\cD(z_0,N^{-1})$. By the
preceding the Riesz mass of $u(\cdot)$  on this disk is at most
$(\log N)^{C_0}$. Hence, we can find a radius $r\asymp \exp(-H)$ so
that
\[
\min_{|z-z_0|=r}\log \big |f_N\bigl(z,\omega,E\bigr)\big | >
NL(\omega,E) - H(\log N)^{C_2}
\]
Now if
\[ \log \big |f_N\bigl(z_0,\omega,E\bigr)\big | < NL(\omega,E) - H(\log N)^{C_2}, \]
then from the maximum principle $f_N\bigl(z_0,\omega,E\bigr)=0$ for
some $|z-z_0|< r$ as claimed.
\end{proof}

Corollary~\ref{cor:2.lexcepzero} should be thought of as a {\em
converse} to the large deviation theorem in some sense; indeed, it
shows that if $\log|f_N|$ is too small at some point, then nearby
there must be a zero. In other words, and not surprisingly, zeros
are responsible for the failure of the large deviation estimates.

The following result allows us to translate separations of an
energy~$E_0$ from the spectrum of $H_N(x,\omega)$ into quantitative
lower bounds on $\log |f_N(x,\omega,E_0)|$. For it we need $V$ to be
real-valued on $\tor$. As usual, $\omega \in \tor_{c,a}$, and we
remind the reader that $C_1,C_2$ etc.\ depend on~$\omega$, see
Definition~\ref{def:constants}. Before proving it, we recall a basic
fact of Hermitian matrices. It will be applied repeatedly in this
paper.

\begin{lemma}
  \label{lem:herm} Suppose $A$ is a Hermitian $n\times n$-matrix.
  Further, let $B$ be another $n\times n$-matrix with $\|A-B\|< \eps$ in
  operator norm. Then
  \[
\dist(\spec(A),\spec(B))<\eps
  \]
\end{lemma}
\begin{proof}
  Suppose $z\in\IC$ satisfies $\dist(z,\spec(A))\ge\eps$.
Since $A$ is Hermitian, we see (from the spectral theorem) that
\[
\|(A-z)^{-1}\|\le \eps^{-1}
\]
Then
  \[
 R(z):=\sum_{n=0}^\infty (A-B)^n (A-z)^{-(n+1)}
  \]
converges as a Neuman series in operator norm.   Moreover, $R(z)
(B-z)= (B-z)R(z) = I$, the identity. Hence
$z\in\IC\setminus\spec(B)$ whence the lemma.
\end{proof}

We can now state another important type of converse of the
large deviation theorem.

\begin{cor}
\label{cor:2.uniexcepzero} Let $V$ be real-valued on $\tor$. Assume
that for sufficiently large $N$, $x_0 \in \tor$, $E_0 \in \IR$ one
has
$$
\bigl(E_0 - \eta, E_0 + \eta\bigr) \cap \spec (H_N(x_0, \omega)) =
\emptyset
$$
with $\eta \le \exp\bigl(-\bigl(\log N\bigr)^{C_1}\bigr)$.  Then
$$
\log \big |f_N\bigl(e(x_0), \omega, E\bigr) \big | > NL(\omega, E_0)
- (\log N)^{C_1} \log \frac{1}{\eta}
$$
for any $|E_0 - E|\le \frac{\eta}{2}$.
\end{cor}
\begin{proof}
Suppose that
\[
\log \big |f_N\bigl(e(x_0), \omega, E_1\bigr) \big | < NL(\omega,
E_0) - (\log N)^{C_1}\, \log \eta^{-1}
\]
for some $|E_1-E_0|\le  \frac{\eta}{2}$. Then there is $z_1\in\IC$
with $|z_1-x_0|\ll \eta$ so that
\[
f_N\bigl(e(z_1), \omega, E_1\bigr) =0
\]
Since $H_N(x,\omega)$ is Hermitian  for $x\in\tor$, it follows from
Lemma~\ref{lem:herm} that the eigenvalues $E_j^{(N)}(\cdot,\omega)$
satisfy
\[
|E_j^{(N)}(z,\omega) - E_j^{(N)}(x_0,\omega)|\le C|x_0-z| \quad
\forall\, z\in\cA_{{\rho_0}/2}
\]
In other words, there is some $E_2$ with $|E_2- E_0|<\eta$ such that
\[
f_N\bigl(e(x_0), \omega, E_2\bigr) =0
\]
However, this contradicts our assumption.
\end{proof}

We now address the important issue of a large deviation estimate with regard
to the $E$ variable.

\begin{lemma}
 \label{lem:moveE} Let $\omega_0\in\tor_{c,a}$ and $x_0\in\tor$. Then there exists $x_1\in\tor$
so that
\begin{align*}
|x_1-x_0| &<\exp(-(\log N)^{C_0}) \\
 \dist\big( \spec(H_N(x_1,\omega_0)),\,\spec(H_N(x_0,\omega_0)) ) &> \exp(-(\log N)^{C_1})
\end{align*}
where $C_0< C_1$.
\end{lemma}
\begin{proof}
Write
$
 \spec(H_N(x_0,\omega_0)) = \{ E_j(x_0,\omega_0)\}_{j=1}^N
$.
By Lemma~\ref{lem:wegner},
\[
 \mes \bigl\{x \in \tor\::\: \min_{1\le j\le N}\dist\bigl(\spec(
H_N(x,\omega)), E_j(x_0,\omega_0) \bigr) < \exp(- (\log
N)^{C_1})\bigr\} \les \exp \bigl(-(\log
N)^{C_0}\bigr)
\]
where $C_0<C_1$, and we are done.
\end{proof}

\begin{lemma}
 \label{lem:E_1lower} Let $\omega_0\in\tor_{c,a}$ and fix $x_0\in\tor$, $E_0\in\IR$. There exists
$|E_1-E_0|<\exp(-(\log N)^{C_0})$ with
\begin{equation}
 \label{eq:E_1lower}
\log|f_N(e(x_0),\omega_0,E_1)| > NL(E_0,\omega_0) - (\log N)^{C_2}
\end{equation}
where $C_2>C_0$.
\end{lemma}
\begin{proof}
 If
\[
 \dist(E_0, \spec(H_N(x_0,\omega_0)) ) > \exp(-(\log N)^{C_1})
\]
then
\[
 \log|f_N(e(x_0),\omega_0,E_0)| > NL(E_0,\omega_0) - (\log N)^{2C_1}
\]
by Corollary~\ref{cor:2.uniexcepzero}. Hence, in this case we can choose $E_1=E_0$.
Now assume that
\[
 \dist(E_0, \spec(H_N(x_0,\omega_0)) ) \le  \exp(-(\log N)^{C_1}),
\]
By the previous lemma we choose $|x_1-x_0|<\exp(-(\log N)^{C_0})$ such that
\begin{equation}\label{eq:spec_distx1}
 \dist\big( \spec(H_N(x_1,\omega_0)),\,\spec(H_N(x_0,\omega_0)) \big) > \exp(-(\log N)^{C_1})
\end{equation}
By self-adjointness, there exists $E_1\in\spec(H_N(x_1,\omega_0))$ with
\[
 |E_1 - E_0|< C\exp(-(\log N)^{C_0})
\]
which, in view of \eqref{eq:spec_distx1} also satisfies
\[
 \dist\big(\spec(H_N(x_0,\omega_0)),\,E_1 \big) >  \exp(-(\log N)^{C_1})
\]
By Corollary~\ref{cor:2.uniexcepzero} we conclude that
\[
 \log|f_N(e(x_0),\omega_0,E_1)| > NL(E_0,\omega_0) - (\log N)^{2C_1}
\]
The lemma follows with $C_2=2C_1$.
\end{proof}

We can now state the large deviation estimate with respect to the $E$-variable.

\begin{prop}
 \label{prop:LDEinE}   Let $\omega_0\in\tor_{c,a}$ and assume that $L(\omega_0,E)>\gamma>0$ for
all $E\in[E',E'']$. Then for large $N$, and all $x_0\in\tor$,
\begin{equation}\label{eq:Elde}
\mes\{E\in[E',E'']\::\:
|\log|f_N(e(x_0),\omega_0,E)|-NL(\omega_0,E)|>H\} \le
C\exp(-H/(\log N)^{C_1})
\end{equation}
for all $H>(\log N)^{2C_1}$.
\end{prop}
\begin{proof} Let $C_0$ be as in the previous lemma.
 Covering $[E',E'']$ by intervals of length $100\exp(-(\log N)^{C_0})$ we see that
it suffices to prove~\eqref{eq:Elde} locally on such an interval. Thus, consider a disk
$\cD(E_0, r_0)$ where $r_0= 100\exp(-(\log N)^{C_0})$. By Lemma~\ref{lem:E_1lower} there exists
$E_1\in \cD(E_0, r_0/100)$ with
\[
 \log|f_N(e(x_0),\omega_0,E_1)| > NL(E_0,\omega_0) - (\log N)^{C_2}
\]
On the other hand, there is the uniform upper bound
\[
 \sup_{E\in \cD(E_0, r_0/100)} \log |f_N(e(x_0),\omega_0,E_1)| \le NL(E_0,\omega_0) + (\log N)^{C_2}
\]
see Corollary~\ref{cor:2.uniflocb}. Now the proposition follows from Cartan's estimate.
\end{proof}

\begin{remark}
\label{rem:ldee_cartan}
Even though \eqref{eq:Elde} was stated for real $E$, one can pass to a version of this estimate in the
complex plane via Cartan's theorem: for all $H>(\log N)^{3C_1}$ there exist disks $\{\cD(\zeta_j,r_j)\}_{j=1}^J$
with $\sum_j r_j < \exp(-H(\log N)^{-2C_1})$, $J\le (\log N)^{C_2}$ and
\[
\{E\in[E',E'']+\cD(0,N^{-1})\::\:
|\log|f_N(e(x_0),\omega_0,E)|-NL(\omega_0,E)|>H\}  \subset \bigcup_j \cD(\zeta_j,r_j)
\]
for large $N$. This follows from Proposition~\ref{prop:LDEinE} by choosing $H=(\log N)^{2C_1}$ (where $C_1$ is large
depending on $(E',E'')$) which insures that there is at least one energy in $(E',E'')$ satisfying~\eqref{eq:Elde}.
Now apply Cartan's theorem as in part (a) of this section.
\end{remark}

We close this subsection with an  important consequence of the previous estimates; it allows
us to  bounds the number of zeros of the determinants with respect to both the~$z$ and~$E$ variables.

\begin{prop}
\label{prop:2.nzeroprop} Let $V$ be analytic on $\cA_{\rho_0}$ and
real-valued on~$\tor$. Let $\omega \in \tor_{c,a}$. Then for any
$x_0 \in \tor$, $E_0 \in \IR$ one has
\begin{align}
\label{eq:null1} \# \bigl\{E \in \IR\::\: f_N\bigl(e(x_0),
\omega,E\bigr) = 0,\;
 |E - E_0| < \exp \bigl(-(\log N)^{C_1}\bigr)\bigr\} &\le (\log N)^{C_1} \\
\# \bigl\{z \in \IC\::\: f_N(z, \omega, E_0) = 0,\; |z - e(x_0)| <
N^{-1}\bigr\} &\le (\log N)^{C_1} \label{eq:null2}
\end{align}
for all sufficiently large $N\ge N(V,\gamma,\rho_0,\omega,E_0)$.
\end{prop}
\begin{proof}
By the uniform upper bound
\[
 \sup \bigl\{ \log \big | f_N(e(x), \omega,E)\big | \::\: x \in \tor,\, E \in \IC, \,
 |E - E_1| < \exp\bigl(-(\log N)^{C_1}\bigr) \bigr\} \le NL_N(\omega,E_1) + (\log
 N)^{C_1}
\]
for any $E_1$. Due to the large deviation theorem with respect to the $E$
variable, see Proposition~\ref{prop:LDEinE},  there exist $x_1, E_1$
such that $|x_0 - x_1| < \exp\bigl(-(\log N)^{2C_1}\bigr)$, $|E_0 -
E_1| < \exp \bigl(-(\log N)^{2C_1}\bigr)$ so that
$$
\log \big | f_N(e(x_1),  \omega, E_1) \big | > NL_N(\omega,E_1) -
(\log N)^{C_1}.
$$
By Jensen's formula \eqref{eq:2.jensb},
$$
\# \bigl\{E\::\: f_N(e(x_1), \omega,E) = 0, |E - E_1| < \exp
\bigl(-(\log N)^{C_1}\bigr) \bigr\} \le 2(\log N)^{C_1}.
$$
Since $\big \| H_N (x_0, \omega) - H_N (x_1, \omega) \big \|
\lesssim \exp \bigl(- (\log N)^{2C_1}\bigr)$ and since $H_N (x_0,
\omega)$ is Hermitian one has
\begin{equation*}
\begin{split}
& \# \bigl\{E: f_N(e(x_0), \omega,E) = 0, |E - E_0| < \exp \bigl(-(\log N)^{2C_1}\bigr) \bigr\}\\
& \le \#\bigl\{E: f_N(e(x_1), \omega,E) = 0, |E- E_1| <
\exp\bigl(-(\log N)^{C_1}\bigr) \bigr\} \le (\log N)^{C_1}\ .
\end{split}
\end{equation*}
That proves~\eqref{eq:null1}.  The proof of \eqref{eq:null2} is
similar. Indeed, due to the uniform upper bound
\[
 \sup \bigl\{ \log \big | f_N(e(x+iy), \omega, E_0)\big | \::\: x \in \tor,\, |y|<2N^{-1}
  \bigr\} \le NL_N(\omega, E_0) + (\log N)^{C_1}.
\]
By the large deviation theorem, there is $x_1$ with $|x_0-x_1|<
\exp(-(\log N)^{C_1/2})$ such that
\[ \log \big | f_N(e(x_1), \omega, E_0)\big | > NL_N(\omega, E_0)-(\log N)^{C_1}\]
Hence, by Jensen's formula \eqref{eq:2.jensb},
\[ \# \bigl\{z\::\: f_N(z, \omega,E) = 0, |z-e(x_1)| < 2N^{-1} \bigr\} \le 2(\log N)^{C_1},\]
and \eqref{eq:null2} follows.
\end{proof}

\bigskip
\noindent {\em (g)~The Weierstrass preparation theorem for Dirichlet
determinants} \label{sec:weier}

\medskip Recall the Weierstrass preparation theorem for an analytic
function $f(z, w_1, \dots, w_d)$ defined in a polydisk
\begin{equation}\label{eq:2.polydisk}
\cP = \cD(z_0, R_0) \times \prod^d_{j=1} \cD(w_{j,0}, R_0),\quad z_0,\ w_{j, 0} \in \IC\qquad
\frac12\ge R_0 > 0\ .
\end{equation}

\begin{prop}
\label{th:2.weier}
Assume that $f(\cdot, w_1, \dots, w_d)$ has no zeros on some circle $\bigl\{z:
|z-z_0| = \rho_0 \bigr\}$, $0 < \rho_0 < R_0/2$, for any $\uw = (w_1, \dots, w_d) \in
\cP_1 = \prod\limits^d_{j=1} \cD(w_{j, 0}, r_1)$ where $0<r_1<R_0$.  Then there exist a polynomial $P(z, \uw)
= z^k +a_{k-1} (\uw) z^{k-1} + \cdots + a_0 (\uw)$ with $a_j(\uw)$ analytic in
$\cP_1$ and an analytic function $g(z, \uw), (z, \uw) \in \cD(z_0, \rho_0) \times \cP_1$
so that the following properties hold:
\begin{enumerate}
\item[(a)] $f(z, \uw) = P(z, \uw) g(z, \uw)$ for any $(z, \uw) \in \cD(z_0, \rho_0) \times
\cP_1$.

\item[(b)] $g(z, \uw) \ne 0$ for any $(z, \uw) \in \cD(z_0, \rho_0) \times \cP_1$

\item[(c)] For any $\uw \in \cP_1$, $P(\cdot, \uw)$ has no zeros in $\IC \setminus \cD(z_0,
\rho_0)$.
\end{enumerate}
\end{prop}
\begin{proof} By the classical Weierstrass argument,
$$
b_p(\uw) := \sum^k_{j=1} \zeta^p_j (\uw) = {1\over 2\pi i} \oint\limits_{|z-z_0|= \rho_0} z^p \
{\partial _z f(z, \uw) \over f(z, \uw)}\, dz
$$
are analytic in $\uw \in \cP_1$.  Here $\zeta_j (\uw)$ are the zeros
of $f(\cdot, \uw)$ in $\cD(z_0, \rho_0)$ counted with multiplicity.
Since the coefficients $a_j(\uw)$ are linear combinations of the
$b_p$, they are analytic in $\uw$.  Analyticity of $g$ follows by
standard arguments.
\end{proof}

Since there is an estimate for the local number of the zeros of the
Dirichlet determinant and also the local number of the  Dirichlet
eigenvalues, one can apply Proposition~\ref{th:2.weier} to  $f_N(z,
\omega,E)$. We need to do this in both the $z$ and the $E$
variables. See Section~6 of \cite{Gol Sch2} for more details. In
what follows recall the convention adopted in
Definition~\ref{def:constants}.

\begin{prop}
\label{prop:2.weierz}
 Given $z_0 \in \cA_{\rho_0/2}$, $E_0 \in \IC$, and $\omega_0 \in \tor_{c,a}$, there
 exists $N_0=N_0(V,\rho_0,a,c,\gamma)$ so that the following holds:
 for any $N\ge N_0$
 there
exists a polynomial
\[  P_N(z, \omega,E) = z^k + a_{k-1} (\omega,E) z^{k-1} + \cdots + a_0(E,
\omega)\]
 with $a_j(\omega,E)$ analytic in $\cD(E_0, r_1)\times \cD(\omega_0, r_1)$, $r_1 \asymp \exp
\bigl(-(\log N)^{C_1}\bigr)$ and an analytic function \[g_N(z,
\omega,E),\quad (z, \omega,E) \in \cP := \cD(z_0, r_0) \times
\cD(E_0, r_1) \times \cD(\omega_0,r_1)\] with $N^{-1}\le r_0 \le 2
N^{-1}$ such that:
\begin{enumerate}
\item[(a)] $f_N(z, \omega,E) = P_N(z, \omega,E) g_N(z, \omega,E)$

\item[(b)] $g_N(z, \omega,E) \ne 0$ for any $(z, \omega,E) \in \cP$

\item[(c)] For any $(\omega,E) \in  \cD(\omega_0, r_1)\times \cD(E_0, r_1) $, the
polynomial
 $P_N(\cdot, \omega,E)$ has no zeros in $\IC \setminus \cD(z_0, r_0)$

\item[(d)] $k = \degg P_N(\cdot, \omega,E) \le (\log N)^{C_0}$.
\end{enumerate}
\end{prop}
\begin{proof}
  With $r_0:= 2N^{-1}$ and $r_1:= \exp(-(\log N)^{C_1})$, we set
  \[
f(\zeta,w_1,w_2) := f_N(z_0+ N^{-1} \zeta, \omega_0+ r_1 w_1, E_0 +
r_1 w_2)\quad \forall\;(\zeta,w_1,w_2) \in \cD(0,1)^3
  \]
  Then by the uniform upper bound $ |f|\le \exp(NL(\omega_0,E_0)+(\log N)^{C_0})=:M $ on
  $\cD(0,1)^3$ and, by the large deviation theorem,
  \[|f(\zeta,0,0)| >  \exp\big(NL(\omega_0,E_0)-(\log N)^{C_0}\big)   \]
for all $|\zeta|=r$ and some $\frac12<r<1$. Moreover, by Cauchy's
estimate
\[
|f(\zeta,0,0) - f(\zeta,w_1,w_2)|\le 2 M(|w_1|+|w_2|)\le \frac{M}{2}
\exp(-2(\log N)^{C_0})\]
for all $|w_1|+|w_2|<
\frac14\exp(-2(\log N)^{C_0})$.
In particular,
\[
f(\zeta,w_1,w_2)\ne0 \quad\forall\;|\zeta|=r,\; |w_1|+|w_2|<
\frac14\exp(-2(\log N)^{C_0})
\]
The proposition follows by applying Proposition~\ref{th:2.weier} and a
rescaling.
\end{proof}

Later we shall need to localize Proposition~\ref{prop:2.weierz} to smaller regions in $E$ and~$z$. 

\begin{cor}
 \label{cor:weierz_loc} Using the notations of the previous proposition, let $0< r_2< e^{-(\log N)^{C_1}}$ be given.
With the same hypotheses, the conclusions of Proposition~\ref{prop:2.weierz}
hold on the smaller poly-disk \[ \cD(z_0,r_2')\times \cD(E_0, r_2'')\times \cD(\omega_0, r_2'')\] where  
\[
 r_2'\asymp r_2,\qquad r_2'' \asymp r_2^{3(\log N)^{C_1}}
\]
\end{cor}
\begin{proof}
 Apply the proposition and let $z_j(E,\omega)$ be the zeros of $P_N(\cdot,\omega,E)$. Then
\[
 P_N(z,\omega,E)= \prod_{j=1}^k (z-z_j(E,\omega))
\]
Select $r_2'$ so that
\[
 \inf_{|z-z_0|=r_2'} |P_N(z,\omega,E)|\ge  \big( r_2/ (\log N)^{C_1} \big)^{(\log N)^{C_1}}
\]
Since $|a_j(\omega,E)|\le 1$, it follows that 
\[
 \inf_{\substack{|E-E_0|\le r_2''\\|\omega-\omega_0|\le r_2''}} \;\;\inf_{|z-z_0|=r_2'} |P_N(z,\omega,E)|\ge  \frac12\big( r_2/ (\log N)^{C_1} \big)^{(\log N)^{C_1}}
\ge r_2^{2(\log N)^{C_1}}
\]
where $r_2''$ is as above. We can now apply Proposition~\ref{th:2.weier} as before. 
\end{proof}

The preparation theorem relative to $E$ is easier since we need it
only in the neighborhood of the unit circle,
 i.e., in the
neighborhood of points $e(x_0)$ with $x_0\in \tor$. In this case,
one can use the fact that $H_N(e(x_0),\omega)$ is Hermitian.

\begin{prop}
\label{prop:2.weierE}
 Given $x_0 \in \tor$, $E_0 \in \IC$, and $\omega_0\in \tor_{c,a}$, there exist a
polynomial
\[ P_N(z, \omega,E) = E^k + a_{k-1} (z, \omega)E^{k-1} + \cdots + a_0(z,
\omega) \] with $a_j(z, \omega)$ analytic in $\cD(z_0, r_1) \times
\cD(\omega_0, r_1)$, $z_0=e(x_0)$,  $r_1 \asymp \exp \bigl(- (\log
N)^{2C_1}\bigr)$ and an analytic function $g_N(z, \omega,E)$, $(z,
\omega,E) \in \cP = \cD(z_0, r_1) \times  \cD(\omega_0, r_1)\times
\cD(E_0, r_1)$ such that
\begin{enumerate}
\item[(a)] $f_N(z, \omega,E) = P_N(z, \omega,E) g_N(z, \omega,E)$

\item[(b)] $g_N(z, \omega,E) \ne 0$ for any $(z, \omega,E) \in \cP$

\item[(c)] For any $(z, \omega) \in \cD(z_0, r_1) \times \cD(\omega_0, r_1)$, the polynomial $P_N(z, \omega, \cdot)$
has no zeros in $\IC \setminus \cD(E_0, r_0)$, $r_0 \asymp \exp
\bigl(-(\log N)^{C_1}\bigr)$

\item[(d)] $k = \degg P_N(z, \omega,\cdot) \le (\log N)^{C_2}$
\end{enumerate}
\end{prop}
\begin{proof} Recall that due to Proposition~\ref{prop:2.nzeroprop} one has
$$
\#\bigl\{E\in \IC: f_N(z_0, \omega_0,E) = 0,\quad |E - E_0| <
\exp\bigl(-(\log N)^{C_1}\bigr) \bigr\} \le (\log N)^{C_1}
$$
Find $r_0 \asymp \exp\bigl(-(\log N)^{C_1}\bigr)$ such that
$f_N(z_0, \omega_0, \cdot)$ has no zeros in the annulus
\[ \bigl\{r_0(1-2N^{-2}) < |E - E_0| < r_0(1+2N^{-2})\bigr\}. \]
  Since $H_N(z_0, \omega_0)$ is self-adjoint,
  $f_N(z,\omega,\cdot)$ has no zeros in the annulus
\[ \bigl\{r_0 (1-N^{-2}) < |E - E_0| < r_0 (1+N^{-2})\bigr\},\]
provided $|z - z_0| \ll r_1:= r_0 N^{-4}$, $|\omega - \omega_0| \ll
r_1 $, see Lemma~\ref{lem:herm}.  The proposition now follows from Proposition~\ref{th:2.weier}.
\end{proof}

\section{The trace of $M_N$ and Hill's discriminant of the periodic problem} \label{sec:trace}

This section establishes large deviation estimates for the trace of $M_N$ as well
as other useful relations involving the trace. The importance of this section, which does not appear in~\cite{Gol Sch2},
lies with periodic boundary conditions: recall that the determinant of the Hamiltonian $H_{[1,N]}$
with periodic boundary conditions equals the trace of the Monodromy matrix $M_N$ up to a constant (the
latter trace is referred to as ``Hill's discriminant''). In our proof of gap formation periodic boundary
conditions play an important technical role, whence the relevance of this section.
Let us recall some  properties of matrices in $SL(2,\IR)$. It
follows from the polar decomposition that for any $M \in SL (2,\IR)$
there are unit vectors $\uu^+_M$, $\uu^-_M$, $\uv^+_M$, $\uv^-_M$ so
that $M\uu^+_M = \| M \| \uv^+_M $, $M\uu^-_M = \| M \|^{-1}
\uv^-_M$. Moreover, $\uu^+_M \perp \uu^-_M$ and $\uv^+_M \perp
\uv^-_M$.

\begin{lemma}\label{lem:trace}
For any $M\in SL(2,\IR)$,
\begin{equation}
\label{eq:tr} \|M^2\| - 4 \le \|M\|\,|\rtr M | \le \|M^2\| + 2.
\end{equation}
\end{lemma}
 \begin{proof}
 Due to the properties of the vectors $ \uu^+ = \uu^+_M$ , $\uu^- =
 \uu^-_M$ one has
\begin{equation}
\label{eq:spur} \rtr M = \|M\| \uv^{+}\cdot\uu^{+} + \|M\|^{-1}
\uv^{-}\cdot \uu^{-}
\end{equation}
On the other hand,
\begin{equation}
\label{eq:square}
 M^2 \uu^{+} = \|M\|^2 (\uu^+\cdot \uv^+) \uv^+ + (\uv^+\cdot \uu^-)\uv^{-}.
\end{equation}
It follows from \eqref{eq:spur} and \eqref{eq:square} that
\[ |\rtr M|\,\|M\| \le \|M\|^2\,|\uu^+\cdot\uv^+| + 1\le \|M^2\| + 2,\]
as well as
\[ |\rtr M|\,\|M\| \ge \|M\|^2\,|\uu^+\cdot\uv^+| - 1 \ge \|M^2\uu^+\| - 2.\]
Finally, using that $\|M\|\ge1$ one checks that $\|M^2 \uu^{-}\| \le
2$, and thus $\|M^2\uu^{+}\|\ge \|M^2\|-2$. Inserting this bound
into the last line finishes the proof.
\end{proof}

The following lemma establishes the large deviation estimate for a product
of mondromy matrices. The technical (albeit, important) twist here is that
we shift the phase in the second factor by a small but {\em fixed} amount. This will be essential for applications
to the trace. Indeed, in view of the previous lemma, in order to prove a large deviation theorem for $\rtr M_N$, say,
  it will be necessary to do the same for $M_N^2$. The latter should behave like $M_{2N}$, but more precisely it is equal to $M_N(x+N\omega+\kappa)M_N(x)$ where $\kappa\equiv -N\omega \;(\mod\;1)$.

\begin{lemma}
\label{lem:tracecond}
Assume that for some $\omega \in \tor_{c,a}$ and $E\in \IR$, one has $L(\omega, E) \geq\gamma > 0$.
Then there exists $\kappa_{0}=\kappa_{0}(V,\omega,\gamma,E)>0$
such that for any  $|\kappa|\le \kappa_0$
\begin{equation}\begin{split}
&\mes\bigl\{x\in\tor\::\: \bigl|\log \big \|M_N\bigl(e(x+N\omega+\kappa),\omega,
E\bigr)M_N\bigl(e(x),\omega,
E\bigr)\big\| - 2NL(E,\omega)\bigr| > H\bigr\} \\
&\le
\label{eq:LDTshift} C\exp\bigl(-H\big /\bigl(\log N\bigr)^{C_2}\bigr)\end{split}\end{equation}
for all $H>0$ and $N\ge2$.
\end{lemma}

\begin{proof} This will be done by induction in $N$; more precisely, we will introduce an increasing integer sequence
$\{N_j\}_{j\ge0}$ so that if \eqref{eq:LDTshift} holds for all $N_j\le N< N_{j+1}$, then it also holds
in the range $N_{j+1}\le N< N_{j+2}$. Clearly, by choosing  $N_0:=N_0(\gamma,V,\omega,E)$ large and
$
 \kappa_0:= \exp(-CN_1)
$
we see that the case $j=0$ can be made to hold for any $N_1$. Next, let $N_{j+1}\le N< N_{j+2}$
and set $n:= [(\log N)^{C_1}]$ where $C_1=2C_0$ with
$C_0$ as in~\eqref{eq:2.ldtm}.  Also, we define $N_{j+1}:=\exp\big(N_j^{\frac{1}{2C_1}}\big)$.
By the large deviation theorem from Section~\ref{sec:basictools} as well
as our inductive assumption (applied with $H=n^{2/3}$, say), there is an avalanche principle expansion of the form
\begin{align*}
 &\log\big \|M_N\bigl(e(x+N\omega+\kappa),\omega,
E\bigr)M_N\bigl(e(x),\omega,
E\bigr) \big\| - \log\big \|M_N\bigl(e(x+N\omega+\kappa),\omega,
E\bigr)\big \| \\
&\quad - \log\big\| M_N\bigl(e(x),\omega,
E\bigr) \big\| \\&= \log \|M_n \bigl(e(x+N\omega+\kappa),\omega,
E\bigr)M_n \bigl(e(x+(N-n)\omega),\omega,
E\bigr)\big \| - \log \|M_n \bigl(e(x+N\omega+\kappa),\omega,
E\bigr)\| \\
& \quad - \log \|M_n \bigl(e(x+(N-n)\omega),\omega,
E\bigr)\big \| + O(e^{-\sqrt{n}})
\end{align*}
for all  $x\in\tor\setminus \cB$ where $\mes(\cB)< e^{-\sqrt{n}}$.
In particular, this implies that
\[
 \log\big \|M_N\bigl(e(x_1+N\omega+\kappa),\omega,
E\bigr)M_N\bigl(e(x_1),\omega,
E\bigr) \big\| > 2NL(\omega,E) - (\log N)^{C_1}
\]
for all $x_1\in\tor\setminus \cB$.  Next, note from the uniform upper bounds in Part (e) of Section~\ref{sec:basictools} that
\[
 \sup_{x\in\tor,\,|y|<N^{-1}} \log \big \|M_N\bigl(e(x+iy+N\omega+\kappa),\omega,
E\bigr)M_N\bigl(e(x+iy),\omega,
E\bigr) \big\| < 2NL(\omega,E) + (\log N)^{C_1}
\]
By averaging, we conclude that
\[
 \Big|\int_0^1 \log \big \|M_N\bigl(e(x+iy+N\omega+\kappa),\omega,
E\bigr)M_N\bigl(e(x+iy),\omega,
E\bigr) \big\|\, dx - 2NL(\omega,E)\Big| < (\log N)^{2C_1}
\]
Furthermore,
by Cartan's theorem, see Lemma~\ref{lem:2.high_cart},  this yields \eqref{eq:LDTshift} for $N$ with $C_2=4C_1$, say.
\end{proof}

We can now state and prove the main result of this section. Note that even the average of $\log|\rtr M_N|$, which
appears as the first statement below, is far from being clear and requires much of the machinery developed so far
in this paper.

\begin{prop}
\label{th:2.avatrace}
 Assume that for some $\omega \in \tor_{c,a}$, $E\in \IR$ one has $L(\omega, E) \geq\gamma > 0$
 and let $\kappa_0>0$ be as in the previous lemma.
Then there exists
$N_{0}=N_0 (V,c,a,E,\gamma)$
such that for any $N\ge N_0$ satisfying $\|N\omega\|<\kappa_0$ the
following properties hold:
\begin{itemize}
\item  $\int_\tor \log|\rtr M_N(e(x),\omega,E)|\, dx = NL(\omega,E)+ O(1) $
\item  {\em Large deviation estimate for the traces:}
\[ \mes \big\{x\in[0, 1] : | \log | \rtr M_N(e(x+iy),\omega,E)|
 -N L(\omega,E)| > H \big\}\le C\exp (- H(\log N)^{-C_0})
\]
for all $|y|<N^{-1}$ and all $H>0$.
\end{itemize}
\end{prop}

\begin{proof} For simplicity, we set $y=0$. We begin with the simple observation that
\begin{equation}
\label{eq:twoNmonodr}
\log\|M_{N}(e(x),\omega,E)^2\| = \log\|M_{N}(e(x+N\omega + \kappa),\omega,E)
M_N(e(x),\omega,E)\|
\end{equation}
where $\kappa\equiv -N\omega\;(\mod\,1)$. By assumption, we can choose $|\kappa|<\kappa_0$ from Lemma~\ref{lem:tracecond}.
We apply the avalanche principle to $M_{N}(e(x+N\omega + \kappa),\omega,E)
M_N(e(x),\omega,E)$. To this end, define for $1\le j\le m$
\[
 A_j=M_{n_j}(e(x+t_j\omega),\omega,E), \; (\log N)^{C_0}<n_j< \sqrt{N},\; t_j=\sum_{i<j} n_i,\; \sum_{i=1}^m n_i =N
\]
as well as  $A_{j+m}(x):= A_j(x+ N\omega+\kappa)=A_j(x)$. We also require that $n_1=n_m$.
Then by the large deviation estimate~\eqref{eq:2.ldtm} and Lemma~\ref{lem:tracecond}, we conclude that
\begin{align}
 \log\|M_{N}(e(x),\omega,E)^2\| &= \sum^{2m-1}_{j=1}
\log\| A_{j+1}(x)A_j(x)\|- \sum^{2m-1}_{j=1} \log\|A_j(x)\|+0(e^{-\sqrt{\un}}) \nn \\
&= 2\Big[\sum^{m-1}_{j=1}
\log\| A_{j+1}(x)A_j(x)\|- \sum^{m-1}_{j=1} \log\|A_j(x)\|\Big] + \label{eq:trace_exp} \\
&\quad +  \log\| A_{1}(x)A_m(x)\| - \log \|A_1(x)\| - \log\|A_m(x)\| +  0(e^{-\sqrt{\un}}) \nn
\end{align}
for any $x\in \tor\setminus \cB$, with $\mes \cB <
\exp(-\un^{1/2})$ where $\un=\min_i n_i$. Interpreting the expression in brackets
as expansion of $\log \|M_N(e(x),\omega,E)\|$, and in view of Lemma~\ref{lem:trace}, we see that
\begin{equation}\label{eq:rtr_MN} \begin{split}
 \log|\rtr M_N(e(x),\omega,E)| - \log\|M_N(e(x),\omega,E)\| & = \log\| A_{1}(x)A_m(x)\| - \log \|A_1(x)\| \\
&\quad - \log\|A_m(x)\| + 0(e^{-\sqrt{\un}})
\end{split}
\end{equation}
for any $x\in \tor\setminus \cB$ with the same $\cB$.
The proposition now follows from Lemma~\ref{lem:tracecond} and the standard large deviation theorem
for the matrices $M_N$.
\end{proof}

For future reference, we remark that the {\em avalanche principle expansion} of $\log|\rtr M_N|$ given
by~\eqref{eq:trace_exp}, and the comparison statement~\eqref{eq:rtr_MN} are of independent interest.
Note that \[ \log|\rtr M_N(e(x),\omega,E)|\le \log \|M_N(e(x),\omega,E)\| + \log 2\]
 In particular, due to Proposition~\ref{prop:2.unifb},
one has the following uniform upper bound for the trace:

\begin{corollary}
\label{cor:2.unifbfortrace} Assume $L(\omega, E)\ge\gamma > 0$, $\omega \in
\tor_{c,a}$. Then
\begin{equation}\nn
\sup\limits_{x\in\tor} \log |\rtr M_N (x,\omega,E) | \le
NL_N(\omega,E)+ C(\log N)^{C_0} \ ,
\end{equation}
for all $N\ge2$, provided $\|N\omega\|\le \kappa_0(V,c,a,\gamma)$.

\end{corollary}
Now just as in Section~\ref{sec:basictools} one has the following
\begin{corollary}
\label{cor:3.liptr}  Assume $L(\omega, E_1)\ge\gamma > 0$, $\omega_1 \in
\tor_{c,a}$. Then for any $\omega\in \tor$, $y_1$, $y\in \IR$, $|y_1|$, $|y| \le 1/N$,
and any $x_1$,$x\in \tor$ one has
\begin{equation}
\label{eq:2.lipnormtrace}
\begin{split}
\bigl| \log {\bigl | \rtr M_N\bigl(e(x+iy), \omega, E\bigr)\bigr|\over
\bigl| \rtr M_N\bigl(e(x_1 + iy_1), \omega_1, E_1\bigr)\bigr|} \bigr| &
\le
\bigl(|E -E_1| + |\omega - \omega_1| + |x - x_1|+ |y - y_1|\bigr)\\
&\qquad {\exp\bigl(NL(y_1,\omega_1, E_1\bigr) + C(\log
N)^{C_0}\bigr)\over \bigl| \rtr M_N\bigl(e(x_1 + iy_1), \omega_1,
E_1\bigr)\bigr|}\ ,
\end{split}
\end{equation}
provided  $\|N\omega\|\le \kappa_0(V,c,a,\gamma)$ and the right-hand
side of~\eqref{eq:2.lipnormtrace} is less than $1/2$.
\end{corollary}

With these results on the trace at our disposal, we now turn to their implications
 for the periodic problem. To fix notation, let
\begin{equation}\label{eq:5.schred}
\bigl[ H_{[1,N]}^{(\pm P)}(x,\omega)\psi \bigr] (n) = -\psi
(n-1)-\psi(n+1) + V(e(x+n\omega))\psi(n)
\end{equation}
be the Schr\"odinger operator on $[1,N]$ with periodic
(respectively, antiperiodic) boundary conditions.
\begin{equation}\label{eq:5.bc}
\psi(0)=\pm\psi(N), \,\, \psi(1)=\pm\psi(N+1)
\end{equation}
Let
\begin{equation}\label{eq:5.eva}
E_1^{(N,\pm P)}(x,\omega)\le E_2^{(N,\pm P)}(x,\omega)\le \dots \le
E_N^{(N,\pm P)}(x, \omega)
\end{equation}
be the eigenvalues of $H_{[1,N]}^{(\pm P)}(x,\omega)$. Recall that
the characteristic determinant
\[
g_N^{(\pm)}(e(x), \omega, E):= \det \bigl(H_{[1,N]}^{(\pm
P)}(x,\omega)-E \bigr)
\]
which takes the form
\begin{equation}\nonumber
g_N^{(\pm)}(e(x),\omega,E)=\det \left[{\begin{array}{cccccc}
V(e(x+\omega)-E & -1 & 0 & \dots & 0 & \mp1 \\
-1 & V(e(x+2\omega)-E & -1 & \dots & 0 & 0 \\
  \ldots &\ldots &\ldots &\ldots &\ldots & \ldots  \\
  \ldots &\ldots &\ldots &\ldots &\ldots & \ldots \\
  \ldots &\ldots &\ldots &\ldots &\ldots & \ldots \\
0 & 0 & \dots & \dots & & -1 \\
\mp 1 & 0 & \dots & \dots & -1 & V(e(x+N\omega)-E \\
\end{array}}\right]
\end{equation}
satisfies
\begin{equation}\label{eq:5.pdethill}
g_N^{(\pm)}(e(x),\omega,E)= h_N(e(x),\omega,E)\mp 2
\end{equation}
where
\begin{equation}\label{eq:5.hill'}
h_N(e(x),\omega,E):=\rtr M_N(e(x),\omega,E)=f_{[1,N]}(e(x),\omega,E)-f_{[2,N-1]}(e(x),\omega,E)
\end{equation}
is Hill's discriminant. Cramer's rule then yields
\begin{equation}\label{eq:5.grep}
(H_{[1,N]}^{(\pm P)}(x,\omega)-E)^{-1}(m, n) =
\displaystyle\frac{f_{[1,m-1]}(e(x),\omega,E)f_{[n+1,N]}(e(x),\omega,E)}{g_N^{(\pm)}(e(x),\omega,E)},
\,\, 1\le m\le n\le N
\end{equation}
The large deviation estimate for $\log|\rtr
M_N(e(x),\omega,E)|$ from above implies the following lemma concerning the spectrum of the periodic problem.
In particular, we obtain an analogue of the Wegner bound.

\begin{lemma}\label{lem:5.basicperiod}
Let $V(e(x))$ be real analytic,
and let $\omega\in\tor_{c,a}$. Assume that $L(\omega,E)\ge
\gamma>0$
 for some $E\in\IR$. Then there exists $N_0=N_0(V,a,c,\gamma,E)$ such that for any $N\ge N_0$
 with  $\|N\omega\|\le \kappa_0(V,c,a,\gamma)$ the following properties hold:
\begin{itemize}
 \item[(a)] For any $H>(\log N)^{2C_0}$ one has
\[
\mes \{x\in\tor: \dist (\spec ( H_{[1,N]}^{(\pm
P)}(x,\omega)),E)<\exp(-H) \} \le \exp(-H/(\log N )^{C_0})
\]

\item[(b)] If for some $x\in\tor$
\[
(E-\eta, E+\eta) \cap \spec \big(H_N^{(\pm P)}(x,\omega)\big) = \emptyset
\]
with $\eta \le \exp(-(\log N)^{2C_0})$ then
\[
\log |g_N^{(\pm P)}(e(x'),\omega, E')|>NL(E,\omega)-(\log
N)^{C_1}\log\frac{1}{\eta}
\]
for any $|x'-x| + |E'-E|<\eta /C$.
\end{itemize}
\end{lemma}
\begin{proof}
The proof of (a) is analogous to the proof of Lemma~\ref{lem:wegner}; the only difference is that
we apply the large deviation theorem for the traces instead of for the determinants. Part~(b) is
analogous to Corollary~\ref{cor:2.uniexcepzero}, at least when $x'=x$,  and we skip the proof.
Finally, to move $x$ to $x'$ one uses Corollary~\ref{cor:3.liptr}.
\end{proof}

To close this section, we prove the following
 large deviation theorem for the traces with respect to the $E$ variable, cf.~Proposition~\ref{prop:LDEinE}.
It will play an important role later when we start counting eigenvalues by means of the Jensen average
machinery which is subject of the following section.

\begin{prop}
 \label{prop:LDEinEtrace}   Let $\omega_0\in\tor_{c,a}$ and assume that $L(\omega_0,E)>\gamma>0$ for
all $E\in[E',E'']$. Then for large $N$, with  $\|N\omega\|\le
\kappa_0(V,c,a,\gamma)$ and all $x_0\in\tor$,
\begin{equation}\label{eq:EldegN}
\mes\{E\in[E',E'']\::\:
|\log|g_N^{(\pm P)}(e(x_0),\omega_0,E)|-NL(\omega_0,E)|>H\} \le
C\exp(-H/(\log N)^{C_1})
\end{equation}
for all $H>(\log N)^{2C_1}$. The same statement applies to $\rtr M_N(e(x_0),\omega_0,\cdot)$.
\end{prop}
\begin{proof}
 This proof is completely analogous to the proof of Proposition~\ref{prop:LDEinE}. Indeed, the previous
Lemma~\ref{lem:5.basicperiod} replaces the Wegner type lemma used there as well as Corollary~\ref{cor:2.uniexcepzero}.
\end{proof}

\section{Zeros, eigenvalues, and the Jensen formula} \label{sec:jensen}

\medskip
This section introduces a key element in our approach
to the problem of determining
the location
of the spectrum and of the spectral gaps. More specifically, we identify the spectral values
and the spectral gaps according to whether $f_N(\cdot,\omega
,E)$ has a sequence of real or  complex zeros in the annulus
$\cA_{\rho_N} = \bigl\{z \in \IC: 1 - \rho_N < |z| < 1 +
\rho_N\bigr\}$ with appropriate $\rho_N$.  We would also like to
single out Lemma~\ref{lem:2.avernstab} below. It guarantees that the gaps
in finite volume stabilize after a finite (and uniformly bounded) number of iterations
in our ``induction on scales'' procedure. An important feature of the machinery developed here
lies with the fact that it applies equally well to the $z$-variable as to the $E$-variable
of~$f_N(z,\omega,E)$. Another important feature is the ``linearity'' of our bounds which means
that the zero count is additive. 

Now for the details, which for the most part already appear in~\cite{Gol Sch2} (with the exception of 
the crucial Lemma~\ref{lem:2.avernstab}).
The Jensen formula states
that for any function $f$ analytic on a neighborhood of
$\cD(z_0,R)$, see~\cite{levin}, \beeq \label{eq:2.jensb}
 \int_0^1 \log |f(z_0+Re(\theta))|\, d\theta - \log|f(z_0)| = \sum_{\zeta:f(\zeta)=0} \log\frac{R}{|\zeta-z_0|}
\eneq provided $f(z_0)\ne0$. In the previous section, we showed how
to combine this fact with the large deviation theorem and the
uniform upper bounds to bound the number of zeros of $f_N$ which
fall into small disks, in both the $z$ and $E$ variables. In what
follows, we will refine this approach further. For this purpose, it
will be convenient to average over $z_0$ in~\eqref{eq:2.jensb}.
Henceforth, we shall use the notation
\begin{align}
\nu_f(z_0, r) &= \# \{z \in \cD(z_0, r): f(z)=0\} \label{eq:2.zerondisk}\\
\cJ(u, z_0, r_1, r_2) &= \mathop{\nint}\limits_{\cD(z_0, r_1)}  \cJ(u, z, r) \, dx\,
dy \label{eq:2.jensenaver'}
\end{align}
\begin{align}
\cJ(u, z, r) &=\mathop{\nint}\limits_{\cD(z, r)} d \xi d \eta\,
[u(\zeta)-u(z)] \label{eq:2.jensaver}
\end{align}
where $z=x+iy, \zeta=\xi+i\eta$. Recall that
$
\cJ(u, z, r) \ge 0 $
for any subharmonic function $u$.  

\begin{lemma}
\label{lem:2.jensencor} Let $f(z)$ be analytic in $\cD(z_0, R_0)$.
Then for any $0<r_2<r_1<R_0-r_2$
\begin{equation*}
\nu_f(z_0, r_1 - r_2) \leq 4 \frac{r_1^2}{r_2^2} \cJ(\log |f|, z_0,
r_1, r_2) \leq \nu_f(z_0, r_1+r_2)
\end{equation*}
\end{lemma}
\begin{proof}
 Jensen's formula yields
\begin{align*}
\cJ(\log |f|, z_0, r_1, r_2) &= \mathop{\nint}\limits_{\cD(z_0,
r_1)} dx\, dy \bigg [ \frac{2}{r^2_2} \mathop\int\limits_0^{r_2} dr
\bigg (r \mathop\sum\limits_{f(\zeta)=0, \zeta \in \cD(z, r)} \log
\Big(\frac {r}{|\zeta -
z|} \Big)\bigg ) \bigg ]\\
&\leq \mathop\sum\limits_{f(\zeta)=0, \zeta \in \cD(z_0, r_1+r_2)}
\frac{1}{\pi r_1^2} \bigg [ \frac{2}{r^2_2}
\mathop\int\limits_0^{r_2} dr \bigg ( r
\mathop\int\limits_{\cD(\zeta, r)} \log \Big(\frac{r}{|z - \zeta |} \Big)
dx\, dy \bigg)
\bigg]\\
&= \frac{1}{4} \frac{r^2_2}{r^2_1} \nu_f (z_0, r_1+r_2),
\end{align*}
which proves the upper estimate for $\cJ(\log |f|, z_0, r_1, r_2)$.
The proof of the lower estimate is similar.
\end{proof}

\begin{corollary}
\label{cor:2.jensenequal} Let $f$ be analytic in $\cD(z_0, R_0)$, $0
< r_2 < r_1 < R_0 - r_2$. Assume that $f$ has no zeros in the
annulus $\cA = \bigl\{r_1 - r_2 \le |z - z_0| \le r_1 +
r_2\bigr\}$.  Then
$$
\nu_f(z_0, r_1) = 4\ {r_1^2\over r_2^2}\ \cJ\bigl(\log |f|, z_0,
r_1, r_2\bigr)
$$
\end{corollary}

\begin{corollary}
\label{cor:2.jensencomp} Let $f(z), g(z)$ be analytic in $\cD(z_0,
R_0)$.  Assume that for some $0<r_2<r_1<R_0-r_2$
\begin{equation*}
|\cJ(\log |f|, z_0, r_1, r_2)-\cJ(\log |g|, z_0, r_1, r_2)|< \frac
{r^2_2}{4r_1^2}
\end{equation*}
Then
$
\nu_f(z_0, r_1-r_2) \leq \nu_g(z_0, r_1+r_2)$, $\; \nu_g(z_0,
r_1-r_2) \leq \nu_f(z_0, r_1+r_2).
$
\end{corollary}

We shall also need  a simple generalization of these estimates to
averages over general domains. More precisely, set
\begin{align}
\nu_f(\cD) &= \# \{z \in \cD: f(z)=0\} \label{eq:2.jensenndom}\\
\cJ(u, \cD, r_2) &= \mathop{\nint}\limits_{\cD} dx\, dy
\mathop{\nint}\limits_{\cD(z, r_2)} d \xi d \eta\, [u(\zeta)-u(z)].
\nn
\end{align}
Given a domain $\cD$ and $r>0$ , set $ \cD (r)= \{z: \dist(z,\cD) <r
\}$. Let $f(z)$ be analytic in $\cD(R)$. Then for any
$0<r_2<r_1<R-r_2$
\begin{equation}
\label{eq:2.jensencordom} \nu_f( \cD(r_1 - r_2)) \leq 2
\frac{\mes(\cD)}{\pi r_2^2} \cJ(\log |f|, \cD( r_1), r_2) \leq
\nu_f( \cD(r_1+r_2))
\end{equation}

\noindent  Let $\cA_{R_1,R_2} := \{z\in \IC\::\: R_1<|z|<R_2 \}$.

\begin{lemma}\label{lem:2.jensenannul}
\begin{equation}\label{eq:2.jensenannulf}
\begin{aligned}
& N^{-1}\cJ\bigl(\log \big |f_N(\cdot, \omega, E)\big |, \cA_{R_1,R_2}, r_2\bigr)=\\
& 4(R_2^2-R_1^2)^{-1} r_2^{-2}\int^{R_2}_{R_1}\, \rho\, d\rho
\int_0^{r_2}\, r\,dr
 \int_0^1\, dy \bigl[L_N(\xi(\rho,r,y), \omega, E ) -
 L_N(\xi(\rho),\omega,E)\bigr]
\end{aligned}
\end{equation}
where $\xi(\rho,r,y)=\log|\rho+re(y)|$, $\xi(\rho)=\log\rho$.
\end{lemma}
\begin{proof}
Due to the definition of $\cJ(u,{\mathcal D}, r_2)$ one has
\begin{align*}
&N^{-1}\cJ(\log|f_N(\cdot,\omega,E)|, \cA_{R_1,R_2},r_2) \nn\\
&= \frac{4\pi N^{-1}}{ |\cA_{R_1,R_2}| r_2^{2}} \int_{R_1}^{R_2}\!
\!\!\!\rho\, d\rho \int_0^{r_2} \!\!r\, dr \biggl \{ \int_0^1\!\!dx\int_0^1\!\! dy
\bigl[ \log|f_N(\rho e(x)+r
e(y),\omega,E)|-\log|f_N(\rho e(x),\omega,E)| \bigr]\!\!\biggr\}\nn \\
&=\frac{4\pi N^{-1}}{ |\cA_{R_1,R_2}| r_2^{2}} \int_{R_1}^{R_2}\!
\!\!\!\rho\, d\rho \int_0^{r_2} \!\!r\, dr \biggl \{ \int_0^1\!\!dx\!\!\int_0^1\!\! dy
\bigl[ \log|f_N(|\rho +r e(y)|e(x),\omega,E)|-\log|f_N(\rho
e(x),\omega,E)|
\bigr]\!\!\biggr\}\nn\\
&= 4(R_2^2-R_1^2)^{-1} r_2^{-2}\int^{R_2}_{R_1}\, \rho\, d\rho
\int_0^{r_2}\, r\,dr
 \int_0^1\!\! dy \bigl[L_N(\xi(\rho,r,y), \omega, E ) -
 L_N(\xi(\rho),\omega,E)\bigr]
\end{align*}
as claimed.
\end{proof}

\noindent Set
\begin{equation}\label{eq:2.avernzero}
\cM_N(\omega, E, R_1, R_2)  := N^{-1} \# \bigl\{z \in \cA_{R_1,
R_2}: f_N\zoe = 0\bigr\}
\end{equation}

\begin{remark}\label{rem:5.totnumzerosupperbound} Recall that
$\log|f_N(z,\omega,E)|\le C(V)N$.  Corollary~\ref{cor:2.jensenequal} and the previous lemma 
therefore imply that 
\[
\cM_N(\omega, E, R_1, R_2)\le C(V)
\]
for any $N$ and $ R_1, R_2$. Furthermore, if $V(e(x))$ is a
trigonometric polynomial of degree $k_0$ then
\[
V(z)=z^{-k_0}P(z)
\]
where $P(z)$ is a polynomial of degree $2k_0$. Hence, with $\omega$
and $E$ being fixed one has
\[
z^{Nk_0}f_N(z,\omega,E)=F_N(z)
\]
where $F_N(z)$ is a polynomial of degree $2Nk_0$. Therefore, in this
case
\[
\cM_N(\omega, E, R_1, R_2)\le 2k_0
\]
which will be crucial for the $2k_0$ bound at the end of Theorem~\ref{th:1.sinai}. 
\end{remark} 

The following lemma allows us to compare
these averages for different scales. Later, this will be the crucial
device that prevents ``pre-gaps'' from collapsing at subsequent
stages of the iteration.

\begin{lemma}\label{lem:2.avernstab} Assume $\gamma = L(\omega, E) > 0$
and fix some small $0<\sigma\ll 1$.  There exist $N_0 = N_0(V,
\omega, \gamma,\sigma)$, $\rho^{(0)} = \rho^{(0)}(V, \omega, \gamma)
> 0$ such that for any $n > N_0$, $N > \exp(\gamma\, n^\sigma)$,
$1- \rho^{(0)} < R_1 < R_2 < 1 + \rho^{(0)}$ one has
\begin{equation}\label{eq:2.avernzel}
\begin{aligned}
\cM_N(\omega, E, R_1 + r_2, R_2 - r_2) & \le \cM_n(\omega, E, R_1 - r_2, R_2 + r_2) + n^{-1/4}\\[5pt]
\cM_n(\omega, E, R_1 + r_2, R_2 - r_2) & \le \cM_N(\omega, E, R_1 -
r_2, R_2 + r_2) + n^{-1/4}
\end{aligned}
\end{equation}
where $r_2 = n^{-1/4}(R_2 - R_1)$ and provided $r_2>\exp(-\gamma\,
n^\sigma/100)$.
\end{lemma}

\begin{proof} Recall that due to avalanche principle expansion one has
\begin{equation*}
\begin{aligned}
& \bigg| \log {\big \|M_n\bigl(e(x + n\omega + iy), \omega, E\bigr)
\big \|\ \big \|M_n\bigl(e(x + iy), \omega, E\bigr) \big \|\over
\big \|M_{2n}\bigl(e(x + iy), \omega, E\bigr) \big \|} -\\
& \log{\big\|M_\ell\bigl(e(x + n\omega + iy)\bigr) \big \|\  \big
\|M_\ell \bigl(e(x + (n - \ell)\omega + iy), \omega, E\bigr) \big
\|\over \big \|M_{2\ell} \bigl(e(x + (n-\ell)\omega + iy\bigr) \big
\|}\bigg | \le \exp\bigl(-\gamma_1\, n^{1/2}\bigr)
\end{aligned}
\end{equation*}
for any $|y| < \rho_0/2$, $x \in \tor \setminus \cB_y$, $\mes \cB_y
< \exp\bigl(-\gamma_1\, n^{1/2}\bigr)$ where $\ell =
\bigl[n^{1/2}\bigr]$, $\gamma_k = L(\omega, E)/2^k$.

That implies in particular
\begin{equation}\label{eq:2.rateconv}
\begin{aligned}
&  L_n(y, \omega, E) - L_{2n}(y, \omega, E)  =\\
& \frac{\ell}{n} (L_\ell(y, \omega, E) - L_{2\ell}(y, \omega, E)) +
O\Big(\exp\big(-\gamma_2 n^{1/2}\big)\Big)
\end{aligned}
\end{equation}

Let $\xi(\rho) = \log \rho$, $\xi(\rho, r, y) = \log |\rho +
re(y)|$, $R_1 < \rho < R_2$, $0 < r < r_2$, $0 \le y \le 1$, as in
Lemma~\ref{lem:2.jensenannul}.  Then, by Lemma~\ref{lem:2.lipsub}
\begin{equation}
\label{eq:lipsch} \bigl|L_{j\ell} \bigl(\xi(\rho, r, y), \omega,
E\bigr) - L_{j\ell}\bigl(\xi(\rho), \omega, E\bigr) \bigr| \le
CR_1^{-1} r \qquad  j = 1, 2
\end{equation}
Recall that for any $N > \exp(\gamma_1\,n^\sigma)$ one has
$$
\bigl|L_N(y, \omega, E) - 2L_{2n}(y, \omega, E) + L_n(y, \omega, E)
\bigr| < \exp\bigl(-\gamma_2 n^\sigma\bigr)\ ,
$$
see \cite{Gol Sch1}.
 Hence, due to \eqref{eq:2.jensencordom} and Lemma~\ref{lem:2.jensenannul}
\begin{align}
&\cM_N(\omega, E, R_1 + r_2, R_2 - r_2) \le {4|\cA_{R_1, R_2}|\over
r_2^2 N}\, \cJ\Bigl(\log \big |f_{N}(\cdot, \omega, E)\big |,
 \cA_{R_1, R_2}, r_2\Bigr) \nn \\
 & =  {4|\cA_{R_1, R_2}|\over
r_2^2 }\, \cJ\Bigl(n^{-1}[\log \big |f_{2n}(\cdot, \omega, E)\big
|-\log \big |f_{n}(\cdot, \omega, E)\big |],
 \cA_{R_1, R_2}, r_2\Bigr)
\label{eq:LS} \\
& \quad + O\big((R_2-R_1)r_2^{-2} \exp(-\gamma_2 n^{\sigma})
\big)\label{eq:error1}
\end{align}
Next, we rewrite the Jensen average in \eqref{eq:LS} using
Lemma~\ref{lem:2.jensenannul}
\begin{align}
& \cJ\Bigl(n^{-1}[\log \big |f_{2n}(\cdot, \omega, E)\big |-\log
\big |f_{n}(\cdot, \omega, E)\big |],
 \cA_{R_1, R_2}, r_2\Bigr) \nn \\
 & = 2\cJ\Bigl(\frac{1}{2n}\log \big |f_{2n}(\cdot, \omega, E)\big |-\frac{1}{n}\log
\big |f_{n}(\cdot, \omega, E)\big |,
 \cA_{R_1, R_2}, r_2\Bigr) \label{eq:error2} \\
 & + \cJ\Bigl(n^{-1}\log \big |f_{n}(\cdot, \omega, E)\big |,
 \cA_{R_1, R_2}, r_2\Bigr) \label{eq:main}
\end{align}
Inserting \eqref{eq:main} into \eqref{eq:LS} leads to the main term
on the right-hand side of~\eqref{eq:2.avernzel}. It is bounded above
by $\cM_n(\omega, E, R_1 - r_2, R_2 + r_2)$ in view
of~\eqref{eq:2.jensencordom}. It remains to bound the error
term~\eqref{eq:error2}. We introduce the short-hand notation
\begin{align*}
&\cS [L_n\bigl(\xi(\rho, r, y),\omega, E\bigr) - L_n
\bigl(\xi(\rho), \omega,
 E\bigr)] \\
 & = {4\pi\over (R_2^2 - R_1^2)r^2_2} \int^{R_2}_{R_1} \rho \, d\rho\
 \int^{r_2}_0 r dr \int^1_0 dy  \bigl[L_n\bigl(\xi(\rho, r, y),\omega, E\bigr) - L_n \bigl(\xi(\rho), \omega,
 E\bigr)\bigr]
 \end{align*}
Hence, the Jensen-average in \eqref{eq:error2} equals, see
\eqref{eq:2.rateconv},
\begin{align*}
& \cS[ L_{2n}\bigl(\xi(\rho, r, y),\omega, E\bigr) -
L_n\bigl(\xi(\rho, r, y),\omega, E\bigr)] - \cS [L_{2n}
\bigl(\xi(\rho), \omega, E\bigr) - L_n \bigl(\xi(\rho), \omega,
 E\bigr)] \\
 &= \frac{\ell}{n} \cS[ L_{2\ell}\bigl(\xi(\rho, r, y),\omega, E\bigr) -
L_\ell\bigl(\xi(\rho, r, y),\omega, E\bigr)] - \frac{\ell}{n}\cS
[L_{2\ell} \bigl(\xi(\rho), \omega, E\bigr) - L_\ell
\bigl(\xi(\rho), \omega,  E\bigr)] \\
& \quad + O\Big(\exp\big(-\gamma_2 n^{1/2}\big)\Big)
\end{align*}
By the Lipschitz bound \eqref{eq:lipsch}, we can further estimate
the absolute value here by
\begin{align*}
& \lesssim  \Big|\frac{\ell}{n} \cS[ L_{2\ell}\bigl(\xi(\rho, r,
y),\omega, E\bigr) - L_{2\ell} \bigl(\xi(\rho), \omega,
E\bigr)]\Big| + \Big| \frac{\ell}{n}\cS [L_\ell\bigl(\xi(\rho, r,
y),\omega, E\bigr) - L_\ell
\bigl(\xi(\rho), \omega,  E\bigr)] \Big|\\
& \quad + O\Big(\exp\big(-\gamma_2 n^{1/2}\big)\Big) \\
& \lesssim n^{-1/2} r_2 + O\Big(\exp\big(-\gamma_2 n^{1/2}\big)\Big)
\end{align*}
So the total error is the sum of this term times ${4|\cA_{R_1,
R_2}|\over r_2^2 }$ plus the error in~\eqref{eq:error1}. In view of
our assumptions on $r_2$ the lemma is proved.
\end{proof}

\section{Eliminating resonances via resultants}
\label{sec:resultant}

In this section we describe the mechanism behind the process of
eliminating ``bad'' frequencies~$\omega$, which is fundamental to
everything we do. ``Bad'' here refers to those $\omega$ which
produce too many too close resonances. More precisely, we will need
to ensure that the zeros of
\[
f_{\ell_1}(\cdot,\omega,E) \text{\ \ and\ \ } f_{\ell_2}(\cdot
e(t\omega),\omega,E)
\]
do not come too close. This requires the elimination of a  set of
energies of small measure and not too large complexity. The
elimination method is based on the natural idea that $\omega$ and
$t\omega$ act almost as independent variables (due to $\omega$ being
the slow
 and $t\omega$ being the fast variable). On a technical level this
 will be accomplished via a Cartan estimate on the {\em resultant} of two
 polynomials (which themselves come from the Weierstrass preparation
 theorem applied to $f_{\ell_1}$, $f_{\ell_2}$).
For those properties of the resultant which we need here, we refer
the reader to Appendix~A.
We would like to draw the reader's attention to the fact that all sets
removed in this section are very small in terms of Hausdorff dimension.
Indeed, the complexity of the bad sets is always less than their measure
raised to an arbitrarily small negative number (at least for $f_N$ with $N$ large).
This is one reason why we are able to eventually remove sets of Hausdorff dimension zero.
Another reason has to do with the second, an completely different, elimination method
used in this paper. It is designed to remove triple resonances and is based on the implicit
function theorem rather than resultants, see Section~\ref{sec:tripleelim}. It should be emphasized,
though, that Section~\ref{sec:tripleelim} must come after this section in the sense that the methods
there can only possibly work after we have removed the frequencies (as well as energies) specified
in this section. This is simply due to the fact that the implicit function theorem requires a
non-degeneracy condition which can be guaranteed only via the results of this section.

\begin{lemma}\label{lem:2.resultanal}
Let $f(z; w) = z^k + a_{k-1} (w) z^{k-1} + \cdots + a_0(w)$, $g(z;
w) = z^m + b_{m-1}(w) z^{m-1} + \cdots + b_0(w)$ be polynomials
whose coefficients $a_i(w)$, $b_j(w)$ are analytic functions defined
in a domain $G\subset \IC^d$.  Then $\Res(f(\cdot, w), g(\cdot, w))$
is analytic in $G$.
\end{lemma}

Our goal here is to separate the zeros of two analytic functions
using the resultants by means of shifts in the argument, see
Section~7 of \cite{Gol Sch2}, in particular Lemma~7.4. This can be
reduced to the same question for polynomials due to the Weierstrass
preparation theorem. Here is a simple observation regarding the
resultant of a polynomial and a shifted version of another
polynomial.

\begin{lemma}\label{lem:2.resultshiftlem}
 Let $f(z) = z^k + a_{k-1} z^{k-1}+ \cdots + a_0$, $g(z) = z^m +
b_{m-1} z^{m-1} + \cdots + b_0$ be polynomials.  Then
\begin{equation}
\label{eq:2.resultshift} \Res\bigl(f(\cdot + w), g(\cdot)\bigr) =
(-w)^n + c_{n-1} w^{n-1} + \cdots + c_0
\end{equation}
where $n = km$, and $c_0, c_1 \cdots$ are polynomials in the $a_i,
b_j$.
\end{lemma}
\begin{proof} Let $\zeta_j$, $1 \le j \le k$ (resp. $\eta_i, 1 \le i \le m$) be
the zeros of $f(\cdot)$ (resp. $g(\cdot)$).  The zeros of $f(\cdot +
w)$ are $\zeta_j - w$, $1 \le j \le k$.  Hence
\begin{equation}\label{eq:res_exp}
\Res\bigl(f(\cdot + w), g(\cdot)\bigr) = \prod_{i, j} (\zeta_j - w
- \eta_i)
\end{equation}
and \eqref{eq:2.resultshift} follows.
\end{proof}

The following lemma gives some information on the coefficients
in~\eqref{eq:2.resultshift}.

\begin{lemma}
\label{lem:bjuw} Let $P_s(z, \uw) = z^{k_s} + a_{s, k_s -1} (\uw)
z^{k_s -1} + \cdots + a_{s,0}(\uw)$, $z \in \IC$, where $a_{s,
j}(\uw)$ are analytic functions defined in some polydisk $\cP =
\prod\limits_i \cD(w_{i,0}, r)$, $\uw = (w_1, \dots, w_d) \in
\IC^d$, $\uw_0 = (w_{1,0}, \dots, w_{d, 0}) \in \IC^d$, $s = 1, 2$.
Set $\chi(\eta, \uw) = \Res\bigl(P_1(\cdot, \uw), P_2(\cdot + \eta,
\uw)\bigr)$, $\eta \in \IC$, $\uw \in \cP$.  Then
\begin{equation}
\label{eq:etauw} \chi(\eta, \uw) = (-\eta)^k + b_{k-1}(\uw)
\eta^{k-1} + \cdots + b_0(\uw)
\end{equation}
where $k = k_1 k_2$, $b_j(\uw)$ are analytic in $\cP$,  $j = 0, 1,
\dots, k-1$.  Moreover, if the zeros of $P_i(\cdot, \uw)$ belong to
the same disk $D(z_0, r_0)$, $i = 1, 2$ then for all $0 \le j \le
k-1$,
\begin{equation}
\label{eq:bjuw} \big | b_j(\uw)\big | \le \binom{k}{k-j}
(2r_0)^{k-j} \le(2r_0 k)^{k-j}
\end{equation}
\end{lemma}
\begin{proof} The relation~\eqref{eq:etauw}  with some coefficients $b_j(\uw)$ follows
from Lemma~\ref{lem:2.resultshiftlem} whereas the
bound~\eqref{eq:bjuw} follows from the expansion~\eqref{eq:res_exp}.
By Lemma~\ref{lem:2.resultshiftlem}, $\chi(\eta, \uw)$ is analytic
in $\IC \times \cP$.  Therefore $b_j(\uw) =
(j!)^{-1}(\partial_\eta)^j \chi(\eta,\uw) \Bigm |_{\eta=0}$ are
analytic $j = 0, 1, \dots, k-1$.
\end{proof}

The following lemma allows for the separation of the zeros of one
polynomial from those of a shifted version of another polynomial.
This will be the main mechanism for eliminating certain ``bad''
rotation numbers $\omega$. The logic is as follows: due to the basic
definition of the resultant,
$$
|\chi(\eta, \uw)| = \prod_{i, j} \bigl|\zeta_{i, 1}(\uw) - \zeta_{j,
2} (\eta, \uw)\bigr|
$$
where $\zeta_{i, 1}(\uw)$, $\zeta_{j, 2}(\eta, \uw)$ are the zeros
of $P_1(\cdot, \uw)$ and $P_2(\cdot + \eta, \uw)$, respectively.
Therefore, if $\big|\zeta_{i,1}(\uw) - \zeta_{j, 2}(\uw) \big | <
\exp(-kH)$ for one choice of~$i$, then $\big |\chi(\eta, \uw)\big |
< \exp(-kH)$. This simple fact allows one to separate the zeros
$\zeta_{i,1}(\uw)$ from the zeros $\zeta_{j,2}(\uw)$ provided $\uw$
falls outside of a set whose measure and complexity is controlled by
Cartan's estimate. More specifically, we obtain this separation by
means of a shift by~$t\omega_1$ in the $z$-slot:

\begin{lemma}
\label{lem:H2O} Let $P_s(z, \uw)$ be polynomials in $z$ as in
Lemma~\ref{lem:bjuw}, $s = 1, 2$. In particular, $\uw\in\cP$ where
$\cP$ is a polydisk of some given radius $r>0$.  Assume that $k_s >
0$, $s = 1, 2$ and set $k=k_1k_2$. Suppose that for any $\uw \in
\cP$ the zeros of $P_s(\cdot, \uw)$ belong to the same disk
$\cD(z_0, r_0)$, $r_0 \ll 1$, $s = 1, 2$.  Let $t > 16k\, r_0\,
r^{-1}$. Given $H \gg 1$ there exists a set \[\cB_{H,t}\subset
\widetilde\cP:=\cD(w_{1,0}, 8kr_0/t) \times \prod\limits^d_{j=2}
\cD(w_{j,0}, r/2)\]  such that $S_{\uw_0, (16kr_0t^{-1}, r,\dots,
r)} (\cB_{H,t}) \in \car_d(H^{1/d}, K)$, $K = CHk$ and for any $\uw
\in \widetilde\cP \setminus \cB_{H,t}$ one has
\begin{equation}
\label{eq:zero_dist} \dist\bigl(\bigl\{\mbox{\rm zeros of
$P_1(\cdot, \uw)$}\bigr\}, \bigl\{\mbox{\rm zeros of $P_2\bigl(\cdot
+ t(w_1 - w_{1,0}), \uw\bigr)\bigr\}$}\bigr) \ge e^{-CHk}
\end{equation}
\end{lemma}
\begin{proof} Define $\chi(\eta, \uw)$ as in Lemma~\ref{lem:bjuw}.  Note that for any
$\uw \in \cP$ one has
$$
|\chi(\eta, \uw)| \ge |\eta|^k \bigl[ 1 - \sum^\infty_{j=1}
\bigl({2r_0 k\over |\eta|}\bigr)^j\bigr] \ge {1\over 2} |\eta|^k
$$
provided $|\eta| \ge 8r_0\, k$.  Furthermore, for any $\uw\in\cP$,
$$
|\chi(\eta, \uw)| \le |\eta|^k \bigl[ 1 + \sum^\infty_{j=1}
\bigl({2r_0 k\over |\eta|}\bigr)^j\bigr] \le 2 |\eta|^k.
$$
provided $|\eta| \ge 8r_0\, k$. Hence, by the maximum principle,
$$
\sup \bigl\{\big | \chi(\eta, \uw) \big |: |\eta| \le 16r_0k\bigr\}
\le 2(16kr_0)^k\ .
$$
Set \[f(\uw) = \chi\bigl(t(w_1 - w_{1,0}), (w_1,w_2, \dots,
w_d)\bigr),\quad w_1 \in \cD(w_{1,0}, 16kr_0/t),\quad (w_2, \dots,
w_d) \in \prod\limits^d_{j=2} \cD(w_{j,0}, r). \] This function is
well-defined because $16kr_0/t<r$ by our lower bound on $t$. By the
preceding,
\[ \sup\limits_{\uw} |f(\uw)| \le 2(16kr_0)^k, \quad \big |f\bigl(w_{1,0}
+ 8kr_0/t, w_{2,0}, \dots, w_{d,0}\bigr)\big | > \frac12(8kr_0)^k.\]
We can therefore apply Lemma~\ref{lem:2.high_cart}  to
\[ \phi=f\circ S_{\uw_0,(16kr_0/t,r,\ldots,r)}^{-1} \text{\ \ with\ \ }
M = \log2 + k \log(16kr_0), \quad m=-\log 2 + k\log(8kr_0)\] on a
polydisk of unit size.   Thus, given $H\gg1$ there exists
$\cB^{(1)}_{H,t} \subset \cP$ such that
$$
S_{\uw_0,(16kr_0t^{-1}, r, \dots, r)} \bigl(\cB^{(1)}_{H,t}\bigr)
\in \car\nolimits_d \bigl(H^{1/d}, K\bigr), \quad K=CkH,
$$
and such that for any \[(w_1, \dots, w_d) \in \cD(w_{1,0}, 8kr_0/t)
\times \prod\limits^d_{j=2} \cD(w_{j, 0}, r/2) \setminus
\cB^{(1)}_{H,t}\] one has $|f(\uw)| > e^{-CHk}$.  Recall that due to
basic properties of the resultant
$$
|f(w)| = \prod_{i, j} |\zeta_{i, 1} (\uw) - \zeta_{j, 2} (\uw) |
$$
where $\zeta_{i, 1} (\uw)$, $\zeta_{j, 2} (\uw)$ are the zeros of
$P_1(\cdot, \uw)$,  and $P_2(\cdot + t(w_1 - w_{1,0}), \uw)$,
respectively. Since $r_0\ll1$, this implies \eqref{eq:zero_dist},
and we are done.
\end{proof}
Lemma~\ref{lem:H2O} of course applies to polynomials $P_s(z)$ that
do not depend on $\uw$ at all. This example is important, and
explains why quantities like $K$ have the stated form.

This method of elimination applies to the Dirichlet determinants
$f_{\ell_1}(\cdot, \omega, E)$ and $f_{\ell_2}(\cdot e(t\omega),
\omega, E)$. We now state a result in this direction. We shall use
the following notation
$$
\cZ(f, \Omega) = \bigl\{z \in \Omega: f(z) = 0\bigr\}
$$
and
$$
\cZ(f,z_0,r_0)=\cZ(f,\cD(z_0,r_0))
$$

\begin{prop}\label{prop:2.zerosepar}
Let $V$ be analytic on $\cA_{\rho_0}$ and real-valued on~$\tor$.
Assume that $L(\omega,E)>\gamma>0$ for all\footnote{One can localize
here to intervals of $\omega$ and $E$.} $\omega, E$. Fix
$a>1$, and $1\gg c>0$ as well as a bounded set $\cS\subset\IC$.
 There exists $\ell_0 = \ell_0(V,\rho_0, a,c, \gamma,\cS)$,
  such that for
any $\ell_1 \ge \ell_2 \ge \ell_0$ the following holds:  Given $t >
\exp\bigl(\bigl(\log \ell_1\bigr)^{C_0}\bigr)$, $H \ge 1$, there
exists a set $\Omega_{\ell_1, \ell_2, t, H} \subset \tor$, with
\begin{equation}\label{eq:OmNbad}
  \begin{aligned}
\mes\bigl(\Omega_{\ell_1, \ell_2, t, H}\bigr) & <  \exp\big((\log\ell_1)^{C_1}\big) e^{-\sqrt H}\\
\compl\bigl(\Omega_{\ell_1, \ell_2, t, H}\bigr) & < t
\exp\big((\log\ell_1)^{C_1}\big) H
\end{aligned}
\end{equation}
such that for any $\omega \in \tor_{c,a} \setminus \Omega_{\ell_1,
\ell_2, t, H}$ there exists a set $\cE_{\ell_1, \ell_2, t, H,
\omega}$ with
\begin{align*}
\mes\bigl(\cE_{\ell_1, \ell_2, t, H, \omega}\bigr) & <  t\exp\big((\log\ell_1)^{C_1}\big)\, e^{-\sqrt H}\\
\compl\bigl(\cE_{\ell_1, \ell_2, t, H, \omega}\bigr) & <  t
\exp\big( (\log\ell_1)^{C_1}\big) H
\end{align*}
such that for any $E \in \cS \setminus \cE_{\ell_1, \ell_2, t, H,
\omega}$ one has
\begin{align}\label{eq:z1z2sep} \dist\bigl(\cZ\bigl(f_{\ell_1}(\cdot, \omega,
E), \cA_{\rho_0}\bigr), \cZ\bigl(f_{\ell_2}(\cdot e(t\omega),\omega,
E),
 \cA_{\rho_0}\bigr)\bigr) &> e^{-H(\log \ell_1)^{C_2}} \\
\dist\bigl(\spec\bigl(H_{\ell_1} (e(x_0), \omega)\bigr)
\setminus\cE_{\ell_1, \ell_2, t, H, \omega}\;,\;
\spec\bigl(H_{\ell_2}(e(x_0 + t\omega),
\omega)\bigr)\bigr) &\ge e^{-H(\log \ell_1)^{3C_2}}.\label{eq:E1E2sep}
\end{align}
Here $C_j$ are also allowed to depend on $\cS$.
\end{prop}
\begin{proof} Fix some choice of $z_0\in\cA_{\rho_0}$,
$E_0\in\cS$, and  $\omega_0\in \tor_{c,a}$. By the Weierstrass
preparation theorem of the previous section we can write
\begin{align*}
f_{\ell_1}(e(z),\omega,E) &= P_1(z,\omega,E) g_1(z,\omega,E) \\
f_{\ell_2}(e(z+t\omega_0),\omega,E) &= P_2(z,\omega,E)
g_2(z,\omega,E)
\end{align*}
for all
\[(z,\omega,E) \in \cP_0:=\cD(z_0,r_0)\times
\cD(\omega_0,r)\times \cD(E_0,r)
\]
where\footnote{We remind the reader that $\asymp$ means
proportional. The constant of proportionality here is allowed to
depend on $V,\rho_0,a,c,\gamma,\cS$.} $r_0\asymp \ell_1^{-1}$,
$r\asymp \exp(-(\log \ell_1)^{C_0})$, and $g_j$ does not vanish
on~$\cP_0$. Moreover, each $P_j(\cdot,\omega,E)$ is a polynomial of
degree $k_j\les (\log \ell_j)^{C_0}$ for all $(\omega,E)\in
\cD(\omega_j,r)\times \cD(E_i,r)$ and all its zeros belong
to~$\cD(z_0,r_0)$.
 Apply Lemma~\ref{lem:H2O}, with $t>(\log
 \ell_1)^{C_1} > 16k_1k_2 r_0r^{-1}$, to the polynomials
 \[
P_1(\cdot,\omega,E),\qquad P_2(\cdot +t(\omega-\omega_0),\omega,E)
 \]
Thus, for any $H\ge1$ there exists $\cB_{H,t}\subset
\cD(\omega_0,8kr_0/t)\times \cD(E_0,r)$ with
\[
\Big\{ \big(t(\omega-\omega_0)/(16kr_0), (E-E_0)/r\big) \::\:
(\omega,E)\in \cB_{H,t} \Big\} \in \car_2(H^{1/2},K),\quad K=CHk,
\]
so that for any $(\omega,E)\in \cD(\omega_0,8kr_0/t)\times
\cD(E_0,r/2)\setminus \cB_{H,t}$ one has
\begin{equation}
  \dist\bigl(\bigl\{\mbox{\rm zeros of $P_1(\cdot,
\omega,E)$}\bigr\}, \bigl\{\mbox{\rm zeros of $P_2\bigl(\cdot +
t(w_1 - w_{1,0}), \omega,E\bigr)\bigr\}$}\bigr) \ge e^{-CHk}
\label{eq:P1P2_sep2}
\end{equation}
By definition of $P_1, P_2$, \eqref{eq:P1P2_sep2} this implies that
\[
\dist\bigl(\cZ\bigl(f_{\ell_1}(\cdot, \omega, E), z_0, r_0\bigr),
\cZ\bigl(f_{\ell_2}(\cdot e(t\omega),\omega, E),
 z_0, r_0\bigr)\bigr) > e^{-H(\log \ell_1)^{C_1}}
\]
Now let $z_0,\omega_0,E_0$ run over a net
\[\cN=\{(z_j,\omega_j,E_j)\}_{j=1}^J\subset \cA_{\rho_0}\times
\tor_{c,a}\times \cS \] so that each point in $\cA_{\rho_0}\times
\tor_{c,a}\times \cS $ comes $(r_0, kr_0/t, r)$--close to one of the
points in~$\cN$ and no two points in~$\cN$ are closer than this
distance. Denoting by $\cB_{H,t}(j)$ the bad set constructed above
for each point in~$\cN$, there exist \[\Omega_j\in \omega_j+
t^{-1}kr_0\car_1(\sqrt{H},K)\] so that for each \[z\in
\cD(\omega_j,kr_0/t)\setminus \Omega_j \] the $z$-slice $ \cE_{j,z}
:= \cB_{H,t}(j)|_z $ belongs to $E_j + r \car_1(\sqrt{H},K) $.  Now
define
\[
\Omega_{\ell_1,\ell_2,t,H} := \bigcup_j \Omega_j
\]
By construction, $ \Omega_{\ell_1,\ell_2,t,H}$
satisfies~\eqref{eq:OmNbad}. Moreover,  for each
$\omega\in\tor_{c,a}\setminus\Omega_{\ell_1,\ell_2,t,H}$ define
\[
\cE_{\ell_1,\ell_2,t,H,\omega} := \bigcup_j \cE_{j,z}
\]
Then
\[
\mes (\cE_{\ell_1,\ell_2,t,H,\omega} ) \le t\exp((\log\ell_1)^{C_1})
e^{-\sqrt{H}},\quad \compl(\cE_{\ell_1,\ell_2,t,H,\omega} )  \le
\exp( (\log\ell_1)^{C_1}) t H^2
\]
as desired. If \eqref{eq:z1z2sep} failed, then there would have to
be $z_1, z_2\in \cA_{\rho_0}$ and $\omega_0\in \tor_{c,a}\setminus
\Omega_{\ell_1,\ell_2,t,H}$ and $E_0\in \cS\setminus
\cE_{\ell_1,\ell_2,t,H,\omega}$ such that
\[
f_{\ell_1}(z_1,\omega_0,E_0)=
f_{\ell_2}(z_2e(t\omega_0),\omega_0,E_0)= 0,\quad |z_1-z_2|<
e^{-H(\log\ell_1)^{C_2}}
\]
Then $(z_1,\omega_0,E_0)\in \cD(z_j,r_0)\times
\cD(\omega_j,kr_0/t)\times \cD(E_0,r)$ for some choice of~$j$. By
construction, $\omega_0\in \cD(\omega_j,kr_0/t)\setminus \Omega_j$
and $E_0 \in \cD(E_0,r)\setminus \cE_{j,z}$ which implies that
\[
|z_1-z_2|\ge e^{-CHk},
\]
see \eqref{eq:P1P2_sep2}. This is a contradiction and we are done with
\eqref{eq:z1z2sep}. For \eqref{eq:E1E2sep}, assume that $f_{\ell_1}(z_1, \omega, E_1) = 0$, $f_{\ell_2}\bigl(z_1
e(t\omega), \omega, E_2 \bigr) = 0$ for arbitrary $z_1=e(x_0)$, and
\beeq
\label{eq:E1E2} |E_1 - E_2| < e^{-H(\log \ell_1)^{3C_2}}, \qquad E_1\in [-C,C]\setminus
\cE_{\ell_1,\ell_2,t,H,\omega},\;\omega\in\tor\setminus\Omega_{\ell_1,\ell_2,t,H}\,.
\eneq
Then, by Corollary~\ref{cor:2.derlocb},
\begin{align*}
 \big | f_{\ell_2}(z_1e(t\omega),  \omega, E_1)\big | &\lesssim
|E_1 - E_2| \exp\bigl(\ell_2 L(\omega, E_1) + (\log \ell_2)^B\bigr) <
\exp\bigl(\ell_2 L(\omega,E_1) - H(\log \ell_1)^{2C_2}\bigr)
\end{align*}
By our choice of $E_1$,  there exists
$z_2$ so that $|z_2 - z_1| < \exp(-100 H(\log \ell_1)^{C_2})$, for which
\[ f_{\ell_2}\bigl(z_2 e(t\omega), \omega, E_1\bigr) = 0,\]
see Corollary~\ref{cor:2.lexcepzero}. But this would contradict~\eqref{eq:z1z2sep}
and we are done.
\end{proof}

\section{Localized  eigenfunctions in finite volume}
\label{sec:anderson}

In this section we apply the results of the previous section to the
study of the eigenfunctions of the Hamiltonian restricted to
intervals on the integer lattice. More precisely, we shall obtain a
finite-volume version of Anderson localization (albeit, at the
expense of removing a small set of energies). This section
corresponds to Section~9 of~\cite{Gol Sch2}.

\begin{lemma}
\label{lem:3.Green}
 Let $\omega \in \tor_{c,a}$,  $E_0 \in \IR$, $L(\omega,E_0)>\gamma > 0$, and $N\ge
 N_0(V,\rho_0,a,c,\gamma,E_0)$. Furthermore, assume that
\begin{equation}
\label{eq:3.fN_unter}
\log \big | f_N(z_0, \omega,E_0) \big | > NL(\omega,E_0) - K/2
\end{equation}
for some $z_0 = e(x_0)$, $x_0 \in \tor$,  $K > (\log N)^{C_0}$. Then
\begin{align}
\label{eq:3.Gjk}
\big | \cG_{[1, N]} (z_0,\omega,E) (j, k)\big |  & \le \exp\bigl(-
{\gamma}(k - j) + K\bigr)\\
\big \| \cG_{[1, N]} (z_0, \omega,E) \big \| & \le \exp(K)
\label{eq:3.normG}
\end{align}
where $\cG_{[1, N]}(z_0, \omega,E_0) = \bigl(H(z_0, \omega) -
E_0\bigr)^{-1}$ is the Green function, $\gamma = L(\omega,E_0)$, $1
\le j \le k \le N$.
\end{lemma}
\begin{proof}
By Cramer's rule and the uniform upper bound of
Proposition~\ref{prop:2.unifb} as well as the rate of convergence
estimate~\eqref{eq:2.rateconv},
\begin{equation}\label{eq:3.4}
\begin{split}
\big | \cG_{[1, N]}(z_0, \omega,E) (j, k) \big | & = \big | f_{j-1}
(z_0, \omega, E_0)\big | \cdot \big | f_{N - k} \bigl(z_0e(k\omega),
\omega, E_0\bigr) \big | \cdot \big | f_N(z_0, \omega, E_0)\big |^{-1}\\
& \le \big | f_N(z_0,
\omega, E_0)\big |^{-1}  \exp \bigl(NL(\omega, E_0) - (k-j) L(\omega, E_0) + (\log N)^C\bigr)
\end{split}
\end{equation}
Therefore, \eqref{eq:3.Gjk} follows from condition
\eqref{eq:3.fN_unter}.  The estimate~\eqref{eq:3.normG} follows from
\eqref{eq:3.Gjk} via Hilbert-Schmidt norms.
\end{proof}

Any solution of the equation
\begin{equation}
\label{eq:3.hamilton}
-\psi(n+1) - \psi(n-1) + v(n)\psi(n) = E\psi(n)\ ,\quad n \in \IZ\ ,
\end{equation}
obeys the relation (``Poisson formula'')
\beeq \label{eq:3.poisson}
\psi(m) = \cG_{[a, b]} (E)(m, a-1)\psi(a-1) + \cG_{[a, b]} (E)(m,
b+1)\psi(b+1),\quad m \in [a, b]. \eneq where $\cG_{[a,b]} (E) =
\bigl(H_{[a,b]} -E\bigr)^{-1}$ is the Green function, $H_{[a,b]}$
is the linear operator defined by \eqref{eq:3.hamilton} for $n \in
[a, b]$ with zero boundary conditions. In particular, if $\psi$ is a
solution of equation \eqref{eq:3.hamilton},  which satisfies a zero
boundary condition at the left (right) edge, i.e.,
\begin{equation}
\nn
\psi(a-1) = 0 \quad \mbox{(resp. $\psi(b+1) = 0$)}\ ,
\end{equation}
then
\begin{align}
\psi(m) & =\cG_{[a,b]}(m, b+1)\psi(b+1) \nn \\
\bigl(\mbox{resp. $\psi(m)$} & = \cG_{[a,b]}(m,a-1) \psi(a-1)\; \bigr)
\nonumber
\end{align}
The following lemma states that after removal of certain rotation
numbers $\omega$ and energies $E$, but uniformly in $x\in\tor$, only
one choice of $n\in[1,N]$ can lead to a determinant
$f_\ell(x+n\omega,\omega,E)$ with $\ell \asymp (\log n)^C$ which is
not large. This relies on the elimination results, see (g) in
Section 2, and is of crucial importance for all our  work.

\begin{lemma}
\label{lem:3.waschno} Fix $a>1,c>0$ and assume that
$L(\omega,E)>\gamma>0$ for all\footnote{This can of course be localized to intervals of $\omega$ and $E$.} $\omega, E$. Given $N\ge
N_0(V,\rho_0,\gamma,a,c)$ large, there exist a constant
$B=B(V,\rho_0,\gamma,a,c)$ and a set $\Omega_N \subset \tor$ with
\[\mes(\Omega_N) < \exp\bigl(-(\log N)^{B}\bigr),\quad
\compl(\Omega_N) < N^2,\] such that for all $\omega \in \tor_{c,a}
\setminus \Omega_N$ there is a set\footnote{The sets
$\Omega_N,\cE_{N,\omega}$ also depend on $V,\rho_0,\gamma,a,c$ but we
omit these parameters from our notation.} $\cE_{N, \omega} \subset
\IR$,
\[ \mes(\cE_{N,\omega}) < \exp\big(-(\log N)^{B}\big),\quad
\compl(\cE_{N,\omega}) <  N^3,\] with the following property: For
any $x \in \tor$ and any $\omega \in \tor_{c,a}\setminus \Omega_N$,
$E \in \IR \setminus\cE_{N,\omega}$ either
\begin{equation}\label{eq:3.star}
\log \big | f_\ell\bigl(e(x+n\omega), \omega, E\bigr)\big | > \ell L(\omega,E) - \sqrt\ell
\end{equation}
for all \begin{equation}   \label{eq:ell_int} (\log N)^{20B}\le \ell
\le 4(\log N)^{20B}
\end{equation}
and all $1 \le n \le N$, or there exists $n_1 = n_1(x, \omega, E)
\in [1, N]$ such that~(\ref{eq:3.star}) holds for all $n \in [1, N]
\setminus [n_1 - Q, n_1 + Q]$,  with
\[Q=\big\lfloor\exp \bigl((\log\log N)^{B_1}\bigr)\big\rfloor
,\quad B_1=B_1(V,\rho_0,\gamma,a,c)\] but not for $n = n_1$. Moreover,
in this case
\begin{equation}\label{eq:3.f1n}
\big | f_{[1,n]} \bigl(e(x), \omega, E\bigr) \big | > \exp
\bigl(nL(\omega, E) - (\log N)^{100B}\bigr)
\end{equation}
for each $1 \le n \le n_1 - Q$ and
\begin{equation}
\label{eq:3.fnN} \big | f_{[n, N]} \bigl(e(x), \omega, E\bigr)\big |
> \exp \bigl((N-n)L(\omega, E) - (\log N)^{100B}\bigr)
\end{equation}
for each $n_1 + Q\le n \le N$.
\end{lemma}
\begin{proof} Let  $B\gamma\ge 2$ (below, we will need to make $B$ large depending on~$a$ as well).
With $a>1$ and $ c>0$ fixed, we let $\Omega_{\ell_1, \ell_2, t, H}$
be as in Proposition~\ref{prop:2.zerosepar} and define
\[\Omega_N := \bigcup\, \Omega_{\ell_1, \ell_2, t, H}\] where
 $H =(\log N)^{3B}$ is fixed, and the union runs over $\ell_1, \ell_2$ as
in~\eqref{eq:ell_int}, and $N > t > \exp \bigl((\log\log
N)^{2C_0}\bigr)$ where $C_0$ is from
Proposition~\ref{prop:2.zerosepar} (thus, take $B_1=2C_0$). For any
$\omega \in \tor_{c,a} \setminus \Omega_N$ define
\[ \cE_{N,\omega} := \bigcup\, \cE_{\ell_1, \ell_2, t, H, \omega}\]
where $\cE_{\ell_1, \ell_2, t, H, \omega}$ is from the proposition
and the union is the same as before.  Now fix
$\omega\in\tor_{c,a}\setminus \Omega_N$,
$E\in\IR\setminus\cE_{N,\omega}$ and suppose \eqref{eq:3.star} fails
somewhere, i.e.,
$$
\log \big | f_{\ell_1}\bigl(e(x + n_1 \omega), \omega, E\bigr)\big | < \ell_1 L(\omega,E) - \sqrt{\ell_1}
$$
for some $1 \le n_1 \le N$ and $\ell_1$ as in~\eqref{eq:ell_int}. By
Corollary~\ref{cor:2.lexcepzero} there exists $z_1$ with $|z_1 - e(x
+ n_1 \omega) | < e^{-\ell_1^{{1\over 4}}}$ and
$$
f_{\ell_1} (z_1, \omega, E) = 0
$$
If
$$
\log \big | f_{\ell_2}\bigl(e(x + n_2 \omega), \omega, E\bigr)\big | < \ell_2 L(\omega,E) -
\sqrt{\ell_2}
$$
for some $\ell_2$ as in~\eqref{eq:ell_int} and $|n_2 - n_1| > \exp
\bigl((\log\log N)^{B_1}\bigr)$, then for some $z_2$, and $t := n_1
- n_2$
$$
f_{\ell_2} \bigl(z_2 e(t\omega), \omega, E\bigr) = 0
$$
with \[ |z_1 - z_2| < e^{-(\log N)^{4B}}\] However, by our choice of
$(\omega, E)$ \[ |z_1-z_2|>\exp\big(-CH(\log \log N)^{C_1}\big) =
\exp\big(- C(\log N)^{3B} (\log \log N)^{C_1}\big)
\]
which is a contradiction for $N\ge N_0$ large.
Thus~\eqref{eq:3.star} holds for all $\ell$ as
in~\eqref{eq:ell_int},  and any $1 \le n \le N$ such that $|n - n_1|
> \exp \bigl((\log\log N)^{B_1}\bigr)$.
This property allows one to apply the avalanche principle
Proposition~\ref{prop:AP} to the determinants appearing
in~\eqref{eq:3.f1n} and~\eqref{eq:3.fnN}. It will suffice to
consider the former with $n\ge (\log N)^{20B}$: in view
of~\eqref{eq:3.star} the conditions of Proposition~\ref{prop:AP}
hold if we choose the factor matrices $A_j$ there to be of length as
in~\eqref{eq:ell_int}.  It yields that
$$
\log \big | f_{[1, n]} \bigl(e(x), \omega, E\bigr)\big | \ge
nL(\omega, E) - C{n\over (\log N)^{5B}} > 0
$$
provided $N_0$ is large.  In fact, we can vary $x$ here:  note that
by Corollary~\ref{cor:2.lipnorm}, if \eqref{eq:3.star} holds at $x$,
then also for all $z \in \cD\bigl(e(x), e^{-\ell}\bigr)$,
$\ell=(\log N)^{20B}$. Repeating the avalanche principle expansion
for those~$z$ yields
\begin{equation}
\label{eq:3.fN_notzer}
f_{[1, n]} (z, \omega, E) \ne 0
\end{equation}
Now suppose
$$
\log \big | f_{[1, n]} \bigl(e(x), \omega, E\bigr) \big | \le
nL(\omega, E) - (\log N)^{100B}
$$
By Corollary~\ref{cor:2.lexcepzero},
$$
f_{[1, n]} (z, \omega, E) = 0
$$
for some $|z - e(x)| < \exp \bigl(-(\log N)^{50B}\bigr)$ provided
$B$ is sufficiently large (depending on~$a$).  This contradicts
\eqref{eq:3.fN_notzer} and we are done.
\end{proof}

\begin{remark}
\label{rem:3.expstab}
It follows from Corollary~\ref{cor:2.lipnorm} that \eqref{eq:3.star} is stable under perturbations of $E$ by an
amount $< e^{-C\ell}$.  More precisely, if \eqref{eq:3.star} holds for $E$, then
\begin{equation}
\nn
\log \big | f_\ell \bigl(e(x + n\omega), \omega, E'\bigr) \big | > \ell L(E', \omega) - 2 \sqrt{\ell}
\end{equation}
for any $E'$ with $|E' - E| <  e^{-C\ell}$.  Inspection of the
previous proof now shows that \eqref{eq:3.f1n} and \eqref{eq:3.fnN}
are therefore also stable under such perturbations.
\end{remark}

The previous lemma yields the following finite volume version of Anderson localization.

\begin{lemma}
\label{lem:3.4} Fix $a>1, c>0$, assume that $L(\omega,E)>\gamma>0$
for all $\omega, E$, and let $\Omega_N$ and $\cE_{N, \omega}$ be as
in the previous lemma with $N\ge N_0$. For any $x , \omega \in
\tor$, let $\bigl\{E_j^{(N)} (x, \omega)\bigr\}_{j=1}^N$ and
$\bigl\{\psi_j^{(N)} (x, \omega, \cdot)\bigr\}_{j=1}^N$ denote the
eigenvalues and normalized eigenvectors of $H_{[1, N]}(x, \omega)$,
respectively. If $\omega \in \tor_{c,a} \setminus \Omega_N$ and for
some $j$, $E_j^{(N)}(x, \omega) \notin \cE_{N, \omega}$, then there
exists a point $\nu_j^{(N)} (x, \omega) \in [1, N]$ (which we call
the center of localization) so that for any $\exp \bigl((\log\log
N)^{B_1}\bigr) \le Q \le N$ and with $\Lambda_Q := [1, N] \cap
\bigl[\nu_j^{(N)}(x, \omega) - 2Q, \nu_j^{(N)} (x, \omega) +
2Q\bigr)$ one has
\begin{enumerate}
\item[{\rm(i)}] $\dist\bigl(E_j^{(N)}(x, \omega), \spec\bigl(H_{\Lambda_Q}(x, \omega)\bigr)\bigr) <
e^{-\gamma Q/4}$

\item[{\rm(ii)}] $\sum\limits_{k \in [1, N] \setminus \Lambda_Q} \big | \psi_j^{(N)} (x, \omega; k)
\big |^2 < e^{-Q\gamma}$, where $\gamma > 0$ is a lower bound for
the Lyapunov exponents.
\end{enumerate}
Here $B,B_1$ are the constants from the previous lemma.
\end{lemma}
\begin{proof} Fix $N\ge N_0$, $\omega \in \tor_{c,a} \setminus \Omega_N$ and $E_j^{(N)} (x, \omega) \notin
\cE_{N, \omega}$.  Let $n_1 = \nu_j^{(N)} (x, \omega)$ be such that
$$
\big | \psi_j^{(N)} (x, \omega; n_1)\big | = \max_{1 \le n \le N}
\big | \psi_j^{(N)} (x, \omega; n)\big |
$$
Fix an  $\ell$ as in~\eqref{eq:ell_int} and suppose that, with $E =
E_j^{(N)} (x, \omega)$, and $\Lambda_0 := [1, N] \cap [n_1 -\ell,
n_1 + \ell]$,
\begin{equation}
\label{eq:3.11}
\log \big | f_{\Lambda_0} (x, \omega, E) \big | > |\Lambda_0| L(\omega,E) - \sqrt{\ell}
\end{equation}
By Cramer's rule, see~\eqref{eq:cramerG},  this would then imply
that
$$
\big | G_{\Lambda_0} (x, \omega, E) (k, j) \big | < \exp \bigl(-
\gamma  \big |k - j\big | + C\sqrt{\ell}\bigr)
$$
for all $k, j \in \Lambda_0$.  But this contradicts the maximality
of $\big | \psi_j^{(N)} (x, \omega; n_1) \big |$ due to
\eqref{eq:3.poisson} and $\ell$ being large.  Hence \eqref{eq:3.11}
above fails, and we conclude from Lemma~\ref{lem:3.waschno} that
$$
\log \big | f_{\Lambda_1} (x, \omega, E) \big | > |\Lambda_1| L(\omega,E) - \sqrt{\ell}
$$
for every $\Lambda_1 = [k - \ell, k + \ell] \cap [1, N]$ provided
$|k - n_1 | > \exp \bigl((\log\log N)^{B_1}\bigr)$.  Since
\eqref{eq:3.11} fails, we conclude that $f_{\Lambda_0} (z_0, \omega,
E) = 0$ for some $z_0$ with $|z_0 - e(x)| < e^{-\ell^{1/4}}$.  By
self-adjointness of $H_{\Lambda_0} (x, \omega, E)$ we obtain
$$
\dist\bigl(E, \spec\bigl(H_{\Lambda_0} (x, \omega)\bigr)\bigr) <
Ce^{-\ell^{{1/4}}}\ ,
$$
as claimed (the same arguments applies to the larger intervals
$\Lambda_Q$ around $n_0$). From \eqref{eq:3.f1n} of the previous
lemma with $n = n_1 - Q$ (if $n_1 - Q < (\log N)^{2C_0}$, then
proceed to the next case) one concludes from Cramer's rule,
see~\eqref{eq:cramerG}, that
\begin{equation}
\label{eq:3.12} \big | G_{[1, n_1 -  Q]} (x, \omega, E) (k, m) \big
| < \exp \bigl(-\gamma |k-m| + (\log N)^{C_0}\bigr)
\end{equation}
for all $1 \le k, m \le n_1 -  Q$.  In particular,
$$
\big | \psi_j^{(N)} (x, \omega; k) \big | < e^{-\frac{\gamma}{2} |
n_1 -
 Q - k|}
$$
for all $1 \le k \le n_1 - 2Q$. The same reasoning applies to
$$
G_{[n_1 + Q, N]} (x, \omega, E)
$$
via \eqref{eq:3.fnN} of the previous lemma, and (ii) follows. For
(i), note that \eqref{eq:3.poisson} and~(ii) imply that
\[
\| (H_{\Lambda_Q}(x,\omega)- E_j^{(N)}(x,\omega)) \psi_j^{(N)}\|\le
e^{-\gamma Q/2}
\]
Since $\|\psi_j^{(N)}\|_{\ell^2(\Lambda_Q)} \ge 1 - e^{\gamma Q}$,
we obtain~(i).
\end{proof}

Since eigenvalues of the Dirichlet problem are simple, we can order
the $E_j^{(N)}(x,\omega)$ according to the convention
\[
E_1^{(N)}(x,\omega)<E_2^{(N)}(x,\omega)<\ldots< E_N^{(N)}(x,\omega)
\]
 The following
corollary deals with the stability of the localization statement of
Lemma~\ref{lem:3.4} with respect to the energy. As in previous
stability results of this type in this paper,  the most important
issue is the relatively large size of the perturbation, i.e.,
$\exp(-(\log N)^C)$ instead of $e^{-N}$, say.

\begin{corollary}
\label{cor:3.Estab} Using the assumptions and terminology of the
previous proposition, let $\omega\in \tor_{c,a} \setminus \Omega_N$,
$E_j^{(N)}(x, \omega) \notin \cE_{N,\omega}$,  and
$\nu_j^{(N)}(x,\omega)$ be associated with~$\psi_j^{(N)}(x, \omega;
\cdot)$ as stated there. If $|E - E_j^{(N)}(x, \omega)| < e^{-(\log
N)^{40B}}$ with $B$ as above, then
\begin{equation}
\label{eq:3.13} \sum^{\nu_j^{(N)}(x,\omega) - Q}_{n=1} \big | f_{[1,
n]} \bigl(e(x), \omega, E\bigr) \big |^2 < e^{-\gamma Q/2} \sum_{n
\in \Lambda_Q} \big | f_{[1, n]}\bigl(e(x), \omega, E\bigr)\big |^2
\end{equation}
where $\Lambda_Q = \bigl[\nu_j^{(N)} (x, \omega) - Q, \nu_j^{(N)}(x, \omega) + Q\bigr] \cap [1,
N]$.  Similarly,
\begin{equation}
\label{eq:3.14} \sum^N_{n = \nu_j^{(N)}(x, \omega)+ Q} \big | f_{[n,
N]}(x, \omega, E)\big |^2 < e^{-\gamma Q/2} \sum_{n \in \Lambda_{Q}}
\big | f_{[n, N]}(x, \omega, E) \big |^2
\end{equation}
Finally, under the same assumptions one has
\begin{equation}
\label{eq:3.15}\begin{split}   &\big | f_{[1, n]} \bigl(e(x),
\omega, E\bigr) - f_{[1, n]}\bigl(e(x), \omega, E_j^{(N)} (x,
\omega)\bigr) \big | \\
&\le \exp\bigl((\log N)^C\bigr) \big | E - E_j^{(N)}\bigl(x,
\omega)\big |\, \big |f_{[1, n]}(e(x), \omega, E_j^{(N)} (x,
\omega)\bigr)\big |
\end{split}
\end{equation}
provided $1 \le n \le \nu_j^{(N)} (x, \omega) - Q$, and similarly for $f_{[n, N]}$.
\end{corollary}
\begin{proof} For each $j$ there exists a constant $\mu_j(x, \omega)$ so that
$$
\psi_j^{(N)}(x, \omega; n) = \mu_j(x, \omega) f_{[1, n-1]}\bigl(x, \omega; E_j^{(N)} (x,
\omega)\bigr)
$$
for all $1 \le n \le N$ (with the convention that $f_{[1, 0]} = 1$).  A
similar formula holds for
\[ f_{[n+1, N]}\bigl(e(x), \omega, E_j^{(N)} (x,
\omega)\bigr).\]
As in the previous proof, this implies
estimate~\eqref{eq:3.12} with $E = E_j^{(N)} (x, \omega)$.  Thus, for all $1 \le n \le
\nu_j^{(N)}(x, \omega) - Q$ one has 
$$
\big | f_{[1, n]}\bigl(e(x), \omega, E_j^{(N)} (x,
\omega)\bigr)\big | < e^{-\gamma|\nu_j^{(N)} (x, \omega) - n|/2}\
\big | f_{[1, \nu_j^{(N)}(x, \omega)]} \bigl(e(x), \omega, E_j^{(N)}
(x, \omega)\bigr)\big |\ ,
$$
which implies \eqref{eq:3.13} for $E = E_j^{(N)}(x, \omega)$, and \eqref{eq:3.14} follows by a similar argument for this $E$.
Corollary~\ref{cor:2.lipnorm} implies that
\begin{align*}
& \big | f_{[1, n]}\bigl(e(x), \omega, E\bigr) - f_{[1, n]}\bigl(e(x), \omega,
E_j^{(N)}(x, \omega)\bigr)\big | \\
& \le \exp \bigl((\log N)^C\bigr) \big | E - E_j^{(N)} (x, \omega)\big
|\ \big | f_{[1, n]} \bigl(e(x), \omega, E_j^{(N)}(x, \omega)\bigr)\big |
\end{align*}
for all $1 \le n \le \nu_j^{(N)}(x, \omega) - Q$,
and \eqref{eq:3.15} follows for all $\big | E - E_j^{(N)}(x, \omega)\big | <
\exp \bigl(-(\log N)^B\bigr)$.
\end{proof}

\section{Quantitative separation of the Dirichlet eigenvalues in finite volume} \label{sec:separation}

At least conceptually, this section provides arguably the most important
single ingredient in the proof of gap formation. It already played a crucial
role in our earlier work~\cite{Gol Sch2}.
Based on the finite volume Anderson localization from the previous section, we
shall now obtain a quantitative separation property of the eigenvalues on finite
volume. Note carefully that {\em localization} does not depend on the off-diagonal terms
in the Hamiltonian - indeed, a diagonal Hamiltonian has $\delta$-function eigenstates
which are perfectly localized. In contrast to this, the separation or even the
simplicity of the eigenvalues of course crucially depend on these off-diagonal terms.
Needless to say, the gaps in the spectrum also hinge on this property and this
section is one of the places where it enters in an essential way. The reader will
easily see this in the proofs of the first two results of this section. The off-diagonal
terms enter there simply through the mechanism of transfer matrices; or in other words, we
have to exploit that we are dealing with a second order difference equation which we can
solve via initial conditions.

We now turn to the details.
In this section it will be convenient for us to work with the
 operators $H_{[-N, N]}(x, \omega)$ instead of
$H_{[1, N]}(x, \omega)$. Abusing our notation somewhat, we use the
symbols $E_j^{(N)}, \psi_j^{(N)}$ to denote the eigenvalues and
normalized eigenfunctions of $H_{[-N,N]}(x, \omega)$, rather than
the eigenvalues and normalized eigenfunctions of $H_{[1, N]}(x,
\omega)$, as in the previous section. A similar comment applies to
$\Omega_N, \cE_{N, \omega}.$

The following proposition states that the eigenvalues $\{
E_j^{(N)}(x,\omega)\}_{j=1}^{2N+1}$  are separated from each other
by at least $e^{-N^\delta}$  provided $\omega\not\in\Omega_N$ and
provided we delete those eigenvalues that fall into a bad set
$\cE_{N,\omega}$ of energies. We remind the reader that
\[ \mes(\cE_{N,\omega})\les \exp(-(\log N)^{B}),\quad \compl(\cE_{N,\omega})\les N^3, \]
 and similarly for $\Omega_N$, see Lemma~\ref{lem:3.waschno}. This
section corresponds to Section~11 of~\cite{Gol Sch2}.

\begin{prop}
\label{prop:4.Ej_sep} Fix $a>1$, $c>0$ and assume that
$L(\omega,E)>\gamma>0$ for all $\omega, E$. Let $\Omega_N$,
$\cE_{\omega,N}$ be as in~Lemma~\ref{lem:3.waschno}. Furthermore,
fix $\delta\in (0,1)$. Then there exists
$N_0=N_0(\delta,V,\rho_0,a,c,\gamma)$ so that for any $N \ge N_0$, any
 $\omega \in\tor_{c,a}\setminus \Omega_N$ and all $x$ one has
\begin{equation}
\label{eq:4.Ej_sep'} \big | E_j^{(N)} (x, \omega) - E_k^{(N)} (x,
\omega) \big | > e^{-N^\delta}
\end{equation} for all
$j, k$ provided $E_j^{(N)}(x, \omega) \notin \cE_{N,\omega}$.
\end{prop}
\begin{proof} Fix $x \in \tor, E_j^{(N)} (x, \omega) \notin
\cE_{N,\omega}$.  Let $Q\asymp \exp \bigl((\log\log
N)^{B_1}\bigr)$, see Lemma~\ref{lem:3.waschno}. By
Lemma~\ref{lem:3.4} there exists
$$
\Lambda_Q := \bigl[\nu_j^{(N)} (x, \omega) - Q, \nu_j^{(N)} (x, \omega)
+ Q\bigr] \cap [-N, N]
$$
so that
\begin{equation}
\label{eq:4.2}
\begin{split}
& \sum_{n \in [-N, N]\setminus \Lambda_Q} \big | f_{[-N, n]} \bigl(e(x),
\omega; E_j^{(N)}(x, \omega)\bigr)\big |^2 \\
& < e^{-Q \gamma} \sum^N_{n=-N} \big |f_{[-N, n]} \bigl(e(x),
\omega; E_j^{(N)}(x, \omega)\bigr)\big|^2\ .
\end{split}
\end{equation}
Here we used that with some $\mu = \const$
$$
\psi_j^{(N)} (x, \omega; n) = \mu \cdot f_{[-N, n-1]}\bigl(e(x), \omega;
E_j^{(N)}(x, \omega)\bigr)
$$
for $ -N \le n \le N$. Note the convention that
\[ f_{[-N,-N-1]}=0, \quad f_{[-N,-N]}=1.\]
One can assume $\nu_j^{(N)}(x, \omega)
\ge 0$ by symmetry.
Using Corollary~\ref{cor:2.lipnorm} and \eqref{eq:3.15}, we conclude that
\begin{equation}
\label{eq:4.3}
\begin{split}
& \sum_{n= -N}^{\nu_j^{(N)}(x, \omega) - Q} \big | f_{[-N, n]} (e(x),
\omega, E) - f_{[-N, n]}\bigl(e(x), \omega, E_j^{(N)} (x, \omega)\bigr)\big
|^2 \\
& \le e^{-\gamma Q} \big |E - E_j^{(N)}(x, \omega) \big |^2 e^{(\log
N)^C} \sum_{n \in \Lambda_Q} \big |f_{[-N, n]} \bigl(e(x), \omega,
E_j^{(N)} (x, \omega)\bigr)\big |^2
\end{split}
\end{equation}
Let $n_1 = \nu_j^{(N)} (x, \omega) - Q -1$.  Furthermore,
\begin{equation}
\label{eq:4.4}
\begin{split}
 &\bigg \| \binom{f_{[-N, n+1]}(e(x), \omega, E)}{f_{[-N, n]}(e(x), \omega, E)}
- \binom{f_{[-N, n+1]} \bigl(e(x), \omega, E_j^{(N)} (x,
  \omega)\bigr)}{f_{[-N, n]}\bigl(e(x), \omega, E_j^{(N)}(x, \omega)\bigr)}
  \bigg \|\\
  &= \bigg \| M_{[n_1+1, n]} (e(x), \omega, E) \binom{f_{[-N, n_1
  +1]}(e(x),
  \omega, E)}{f_{[-N, n_1]}(e(x), \omega, E)} \\
  & \qquad\qquad- M_{[n_1 +1,
  n]}\bigl(e(x),
  \omega, E_j^{(N)}(x, \omega)\bigr) \binom{f_{[-N, n_1+1]}\bigl(e(x),
  \omega, E_j^{(N)}\bigr)}{f_{[-N, n_1]}\bigl(e(x), \omega,
  E_j^{(N)}\bigr)}\bigg \|\\
& \le e^{C(n - n_1)} e^{-\gamma Q/2} \big |E - E_j^{(N)} (x, \omega)
\big | e^{(\log N)^C} \Bigl(\sum_{n \in \Lambda_Q} \big | f_{[-N,
n]} \bigl(e(x), \omega, E_j^{(N)}(x, \omega)\bigr) \big
|^2\Bigr)^{{1\over
2}}\ .\\
  \end{split}
  \end{equation}
Now suppose there is $E_k^{(N)}(x, \omega)$ with $\big
|E_k^{(N)}(e(x), \omega) - E_j^{(N)} (x, \omega) \big | <
e^{-N^\delta}$ where the small $\delta > 0$ is arbitrary but fixed.
Then (\ref{eq:4.3}), (\ref{eq:4.4}) imply that
\begin{equation}
\label{eq:4.5}
\begin{split}
& \sum^{\nu_j^{(N)}(x, \omega) + Q}_{n= -N} \big | f_{[-N, n]}
\bigl(e(x),
\omega, E_j^{(N)} (x, \omega)\bigr) - f_{[-N, n]}\bigl(e(x), \omega,
E_k^{(N)} (x, \omega)\bigr)\big |^2\\
& < e^{-{1\over 2} N^\delta} \sum_{n \in \Lambda_Q} \big |f_{[-N, n]}
\bigl(e(x), \omega, E_j^{(N)}(x, \omega)\bigr)\big |^2\ ,
\end{split}
\end{equation}
provided $N^\delta > \exp\bigl((\log\log N)^{B_1}\bigr)$. Let us
estimate the contributions of $\bigl[\nu_j^{(N)} (x, \omega) + Q,
N\bigr]$ to the sum terms in the left-hand side of (\ref{eq:4.5}).

For both $E = E_j^{(N)}$ and $E_k^{(N)}$ one has
$$
f_{[-N, n]}(e(x), \omega, E) = G_{[\nu_j^{(N)} (x, \omega) + {Q\over 2},
N]} (e(x), \omega, E) \bigl(n ,\nu_j^{(N)}(x, \omega) + {Q\over 2}\bigr)
f_{[-N, \nu_j^{(N)}(x, \omega) + {Q\over 2} -1]}(e(x), \omega, E)
$$
due to the zero boundary condition at $N +1$, i.e.,
$$
f_{[-N, N]} \bigl(e(x), \omega, E_j^{(N)} (x, \omega)\bigr) = f_{[-N, N]}
\bigl(e(x), \omega, E_k^{(N)} (x, \omega) \bigr) = 0\ .
$$
Therefore, 
\begin{equation}
\label{eq:4.6}
\sum^N_{n = \nu_j^{(N)} + Q} \big | f_{[-N, n]}(e(x), \omega, E)\big |^2
\le e^{-{\gamma Q\over 4}} \sum_{k \in \Lambda_Q}\big |f_{[-N, k]}(e(x),
\omega, E)\big |^2
\end{equation}
again for both $E = E_j^{(N)}(x, \omega)$ and $E = E_k^{(N)}(x,
\omega)$.  Finally, in view of (\ref{eq:4.5}) and (\ref{eq:4.6}),
\begin{equation}
\label{eq:4.7}
\begin{split}
 &\sum^N_{n=-N} \big | f_{[-N, n]}\bigl( e(x), \omega, E_j^{(N)} (x,
\omega)\bigr) - f_{[-N, n]}\bigl(e(x), \omega, E_k^{(N)} (x,
\omega)\bigr)\big |^2\\
& < e^{-{\gamma Q\over 4}} \Bigl[\sum_{n \in \Lambda_Q} \big |
f_{[-N, n]} \bigl(e(x), \omega, E_j^{(N)} (x, \omega)\bigr)\big |^2\\
&\quad +
\sum_{n \in \Lambda_Q} \big | f_{[-N, n]}\bigl(e(x), \omega, E_k^{(N)} (x,
\omega)\bigr) \big |^2 \Bigr ]
\end{split}
\end{equation}
By orthogonality of $\bigl\{f_{[-N, n]}\bigl(e(x), \omega, E_j^{(N)}
(x, \omega)\bigr)\bigr\}^N_{n=-N}$ and $\bigl\{ f_{[-N, n]}
\bigl(e(x), \omega, E_k^{(N)} (x, \omega)\bigr)\bigr\}^N_{n=-N}$ we
obtain a contradiction from (\ref{eq:4.7}).
\end{proof}

In order to capture the mechanism behind Figure~2 later in this paper, we
will need to achieve the separation property of the previous proposition without
removing any energies. Rather, we will be using some apriori information which
has the same effect as requiring the energies to be ``good''. This is done
in the following result which is a corollary of the preceding proof rather than of the
statement itself.

\begin{cor}
   \label{cor:sep}  Assume $L(\omega,E)>\gamma>0$ for all $\omega,E$ and let $(x_0,\omega)\in\tor\times \tor_{c,a}$ be arbitrary
   where $a>1$ and $c>0$ are fixed. Moreover, fix constants $1
   \ge \delta>
   \eps>0$ and let $N\ge N_0(\delta,\eps,V,a,c,\gamma)$ be
sufficiently large.
   Suppose there is $\Lambda=[N',N'']\subset [-N,N]$ satisfying \[(\log N)^{2C_0}\le |\Lambda|\le
   N^\ve, \qquad
   100(\log N)^{2C_0}<N'<N''<N-100(\log N)^{2C_0},\] and
such that for some pair $E_1, E_2\in\IR$
 \[f_k(ze(s\omega),\omega,E)\ne0
 \quad\forall\, z\in\cD(e(x_0),r_0),\;\forall\, E\in\cD(E_1,r_0)\cup \cD(E_2,r_0),\;\forall\, s\in[-N,N]\setminus\Lambda\]
 and  all choices of $(\log N)^{C_0} \le k \le 100(\log N)^{C_0} $  with $r_0:=\exp(-(\log N)^{C_0/2})$.
If $E_1, E_2$ are eigenvalues of $H_{[-N,N]}(x_0,\omega)$, then
\[ |E_1-E_2|>e^{-N^\delta}\]
Furthermore, suppose $\psi_{j}(\cdot)$ are normalized eigenfunctions
of the Dirichlet problem on~$[-N,N]$
\[ H_{[-N,N]}(x_0,\omega)\psi_j = E_j\, \psi_j \]
Then \begin{equation}   \label{eq:psij_decay} \big |\psi_j(n) \big |
\le \exp\bigl(-\gamma \dist(n,\Lambda) /2\bigr),\quad j=1,2
\end{equation}
for all $n\in [-N,N]$.
\end{cor}
\begin{proof}
Let $\bar N':= N'+10(\log N)^{C_0}$ and $\bar N'' := N''- 10(\log
N)^{C_0}$. We apply the avalanche principle
Proposition~\ref{prop:AP} to $f_{[-N,\bar N']}(x_0,\omega,E)$ and
$f_{[\bar N'',N]}(x_0,\omega,E)$ with arbitrary
$E\in\cD(E_1,r_0)\cup \cD(E_2,r_0)$. It will suffice to consider the
former. If for some choice of $k$ as above and $-N\le s\le N'$
\[ \log|f_k(x_0 e(s\omega),\omega,E)|< kL(\omega,E)-{k^{\frac34}}\]
then by Corollary~\ref{cor:2.lexcepzero} one has $f_k(z_0
e(s\omega),\omega,E)=0$ for some $z_0\in \cD(e(x_0), e^{-\sqrt{k}})$
contradicting our hypothesis. Hence, Proposition~\ref{prop:AP}
implies that
\[
f_{[-N,\bar N']}(x_0,\omega,E) \ne 0\quad \forall \,
E\in\cD(E_1,r_0)\cup \cD(E_2,r_0)
\]
and similarly for $f_{[-N,\bar N']}(x_0,\omega,E)$. In fact, by the
same argument,
\[
f_{[-N,\bar N']}(z,\omega,E) \ne 0\quad \forall \,
z\in\cD(e(x_0),r_0),\;\forall\, E\in\cD(E_1,r_0)\cup \cD(E_2,r_0)
\]
Now suppose that for some choice of $E\in\cD(E_1,r_0)\cup
\cD(E_2,r_0)$,
\[
\log|f_{[-N,\bar N']}(x_0,\omega,E)| < (\bar N'-N)L - (\log
N)^{2C_0}
\]
Then Corollary~\ref{cor:2.lexcepzero} implies that $f_{[-N,\bar
N']}(z,\omega,E)=0$ where $|z-e(x_0)|< \exp(-(\log N)^{C_0})$, a
contradiction to our choice of~$r_0$. Thus,
\begin{align*}
\log|f_{[-N,\bar N']}(x_0,\omega,E)| &> (\bar N'-N)L - (\log
N)^{2C_0}\\
\log|f_{[\bar N'',N]}(x_0,\omega,E)| &> (N-\bar N'')L - (\log
N)^{2C_0}
\end{align*}
which in its turn imply the Green function bounds
\begin{align*}
  |G_{[-N,\bar N']}(x_0,\omega,E)(p,q)| &\le \exp\big(-\gamma |p-q|
  +  (\log
N)^{2C_0}\big) \quad\forall\, p,q\in [-N,\bar N']\\
|G_{[\bar N'',N]}(x_0,\omega,E)(p,q)| &\le \exp\big(-\gamma |p-q|
  +  (\log
N)^{2C_0}\big) \quad\forall\, p,q\in [\bar N'',N]
\end{align*}
and some constant $C_0=C_0(a)$, see Cramer's rule
\eqref{eq:cramerG}. These bounds prove~\eqref{eq:psij_decay} via the
Poisson formula~\eqref{eq:3.poisson}. Furthermore, inspection of the
proofs of Corollary~\ref{cor:3.Estab} and
Proposition~\ref{prop:4.Ej_sep} shows that they apply verbatim to
the situation at hand (the only difference here is that
$e^{-N^\delta}$ needs to beat $e^{N^{\vep}}$). In particular,
$|E_1-E_2|>e^{-N^\delta}$ for large~$N$ as desired.
\end{proof}

The eigenvalues $E_j^{(N)}(x,\omega)$ of the Dirichlet problem on
$[-N,N]$ are real-analytic functions of $x\in\tor$ and can therefore
be extended analytically to a complex neighborhood of $\tor$.
Moreover, by simplicity of the eigenvalues of the Dirichlet problem,
the graphs of these functions of $x$ do not cross.
Proposition~\ref{prop:4.Ej_sep} makes this non-crossing
quantitative, up to certain sections of the graphs where we lose
control. These are the portions of the graph that intersect
horizontal strips corresponding to energies in~$\cE_{N,\omega}$. The
quantitative control provided by~\eqref{eq:4.Ej_sep} allows us to
give lower bounds on the radii of the disks to which the functions
$E_j^{(N)}(x,\omega)$ extend analytically.

\begin{corollary}
\label{cor:4.2} Fix $a>1$, $ c>0$ as well as $\delta\in(0,1)$,
assume $L(\omega,E)>\gamma>0$ for all $(\omega,E)$,  and let
$\Omega_N$, $\cE_{N,\omega}$ be as in Lemma~\ref{lem:3.waschno}.
There exists a large integer $N_0(V,\gamma,a,c,\delta)$ such that
for all $N\ge N_0$ the following holds: assume $f_N(x_0,\omega_0,
E_0) = 0$ for some $x_0\in\tor$,
$\omega_0\in\tor_{c,a}\setminus\Omega_N$, and $E_0\in\IR\setminus
\cE_{N, \omega_0}$.  Then (with $\omega_0$ fixed)
\begin{equation}
\label{eq:4.b0z}
f_N(z,\omega_0, E) = \bigl(E - b_0(z)\bigr) \chi(z, E)
\end{equation}
for all $z \in \cD(x_0, r_0)$, $E \in \cD(E_0, r_1)$ where  $r_1 =
e^{-N^\delta}$, $r_0 =C^{-1} r_1$. Moreover, $b_0(z)$ is analytic on
$\cD(x_0, r_0)$, $\chi(z, E)$ is analytic and nonzero on $\cD(x_0,
r_0) \times \cD(E_0, r_1)$, $b_0(x_0) = E_0$.
\end{corollary}
\begin{proof} By Proposition~\ref{prop:4.Ej_sep}, $f_N(x_0,\omega_0, E) \ne 0$ if $E \in \cD(E_0,
r_1)$, $E \ne E_0$.  Since $H_N(x_0,\omega_0)$ is self adjoint and
\[ \big\| H_N(z,\omega_0)
- H_N(x_0,\omega_0)\bigr\| \lesssim |z - x_0|,
\] it follows that $f_N(z,\omega_0, E) \ne 0$
for any $|z - x_0| \le C^{-1} r_1$, $\frac{r_1}{2}< |E - E_0| <
{3\over 4}r_1$. The representation~\eqref{eq:4.b0z} is now obtained
by the same arguments that lead to the  Weierstrass preparation
theorem, see Proposition~\ref{th:2.weier}.
\end{proof}

We shall also require quantitative control on the function $\chi$. Let $z_0:=e(x_0)$.

\begin{corollary}
\label{cor:7.3A}  Using the notations of the previous corollary one has
$$
f_N(z, \omega_0, E) = (E-b_0 (z) ) \chi(z, E)
$$
where $\chi(z, E)$ is analytic in $\cD(z_0, r_0) \times
\cD(E_0, r_0)$ and obeys the bound
$$
NL(E_0,\omega_0) - N^{2\delta} \leq \log | \chi(z,E)| \leq NL(E_0,\omega_0) +
N^{2\delta}
$$
for any $(z, E) \in \cD(z_0, r_0/2) \times \cD(E_0, r_0/2)$.
\end{corollary}

\begin{proof} Due to the uniform upper estimates on $\log|f_N|$ one has
$$
\big| E - b_0(z) \big| \big|\chi(z,E)\big| \leq \exp(N L (E_0,\omega_0)
+ (\log N)^A )
$$
for any $(z,E) \in \cD(z_0, r_0)\times\cD(E_0, r_0)$. Take
arbitrary $z_1 \in \cD(z_0, r_0)$, $E_1 \in \cD(E_0, r_0/2)$. We distinguish
two cases: (a) If
$$
\big| E_1 - b_0(z_1) \big| \geq r_0/4
$$
then
\begin{align}
\label{eq:7.kappaup} \big| \chi(z_1, E_1) \big| &\leq 4 \exp (N L(E_0,\omega_0) +
N^\delta),\qquad  \log |\chi(z_1, E_1) | &\leq N L (E_0,\omega_0) +
2N^\delta
\end{align}
(b) Otherwise,
$$
|E-b_0(z_1)| \geq r_0/4
$$
for any $|E-E_1| = r_0/2$. Hence
$$
|\chi(z_1, E) | \leq 4 \exp(N L(E_0,\omega_0) + N^\delta)
$$
for any $|E-E_1|=r_0/2$ in this case. The maximum
principle implies~\eqref{eq:7.kappaup}. Thus~\eqref{eq:7.kappaup}
holds for any $z_1 \in \cD(z_0, r_0)$, $E_1 \in \cD(E_0,
r_0/2)$ which proves
the upper estimate for $\log | \chi(z, E)|$. Furthermore,  $|b_0(z)| \lesssim 1$ for any $z \in \cD(z_0, r_0)$.
Hence
$$
\log |f_N(z, E)| \leq \log |\chi(z, E)| + C
$$
It follows from the large deviation estimate  that
given $E_1 \in \cD(E_0, r_0)$ there exists $z_1 \in \cD(z_0,
r_0/4)$ such that
$$
\log |\chi(z_1, E_1)| > NL(E_0,\omega_0) -N^{2\delta}
$$
Therefore, the lower bound for $\log |\chi(z, E)|$ follows from the
Harnack inequalities.
\end{proof}

\section{Evaluating Jensen averages via  the Harnack inequality
}\label{sec:harnakjen}

This section is part of the ``zero counting machinery'' from~\cite{Gol Sch2}.
This machinery is of crucial importance to the gap development, as we shall see later.
Technically speaking, the goal of this section is to develop estimates for $\log\|M_N\|$ that
are analogous to those valid for $\log|f(z)|$ where $f$ is analytic.
Special attention will be paid to the location of the zeros of the entries of $M_N$.
Similar considerations appear in~\cite{Gol Sch2}, and as in that paper it will
be convenient to work in the following degree of generality:

\begin{defi}\label{def:ALDT}
Let $M(z)$ be a $ 2\times 2$ matrix-function defined
in a disk $\cD(z_0, r_0) \subset \IC$, $r_0 \ll 1$. Thus, let
\begin{equation}
\label{eq:14.11} M(z) = \begin{bmatrix}
a_{11}(z) & a_{12}(z)\\
a_{21}(z) & a_{22}(z)
\end{bmatrix}
\end{equation}
$\det M(z) = 1$. We say that $M(z)$ satisfies an {\em abstract large deviation estimate} (ALDE) provided
the following holds:  let $100\le  K:= \sup \bigl\{\|M(z)\|: z
\in \cD(z_0, r_0) \bigr\} < \infty$. Then for any $H \ge
(\log\log K)^C$ one has
\begin{equation}
\log \bigl |a_{ij}(z)\big | > \log K - H \tag{ALDE}
\end{equation}
for any entry $a_{ij}$ which is not identically zero and all
\[ z \in \cD(z_0, r_0) \setminus \cB, \quad \cB =
\bigcup\limits^J_{j=1} \cD(\zeta_j, r),\quad r= r_0 \exp
\Bigl({-H\over (\log\log K)^C}\Bigr),\] and $J \lesssim (\log\log
K)^A$.
\end{defi}

In our applications, $A,C$ will be constants as in Definition~\ref{def:constants}.
This definition is set up to be scaling invariant, which is useful throughout this section.
We begin by recalling the following version of Harnack's inequality:

\begin{lemma}
\label{lem:9.harn}
Let $f(x)$ be analytic in $\cD(z_0, r_0)$ and
non-vanishing in  $\cD(z_0, r_1)$ with $0<r_1\le r_0$. Assume that
\[ 100 \le K:= \sup \bigl\{|f(z)|:
z \in \cD(z_0, r_0) \bigr\} < \infty
\]
Assume also that\begin{equation} \label{harcond} |f(z_0)| \ge K^{-1}
\end{equation}
Then, with some absolute constant $C$,
\[ |f(\zeta)| \le C|f(z)|
\]
for any $z,\zeta \in \cD(z_0, r_2), r_2=(1 + \log K)^{-2}r_1$.
\end{lemma}
\begin{proof}
The function $u(z):= \log K - \log |f(z)|$ is harmonic and non-negative
in $\cD(z_0, r_1)$. Applying  Harnack's inequality to it in $\cD(z_0, r_1)$
yields
\begin{align*}
& [1-2(1 + \log K)^{-2}](\log K - \log|f(z_0)|)  \\
&  \le
\log K - \log|f(z)| \le [1+ 3(1 + \log K)^{-2}](\log K
-\log|f(z_0)|)
\end{align*}
for any $z\in \cD(z_0, r_2)$. Hence, using \eqref{harcond}, this implies that
\begin{align*}
 -2 - \log|f(z_0)| \le - \log|f(z)| \le 2 -\log|f(z_0)|
 \end{align*}
 for any $z\in \cD(z_0, r_2)$, and the lemma follows with $C=e^4$.
\end{proof}

Next, we turn to matrices as in Definition~\ref{def:ALDT}.

\begin{lemma}
\label{lem:2.30}
Fix some $0<\delta_1<\frac{1}{10}$. Then for $K$ sufficiently large depending on $\delta_1$
the following holds: suppose that one of the entries $a_{ij}(z)$ has no
zeros in $\cD(z_0, r_1)$, with $\exp(-(\log K)^{2\delta_1})r_0 \le
r_1\le \exp(-(\log K)^{\delta_1}) r_0$.  Then
\begin{equation}
\label{eq:2.quad} \Big | \log {\big \|M(z)\big \|\over \big \|
M(z_0)\big \|} \Big | \le \exp((\log K)^{5\delta_1})|z - z_0|
r_0^{-1}
\end{equation}
for any $|z -z_0|\le r_1$ for which the right-hand side of \eqref{eq:2.quad} is $<\frac12$, say.
\end{lemma}
\begin{proof}
We
 first claim that
\begin{equation}
  \label{eq:claim23}
  \log \|M(z_0)\| \ge \log K - (\log K
)^{5\delta_1}
\end{equation}
Let $a_{i_0j_0}(z)$ be an entry which has no zeros in  $\cD(z_0,
r_1)$. Due to condition (ALDE) there exists $z_1$ with $|z_1 -z_0|\le r_0
\exp(-(\log K)^{3\delta_1})$ such that $|a_{i_0j_0}(z_1)| \ge
K\exp(-(\log K)^{4\delta_1})
> K^{-1} $. Since $a_{i_0j_0}(z)$ is analytic and does not vanish in
$\cD(z_1,r_1/2)$, one can apply Lemma~\ref{lem:9.harn} with $z_1$ in
the role of $z_0$. Therefore,
\begin{align*}
 \log|a_{i_0j_0}(z_1)| - C
 & \le
\log|a_{i_0j_0}(z_0)|
\end{align*}
Hence,
\begin{equation}
\nn\log|a_{i_0j_0}(z_0)| \ge \log K - (\log K
)^{5\delta_1}
\end{equation}
which implies \eqref{eq:claim23}.
 Next, we note that
for any $|z-z_0|<r_0$,
\begin{align*}
\|M(z)- M(z_0)\| &\le |z-z_0| \sup_{|\zeta-z_0|\le r_0} \|M'(\zeta)
\| \le 2|z-z_0| r_0^{-1} \sup_{|\zeta-z_0|\le
2r_0}\|M(\zeta)\| \\
&\le \exp((\log K)^{5\delta_1})|z-z_0| r_0^{-1} \|M(z_0)\|
\end{align*}
where we used \eqref{eq:claim23} to pass to the final inequality.
This implies that
\[
\bigl| \frac{\|M(z)\|}{\|M(z_0)\|} -1 \bigr | \le\exp\big((\log
K)^{5\delta_1}\big) |z-z_0| r_0^{-1}
\]
for all $|z-z_0|\ll r_1$, which is the same as~\eqref{eq:2.quad}.
\end{proof}

Next, we consider the case when all entries $a_{ij}(z)$ have zeros in
$\cD(z_0, r_0)$.  Assume that for some $\zeta_0 \in
\cD(z_0, r_0/4)$ the following conditions\footnote{We will later verify conditions (a) and (b) above for the
Dirichlet determinants by means of Proposition~\ref{prop:4.Ej_sep}.} are valid\footnote{In this section $\rho_0$ is used with
a different meaning than as before; however, there is no danger of confusion.}:
\begin{enumerate}
\item[{\rm{(a)}}] each entry $a_{ij}(z)$  has exactly one zero
 in $\cD(\zeta_0, \rho_0)$, where
\[ \exp(-(\log K)^{2\delta_0})r_0 \le \rho_0 \le \exp(-(\log K)^{\delta_0}) r_0
\]
We denote this unique zero by $\zeta_{ij}$.
\item[{\rm{(b)}}] no entry $a_{ij}(z)$
has any zeros in $\cD(\zeta_0, \rho_1) \setminus \cD(\zeta_0,
\rho_0)$, where
\[\exp(-(\log K)^{\delta_1})r_0 \le \rho_1 \le r_0\]
with $0<10 \delta_1 \le \delta_0 \ll 1$.
\end{enumerate}

$K$ will need to be large depending on $\delta_0,\delta_1$.

\begin{lemma}
\label{lem:14.8} The function $b_{ij}(z):=r_0(z-\zeta_{ij})^{-1}a_{ij}(z)$ is analytic in $\cD(z_0, r_0)$
and non-vanishing in
$\cD(\zeta_0, \rho_1)$. Set $\widetilde M(z) := \bigl\{b_{ij}(z)\bigr\}_{1
\le i, j \le 2}$. Then
\begin{equation}\label{eq:9.2bnozero}
\begin{split}
T
 := \sup \bigl\{ \| \widetilde M(z)\| : z \in \cD(z_0,
r_0/2)\bigr\}  \le 5K
\end{split}
\end{equation}
Furthermore,
\begin{equation}\label{eq:9.27a} \log \|\widetilde M(z)\| \ge
\log K - (\log K )^{3\delta_1}
\end{equation}
for any $z \in \cD(\zeta_0, \rho_2)$ with $\rho_2=\rho_1 (\log K)^{-2}$, and
\begin{equation}\label{eq:frombijtoM}
 \log\|M(z)\| \ge \log K + \log (|z - \zeta_0|r_0^{-1}) - (\log K )^{4\delta_1}
\end{equation}
for any $z\in \cD(\zeta_0, \rho_2) \setminus \cD(\zeta_0,2 \rho_0)$.
\end{lemma}
\begin{proof}  For every $|z-z_0|=r_0/2$  and large $K$,
\[
 |z-\zeta_{ij}|\ge |z-z_0|-|z_0-\zeta_0|-|\zeta_0-\zeta_{ij}| \ge r_0/4-\rho_0\ge r_0/5
\]
which implies via the maximum principle that
\[
 |b_{ij}(z)|\le 5|a_{ij}(z)|\le 5K \quad \forall \;  |z-z_0|\le r_0/2
\]
as claimed. Next, it follows from (ALDE) that for some $\zeta_1\in \cD(\zeta_0,\rho_1/2)$ one has
\[
\log |b_{ij}(\zeta_1)|\ge  \log K - (\log K)^{2\delta_1}
\]
Applying Lemma~\ref{lem:9.harn} above with $\zeta_1$ in the role of $z_0$ and $\rho_1/2$ in the role
of $r_1$ implies that
\[
 |b_{ij}(z)|  \ge  \log K - (\log K)^{3\delta_1}
\]
for all $z\in\cD(\zeta_0,\rho_2)$ where $\rho_2$ is as above.
Finally, \eqref{eq:frombijtoM} follows from~\eqref{eq:9.27a} since
\begin{align*}
\log(|z-\zeta_{ij}|r_0^{-1}) &= \log(|z-\zeta_0| r_0^{-1}) - \log(|z-\zeta_0| |z-\zeta_{ij}|^{-1})  \\
&\ge \log(|z-\zeta_0| r_0^{-1}) - \log(1+\rho_0 |z-\zeta_{ij}|^{-1})  \ge \log(|z-\zeta_0| r_0^{-1}) -\log 2
\end{align*}
since $|z-\zeta_{ij}|\ge |z-\zeta_0|-|\zeta_0-\zeta_{ij}|\ge \rho_0$ for all $z\in \cD(\zeta_0, \rho_2) \setminus \cD(\zeta_0,2 \rho_0)$.
\end{proof}

Next, we obtain the analogue of Lemma~\ref{lem:2.30} for the case of $M(z)$ as in (a), (b) above.

\begin{lemma}
\label{lem:9harnakwithzeros}
For any $z,w\in \cD(\zeta_0,\rho_4)\setminus\cD(\zeta_0,\rho_3)$ one has
\begin{equation}
\begin{split}
\label{eq:0.40'}\Big |\log {\| M(z)\|\over \|M(w)\|} - \log
{|z-\zeta_0|\over |w - \zeta_0|}\Big | \le C\exp(-(\log K
)^{\delta_2})
\end{split}
\end{equation}
where $\rho_4=\exp(-(\log K)^{2\delta_2})\rho_1 $ and $\rho_3=\exp((\log K)^{\delta_2})\rho_0$, where $1>\delta_0\ge \delta_2\ge  3\delta_1>0$ and $K$ is large depending on these parameters.
\end{lemma}
\begin{proof}
Set
\[
\hat a_{ij}(z) := (z- \zeta_{0})r_0^{-1}\, b_{ij}(z),\; i, j = 1, 2, \quad
\widehat M(z) := \bigl\{\hat a_{ij}(z)\bigr\}_{1 \le i, j \le 2}
\]
Then
\begin{align*}
 \log {\| M(z)\|\over \|M(w)\|} - \log
{|z-\zeta_0|\over |w - \zeta_0|} &= \log \Big\{ \frac{\| M(z)\|}{\| \widehat M(z)\|} \frac{\| \widehat M(w)\|}{\|  M(w)\|}  \Big\}
 + \log \frac{\|\widetilde M(z)\|}{\|\widetilde M(w)\|}
\end{align*}
Since
\[
 \frac{\| M(z)\|}{\| \widehat M(z)\|}  \le 1 + \frac{\| M(z)-\widehat M(z) \|}{\| \widehat M(z)\|}
\le 1 + \max_{1\le i,j\le 2} \frac{|\zeta_0-\zeta_{ij}|}{|z-\zeta_0|} \le 1+\frac{\rho_0}{\rho_3}
\]
and similarly,
\[
 \frac{\| M(z)\|}{\| \widehat M(z)\|}  \ge 1 - \frac{\| M(z)-\widehat M(z) \|}{\| \widehat M(z)\|}
\ge 1 - \max_{1\le i,j\le 2} \frac{|\zeta_0-\zeta_{ij}|}{|z-\zeta_0|} \ge 1 - \frac{\rho_0}{\rho_3}
\]
we see that
\[
 \Big| \log \Big\{ \frac{\| M(z)\|}{\| \widehat M(z)\|} \frac{\| \widehat M(w)\|}{\|  M(w)\|}  \Big\} \Big | \le C \frac{\rho_0}{\rho_3}
\]
whereas from Lemma~\ref{lem:2.30} one has
\[
 \Big |\log \frac{\|\widetilde M(z)\|}{\|\widetilde M(w)\|} \Big | \le C\exp\big((\log K)^{3\delta_1}\big) |z-w|r_0^{-1}
\]
for any $z,w\in \cD(\zeta_0,\rho_2)$ where $\rho_2$ is as above.  In conclusion,
\begin{align*}
 \Big| \log {\| M(z)\|\over \|M(w)\|} - \log
{|z-\zeta_0|\over |w - \zeta_0|} \Big| &\le C\exp\big((\log K)^{3\delta_1}\big) \rho_4 r_0^{-1} +  C \frac{\rho_0}{\rho_3}
\end{align*}
which implies the lemma.
\end{proof}

Next, we apply these results to the propagator matrices $M_N$. Of particular importance
for the zero count are
 the {\em Jensen
averages}
\begin{align}
\cJ(u, z, r) &= \mathop{\nint}\limits_{\cD(z, r)}
[u(\zeta)-u(z)] \,  d \xi d \eta\,\label{eq:2.jensenaver}\\
\cJ(u, \cD, r_2) &= \mathop{\nint}\limits_{\cD} \cJ(u, z,
r_2)\,dx\, dy  \label{eq:9.jensenaver}
\end{align}
where $\cD$ is an arbitrary bounded domain in the second line (as usual, $\nint$ means
average).  Recall that we are assuming that $V$ is analytic on some annulus $\cA_{\rho_0}$.

\begin{prop}
\label{prop:14.10} Let $0<10\delta_1<\delta_0$ be fixed small parameters, $\omega\in\tor_{c,a}$, and $z_0\in \cA_{\rho_0/2}$.
There exists a positive integer $N_0(\delta_0,\delta_1,c,a,\gamma, V,E)$ so that the following holds:

$(1)$ Suppose that one of the Dirichlet determinants
\[ f_{[1,
N]}(\cdot, \omega, E),\; f_{[1, N-1]}(\cdot, \omega, E),\; f_{[2,
N]}(\cdot, \omega, E),\; f_{[2, N-1]}(\cdot, \omega, E)
\]
has no zeros in $\cD(z_0, r_1)$, $\exp(-N^{2\delta_1})\le r_1 \le
\exp(-N^{\delta_1})$.  Then
\begin{equation}
 \Big | \log {\big \|M_N(z, \omega, E)\big \|\over \big \| M_N(z_0,
\omega, E)\big \|}  \Big | \le C\exp(N^{5\delta_1})|z - z_0|
\end{equation}
for any $z\in\cD(z_0, r_1/2)$. In particular, with $\cJ(u, z, r)$ defined as
in~\eqref{eq:2.jensenaver} one has
\begin{equation}
 \label{eq:2.28}  \cJ \bigl(\log \|M_N(\cdot, \omega, E)\|, z_0,
r\bigr) \le r\exp(N^{5\delta_1})
\end{equation}
for any $e^{-\sqrt{N}}\le r \le r_1/2$.

\smallskip
 $(2)$
 Assume that for some
$\zeta_0\in\cA_{\rho_0/8} $ the following conditions are valid:
\begin{enumerate}
\item[{\rm{(a)}}] each
determinant
\[ f_{[1,
N]}(\cdot, \omega, E),\; f_{[1, N-1]}(\cdot, \omega, E),\; f_{[2,
N]}(\cdot, \omega, E),\; f_{[2, N-1]}(\cdot, \omega, E)
\]
has exactly one zero in $\cD(\zeta_0, \rho_0)$, where
$\exp(-N^{2\delta_0}) \le \rho_0 \le \exp(-N^{\delta_0})$

\item[{\rm{(b)}}] none of these determinants
has any zeros in $\cD(\zeta_0, \rho_1) \setminus \cD(\zeta_0,
\rho_0)$, where $\exp(-N^{\delta_1})\le \rho_1 \le \rho_0/10$.

\end{enumerate}
 Then
 \begin{equation}
\begin{split}
\label{eq:0.40}\Big |\log {\| M_N(z,\omega,E)\|\over \|M_N(\zeta,
\omega, E)\|} - \log {|z-\zeta_0|\over |\zeta - \zeta_0|} \Big| \le
\exp(-N^{\delta_2})
\end{split}
\end{equation}
for any
\[
 \rho_0 \exp(N^{\delta_2})  \le |z-\zeta_0|, |w-\zeta_0| \le \rho_1 \exp(-N^{2\delta_2})
\]
where
$\frac14>\delta_0 \ge\delta_2 \ge 3\delta_1>0$. Furthermore,
\begin{equation}
\label{eq:2.28witzeros} \Big | \cJ \bigl([ \log \|M_N(\cdot, \omega,
E)\| - \log|\cdot -\zeta_0|], z, r \bigr) \Big | \le
\exp(-N^{2\delta_2})
\end{equation}
for any $z\in\cD(\zeta_0, \rho_1/2)$ and any $\rho_0 \exp(N^{\delta_2})  \le r \le \rho_1 \exp(-N^{2\delta_2})$.
  Statements similar to $(1)$, $(2)$ hold
with respect to zeros in the $E$-variable with $z=e(x)$ arbitrary but
fixed.
\end{prop}
\begin{proof}  Part (1) follows from Lemma~\ref{lem:2.30} with a choice of $r_0=\rho_0/2$.
Part (2) follows from Lemma~\ref{lem:9harnakwithzeros}. In both cases, the (ALDE) holds due
to the large deviation theorem for $M_N$ and $f_N$, see Propositions \ref{prop:2.ldtst} (for the LDE in $x$)
and~\ref{prop:LDEinE} and Remark~\ref{rem:ldee_cartan} (for the LDE in $E$).
\end{proof}

\section{A multi-scale approach to counting zeros of Dirichlet determinants} \label{sec:zerocount}

The basic question motivating this section is as follows: suppose
\[
[1,N) = \bigcup_{j=1}^{J-1} [n_0,n_{j+1}), \quad
1=n_0<n_1<\ldots<n_J=N
\]
Can we relate the number of zeros of $f_{[1,N)}(\cdot,\omega,E)$ in
$\cD(z_0,r)$ to the sum of the numbers of zeros of
$f_{[n_j,n_{j+1})}(\cdot,\omega,E)$ in $\cD(z_0,r)$?  This seems
like a very farfetched question; indeed, since typically
\begin{equation}
  \label{eq:fake_prod}  f_{[1,N)}(\cdot,\omega,E)\ne \prod_{j=0}^{J-1}
f_{[n_j,n_{j+1})}(\cdot,\omega,E)
\end{equation}
there is no reason to assume that the zeros on the left-hand side
are in any way related to the zeros on the right-hand side.
Nevertheless, we shall see in this section that under certain
conditions (which will still be flexible enough for our purposes)
such an addition theorem does hold for the number of zeros. The
basic tools here are the avalanche principle and the Jensen averages
from Section~\ref{sec:basictools}. The former will give us something akin
to~\eqref{eq:fake_prod}, whereas the latter allows for an effective
zero count based on averaging. Averaging here is particularly
important as it ``washes out'' a set of exceptional phases that need
to be removed for the avalanche principle to hold.
Results similar to those of this section can be found
 in Sections~12 and
13 of~\cite{Gol Sch2}.
  We describe how to combine
Proposition~\ref{prop:14.10} with the avalanche principle expansion
to count precisely the number of the zeros of Dirichlet
determinants. The following definition is very important in this
regard. For the following definition recall that
$\cZ(f,z,r)=\{\zeta\in\cD(z,r)\::\: f(\zeta)=0\}$.

\begin{defi}
\label{def:2.adj}
 Let $\ell\ge 1$ be some integer, and $s\in\ZZ$. Fix $(\omega,E)$ as well
 as some disk
 $\cD(z_0,r_0)$. We say that $s$ is {\em adjusted}
to  $(\cD(z_0,r_0),\omega,E)$ at scale $\ell$ if for all $\ell\le
k\le100\ell$
\[ \cZ(f_{k}(\cdot e((s+m)\omega),\omega,E),z_0,r_0)=\emptyset \qquad \forall\;|m|\le 100\ell. \]
\end{defi}
First, an easy but useful observation: the determinants appearing in
the definition of ``adjusted'' automatically satisfy a large
deviation type estimate.

\begin{lemma}
  \label{lem:adjLDT} Let $\omega_0\in\tor_{c,a}$ and $E_0\in\IC$.
  There exists $\ell_0=\ell_0(V,\rho_0,a,c)$ so that if $s$ is
  adjusted to $(\cD(z_0,r_0),\omega_0,E_0)$ at scale $\ell\ge \ell_0$ with $r_0> e^{-\sqrt{\ell}}$, then
  \[
\log|f_k(z e((s+m)\omega_0),\omega_0,E_0)| > kL(\omega_0,E_0) -
k^{\frac34} \quad \forall\; |z-z_0|< r_0/2
  \]
  for all $|m|\le 100\ell$ and $\ell\le
k\le100\ell$.
\end{lemma}
\begin{proof}
  Suppose not. Then
  \[
\log|f_k(z_1 e((s+m)\omega_0),\omega_0,E_0)| < kL(\omega_0,E_0) -
k^{\frac34}
  \]
  for some choice of $|m|\le 100\ell$, $\ell\le
k\le100\ell$, and $|z_1-z_0|< r_0/2$. By
Corollary~\ref{cor:2.lexcepzero}, there exists
\[
|z_2-z_1|< \exp\big(-k^{\frac34}/(\log k)^{C_0}\big) < r_0/2
\]
such that $f_k(z_2 e((s+m)\omega_0),\omega_0,E_0)=0$.  But this
contradicts Definition~\ref{def:2.adj}.
\end{proof}

We shall now prove that the notion of ``adjusted'' allows for an
affirmative answer to the ``additivity of the zero count'' question
stated at the beginning of this section. In the following
proposition, the constants implicit in the $\ll$ and $\les$
notations are absolute.
 For the $\nu_f$ notation,
see~\eqref{eq:2.zerondisk}. Also, as usual $V$ is analytic on
$\cA_{\rho_0}$ for some $\rho_0>0$.
\begin{prop}
  \label{prop:add_zeros} Let $a>1, c>0$ and fix $\omega_0\in \tor_{c,a}$. Assume that
$L(\omega_0,E_0)>\gamma>0$ where $E_0\in\IC$ is arbitrary but fixed.
There exists a large integer $N_0=N_0(V,\rho_0,\gamma,a,c,E_0)$ such
that for any $N\ge N_0$ the following holds.  Let $\ell$ be an
integer $(\log N)^A \le \ell$ where
$A=A(V,\rho_0,\gamma,a,c,E_0)$ is a large constant.  Suppose
that with some $n_0:=1<n_1<n_2 <\ldots< n_J <n_{J+1}:=N$,
  \[
[1,N] = \bigcup_{j=1}^J [n_{j-1},n_{j}) \cup [n_{J},N],\quad
m:=\min_{0\le j\le J} (n_{j+1}-n_j) > 10\ell
  \]
We assume that  $n_j$ is adjusted to $(\cD(z_0, r_1),\omega_0,E_0)$ at
scale $\ell$ for each $0\le j\le J+1$, where $z_0=e(x_0)$,
$x_0\in\tor$, and $e^{-\sqrt{\ell}}<r_1<\rho_0$. Let $\Lambda_j:=[n_j,n_{j+1})$ with
$0\le j\le J-1$ and $\Lambda_{J}:=[n_J,N]$.
Then, with $e^{-\sqrt{m}} < r_0 < N^{-1}r_1$ and $r_2 = C^{-1} r_0$,
\[
 \cJ(\log|f_{[1,N]}(\cdot,\omega_0,E_0)|, z_0, r_0, r_2) = \sum_{j=1}^J \cJ(\log|f_{\Lambda_j}(\cdot,\omega_0,E_0)|, z_0, r_0, r_2) + O(N\ell^{-1} \big(r_0
r_1^{-1}+e^{-\ell^{\frac12}}\big))
\]
Furthermore, suppose also that for
all $0\le j\le J$ one
has \begin{equation}   \label{eq:zero_freeann}
\cZ\bigl(f_{\Lambda_j}\bigl(\cdot, \omega_0, E_0\bigr), \cD(z_0,
3r_0/2)\setminus \cD(z_0, r_0/2)\bigr)=\emptyset
\end{equation} Then
$$
\nu_{f_{[1,N]}(\cdot, \omega_0, E_0)}(z_0, r_0) = \sum^{J}_{j=0}
\nu_{f_{\Lambda_j}(\cdot , \omega_0, E_0)} (z_0, r_0)
$$
Finally, if every $1\le s\le N$ is adjusted to $(\cD(z_0,r_1),\omega_0,E_0)$
at scale~$\ell$, then
$\nu_{f_N(\cdot,\omega_0,E_0)}(z_0,N^{-1}r_1)=0$.
\end{prop}
\begin{proof}
We begin by noting that
\[
\begin{bmatrix} f_N(z,\omega_0,E_0) & 0\\ 0
  & 0\end{bmatrix} = \begin{bmatrix} 1 & 0\\ 0
  & 0\end{bmatrix} M_N(z,\omega_0,E_0) \begin{bmatrix} 1 & 0\\ 0
  & 0\end{bmatrix}
\]
The idea is to apply the avalanche principle to the matrix on the
right-hand side by writing it as a product of monodromy matrices
corresponding to the $\Lambda_j$. However, we need to connect  any
two such adjacent matrices by (much shorter) ones of length~$\ell$.
Hence, we let
\[
\Lambda_j':= [n_j+2\ell, n_{j+1}-2\ell],\quad A_j(z) :=
M_{\Lambda_j'}(z,\omega_0,E_0), \quad 0\le j \le J
\]
where $M_\Lambda$ for an interval $\Lambda\subset\ZZ$ denotes the
monodromy matrix corresponding to~$\Lambda$. Next, we define
\begin{align*}
   B_{j,1}(z) &:= M_{(n_j-2\ell,n_j-\ell]}(z,\omega_0,E_0),\quad
 B_{j,2}(z):= M_{(n_j-\ell,n_j]}(z,\omega_0,E_0), \\
B_{j,3}(z) &:= M_{(n_j,n_j+\ell)}(z,\omega_0,E_0),\quad B_{j,4}(z):=
M_{[n_j+\ell,n_j+2\ell)}(z,\omega_0,E_0)
\end{align*}
for $1\le j\le J$  and
\[
B_{0,3}(z) := \begin{bmatrix} 1 & 0\\ 0
  & 0\end{bmatrix} M_{[1,\ell)}(z,\omega_0,E_0),\quad B_{J+1,2}(z) :=
  M_{(N-\ell,N]}(z,\omega_0,E_0) \begin{bmatrix} 1 & 0\\ 0
  & 0\end{bmatrix}
\]
as well as $B_{0,4}(z) := M_{[\ell,2\ell)}(z,\omega_0,E_0)$,
$B_{J+1,1}(z) := M_{(N-2\ell,N-\ell]}(z,\omega_0,E_0)$.  Then,
\[
\begin{bmatrix} f_N(z,\omega_0,E_0) & 0\\ 0
  & 0\end{bmatrix} = B_{0,3}(z)B_{0,4}(z)A_0(z)\prod_{j=1}^{J} \Big(B_{j,1}(z)B_{j,2}(z)B_{j,3}(z)B_{j,4}(z)A_j(z)\Big)
   B_{J+1,1}(z)B_{J+1,2}(z)
\]
Since each $n_j$ is adjusted to $(\cD(z_0, r_1),\omega_0,E_0)$ at
scale $\ell$, Lemma~\ref{lem:adjLDT} implies that each $B_{j,i}$
satisfies an estimate of the form
\[
\log\|B_{j,i}(z)\| > \ell L(\omega_0,E_0) - \ell^{\frac34}
\quad\forall\; z\in \cD(z_0, r_1/2)
\]
The point here is that we avoid the removal of sets of measure
$e^{-\sqrt{\ell}}$ in the $z$-variable coming from the large
deviation theorem; such sets would be unnecessarily large. On the
other hand, the large deviation theorem applied to each $A_j$
implies that there exists $\cB:=\cB_{N, \omega_0, E_0}\subset\IC$
with $\mes
  (\cB) \le \exp\bigl(-m^{3/4}\bigr)$  so that
  for any $z \in \cD(z_0,r_1/2)\setminus \cB$
\[
\log\|A_j(z)\| > |\Lambda_j|L(\omega_0,E_0) - |\Lambda_j|^{\frac78}
\]
Therefore, for all $z \in \cD(z_0,r_1/2)\setminus \cB$, one has the
avalanche  principle expansion
\begin{equation}\label{eq:2.29}\begin{aligned}
  \log \big | f_N(z, \omega_0, E_0)\big | &=  \log
   \|B_{0,1}(z) B_{0,2}(z) \| + \log \|B_{0,2}(z)A_0(z)\| + \log\| A_0(z) B_{1,1}(z)\|+\ldots \\
& +\log\|B_{n,4}A_n(z)\|+\log\|A_n
  B_{n+1,1}(z)\| + \log \|B_{n+1,1}(z) B_{n+1,2}(z)\|\\
& -  \Big(\log\|B_{0,2}(z)\|+  \log  \| A_0(z) \| +\ldots + \log\| B_{1,1}(z)\| + \ldots +
  \log\|B_{n,4}(z)\| +\\
& + \log\|A_n(z)\| + \log \|B_{n+1,1}(z)\|+\log \|B_{n+1,2}(z)\|\Big) +
  O(e^{-\ell^{1/2}})
\end{aligned}
\end{equation}
Next, we apply the avalanche principle again, this time to each of
the $f_{\Lambda_j'}(z,\omega_0,E_0)$, $0\le j\le J$. Thus, we write
\[
\begin{bmatrix} f_{\Lambda_j}(z,\omega_0,E_0) & 0\\ 0
  & 0\end{bmatrix} = \begin{bmatrix} 1 & 0\\ 0
  & 0\end{bmatrix} M_{\Lambda_j}(z,\omega_0,E_0) \begin{bmatrix} 1 & 0\\ 0
  & 0\end{bmatrix} =  \widetilde B_{j,3}(z) B_{j,4}(z)
  A_j(z) B_{j+1,1}(z) \widetilde B_{j+1,2}(z)
\]
where
\[
\widetilde B_{j,3}(z) = \begin{bmatrix} 1 & 0\\ 0
  & 0\end{bmatrix} M_{[n_j,n_j+\ell)}(z,\omega_0,E_0), \quad \widetilde
  B_{j,2}(z) := M_{[n_{j}-\ell,n_j )}(z,\omega_0,E_0)\begin{bmatrix} 1 & 0\\ 0
  & 0\end{bmatrix}
\]
Thus, for all $z\in \cD(z_0,r_1/2)\setminus \cB$,
\begin{align*}
 & \log|f_{\Lambda_j}(z,\omega_0,E_0)| = \log \|\widetilde B_{j,3}(z)
  B_{j,4}(z)\| + \log \| B_{j,4}(z) A_j(z)\| + \log\|A_j(z)
  B_{j+1,1}(z)\|  \\
  & + \log\|B_{j+1,1}(z) \widetilde B_{j+1,2}(z)\|- \Big( \log \|B_{j,4}(z)\| + \log \|  A_j(z)\| + \log
  \| B_{j+1,1}(z)\| \Big) + O(e^{-\sqrt{\ell}})
\end{align*}
where we again applied Lemma~\ref{lem:adjLDT}. Summing these
expressions over~$j$, and subtracting the result
from~\eqref{eq:2.29} shows the following: for all $z\in
\cD(z_0,r_1/2)\setminus \cB$
\begin{equation}
  \label{eq:difffN}  \log \big | f_N(z, \omega_0, E_0)\big | - \sum_{j=0}^n \log \big |
f_{\Lambda_j}(z, \omega_0, E_0)\big | = \sum_{k=1}^K \pm \log \|W_k
(z)\| + O(e^{-\sqrt{\ell}})
\end{equation}
with  $K\le N$ and $W_k$ being a $2\times2$-matrix each entry of
which is either identically zero or a determinant
$f_{\Lambda}(z,\omega_0,E_0)$ with $\Lambda\subset[1,N]$ being an
interval of length proportional to~$\ell$; the point here is that we
set up the avalanche principle in such a way that the bulk terms
containing $A_j(z)$ exactly cancel in~\eqref{eq:difffN}. Moreover,
every $W_k$ contains at least one nonzero entry, which necessarily
is a determinant satisfying the conditions of
Definition~\ref{def:2.adj}, i.e., it does not vanish
on~$\cD(z_0,r_1)$. Thus, by Proposition~\ref{prop:14.10}, and with
$r_2:=r_0/10$,
\begin{equation}\label{eq:2.30}
 \Big | \cJ\bigl(\log|f_N(\cdot, \omega_0, E_0) |, z_0, r_0,
r_2\bigr)  - \sum_{j=0}^{J} \cJ(\log|f_{\Lambda_j}(\cdot, \omega_0,
E_0) |, z_0, r_0, r_2) \Big | \les N\ell^{-1} \big(r_0
r_1^{-1}+e^{-\ell^{\frac12}}\big)
\end{equation} We used here that $J\le N\ell^{-1}$.
It is important to note that $r_2$ is very large compared to the
measure of~$\cB$. Hence, when applying the averaging operator~$\cJ$
the exceptional set~$\cB$ only produces an additive error of the
form~$e^{-\ell^{\frac12}}$. By \eqref{eq:2.30}
and~Corollary~\ref{cor:2.jensenequal},
\[
\nu_{f_N (\cdot, \omega_0, E_0)} (z_0, r_0-r_2)\le \nu_g (z_0, r_0 +
r_2),\qquad \nu_g(z_0, r_0 -r_2) \le \nu_{f_N (\cdot, \omega_0,
E_0)} (z_0, r_0 + r_2),
\]
\noindent where
\[
g(z):=\prod^{J}_{j=0} f_{\Lambda_j}\bigl(z, \omega_0, E_0\bigr)
\]
\noindent Replacing $r_0$ by $(r_0 \pm r_2)$, one obtains similarly
\[
\nu_{f_N (\cdot, \omega_0, E_0)} (z_0, r_0) \le \nu_g (z_0, r_0 +
2r_2), \qquad \nu_g (z_0, r_0 - 2r_2)\le \nu_{f_N (\cdot, \omega_0,
E_0)} (z_0, r_0)
\]
\noindent Due to the assumptions of the lemma
\[
\nu_g (z_0, r_0+2r_2)=\nu_g(z_0, r_0-2r_2)
\]
\noindent and the assertion follows.
\end{proof}
Heuristically speaking,  the fact that edges of the main intervals
$\Lambda_j$ are adjusted in Proposition~\ref{prop:add_zeros} ensures
that these intervals are ``independent'' of each other. This is of
course crucial for the zero count to work. One can think of it in
this way: each zero $z_{jk}$ of
$\det(H_{\Lambda_j}(\cdot,\omega_0)-E)$ produces an eigenstate
$\psi_{jk}$ on~$\Lambda_j$ with energy~$E_0$ that is localized
strictly inside~$\Lambda_j$ (since by adjustedness the Green
function decays exponentially close to the edges of $\Lambda_j$).
For this reason, we can simply extend each of these eigenstates
$\psi_{jk}(z)$ to all of~$[1,N]$ by setting them equal to zero
outside of~$\Lambda_j$. This creates an approximate eigenstate of
$H_{[1,N]}(x_0,\omega_0)$ with energy~$E_0$ and therefore a  zero
$z'$ of $\det(H_{[1,N]}(\cdot,\omega_0)-E)$ which is very close
to~$z_{jk}$. Note that in this process the
assumption~\eqref{eq:zero_freeann} plays a very important role: it
guarantees that we do not ``pull in'' any zero from outside of
$\cD(z_0,r_0)$ --- this may certainly occur under the process we
just described. This discussion shows that the assumptions of
Proposition~\ref{prop:add_zeros} are essentially optimal.
We now consider the exact same questions but with regard to zeros
in~$E$ rather than~$z$.
Just as in the case of zeros in $z$ there is a notion of
``adjusted'' for zeros in~$E$.

\begin{defi}\label{def:2.eaddjust}
 Let $\ell\ge 1$ be some integer, and $s\in\ZZ$. Fix $(z,\omega)$. We say that $s$ is {\em adjusted}
to  $(z,\omega,\cD(E_0,r_0))$ at scale $\ell$ if for all $\ell\le
k\le100\ell$
\[ \cZ(f_{k}(z e((s+m)\omega),\omega,\cdot),E_0,r_0)=\emptyset \qquad \forall\;|m|\le 100\ell. \]
\end{defi}

Using self-adjointness of $H_\Lambda(x,\omega)$, for $x\in\IR$ and
real-valued $V$ we can say that the notion of ``adjusted'' is symmetric with
respect to $z$ and~$E$. This will use the crucial Corollary~\ref{cor:2.lexcepzero}.

\begin{lemma}
  \label{lem:move_z}
  Let $V$ be real-valued on~$\tor$ and let $z_0=e(x_0)$ where $x_0\in\IR$.
There exists a constant $C(V)$ with the following property:
  suppose that $s$ is  adjusted
to  $(z_0,\omega_0,\cD(E_0,r_0))$ at scale $\ell$. Then $s$ is
adjusted to  $(z,\omega_0,\cD(E_0,r_0))$ at scale $\ell$ for every
$z\in \cD(z_0,C(V)^{-1} r_0)$. In other words, $s$ is adjusted to
$(\cD(z_0,C(V)^{-1} r_0),\omega_0,E)$ at scale $\ell$ for every
$E\in\cD(E_0,r_0)$. In particular,
\[
\log|f_k(e(z+(s+m)\omega_0),\omega_0,E)| > kL(\omega_0,E) -
k^{\frac34} \quad\forall\; |E-E_0|< r_0/2,\; \forall\,|z-z_0|<
C(V)^{-1} r_0
\]
and for all $|m|\le 100\ell$, $\ell\le k\le100\ell$.
Conversely, if $s$ is adjusted to $(\cD(z_0, r_0),\omega_0,E_0)$ at scale $\ell$, then
$s$ is adjusted to $(z_0,\omega_0,\cD(E_0,r_1))$ at scale $\ell$, where $\log (r_1^{-1}) =\log(r_0^{-1}) (\log \ell)^{C_2}$.
\end{lemma}
\begin{proof}
We need to show that for all $\ell\le k\le100\ell$
\[ \cZ(f_{k}(z e((s+m)\omega_0),\omega_0,\cdot),E_0,r_0/2)=\emptyset \qquad \forall\;|m|\le 100\ell \]
and all $|z-z_0|\le C(V)^{-1} r_0$.
  Suppose $f_{k}(z e((s+m)\omega_0),\omega_0,E)=0$ for some choice
  of $|z-z_0|\le C(V)^{-1} r_0$, $E\in\cD(E_0,r_0/2)$, and $k,m$ in
  the admissible ranges. Since for large enough $C(V)$
  \[
\big\| H_{[1,k]}(z e((s+m)\omega_0),\omega_0) - H_{[1,k]}(z_0
e((s+m)\omega_0),\omega_0) \big\|\le C_1(V) |z-z_0| < r_0/2
  \]
  and $H_{[1,k]}(z_0
e((s+m)\omega_0),\omega_0)$ is Hermitian, we conclude by
  Lemma~\ref{lem:herm} that the
latter operator would need to have an eigenvalue in the interval
$(E-r_0/2,E+r_0/2)$. Since this is included in $\cD(E_0,r_0)$, we
arrive at a contradiction to Definition~\ref{def:2.eaddjust}. The
final statement follows from Lemma~\ref{lem:adjLDT}.

For the converse, note that if $s$ is adjusted to $(\cD(z_0, r_0),\omega_0,E_0)$ at scale $\ell$, then
for any $\ell\le k\le 100\ell$
\[
 \log|f_k(z_0,\omega_0,E_0)|> k L(E_0,\omega_0) - \log(r_0^{-1}) (\log k)^{C_0}
\]
by Corollary~\ref{cor:2.lexcepzero}.  Next, by Corollary~\ref{cor:2.lipnorm},
$f_k(z_0,\omega_0,E)\ne0$ for all $|E-E_0|< r_1$ where $r_1$ is as in the statement of the lemma.
\end{proof}

\noindent We are now in a position to state the analogue of
Proposition~\ref{prop:add_zeros} with regard to the $E$-variable. Recall that the proof of that
result used the large deviation estimate. Here, we shall do the same
but with regard to the matrix function $E\mapsto
M_N(e(x),\omega,E)$. The large deviation estimate for this purpose is provided by Proposition~\ref{prop:LDEinE}.

\begin{prop}
  \label{prop:add_zerosE} Let $V$ be real-valued and  $a>1, c>0$ and fix $\omega_0\in \tor_{c,a}$. Assume that
$L(\omega_0,E_0)>\gamma>0$ where $E_0\in\IC$ is arbitrary but fixed.
There exists a large integer $N_0=N_0(V,\rho_0,\gamma,a,c,E_0)$ such
that for any $N\ge N_0$ the following holds. Let $\ell$ be an
integer such that $(\log N)^A\le  \ell$ where
$A=A(V,\rho_0,\gamma,a,c,E_0)$ is a large constant. Suppose
that with some $n_0:=1<n_1<n_2 <\ldots< n_J <n_{J+1}:=N$,
  \[
[1,N] = \bigcup_{j=1}^J [n_{j-1},n_{j}) \cup [n_{J},N],\quad
m:=\min_{0\le j\le J} (n_{j+1}-n_j) > 10\ell
  \]
Suppose moreover that  $n_j$ is adjusted to $(z_0,\omega_0,\cD(E_0,r_1))$ at
scale $\ell$ for each $0\le j\le J+1$, where $z_0=e(x_0)$,
$x_0\in\tor$, and $e^{-\ell^\frac14}<r_1<\exp(-(\log\ell)^{C_0})$.
Let
$\Lambda_j:=[n_j,n_{j+1})$ with $0\le j\le J-1$ and
$\Lambda_{J}:=[n_J,N]$. Let  $e^{-m^{\frac14}} <
r_0 < N^{-1}r_1$ be arbitrary. Then, with $r_2 = C^{-1} r_0$,
\[
 \cJ(\log|f_{[1,N]}(z_0,\omega_0,\cdot)|, E_0, r_0, r_2) = \sum_{j=1}^J \cJ(\log|f_{\Lambda_j}(z_0,\omega_0,\cdot)|, E_0, r_0, r_2) + O(N\ell^{-1} \big(r_0
r_1^{-1}+e^{-\ell^{\frac12}}\big))
\]
Furthermore, suppose  that for all $0\le j\le J$  one has
$$\cZ\bigl(f_{\Lambda_j}\bigl(z_0, \omega_0, \cdot\bigr),
\cD(E_0, 3r_0/2)\setminus \cD(E_0, r_0/2)\bigr)=\emptyset$$ Then
$$
\nu_{f_{[1,N]}(z, \omega_0,\cdot)}(E_0, r_0) = \sum^{J}_{j=0}
\nu_{f_{\Lambda_j}(z , \omega_0,\cdot)} (E_0,r_0) \qquad \forall\;
|z-z_0|< C(V)^{-1} r_0
$$
Finally, if every $1\le s\le N$ is adjusted to $(z_0,\omega_0,\cD(E_0,r_1))$
at scale~$\ell$, then $\nu_{f_N(\cdot,\omega_0,E_0)}(z,N^{-1}r_1)=0$
for all $|z-z_0|< C(V)^{-1}N^{-1}r_1$.
\end{prop}
\begin{proof} This is very similar to the proof of
  Proposition~\ref{prop:add_zeros}.
  More precisely, running the argument of Proposition~\ref{prop:add_zeros}
  in the variable~$E$ rather than~$z$ yields the following:
$$
\nu_{f_{[1,N]}(z, \omega_0,\cdot)}(E_0, ur_0) = \sum^{J}_{j=0}
\nu_{f_{\Lambda_j}(z , \omega_0,\cdot)} (E_0,ur_0)
$$
for all $x\in\tor$
where
$\frac45<u<\frac65$. Therefore, using Lemma~\ref{lem:herm} yields
that
\[
\nu_{f_{[1,N]}(z, \omega_0,\cdot)}(E_0, r_0) = \sum^{J}_{j=0}
\nu_{f_{\Lambda_j}(z , \omega_0,\cdot)} (E_0,r_0)
\]
for all $|z-z_0|< C(V)^{-1} r_0$ as claimed.
\end{proof}

In applications we will need to chose the $n_j$ to be adjusted.
This is can be done via the following results.

\begin{lemma}
\label{lem:181}  Let $
\omega_0 \in
\tor_{c,a} $,  $x_0 \in \tor$, $E_0 \in \IR$, and
 $n_0\in\IZ$.
Given $\ell\gg1$ and $r_1 = \exp(-(\log \ell)^C)$,
there exists
\beeq
\label{eq:m1m0}
n_0'\in [n_0- {\ell}^6 ,n_0+ {\ell}^6]
\eneq
such that with $z_0 = e(x_0)$,
\beeq\label{eq:goodell}
f_k \bigl(\cdot e(n\omega_0), \omega_0, E_0\bigr)\ \mbox{has no zero in
$\cD(z_0, r_1)$}
\eneq
for any $\big |n - n_0'\big | \le 100\ell$ and $\ell \le k\le 100\ell$.  In other words,
each  $n_0'$ is adjusted to $(z_0,\omega_0,\cD(E_0,r_2))$ with
$r_2 = \exp(-(\log \ell)^{2C})$.
\end{lemma}
\begin{proof}
Suppose this fails. Then there exists a sequence $\{k_j\}_{j=1}^J\subset [\ell,100\ell]$ with
$J\ge \ell^4$ as well as an increasing sequence $\{n_j\}_{j=1}^J\subset [n_0- {\ell}^6 ,n_0+ {\ell}^6]$ so that
\[
\cZ(f_{k_j}(\cdot e(n_j\omega_0),\omega_0,E_0)) \cap \cD(z_0,r_1)\ne\emptyset
\]
for each $1\le j\le J$. Since there are at most $100\ell$ choices for $k_j$, there exists some $j_0$ in this range such that
\[
\cZ(f_{k_{j_0}}(\cdot e(m_i\omega_0),\omega_0,E_0)) \cap \cD(z_0,r_1)\ne\emptyset
\]
for some increasing sequence $\{m_i\}_{i=1}^{J'}\subset  [n_0- {\ell}^6 ,n_0+ {\ell}^6]$ where $J'\ge \ell^2$.
Since $f_{k_{j_0}}(\cdot ,\omega_0,E_0)$ has at most $C(V)\ell$ zeros, it follows that there exists $z_1\in \cD(z_0,r_1)$
as well as $m\in[1,2\ell^6]$ with the property that $z_1 e(m\omega_0)\in \cD(z_0,r_1)$. However, this
contradicts the Diophantine property of $\omega_0$.
\end{proof}

This lemma gives  us a lot of room to find adjusted sequences.

\begin{corollary}
\label{cor:18.2}
Assume that
$\omega_0 \in \tor_{c,a}$.
Given a disk $\cD(z_0, r_1)$, $r_1 \asymp
\exp(-(\log \ell)^A)$ and an increasing sequence $\{\tilde n_j\}_{j=1}^{j_0}$
such that
$\tilde n_{j+1}-\tilde n_j >\ell^7$ for $1\le j<j_0$,
there exists an increasing sequence $\{n_j\}_{j=1}^{j_0}$
which is adjusted to $\cD(z_0, r_1)$ at scale $\ell$ and such that
\begin{equation}
\label{eq:184} \big | n_j - \tilde n_j\big | < \ell^6,\quad 1\le j\le j_0.
\end{equation}
\end{corollary}
\begin{proof}
Simply apply the previous lemma to each $\tilde n_j$.
\end{proof}

We now show how one can apply this zero count to improving the bound on the
separation between the zeros as in Proposition~\ref{prop:2.zerosepar}.
The point is that due to passing to a smaller scale we will be able to substantially
reduce the size of~$t$ in Proposition~\ref{prop:2.zerosepar}.

\begin{prop}
\label{prop:4.elimination} Assume that $L(\omega, E) \geq \gamma
> 0$ for all\footnote{One can localize here as usual.} $\omega, E$. Given $c>0$, $a>1$, and $A>1$ there exists $N_0
= N_0(V, c, a, \gamma, A)$ such that for
any $N \geq N_0$ and $T\ge2N$ there exist ${\Omega}_{N,T} \subset \tor$ and
${\cE}_{N, \omega,T} \subset \IR$ with
$$
\mes ({\Omega}_{N,T} ) \leq T\exp(-(\log N)^A),\;\; \compl (\Omega_{N,T})
\leq T^2\, N
$$

$$
\mes ({\cE}_{N, \omega,T}) \leq T\exp(-(\log N)^A),\;\; \compl (
{\cE}_{N, \omega,T}) < T^2\, N
$$ and with the following property:
 for any $\omega \in \tor_{c,a}
\setminus \Omega_{N,T} $, $z_0=e(x_0)$, $E_0 \in \IR \setminus {\cE}_{N, \omega,T}$ and any $N'$, $t$
which satisfy the following conditions
\begin{itemize}
\item[(i)] $(\log \min(N,N'))^{C_0}\ge \log \max(N,N')$
\item[(ii)] $2N\le t\le T$,
\end{itemize}
one has
$$
\cZ(f_N(z_0, \omega, \cdot), \cD(E_0, r_0))\cap \cZ (f_{N'} (z_0 e(t\omega),
\omega, \cdot), \cD(E_0, r_0)) = \emptyset
$$
where $r_0 :=\exp (-(\log N)^{A})$.

\end{prop}

\begin{proof} Let $\ell$ be an integer, $\ell=(\log N)^{4A}$. Let $\Omega_{\ell_1, \ell_2, t', H}$,
and $\cE_{\ell_1, \ell_2, t', H,
\omega}$ be as in Proposition~\ref{prop:2.zerosepar}. Set $H=\ell^{1/4}$,
\begin{align*}
{\Omega}_{N,T} &=\bigcup_{N\le t'\le T} \quad \bigcup_{\ell\le \ell_1,\ell_2 \le 100\ell} \Omega_{\ell_1, \ell_2, t', H}\\
\\ {\cE}_{N, \omega,T} &= \bigcup_{N\le t'\le T}\quad \bigcup_{\ell\le \ell_1,\ell_2 \le 100\ell} \cE_{\ell_1, \ell_2, t', H,
\omega}
\end{align*}
Assume that $\cZ(f_N(z_0, \omega, \cdot), \cD(E_0, r_0))\neq \emptyset$. Due to the last part of Proposition~\ref{prop:add_zerosE}
there exists $1\le s\le N$ which is not adjusted to $(z_0,\omega_0,\cD(E_0,r_1))$
at scale~$\ell$, where $r_1:= \exp(-\ell^{7/24})$; in other words,
 \[\cZ(f_{\ell_1}(z_0e(s\omega), \omega, \cdot), \cD(E_0, r_1))\neq \emptyset\] for some
 $\ell\le$ $\ell_1\le 100\ell$ . Then, due to Proposition~\ref{prop:2.zerosepar}, provided
$\omega\in\tor_{c,a}\setminus \tilde{\Omega}_{N,T}$ and $E_0\in \IR\setminus \tilde{\cE}_{N, \omega,T}$
 one has \[ \cZ(f_{\ell_2}(z_0e(s'\omega), \omega, \cdot), \cD(E_0, r_1))=\emptyset\]
for any $2N\le s' \le N'$.  We used here that $t':=s'-s> N$, and $N> \exp((\log \log \ell )^{C_0})$,
where $C_0$ is the same as int he statement of Proposition~\ref{prop:2.zerosepar}. That suffices for Proposition~\ref{prop:2.zerosepar}.
Hence \[ \cZ(f_{N'}(z_0e(t\omega), \omega, \cdot), \cD(E_0, r_0))=\emptyset \]
due to last part of Proposition~\ref{prop:add_zerosE}.
\end{proof}

Next, we use this improvement in the size of~$t$ to reduce the size of the window of localization
in Section~\ref{sec:anderson}. The gain here is due to a ``induction on scales'' which
enters into the proof of the previous proposition through the zero count used there (Proposition~\ref{prop:4.elimination}).
 We shall use the notations of that section as for example~$\nu_j^{(N)}(x,\omega)$.

\begin{cor}\label{cor:Qreduce}
Assume that $L(\omega, E) \geq \gamma
> 0$ for all $\omega, E$. Given $c>0$, $a>1$, $A$ there exists $N_0
= N_0(V, \gamma, a, c,A)$ such that for
any $N \geq N_0$ there exist ${\Omega}_{N} \subset \tor$,  $\Omega'_{N} \subset \tor$, $\cE_{N,\omega}\subset\IR$, $\cE'_{N,\omega}\subset\IR$ with
\begin{equation}\label{eq:level1_est}\begin{aligned}
\mes ({\Omega}_{N} ) &\leq \exp(-(\log N)^A),\;\; \compl (\Omega_{N})
\leq N^4\\
\mes (\Omega'_{N} ) &\leq \exp(-(\log \log N)^A),\;\; \compl (\Omega'_{N})
\leq \exp((\log \log N)^{A/2})\\
\mes ({\cE}_{N, \omega}) &\leq \exp(-(\log N)^A),\;\; \compl (
{\cE}_{N, \omega}) < N^4\\
\mes (\cE'_{N, \omega}) &\leq \exp(-(\log \log N)^A),\;\; \compl (
\cE'_{N, \omega}) < \exp((\log \log N)^{A/2})
\end{aligned}\end{equation}
satisfying the following properties: for any $\omega\in\tor_{c,a}\setminus(\Omega_N\cup\Omega'_{N})$ and any $x\in\tor$, any  $\ell^2$-normalized eigenfunction
$\psi_j^{(N)}(x,\omega)$ of $H_{[-N,N]}(x,\omega)$ with
associated eigenvalue $E_j^{(N)}(x,\omega)\in\IR\setminus (\cE_{N,\omega}\cup\cE'_{N,\omega})$ satisfies
\[
|\psi_j^{(N)}(x,\omega)(n)|\le C\exp(-\gamma \dist(n,\Lambda_j)/2)
\]
for all $n\in[-N,N]$ where $\Lambda_j:=[\nu_j^{(N)}(x,\omega)-\ell,\nu_j^{(N)}(x,\omega)+\ell]\cap[-N,N]$ where $\ell=(\log N)^{4A}$.
\end{cor}
\begin{proof}
Inspection of the proofs in Section~\ref{sec:anderson} shows that the size of the window of localization
is determined by the size of the shift~$t$ that assures separation of the zeros as in Proposition~\ref{prop:4.elimination}.
Note that we apply that proposition on scale~$\ell$ rather than $N$; the point here is that we then take $T=\exp((\log\log N)^{B_1})$
which is the size of the localization window guaranteed by Proposition~\ref{lem:3.4} from Section~\ref{sec:anderson}. As long
as we choose $A\gg B_1$ the corollary immediately follows.
\end{proof}

\begin{defi} \label{def:level1}
Using the notations of the previous corollary, we set
\[
\Omega_N^{(1)}:= \Omega_N\cup\Omega_{N}',\quad \cE_{N,\omega}^{(1)}:= \cE_{N,\omega}\cup\cE'_{N,\omega}
\]
In what follows we shall use this notation for sets satisfying the estimates from \eqref{eq:level1_est}.
We shall also need to go down one more level: thus, set
\[
\Omega_N^{(2)}:= \Omega_N^{(1)}\cup\Omega_{\ell}^{(1)},\quad \cE_{N,\omega}^{(2)}:= \cE_{N,\omega}^{(1)}\cup \cE_{\ell,\omega}^{(1)}
\]
where $\ell =(\log N)^{4A}$ as in Corollary~\ref{cor:Qreduce}.
\end{defi}

\section{On the parametrization of the Dirichlet eigenfunctions}\label{sec:parameter}

In this section we describe the graphs of the Rellich
parametrization of the eigenvalues and eigenfunctions. It should be
thought of as a preliminary ingredient in the construction of
Sinai's function $\Lambda$, see Section~\ref{sec:sinai}. In this
section we shall assume for simplicity that
\[
L(\omega,E)>\gamma>0\quad \forall (\omega,E)\in\tor\times \IR
\]
This assumption can of course be localized to a rectangle
$(\omega',\omega'')\times (E',E'')$. In addition, we will fix $a>1$,
$c>0$ and consider $\tor_{c,a}$.
\begin{prop} \label{prop:11.Ej_sep}
Given $0<\delta<1$ there exist large constants
$N_0=N_0(\delta,V,\gamma,a,c)$, and  $A=A(\delta,V,\gamma,a,c)$
 such that for any $N
\ge N_0$, and any $ (\log N)^A=  \ell$ there exist
${\Omega}_{N}^{(1)}$, $ {\cE}_{N, \omega}^{(1)}$  as in
Definition~\ref{def:level1}
 such that for any
 $\omega \in\tor_{c,a}\setminus {\Omega}_{N}^{(1)}$ and all $x\in\tor$ one has
\begin{equation}
\label{eq:4.Ej_sep} \big | E_j^{(N)} (x, \omega) - E_k^{(N)} (x,
\omega) \big | > \exp(-\ell^{\delta})
\end{equation} for all
$j\ne k$ provided $E_j^{(N)}(x, \omega) \notin
{\cE}_{N,\omega}^{(1)}$.
\end{prop}
\begin{proof}
This follows from the proof of the eigenvalue separation from
Section~\ref{sec:separation} in combination with the reduction of the size of
the localization window which was obtained in
Corollary~\ref{cor:Qreduce}.
\end{proof}
We emphasize that $\delta A\gg1$. Thus, the separation achieved here
is always much smaller than~$N^{-1}$.

\begin{corollary}\label{cor:11.Ejsep} Using the notations of Definition~\ref{def:level1}, assume that $\omega\in\tor_{c,a}\setminus\Omega_N^{(2)}$ and  $E_{j_k}(x, \omega) \notin$
 $\tilde {\cE}_{N,\omega}^{(2)}$, $k=1,2$ for some $x\in\tor$. If for $j_1\ne j_2$
\[ |\nu_{j_1}^{(N)} (x, \omega)- \nu_{j_2}^{(N)} (x,
\omega)|\le \ell
\]
then \[ |E_{j_1}^{(N)} (x, \omega)- E_{j_2}^{(N)} (x,
\omega)|\ge \exp(-(\log \ell)^A)
\]
where $A$ is a large parameter as in the definition.
\end{corollary}
\begin{proof} Since $|\nu_{j_1}^{(N)} (x, \omega)- \nu_{j_2}^{(N)} (x,
\omega)|\le \ell$,
$$|\psi_{j_s}^{(N)} (x, \omega, n)|\le \exp(-\gamma |n-\nu_{j_1}^{(N)} (x, \omega)|/2)$$
for $|n-\nu_{j_1}^{(N)} (x, \omega)|\ge C\ell$, $s=1,2$, due to
Corollary~\ref{cor:Qreduce} where $C\gg 1$. Hence,
 \begin{equation}\label{eq:Elocal_split}\dist\big[E_{j_s}^{(N)} (x_0, \omega), \spec\big( H_{[\nu_{j_1}^{(N)} (x,
\omega)- C\ell, \nu_{j_1}^{(N)} (x, \omega) + C\ell]} (x,
\omega)\big) \big]\le \exp(-\gamma \ell)\end{equation} and the corollary follows from the previous
Proposition~\ref{prop:11.Ej_sep}. Indeed,
that proposition, applied to \[ H_{[\nu_{j_1}^{(N)} (x,
\omega)- C\ell, \nu_{j_1}^{(N)} (x, \omega) + C\ell]} (x,
\omega)\]  guarantees a splitting between the eigenvalues $E_{j_s}^{(N)}$
by an amount $\exp(-(\log \ell)^{A\delta})$.
Since $A\delta\gg1$ anyway, we simply set $\delta=1$ and we are done.
Note that \eqref{eq:Elocal_split} guarantees that the restriction
to the interval $[\nu_{j_1}^{(N)} (x,
\omega)- C\ell, \nu_{j_1}^{(N)} (x, \omega) + C\ell]$ only affects the estimate by an exponentially
small amount $e^{-\gamma\ell}$ which is acceptable.
\end{proof}

We will now investigate -- and single out --
 those portions of the graphs of the eigenvalues $E_j^{(N)}(x,\omega)$ which have
controlled slopes and controlled separations from the other eigenvalues.
 In order to obtain reasonable complexity bounds (i.e., to efficiently limit the number
 of these portions), we replace the function
$V(e(x)),\,\, x \in \tor$ by an approximating algebraic polynomial.
Since $V(z)$ is analytic in $\cA_{\rho_0}$ the Fourier coefficients
$\hat{v}(n)$ of $v(x) := V(e(x))$ satisfy
\[
|\hat{v}(n)| \le B_0 \exp\big(-\frac{\rho_0}{2} |n| \big)
\]
where $B_0 = \max \{|V(z)| : z \in \cA_{\rho_0/2} \}$. Replacing the
exponent $e(nx)$ by their Taylor polynomials leads to  the
following statement.

\begin{lemma}\label{lem:7.vapproxi} Given $0<\sigma<1$, and $T \ge 1$ there exists a polynomial
$\wt{V}(x), x \in \IR$ with real coefficients such that
\begin{enumerate}
\item[(1)] $ \max\limits_{|x|\le T} |V(e(x)) - \wt{V}(x)| \le \sigma$

\item[(2)] $ \deg \wt{V}(x) \le 4T(K_0 + \log\sigma^{-1}), \,\,\, K_0 = K_0(V)$
\end{enumerate}
\end{lemma}

\noindent For the remainder of this section, we set $\sigma := \sigma_N = \exp(-N^2)$, and $ T:=T_N =
N+1$ and we denote  the corresponding
polynomial from Lemma~\ref{lem:7.vapproxi} by $\wt{V}_N(x)$. Let $\tilde H_{[-N, N]}(x, \omega)$ be the
Schr\"odinger operator on $[-N, N]$ with Dirichlet boundary
conditions and with potential $\wt{V}_N(x+n\omega)$, $ - N \le n \le N$. Furthermore, let
\[
\wt{E}_1^{(N)}(x, \omega) < \wt{E}_2^{(N)}(x, \omega) < \dots <
\wt{E}_{2N+1}^{(N)}(x,\omega)
\]
be the eigenvalues of $\tilde H_{[-N, N]}(x, \omega)$.

\begin{lemma}\label{lem:7.spapproxi} With the previous notation, one has the following estimates:
\begin{enumerate}
\item[(1)]  $\| H_{[-N, N]}(x, \omega) - \wt{H}_{[-N, N]}(x, \omega) \| \le \sigma$,
$\quad\forall\; x,\omega \in [-1, 1]$

\item[(2)]  For any $x, \omega \in [-1, 1]$, one has $1\le j\le 2N+1$ there exists $1\le j_1 \le 2N+1$
such that
\[
|E_j^{(N)}(x,\omega) - \wt{E}_{j_1}^{(N)}(x,\omega) | \le \sigma,
\]
and vice versa.

\item[(3)]  If for some $x, \omega$
\[
\min\limits_{1\le j< 2N} \bigl( E_{j+1}^{(N)}(x, \omega) -
E_j^{(N)}(x, \omega) \bigr) \ge \exp(-N^\delta),
\]
then one has $j_1=j$ and
\begin{align*}
|E_j^{(N)}(x, \omega) - \wt{E}_j^{(N)}(x,\omega)  | &\le \sigma, \,\,
j \in [1, 2N+1] \\
\min\limits_{1\le j\le 2N} (\wt{E}_{j+1}^{(N)}(x,\omega) -
\wt{E}_j^{(N)}(x,\omega) ) &\ge \frac{1}{2} \exp(-N^\delta)
\end{align*}
\end{enumerate}
\end{lemma}

\begin{proof}  (1) follows from the estimate (1) of Lemma~\ref{lem:7.vapproxi}.
Properties (2) and (3) now follow from assertion (1) due to basic facts about perturbations of self-adjoint operators.
\end{proof}

Set
\[
\tilde {f}_N(e(x),\omega,E) := \det(\tilde H_{[-N, N]}(x,\omega) -E)
\]
Note that due to property (1) of Lemma~\ref{lem:7.vapproxi},
$\tilde {f}_N(x,\omega,E)$ is a polynomial in $x, \omega, E$ with
\begin{equation}\label{eq:7.fdeg}
\deg\tilde {f}_N \le N^4, \,\,\, \text{\ \ for \ \ } N \ge N_0(V, c,
a, \gamma)
\end{equation}
Given an arbitrary interval $[\uE, \oE]$ set
\begin{equation}\label{eq:btilde}\begin{aligned}
\cB_{N,\omega}(\uE, \oE) &= \{ x \in [0,1] : \spec \big(
H_{[-N,N]}(x,\omega) \big)\cap (\uE, \oE) \neq \emptyset\}\\
\wt{\cB}_{N,\omega}(\uE, \oE) &= \{ x \in [0,1] : \spec \big(
\wt{H}_{[-N,N]}(x,\omega)\big) \cap (\uE, \oE) \neq \emptyset\}
\end{aligned}
\end{equation}
These sets have a number of simple properties:

\begin{lemma}\label{lem:7.algewgner}
The sets introduced in \eqref{eq:btilde} satisfy the following properties:
\begin{enumerate}
\item[(1)]
$ \wt{\cB}_{N,\omega}(\uE, \oE)  \subset \cB_{N,\omega}(\uE - \sigma, \oE+\sigma)
 \subset \wt{\cB}_{N,\omega}(\uE
- 2\sigma, \oE+2\sigma)$

\item[(2)]  $\mes \wt{\cB}_{N,\omega}(\uE,\oE) \le \exp(- H/(\log N)^C)$, where $e^{-H}= \min(1/2,\oE-\uE)$

\item[(3)] $\wt{\cB}_{N,\omega}(\uE,\oE)$ consists of a union of at most $O(K_0 N^4)$ closed intervals, where $K_0$ is as in Lemma~\ref{lem:7.vapproxi}
\end{enumerate}
\end{lemma}

\begin{proof} Property (1) is due to fact (2) of Lemma~\ref{lem:7.spapproxi}.
 (2) follows from (1) due to Lemma~\ref{lem:wegner} (the analogue of Wegner's estimate).
 The edges of the maximal intervals in the complement of
  $\wt{\cB}_{N,\omega}$ are the roots of the equations
\[
\tilde {f}_N(x, \omega, \uE) = 0 \text{\ \ or \ \ } \tilde {f}_N(x,
\omega, \oE) = 0
\]
Since $\deg\tilde {f}_N \les K_0N^4$, property (3) follows.
\end{proof}

For the rest of this section,  we fix $\omega \in \tor_{c,a} \backslash
\Omega_N$, where $\Omega_N, \cE_{N,\omega} $ are the sets from
Corollary~\ref{cor:sep}  and
Proposition~\ref{prop:4.Ej_sep}..

\begin{corollary}\label{cor:7.spapproxisp} There exist intervals $[\xi'_k, \xi''_k]
 \subset [0,1]$, $k=1,\ldots,k_0$ such that the following conditions hold:
\begin{enumerate}
\item[(1)]  $|E_j^{(N)}(x,\omega) - \wt{E}_j^{(N)}(x, \omega) | \le \sigma$ for any
$x \in \mybigcup_k [\xi'_k, \xi''_k]$ and any $j=1, \ldots, 2N+1$

\item[(2)]
\begin{align*}\min\limits_{1\le j\le 2N} \bigl({E}_{j+1}^{(N)}(x,\omega) - {E}_j^{(N)}(x,\omega) \bigr)
\ge \frac{1}{4} \exp(-N^\delta)\\
 \min\limits_{1\le j\le 2N} \bigl(\wt{E}_{j+1}^{(N)}(x,\omega) - \wt{E}_j^{(N)}(x,\omega) \bigr)
 \ge \frac{1}{4} \exp(-N^\delta)
\end{align*}
  for any $x\in \mybigcup_k[\xi'_k, \xi''_k]$

\item[(3)]  $\mes([0,1] \backslash \mybigcup_k[\xi'_k, \xi''_k] ) \le \exp(-\frac{1}{2}(\log N)^B)$

\item[(4)]  $k_0 \les N^7$
\end{enumerate}
\end{corollary}
\begin{proof}
Set
\begin{align*}
\cE^{(+)}_{N,\omega} &:= \big\{ E\in\IR \: :\: \dist(E, \cE_{N,\omega}) < \sigma_N
\big\} \\
\wt{\cB}_{N,\omega}&:= \Big\{ x \in [0,1] \::\: \spec \big( \tilde {H}_{[-N,N]}(x,
\omega)\big) \bigcap \cE^{(+)}_{N,\omega} \neq \emptyset \Big\}
\end{align*}
Note that
\[
\wt{\cB}_{N,\omega} = \bigcup_{[\uE,\oE]\subset \cE^{(+)}_{N,\omega}}  \wt{\cB}_{N,\omega}(\uE, \oE)
\]
where the union here runs over the maximal subintervals of $\cE^{(+)}_{N,\omega}$.
One the one hand, due to the properties of $\cE_{N,\omega}$, the set $\cE^{(+)}_{N,\omega}$ can be
covered by at most $N^3$ intervals $[E'_k, E''_k]$ with
\[
\sum\limits_k (E_k'' - E_k') \les \exp(-(\log N)^B)
\]
On the other hand, by Lemma~\ref{lem:7.algewgner}, it follows that for each such $k$, the set
$\wt{\cB}_{N,\omega}(\uE, \oE)$ is the union of at most $\les N^4$ intervals on the $x$-axis.
Hence, $\wt{\cB}_{N,\omega}$ is the union of at most $\les N^7$ intervals. The maximal intervals in the complement
of $\wt{\cB}_{N,\omega}$ are now defined to be $[\xi'_k,\xi''_k]$ with $1\le k\le k_0$.
The corollary follows by combining Lemma~\ref{lem:7.spapproxi}
and Lemma~\ref{lem:7.algewgner}.
\end{proof}

We also record the following  standard fact about
the perturbations of analytic matrix functions $M(z)$ which take values in the Hermitian matrices. We state it
for the case of $H_{[-N,
N]}(x,\omega)$ with $N$ large.

\begin{lemma}\label{lem:7.analytconti} Assume that for some $x_0, \omega_0, j_0$
\begin{equation}\label{eq:EjEj0}
\min\limits_{j \neq j_0} | E_j^{(N)}(x_0, \omega_0) -
E_{j_0}^{(N)}(x_0, \omega_0) | \ge \sigma^{(0)} > 0
\end{equation}
Then there exists an analytic function $E_{j_0}^{(N)}(z, \omega),
(z, \omega) \in \cD(x_0, r_0) \times \cD(\omega_0, r_0), r_0 =
\sigma^{(0)}_2 / N^2$ such that
\[
\spec \big( H_{[-N,N]}(z, \omega) \big) \cap \cD(E_{j_0}^{(N)},
{\sigma^{(0)}}/{2}) = \{ E_{j_0}^{(N)}(z, \omega)\}
\]
for any $(z,\omega) \in \cD(x_0, r_0) \times \cD(\omega_0, r_0)$.
Furthermore, suppose \eqref{eq:EjEj0} holds for all $j_0$. If $\frac{1}{2}\sigma^{(0)} > \sigma_N$, then for each $j$
there
exists an analytic function $\wt{E}_j^{(N)}(z,\omega), (z,\omega)
\in \cD(x_0, r_0) \times \cD(\omega_0, r_0)$ such that
\[
 \spec (\wt{H}_{[-N, N]}(z,\omega)) \cap \cD(\wt{E}_j^{(N)}(x_0, \omega_0), {\sigma^{(0)}}/{2}) =
\{ \wt{E}_{j}^{(N)}(z,\omega)\}
\]
for any $(z,\omega) \in \cD(x_0, r_0) \times \cD(\omega_0, r_0)$.
Finally, for each $j$,
\begin{align*}
|E_j^{(N)}(z,\omega) - \wt{E}_j^{(N)}(z, \omega)| \le 2\sigma_N\\
|\partial^\alpha E_j^{(N)}(z,\omega) -
\partial^\alpha\wt{E}_j(z,\omega)| \le
2\alpha!(\frac{r_0}{2})^{-|\alpha|}\sigma_N
\end{align*}
for any $(z,\omega) \in \cD(x_0, r_{0}/2) \times \cD(\omega_0,
r_0/2)$.
\end{lemma}

Combining Lemma~\ref{lem:7.analytconti} with
Corollary~\ref{cor:7.spapproxisp} one obtains the following.

\begin{corollary}\label{cor:7.dirapproxi} Using the notations of Corollary~\ref{cor:7.spapproxisp} one has
\[
|\partial_xE_j^{(N)}(x,\omega)
-\partial_x\wt{E}_j^{(N)}(x,\omega)|\le \sqrt{\sigma_N}
\]
for any $x \in \mybigcup_k[\xi'_k, \xi''_k]$ and any
$j=1,\ldots,2N+1$.
\end{corollary}

\noindent Note that each $\wt{E}_j^{(N)}(\cdot, \omega)$ is an
algebraic function (see Appendix~A). That implies the following
statement.

\begin{lemma}\label{lem:7.ealgebraic} For each $j=1, \ldots, 2N+1$ and each $\tau > 0$
there exist disjoint intervals $[\eta'_{j,m}(\tau),
\eta''_{j,m}(\tau)]$, $ m=1,2,\ldots,m_0$, $ m_0 \le N^9$
such that
\begin{enumerate}
\item[(i)] $|\partial_x\wt{E}_j^{(N)}(x,\omega)|>\tau$
for any $x\in [-1,1]\backslash\mybigcup_m[\eta'_{j,m}(\tau), \eta''_{j,m}(\tau)]$

\item[(ii)] $|\partial_x\wt{E}_j^{(N)}(x,\omega)|\le\tau$ for any $x\in\mybigcup_m[\eta'_{j,m}(\tau), \eta''_{j,m}(\tau)]$

\end{enumerate}
\end{lemma}
\begin{proof}
The degree of $\tilde f_N$ is $\les N^4$, see~\eqref{eq:7.fdeg}. Since $\tilde f_N(x,\omega,\wt{E}_j^{(N)}(x, \omega))=0$,
we see that
\[
\partial_x \wt{E}_j^{(N)}(x, \omega) = - \partial_E \tilde f_N(x,\omega,E)/ \partial_x \tilde f_N(x,\omega,E) \Big|_{E= \wt{E}_j^{(N)}(x, \omega)}
\]
Hence, by Bezout's theorem in the appendix it follows that the equation $\partial_x \wt{E}_j^{(N)}(x, \omega)=\pm\tau$ has at most $2N^8$ solutions,
whence the result.
\end{proof}

\begin{lemma}\label{lem:7.span}
 Using the notations of Lemma~\ref{lem:7.ealgebraic}  define
\begin{align*}
\uE_j(\tau,m) &:=\min\{\wt{E}_j^{(N)}(x,\omega) :
x\in[\eta'_{j,m}(\tau), \eta''_{j,m}(\tau)]\}\\
\oE_j(\tau,m) &:=\max\{\wt{E}_j^{(N)}(x,\omega) :
x\in[\eta'_{j,m}(\tau), \eta''_{j,m}(\tau)]\}
\end{align*}
$m=1,\dots,m_0$. Then
 one has
\begin{align}\label{eq:7.espan}
\oE_j(\tau,m)-\uE_j(\tau,m) &\le 2\tau
\\\label{eq:7.xspan}
\eta''_{j,m}(\tau) - \eta'_{j,m}(\tau) &\le \exp(- \log\tau^{-1} /
(\log N)^C )
\end{align}
\end{lemma}

\begin{proof} The estimate~\eqref{eq:7.espan} follows from part (ii) of Lemma~\ref{lem:7.ealgebraic}.
 The bound~\eqref{eq:7.xspan} follows from Relation~\eqref{eq:7.espan}
 due to Lemma~\ref{lem:wegner} (which is the analogue of Wegner's estimate).
\end{proof}

Now we obtain the main result of this section. It allows one to control the graphs of
the eigenvalues in terms of slopes and separation properties up to the removal of certain
sets.

\begin{prop}\label{prop:11.4rellich}
Assume that $L(\omega, E) \ge \gamma > 0$ for all $\omega \in
(\omega',\omega'')$ and all $E \in (E', E'')$. Given $\delta \ll 1\ll A$
there exists $N_0=N_0(V, c, a, \gamma, \delta, A)
$,  such that for any $N \ge
N_0$, and
arbitrary $\exp(-N^{\delta}) < \tau < \exp(-(\log N)^A )$ there exist
$\cB_{N,\omega}$, $\cB'_{N,\omega}$
 $\cE_{N, \omega}(\tau),\cB_{N,\omega}(\tau)$ such
that, with $\Omega_N^{(1)}$ and $\cE_{N,\omega}^{(1)}$ as in Definition~\ref{def:level1},
\begin{itemize}
\item[(1)]
\begin{equation}\label{eq:11.level1_est}\begin{aligned}
\mes ({\cB}_{N, \omega}) &\leq \exp(-(\log N)^A),\;\; \compl (
{\cB}_{N, \omega}) < N^9\\
\mes (\cB'_{N, \omega}) &\leq \exp(-(\log \log N)^A),\;\; \compl (
\cB'_{N, \omega}) < \exp((\log \log N)^{A/2})
\end{aligned}\end{equation}
as well as
\begin{equation}\label{eq:11.leveltau_est}\begin{aligned}
  \mes(\cE_{N, \omega}(\tau)) \leq \exp(-(\log \tau^{-1})(\log N)^{-C} ), \,\, \compl (
  \cE_{N, \omega}(\tau) ) \leq N^{C}\\ \mes (\cB_{N, \omega}(\tau))\leq \exp(-(\log \tau^{-1})(\log N)^{-C} ), \,\,
\compl (\cB_{N, \omega}(\tau)) \leq N^{C}
\end{aligned}\end{equation}
\item[(2)] If for some $\omega\in \tor_{c,a}\setminus \Omega^{(1)}_N$, and $x \in \tor$ one has
\[
\spec \big(H_{[-N, N]}(x, \omega) \big)\cap \bigl( (E', E'') \setminus ({\cE}^{(1)}_{N, \omega}\cup \cE_{N,
\omega}(\tau)) \bigr) \neq \emptyset,
\]
then $x \in \tor \backslash ({\cB}_{N, \omega}\cup \cB'_{N, \omega}\cup \cB_{N, \omega}(\tau))$

\item[(3)]For  $\omega\in \tor_{c,a}\setminus \Omega_N$ one has
\begin{equation}
|E_{j_1}^{(N)}(x, \omega) - E_{j_2}^{(N)}(x,\omega) | \ge \exp(-N^{\delta})
\label{eq:sep_weak}
\end{equation}

 for any $j_1\neq j_2$, provided  $E_{j_s}^{(N)}(x,\omega)\in  (E', E'')\setminus
{\cE}_{N, \omega}$, $s=1,2$. For  $\omega\in \tor_{c,a}\setminus \Omega^{(1)}_N$ one has
\begin{equation}
 \label{eq:sep_strong}
|E_{j_1}^{(N)}(x, \omega) - E_{j_2}^{(N)}(x,\omega) | \ge \exp(-(\log N)^A)
\end{equation}
 for any $j_1\neq j_2$, provided  $E_{j_s}^{(N)}(x,\omega)\in (E', E'')\setminus
\cE^{(1)}_{N, \omega}$, $s=1,2$
\item[(4)] For $\omega\in \tor_{c,a}\setminus \Omega^{(1)}_N$  one has $|\partial_xE_j^{(N)}(x,\omega) | \ge \tau$ ,
provided  \[E_j^{(N)}(x,\omega) \in  (E', E'') \setminus
({\cE}^{(1)}_{N, \omega}\cup \cE_{N,
\omega}(\tau)) \]
\end{itemize}
  \end{prop}

\begin{proof} In this proof we do not distinguish between $\tilde H_{[-N,N]}$ and $H_{[-N,N]}$. This is justified
by the results of this section which demonstrate that the small ``fattening'' parameter $\sigma_N$ can be ignored.
 Set
\begin{align*} \cB_{N,\omega} &:=\{x\in \tor : \spec \big( H_N(x,\omega) \big)\cap \cE_{N,\omega}\neq \emptyset \} \\
 \cB'_{N,\omega} &:=\{x\in \tor :\spec \big( H_N(x,\omega) \big)\cap \cE'_{N,\omega}\neq \emptyset\}
 \end{align*}
 where $\cE_{N,\omega}^{(1)}=\cE_{N,\omega}\cup\cE_{N,\omega}'$ as in Definition~\ref{def:level1}.
 Then conditions \eqref{eq:11.level1_est} hold for  $\cB_{N,\omega}$,
 $\cB'_{N,\omega}$ due to Lemma~\ref{lem:wegner} (the analogue of Wegner's estimate).
Using the notations of Lemma~\ref{lem:7.span}
set \[\cB_{N,\omega}(\tau) : =
\mybigcup_{j,m}(\oE_j(2\tau,m), \uE_j(2\tau,m)),\quad \cB_{N,\omega}(\tau) : =
\mybigcup_{j,m}(\eta'_{j,m}(2\tau), \eta''_{j,m}(2\tau)).\] Then
conditions \eqref{eq:11.leveltau_est} follow from
Lemma~\ref{lem:7.span}. Property $(2)$ holds due to the definition of the sets involved,
whereas $(3)$ is due to Proposition~\ref{prop:11.Ej_sep}.
Finally, property $(4)$ is due to Lemma~\ref{lem:7.ealgebraic}.
\end{proof}

\begin{remark}\label{rem:7.lolaliz} For future reference we note that in addition to
$(1)$--$(4)$ of Proposition~\ref{prop:11.4rellich},  the following
property holds due to Corollary~\ref{cor:Qreduce}:
 For any $\omega\in\tor_{c,a}\setminus \Omega_N^{(1)}$ and any $x\in\tor$, any  $\ell^2$-normalized eigenfunction
$\psi_j^{(N)}(x,\omega)$ of $H_{[-N,N]}(x,\omega)$ with
associated eigenvalue $E_j^{(N)}(x,\omega)\in\IR\setminus \cE^{(1)}_{N,\omega}$ satisfies
\[
|\psi_j^{(N)}(x,\omega)(n)|\le C\exp(-\gamma \dist(n,\Lambda_j)/2)
\]
for all $n\in[-N,N]$ where $\Lambda_j:=[\nu_j^{(N)}(x,\omega)-\ell,\nu_j^{(N)}(x,\omega)+\ell]\cap[-N,N]$
 with $\ell=(\log N)^{4A}$.
\end{remark}

\section{Segments of eigenvalue parametrizations and their translations}\label{sec:segments}

In this section we discuss the segments of functions
$E_j^{(N)}(\cdot,\omega)$ and establish a self-similar structure of
these functions with regard to the shift by~$\omega$. The later
property is based on the finite volume localization of
Section~\ref{sec:anderson}, and should be considered as a precursor to the
co-variant property of Sinai's function from Section~\ref{sec:sinai}.
Let $\cB_{N, \omega}, \cB_{N, \omega}(\tau)$ be  as in
 Proposition~\ref{prop:11.4rellich}. Set $\tau_N=\exp(-(\log
N)^{B})$, $B\gg A$, and
\[\bar\cB_{N, \omega} = \cB_{N, \omega}^{(1)} \cup \cB_{N, \omega}(\tau_N),\qquad \cB_{N, \omega}^{(1)} := \cB_{N,\omega}\cup \cB_{N,\omega}' \]
There exist intervals $[\xi'_k, \xi''_k], \xi''_k < \xi'_{k+1},
\,\, k = 1, \dots, k_0, \,\, k_0 \le N^{C}$ such that
\[
\tor \backslash \bar\cB_{N, \omega} = \mybigcup_{1\le k \le
k_0}[\xi'_k, \xi''_k]
\]
Set
\[ \cB''_{N,\omega} := \tor \backslash \bigcup\bigl\{[\xi'_k, \xi''_k] \,:\, \xi''_k - \xi'_k > \exp(-(\log N)^{2A}) \bigr\}\]
To point here is that we add all very short  ``good intervals'' into the bad set. This does not increase the measure of the bad set too much. Indeed,
\[ \mes (\cB''_{N, \omega}) \leq \exp(-(\log N)^A), \,\,\, \compl ( \cB''_{N, \omega}) \leq N^{C} \]
We denote by $[\ux_k, \ox_k], \,\,\, k = 1, 2, \dots, k_1, \,\,\,
\ox_k < \ux_{k+1}$ the maximal intervals of $\tor \backslash
\cB''_{N, \omega}$. Recall that due to
Proposition~\ref{prop:11.4rellich}
\begin{equation}\label{eq:EjtauN}
|\partial_xE_j^{(N)}(x, \omega) | > \tau_N \text{\ \ for any \ \ } x
\in \mybigcup_{1\le k \le k_1} [\ux_k, \ox_k]
\end{equation}
We summarize the properties of the intervals $[\ux_k, \ox_k]$
(including those mentioned in Remark~\ref{rem:7.lolaliz}) in the
following lemma.

\begin{lemma}\label{lem:8.maxintervals} Let $\omega \in \tor_{c, a} \backslash \Omega_N^{(1)}$. There exists intervals
$[\ux_k, \ox_k]$ $($depending on $\omega)$, $k=1,\ldots,k_1$,
$\ux_k < \ox_{k+1}$ such that
\begin{enumerate}
\item [(1)]  $|\partial_x E_j^{(N)}(x, \omega)|\ge \exp(-(\log N)^{B})$ for any
$x \in \mybigcup_{1 \le k \le k_1} [\ux_k, \ox_k]$ and any $1\le j\le 2N+1$

\item[(2)]  for any $x \in \mybigcup_{1\le k\le k_1}[\ux_k, \ox_k]$ and any $j$ there exists
$\nu_j^{(N)}(x, \omega) \in [-N, N]$ such that
\[ |\psi_j^{(N)}(x, \omega, n) | \leq \exp\big(-\frac{\gamma}{2}|n-\nu_j^{(N)}(x, \omega)| \big) \]
provided $|n-\nu_j^{(N)}(x, \omega)|>N^{\frac{1}{2}}$

\item [(3)] $\mes \bigl(\tor \backslash \mybigcup_{1\le k\le k_1} [\ux_k, \ox_x] \bigr) \le 3\exp(-(\log N)^A)$

\item[(4)] $k_1 \le N^C$

\item[(5)]  $\ox_k - \ux_k \ge \exp(-(\log N)^{2A}), \,\,\, k=1, 2, \dots, k_1$
\end{enumerate}
\end{lemma}

\begin{defi}\label{def:8.1} We refer to each triplet
 $\{E_j^{(N)}(x, \omega), \ux_k, \ox_k \}, \,\,\, j=1,\dots,2N+1, \,\,\, k=1,\dots,k_1$
 as a segment of Rellich's parametrization, or just a segment of $E_j^{(N)}$. If $\partial_xE_j^{(N)}(x,\omega)>0$, respectively $\partial_xE_j^{(N)}(x,\omega) <0$, 
  for any $x\in[\ux_k, \ox_k]$ then $\{E_j^{(N)}(x,\omega),\ux_k, \ox_k\}$ is called
   a positive--slope, respectively negative--slope, segment. We also refer to these segments
as $I$-segments provided the open interval $I\subset\IR$ is the range of the segment, i.e.,
$
 I = \{E_j^{(N)}(x, \omega)\::\: \ux_k < x <  \ox_k \}
$.
\end{defi}

\begin{remark}\label{rem:12.analytcontinu} Let $\{E_j^{(N)}(x, \omega), \ux_k, \ox_k
\}$ be a segment. Recall that due to
Proposition~\ref{prop:11.4rellich} one has in particular
\[
|E_j^{(N)}(x, \omega)-E_{j_1}^{(N)}(x, \omega)|\ge \exp(-N^{\delta})
\]
for any $j_1\neq j$. Since $V$ is analytic one infers from this and  standard perturbation
theory of Hermitian matrices  that the function
$E_j^{(N)}(\cdot, \cdot)$ admits an analytic continuation to the
polydisk
\[\cP:=\{(z,w)\in \IC^2 \;:\; |z-e(x)| < r_0, |w-\omega| < r_0\}
\]
where $r_0:=C(V)^{-1}\exp(-N^{\delta})$. Moreover,
\[
\sup_{\cP} |E_j^{(N)}(z, w)| \le C(V)
\]
Note that the same bound holds with the stronger estimate~\eqref{eq:sep_strong} instead of the weaker one~\eqref{eq:sep_weak}, albeit
with $\delta>0$ arbitrarily small (simply because of the stronger bound). 
For many applications below the weaker bound suffices, and we have chosen to use it for the remainder of the section. The reader
will have no difficulty replacing it with the stronger one~\eqref{eq:sep_strong} whenever the need arises. 
\end{remark}

We now turn to translations of the segments. For this purpose we are going to
impose the  condition
\begin{equation}\label{eq:8.regseg}
-N+N^{\frac{1}{2}} \le \nu_j^{(N)}(x,\omega) \le N-N^{\frac{1}{2}}
\end{equation}
which guarantees that the window of localization is separated from the boundary of the box $[-N,N]$.

\begin{lemma}\label{lem:4.transl} Using the notations of Proposition~\ref{prop:4.Ej_sep} assume that
\[ \dist\bigl(E_j^{(N)}(x,\omega), \cE_{N,\omega}\bigr)
> 2\exp\bigl(-N^{\delta}\bigr),\quad -N + N^{1/2} < \nu_j^{(N)}\xo < N - N^{1/2}\]
Then for any $k$ such that $-N + N^{1/2}/2 < \nu_j^{(N)}\xo + k < N
- N^{1/2}/2$ there exists a unique $E^{(N)}_{j_k}(x + k\omega,
\omega)\in \spec \big( H_{[-N, N]}(x + k\omega, \omega) \big)$ such that
\begin{align}
& \bigl|E_j^{(N)}\xo - E^{(N)}_{j_k}(x + k\omega, \omega)\bigr|< \exp\bigl(-\gamma_2 N^{1/2}\bigr)\
,\label{eq:4.14}\\
& \bigl|\partial_x E_j^{(N)}\xo -
\partial_x E^{(N)}_{j_k}(x + k\omega, \omega)\bigr|  < \exp\bigl(-\gamma_2 N^{1/2}\bigr)\
,\label{eq:4.14'}
\end{align}
as well as
\begin{align}
& E^{(N)}_{j_k} (x + k\omega, \omega) \notin \cE_{N,\omega}\ ,\label{eq:4.15}\\
& \bigl|\nu^{(N)}_{j_k} (x + k\omega, \omega) - \bigl(\nu_j^{(N)} \xo + k\bigr)\bigr| \le N^{1/2}/4\ ,\label{eq:4.16}\\
& -N + N^{1/2}/4 < \nu^{(N)}_{j_k}(x + k\omega, \omega) < N - N^{1/2}/4\ ,\label{eq:4.17}\\
&\sum_{|m+k - \nu^{(N)}_j\xo| \le N^{1/2}/4} \bigl|\psi^{(N)}_{j_k}
(x + k\omega, m) - \psi_j^{(N)}(x, m + k) \bigr|^2 <
\exp\bigl(-\gamma_3 N^{1/2}\bigr)\label{eq:4.18}
\end{align}
where $\gamma_t = 2^{-t+1}\gamma_1$.
\end{lemma}

\begin{proof} Note that
\begin{equation}\label{eq:4.19}
\begin{aligned}
 H_{[-N, N]}(x + k\omega, \omega) \bigl(\psi^{(N)}_j (x, \omega, \cdot + k)\bigr)(m) & =
 H_{[-N, N]}\xo \bigl(\psi^{(N)}_j(x,\omega .)\bigr)(m +k)\\
& = E_j^{(N)}\xo \psi^{(N)}_j (x, \omega, k+m)
\end{aligned}
\end{equation}
provided $-N < m+k < N$ and $-N<m<N$.  Recall also that
$\bigl|\psi^{(N)}_j(x, \omega, \pm N)\bigr| \le \exp\bigl(-\gamma_3
N^{1/2}\bigr)$, since $-N + N^{1/2} < \nu_j^{(N)} \xo < N - N^{1/2}$
Hence
\begin{equation}\label{eq:4.20}
\Big \| \bigl(H_{[-N, N]}(x + k\omega, \omega) - E_j^{(N)}\xo\bigr)
\psi_j^{(N)} (x, \omega, \cdot + k) \Big \| < \exp\bigl(-\gamma_4
N^{1/2}\bigr)\ .
\end{equation}
Therefore, there exists $E_{j_k}^{(N)} (x + k\omega, \omega) \in
\bigl(E_j^{(N)}\xo - \exp\bigl(-\gamma_5 N^{1/2}\bigr), E_j^{(N)}\xo
+ \exp\bigl(-\gamma_5 N^{1/2}\bigr) \bigr)$.  Moreover, due to our
assumptions on $E_j^{(N)}\xo$, one has $E_{j_k}^{(N)} (x + k\omega,
\omega) \notin \cE_{N,\omega}$.  Hence,
\begin{equation}\label{eq:4.21}
\bigl| E_{j_k}^{(N)} (x + k\omega, \omega) - E_{j'}^{(N)}(x +
k\omega, \omega) \bigr| > \exp(- N^{\delta})
\end{equation}
for any $j' \ne j_k$.  Relations \eqref{eq:4.19}--\eqref{eq:4.20}
combined imply \eqref{eq:4.14}. Relations \eqref{eq:4.15},
\eqref{eq:4.16} follow from \eqref{eq:4.14}. The
estimate~\eqref{eq:4.18} follows from \eqref{eq:4.20}
and~\eqref{eq:4.21} via the spectral theorem for Hermitian matrices.
Finally, \eqref{eq:4.14'} follows from the well-known ``Feynman formula'' (or first order eigenvalue
perturbation formula)
\[ \partial_x E_j^{(N)}(x,\omega) = \sum_{\ell=-N}^N V'(x+\ell
\omega) \big|\psi_j^{(N)}(x,\omega,\ell)\big|^2\] and the preceding
estimates.
\end{proof}

We now illustrate how to relate the localized eigenfunctions of
consecutive scales.

\begin{lemma}\label{lem:4.3} Using the notations of Proposition~\ref{prop:4.Ej_sep}
 assume that $\omega \in \tor_{c,a} \setminus \bigl(\Omega_N \cup \Omega_{N'}\bigr)$,
 where $N' \asymp \exp\bigl((\log\log N)^{C_1}\bigr)$, $C_1 \gg C$, and with
  $Q = \exp\bigl((\log\log N)^C\bigr)$.
If \[E_j^{(N)}(x,\omega) \notin \cE_{N,\omega},\qquad
\dist\bigl(E_j^{(N)} (x,\omega), \cE_{N',\omega}\bigr) >
\exp\bigl(-(N')^{1/2}\bigr),\] then there exists $\nu \in \IZ$,
$|\nu - \nu_j^{(N)} (x,\omega)| \le Q$ and
\begin{equation*}
 E_{j'}^{(N')}(x + \nu\omega, \omega) \in
\bigl(E_j^{(N)}(x,\omega) - \exp(-\gamma_1 N'), E_j^{(N)}(x,\omega)
+ \exp(-\gamma_1 N')\bigr)\ ,
\end{equation*}
where $\gamma_1 = c\gamma$, $\gamma = \inf L(E, \omega)$.  Moreover,
the corresponding normalized eigenfunctions \[\psi_j^{(N)}(x,\omega,
k),\qquad \psi_{j'}^{(N')}(x + \nu\omega, \omega, k - \nu)\] satisfy
\begin{equation}
\label{eq:4.9}
 \sum_{k \in [\nu - N', \nu + N']} \bigl|\psi_j^{(N)}
 (x, \omega, k) - \psi_{j'}^{(N')}(x + \nu\omega, \omega, k -\nu)\bigr|^2
 \le \exp(-\gamma_1 N')\ .
\end{equation}
\end{lemma}
\begin{proof} Assume first $-N + N' < \nu_j^{(N)}(x, \omega) < N - N'$.
Then with $\nu = \nu_j^{(N)} (x, \omega)$ one has:
\begin{align}
& \label{eq:4.10}\bigl\|\bigl(H_{[\nu - N', \nu + N']} (x, \omega) -
E_j^{(N)} (x, \omega)\bigr) \psi_j^{(N)}(x, \omega, \cdot)
\bigr\| \le \exp(-\gamma N'/4)\ ,\\
& \label{eq:4.11} 1 - \sum_{k \in [\nu - N', \nu + N']}
\bigl|\psi_j^{(N)}(x, \omega, k)\bigr|^2 < \exp(-\gamma N'/4)
\end{align}
due to Proposition~\ref{prop:4.Ej_sep}.  Hence, there exists
\[E_{j'}^{(N')} (x + \nu\omega, \omega) \in \bigl(E_j^{(N)} (x,
\omega) - \exp(-\gamma_1 N'), E_j^{(N)} (x,\omega) + \exp(-\gamma_1
N')\bigr).\]  Moreover, due to assumptions on $E_j^{(N)}(x,
\omega)$, one has $E_{j'}^{(N')}(x + \nu\omega, \omega) \notin
\cE_{N',\omega}$. Hence,
\begin{equation}\label{eq:4.12}
\bigl|E_{j'}^{(N')}(x + \nu\omega, \omega) - E_k^{(N')} (x +
\nu\omega, \omega) \bigr| > \exp\bigl(-(N')^\delta\bigr)
\end{equation}
for any $k \ne j'$.  Then (\ref{eq:4.10})--(\ref{eq:4.12}) combined
imply (\ref{eq:4.9}) (expand in the orthonormal basis
$\{\psi^{(N')}_k\}_k$). If $\nu_j^{(N)}(x, \omega) \le - N + N'$
(resp., $\nu_j^{(N)}(x, \omega) \ge N - N'$), then
(\ref{eq:4.10})--(\ref{eq:4.12}) are valid with $\nu_j^{(N)}(x,
\omega) = - N + N'$ (resp., with $\nu_j^{(N)}(x, \omega) = N
- N'$).
\end{proof}

Recall that $\nu_j^{(N)}(x,\omega)$ is stable under perturbations of $H_{[-N,N]}(x,\omega)$  of
magnitude $\exp(-(\log N)^C)$. Since some of the interval $[\ux_k,
\ox_k]$ are definitely of a larger size,
condition~\eqref{eq:8.regseg} can hold on some part of $[\ux_k,
\ox_k]$ and fail on another one. For that reason we consider also
all the triplets
\[ \big\{ E_j^{(N)}(x), \ux, \ox \big\} \]
with $[\ux, \ox]$ being an arbitrary subsegment of some $[\ux_1,
\ox_1], \,\, k=1,2,\dots,k_1$.

\begin{defi}\label{def:8.regsegdef}
 A segment $\{ E_j^{(N)}(x,\omega), \ux, \ox \}$
is called {\em regular} if the condition
\begin{equation}\label{eq:zentral} -N+N^{\frac{1}{2}} \le \nu_j^{(N)}(x,\omega) \le N-N^{\frac{1}{2}}
\end{equation}
holds for every $x\in[\ux,\ox]$.
\end{defi}

The condition \eqref{eq:zentral}, which requires the associated eigenfunction to be properly localized
inside of the base interval $[-N,N]$, will play a central role in this paper. It ensures {\em stability}
of resonances as one passes from one scale to the next.

\section{Formation of Regular Spectral Segments}\label{sec:formation}

\noindent The section is devoted to the comparison of the zeros in energy of each entry
of $M_N(e(x),\omega,E_0)$. Recall that these entries are are determinants $f_{[a,N-b]}(e(x),\omega,E)$ where $a,b=1,2$.
The main theme will be to single out a ``good case'' characterized by each determinant having
a zero very close to~$E_0$ (we refer to $E_0$ as a {\em unconditional spectral value at scale} $N$ in that case, see Definition~\ref{def:9.uncondint} below).
The significance of this idea lies with the induction in the scale; indeed, in
passing from $[1,N]$ to $[-\bar N,\bar N]$ with $\bar N= N^C$ we shall see that if $E_0$ is an unconditional
spectral value at scale $N$, then it remains very close to the spectrum at the larger scale~$\bar N$. Furthermore,
this will be the crucial vehicle for constructing {\em regular} segments as introduced in Definition~\ref{def:8.regsegdef},
see Proposition~\ref{prop:9.uncondreg} below which is the main result of this section.
With the notations of Section~\ref{sec:parameter}, define
\begin{align*}
\tilde \Omega_N &:= \Omega_N\cup \Omega_{N-1}\cup \Omega_{N-2} \\
\tilde\cB_{N,\omega} &:= \cB_{N,\omega} \cup \cB_{N-1,\omega} \cup
\cB_{N-2,\omega} \\
\tilde \Omega_N^{(1)} &:= \Omega_N^{(1)}\cup \Omega_{N-1}^{(1)}\cup \Omega_{N-2}^{(1)} \\
\tilde\cB_{N,\omega}^{(1)} &:= \cB_{N,\omega}^{(1)} \cup \cB_{N-1,\omega}^{(1)} \cup
\cB_{N-2,\omega}^{(1)}
\end{align*}
 where $ \Omega_N $,  $ \cB_{N,\omega}$,  $ \Omega_N^{(1)} $,  $ \cB_{N,\omega}^{(1)}$ are
the same as in Proposition~\ref{prop:11.4rellich} and $E_j^{(N)}$ be
as in the previous section. In the following lemma, we begin with
the comparison of the spectra of the entries, as indicated in the
previous paragraph.

\begin{lemma}
\label{lem:9.iuncond}  Using the notation from above, one has the following:
\begin{enumerate}
\item Let $\omega \in \TT_{c,a} \setminus \Omega_N$,
and  $x_0 \in \TT\setminus\tilde\cB_{N,\omega}$. Then
 \[\begin{split}&\min \big[\dist ( \spec \big (H_{[1,N-1]} (x_0, \omega) \big), E_{j_0}^{(N)}(x_0, \omega)),
  \dist(\spec\big( H_{[2,N]}(x_0, \omega) \big),E_{j_0}^{(N)}(x_0, \omega))\big] \\
&\leq r_0 =
\exp(-N(\log N)^{-C_0})\end{split}\] provided $N$ is large.
\item Furthermore, let $r_0 = e^{-N^\delta}$  where $0<\delta\ll1$  is arbitrary but fixed and assume that
$N\ge N_0(V,a,c,\gamma, \delta)$. If
 \begin{equation} \label{eq:9.uncondit}
 \max\big[\dist ( \spec \big(H_{[1,N-1]} (x_0, \omega) \big), E_{j_0}^{(N)}(x_0, \omega)),
 \dist(\spec \big(H_{[2,N]}(x_0, \omega) \big),E_{j_0}^{(N)}(x_0, \omega))\big] \leq
r_0,
\end{equation}
 then
\begin{equation}\label{eq:9.uncondit2}\dist(\spec \big(H_{[2,N-1]}(x_0, \omega) \big),E_{j_0}^{(N)}(x_0, \omega))\leq r_1
\end{equation} where $r_1=e^{-N^{\delta/2}}$.
\item  Finally, assume that $\omega\in\tor_{c,a}\setminus \tilde \Omega_N^{(1)}$ and $x_0\in \tor\setminus\tilde \cB_{N,\omega}^{(1)}$.  Let $r_0=\exp(-(\log N)^A)$
with some constant $A$. Then
 \eqref{eq:9.uncondit} implies \eqref{eq:9.uncondit2} with $r_1=e^{-(\log N)^{A/2}}$ for $N\ge N_0(V,c,a,A)$.
\end{enumerate}
\end{lemma}
\begin{proof}
 Let $E_0 =
E_j^{(N)}(x_0, \omega)$.  Then $f_N(e(x_0), \omega, E_0)=f_{[1,N]}(e(x_0), \omega, E_0)=0$. Since
\begin{equation} \label{eq:9.unimod} -f_{[1,N]}(e(x_0), \omega, E_0)f_{[2,N-1]}(e(x_0), \omega, E_0) +
f_{[1,N-1]}(e(x_0), \omega ,E_0)f_{[2,N]}(e(x_0), \omega,
E_0)= 1
\end{equation}
one obtains
\[ f_{[1,N-1]}(e(x_0), \omega, E_0)f_{[2,N]}(e(x_0), \omega, E_0)=1
\]
In particular,
\[ \min(|f_{[1,N-1]}(e(x_0), \omega, E_0)|,|f_{[2,N]}(e(x_0), \omega, E_0)|)\leq 1
\]
Assume, for instance, that
\[ |f_{[1,N-1]}(e(x_0), \omega, E_0)| \leq 1
\]
Then, due to Corollary \ref{cor:2.uniexcepzero} with $\eta =
\exp(-N(\log N)^{-C_0})$, one has
\begin{equation}\label{eq:9.speta}
(E_0-\eta,E_0+\eta) \cap \spec \big(H_{[1,N-1]}(x_0, \omega) \big) \neq
\emptyset
\end{equation}
This proves (1).  Assume now that in addition to
\eqref{eq:9.speta} one has
\begin{equation}\label{eq:9.spetao} (E_0-r_0,E_0+r_0) \cap \spec \big(H_{[2,N]}(x_0, \omega) \big) \neq \emptyset
\end{equation}
where $r_0=e^{-N^\delta}$ where $0<\delta$ is small and fixed.
Let $E_j^{[a,N-b]}(x, \omega), j=1,2,\ldots$ stand for the
eigenvalues of $H_{[a,N-b]}(x, \omega)$, $a=1,2$, $b=0,1$.  Due to
\eqref{eq:9.speta} and \eqref{eq:9.spetao},
 there exist $E_{j_1}^{[1,N-1]}(x_0, \omega),\, E_{j_2}^{[2,N]}(x_0, \omega)$ such that
\[|E_{j_1}^{[1,N-1]}(x_0, \omega) - E_0|,\; |E_{j_2}^{[2,N]}(x_0, \omega) - E_0| < r_0
\]
Due to Corollary \ref{cor:4.2}, one now has
\begin{align}\label{eq:9.factor}
 f_{[1,N]}(e(x_0), \omega, E) &=(E-E_{j_0}^{[1,N]}(x_0, \omega))\chi_0(e(x_0), \omega, E) \\
 \label{eq:9.factorr}
f_{[1,N-1]}(e(x_0), \omega, E) &=(E-E_{j_1}^{[1,N-1]}(x_0,
\omega))\chi_1(e(x_0), \omega, E) \\
\label{eq:9.factorl}
f_{[2,N]}(e(x_0, \omega, E) &=(E-E_{j_2}^{[2,N]}(x_0,
\omega))\chi_2(e(x_0), \omega, E)
\end{align}
where $\chi_k(z, \omega, E)$ is analytic in
$\cD(e(x_0),r_2)\times\cD(E_0,r_2)$, $r_2\asymp\exp(-N^{\delta/4})$,
with $\omega$ being fixed,
 $\chi_k(z, \omega, E)\neq 0$ for any $(z, E)\in\cD(e(x_0),r_2)\times\cD(E_0,r_2)$, $k=0,1,2$.
Moreover,
\begin{equation}\label{eq:9.kappadiv}
{NL}(E_0,\omega)-N^{\delta/3} < \log|\chi_k(z, \omega, E)| <
\mathop{NL}(E_0,\omega)+N^{\delta/3}
\end{equation}
for any $(z,E) \in \cD(e(x_0),r_2/2)\times\cD(E_0,r_2/2)$.
It
follows from \eqref{eq:9.unimod} and \eqref{eq:9.factor} --
\eqref{eq:9.kappadiv} that
\begin{equation}
\begin{split} |f_{[2,N-1]}(e(x_0), \omega, E)|\leq |E-E_{j_0}^{[1,N]}(x_0, \omega))|^{-1}\exp(-\mathop{NL}(E_0,\omega)
+2N^{\delta/3}) + \\
|\theta(E)| \, |E-E_{j_0}^{[1,N]}(x_0,
\omega)|^{-1}\, |E-E_{j_1}^{[1,N-1]}(x_0, \omega)|\, |E-E_{j_2}^{[2,N]}(x_0,
\omega)|\exp(\mathop{NL}(E_0,\omega))
\end{split}
\end{equation} for any $|E-E_0| < r_2/2$, where \[ \theta(E):= e^{-NL(E_0,\omega)} \chi_1
\chi_2/\chi_0\] and satisfies the bound
\[ \exp(-3N^{\delta/3})\leq|\theta(E)|\leq \exp(3N^{\delta/3})
\]
Clearly there exists $|E_1-E_0|\le 2r_0$ such that
\[ |E_1-E_{j_0}^{[1,N]}(x_0, \omega)|^{-1}\,|E_1-E_{j_1}^{[1,N-1]}(x_0, \omega)||E_1-E_{j_2}^{[2,N]}(x_0,\omega)|\le 2r_0
\]
So,
\[ |f_{[2,N-1]}(e(x_0), \omega, E_1)|\leq\exp(\mathop{NL}(E_0,\omega)-N^{\delta}+10N^{\delta/3})
\]
Due to Corollary~\ref{cor:2.uniexcepzero}, one has with $\eta
= \exp(-N^{\delta/2})$,
\[ (E_0-\eta,E_0+\eta)\cap\spec \big(H_{[2,N-1]}(x_0,\omega) \big)\neq\emptyset
\]
and we are done with case (2). Case~(3) is completely analogous.
\end{proof}

The following definition introduces the crucial notion of an {\em
unconditional spectral value}.

\begin{defi}\label{def:9.uncondint}
Using the notation of Lemma \ref{lem:9.iuncond} assume that with $r_0 = e^{-N^\delta}$, 
\begin{equation}\label{eq:9.uncon}
\max\big[\dist(\spec \big( H_{[1,N-1]}(x_0, \omega) \big) ,E_{j_0}^{(N)}(x_0,
\omega)),\dist(\spec \big( H_{[2,N]}(x_0, \omega) \big) , E_{j_0}^{(N)}(x_0,
\omega))\big]\leq r_0 
\end{equation}
In this case one says that $E_0:=E_{j_0}^{(N)}(x_0, \omega)$ is an
$r_0$-{\em unconditional spectral value} of $H_{[1,N]}(x_0,
\omega)$. Let $\{E_{j_0}^{(N)}(x, \omega),\ux,\ox\}$ be a segment of
a Rellich graph. Assume that for any $x\in [\ux,\ox]$,
$E_{j_0}^{(N)}(x, \omega)$ is an $r_0$-unconditional spectral value
of $H_{[1,N]}(x, \omega)$. We call this segment an
$r_0$-{\em unconditional spectral segment} of the Hamiltonian $H_{[1,N]}(\cdot,
\omega)$.
\end{defi}

Note that by (2) of Lemma~\ref{lem:9.iuncond}, each entry has a zero in energy which is
close to~$E_{j_0}^{(N)}(x_0,
\omega)$ in case the latter is an unconditional spectral value.
The importance of the  unconditional spectral values lies with the fact that they are stable
with regard to the induction on scales procedure. In other words, the unconditional energies
at a small scale will turn out to belong to the spectrum (or rather, be close to it) at the
next larger scale. This process can also be reversed: we will show later that the conditional spectral
values die out when we pass to the next larger scale. The corresponding analysis appears later in this section
when we discuss
property (NS) which stands for ``non-spectral''.

\begin{corollary}\label{cor:9.uncondmod}
Under the assumptions of Lemma~\ref{lem:9.iuncond} one has the following:
There exists $N_0=N_0(V,c,a,\gamma,\delta)$, or respectively, $N_0=N_0(V,c,a,\gamma,A)$
so that the following holds for all $N\ge N_0$:
Assume that $E_0=E^{(N)}_{j_0}(x_0,\omega)$ is an $r_0$-unconditional spectral value of $H_{[1,N]}(x_0, \omega)$.
Then
\begin{equation}\label{eq:9.normuncond'}
\log\lVert M_N(e(x_0),
\omega, E_0)\rVert \leq \mathop{NL}(E_0,\omega)- N^{\delta/2}
\end{equation}
or, respectively,
\begin{equation}\label{eq:9.normuncond''}
\log\lVert M_N(e(x_0),
\omega, E_0)\rVert \leq \mathop{NL}(E_0,\omega)- (\log N)^{A/2}
\end{equation}
Conversely, assume that
\eqref{eq:9.normuncond'} $($respectively, \eqref{eq:9.normuncond''}$)$ holds. Then $E_0$ is an $r'_0$-unconditional
spectral value of $H_{[1,N]}(x_0, \omega)$, with $ r'_0:=\exp(-N^{\delta/3}) $ $($respectively, $r'_0:=\exp(-(\log N)^{A/3}) $ $)$.
\end{corollary}
\begin{proof}
Inspection  of the proof of Lemma~\ref{lem:9.iuncond} establishes the direct implication.
Assume that \eqref{eq:9.normuncond'} holds. Then
\[
\log |f_{[a,N-b]}(e(x_0),
\omega, E)| \leq \mathop{NL}(E_0,\omega)- N^{\delta/2}
\]
for any $a=1,2$, and $b=0,1$. Due to Corollary~\ref{cor:2.uniexcepzero}, one has with $\eta
= \exp(-N^{\delta/3})$,
\[ (E_0-\eta,E_0+\eta)\cap\spec \big(H_{[a,N-b]}(x_0,\omega)\big)\neq\emptyset
\]
for any $a=1,2$, and $b=0,1$. Due to  Definition~\ref{def:9.uncondint} this means that $E_0$ is an
$r'_0$-unconditional spectral value for $H_{[1,N]}(x_0,\omega)$. The case of \eqref{eq:9.normuncond''} is similar.
\end{proof}

Recall that for any $x\in\TT$, $E \in \IR$, and any integers
$N'_i,N''_i$, $i=1,2$ such that $N_1'<N_2'<N_2''<N_1''$ one has
\begin{equation}\label{eq:13.monddrel}
\big |\log\lVert M_{[N'_1,N''_1]}(e(x), \omega, E)\rVert-\log\lVert
M_{[N'_2,N''_2]}(e(x), \omega, E)\rVert\big | \leq
C(V,E)\max(N'_2-N'_1,N''_1-N''_2)
\end{equation}
Combining \eqref{eq:13.monddrel} with Corollary~\ref{cor:9.uncondmod} implies the following
important {\em stability property} implied by the concept of unconditional spectral values.
Note that this {\em fails} if at least one entry does not have a zero close to $E_0$.

\begin{corollary}\label{cor:9.uncondedges}
Using the notations of Lemma~\ref{lem:9.iuncond} assume that $E_0$ is an $r_0$-unconditional spectral value of $H_{[1,N]}(\cdot, \omega)$.  Then
for any $-(\log N)^C< N' \leq 1 < N \leq N'' <
N+(\log N)^C$ one has
\[\dist(\spec \big(H_{[N',N'']}(x_0, \omega) \big), E_0) < r'_0
\]
where $r'_0:=\exp(-N^{\delta/4})$ , respectively $r'_0:=\exp(-(\log N)^{A/4})$.
\end{corollary}

To continue the analysis of unconditional spectral segments we will
use the avalanche principle expansion as in
Proposition~\ref{prop:add_zerosE} for the following logarithm
\begin{equation}\label{eq:9.nbar}
\log\lvert f_{[-\oN,\oN]}(e(x_0), \omega, E)\rvert
\end{equation}
where $N^2 \le \oN \le N^{10}$,
 $E\in\cD(E_0,r_0)$, and $N,x_0,E_0$ are as in Corollary~\ref{cor:9.uncondedges}.
We arrange the expansion so that one of the
intervals, say $\Lambda_{k_0}=[n_{k_0},n_{k_0+1}]$, will obey the following condition:
\begin{equation}\label{eq:9.expandiv}
-(\log N)^C < n_{k_0} < 1 < N < n_{k_0+1} < N + (\log N)^C
\end{equation}
 Due to Corollary~\ref{cor:18.2}, one can assume that
$n_1$ etc.~are
 adjusted at scale $\ell:=(\log N)^{C/6}$ relative to $\cD(e(x_0),r_1)\times\cD(E_0,r_1)$, $r_1:=\exp(-(\log N)^{A})$.
 Finally, one can assume that
\begin{equation}\label{eq:9.minsize}
N\le \min_k(n_{k+1}-n_k)\le  N + (\log N)^{C}
\end{equation}
Set
\begin{align*}
\bar\cB_{N,\omega} &:=\mybigcup_{N \leq N' \leq N+(\log N)^{C}} \tilde {\cB}_{N',\omega}\\
{\bar\cB_{N,\omega}}^{(1)} &:=\bar\cB_{N,\omega}\cup  \bigcup_{N\le N'\le N+(\log N)^C}{\cB}^{(1)}_{N',\omega}
\end{align*}
The notations of this paragraph will be used in the following lemma and proposition.

\begin{lemma}\label{lem:9.uncondcount}Using the notations of Corollary~\ref{cor:9.uncondedges} assume that $E_0$ is an $r_0$-unconditional spectral value of $H_{[1,N]}(\cdot, \omega)$.
Assume that $x_0\in \tor\setminus \bar\cB_{N,\omega}$ $($respectively $x_0\in \tor\setminus {\bar\cB_{N,\omega}}^{(1)}$ $)$.  Then for any $|x-x_0|\le
C^{-1}r_0$ each entry of the matrix
\[ M_{[n_{k_0},n_{k_0+1}]}(e(x), \omega, \cdot)\]
has exactly one zero in the disk $\cD(E_0,r_0')$.   Furthermore,
 \begin{equation} \label{eq:9evexist}
\nu_{f_{[-\oN,\oN]}(e(x), \omega, \cdot)}(E_0,r_0'')\geq 1
\end{equation}
where $r_0''=(r'_0)^{1/4}$.
\end{lemma}

\begin{proof}
The first part of the statement follows from
Corollary~\ref{cor:9.uncondedges} and
Proposition~\ref{prop:11.4rellich}. Due to
Proposition~\ref{prop:add_zerosE} one has
\begin{equation}\label{eq:13.jenNbarN}
  \cJ(\log|f_{[-\oN,\oN]}(z_0,\omega_0,\cdot)|, E_0, r_1, r_2) \ge \cJ(\log|f_{[n_{k_0}, n_{k_0+1}]}(z_0,\omega_0,\cdot)|, E_0, r_1, r_2) - O(\sqrt {r_1})
\end{equation}
for any $(r_0')^{1/2}\le r_1 \le (r_0')^{1/3}$, $r_2=r_1/4$.
Since $x_0\in \tor\setminus \bar\cB_{N,\omega}$ $($respectively $x_0\in \tor\setminus {\bar\cB_{N,\omega}}^{(1)}$ $)$ one can pick
$r_1$ so that $f_{[n_{k_0}, n_{k_0+1}]}(z_0,\omega_0,\cdot)$ has no zeros in $\cD(E_0,2r_1)\setminus\cD(E_0,r_1)$. Then
\begin{equation}\label{eq:13.zeroNbarN}
4  r_1^2 r_2^{-2} \cJ(\log|f_{[n_{k_0}, n_{k_0+1}]}(z_0,\omega_0,\cdot)|, E_0, r_1, r_2)\ge 1
\end{equation}
\eqref{eq:13.jenNbarN} combined with \eqref{eq:13.zeroNbarN} imply \eqref{eq:9evexist} via Lemma~\ref{lem:2.jensencor}.
\end{proof}

  \noindent Using the notations of Definition~\ref{def:9.uncondint}, assume that $E_0$ is an $r_0$-unconditional
spectral value of $H_{[1,N]}(x_0,\omega)$. Assume also that $x_0\in \tor\setminus \bar\cB_{N,\omega}$
 $($respectively, $x_0\in \tor\setminus {\bar\cB_{N,\omega}}^{(1)}$ $)$.
Due to  Lemma~\ref{lem:8.maxintervals} there exists a segment
$\{E_{j_0}^{(N)}(x, \omega),\ux,\ox\}$ such that the following
conditions hold:
\[(i)\quad x_0\in [\ux,\ox],\quad (ii)\quad |\ux-\ox|\ge
\exp(-(\log r_0^{-1})^B),\quad (iii)\quad |\partial_x E_{j_0}^{(N)}|\ge
\exp(-(\log r_0^{-1})^B)\]  where $1\ll B $. Set
$\ux_1:=\max(x_0-C^{-1}r_0,\ux)$, $\ox_1:=\min(x_0+C^{-1}r_0,\ox)$.
The following proposition is the main statement concerning the formation
of regular spectral segments.

\begin{prop}\label{prop:9.uncondreg}
Using the above notations assume for instance that
$\{E_{j_0}^{(N)}(x, \omega),\ux,\ox\}$ is a positive slope segment.
Then for any $x_1 \in [\ux_1,\ox_1]\setminus \cB_{\oN,\omega}$
$($respectively, $x_1 \in [\ux_1,\ox_1]\setminus\
\cB^{(1)}_{\oN,\omega}$ $)$ and any $a\in [-\frac {1}{8}\oN ,\frac
{7}{8}\oN]$ there exists a positive slope regular spectral segment
$\{E_{j'}^{(\oN)}(x, \omega),\ux',\ox'\}$ of $H_{[-\oN,\oN]}(x,
\omega)$ such that the following conditions hold:
\begin{itemize}
\item[(1)] $x_1\in [\ux',\ox']$,
\item[(2)] $|\ux'-\ox'|\ge \exp(-(\log r_0^{-1})^{B})$,
\item[(3)] $|E_{j'}^{(\oN)}(x, \omega)-E_{j_0}^{(N)}(x,
\omega)| \le (r')_0^{1/2}$,
\item[(4)] $|\partial_x E_{j'}^{(\oN)}|\ge
\exp(-(\log r_0^{-1})^B)$,
\item[(5)] $a-N\le \nu^{(\oN)}_j(x, \omega) \le
a+2N$.
\end{itemize}
\end{prop}
\begin{proof}
Let $a=0$. Due to Lemma \ref{lem:9.uncondcount} one has
\[ \dist\big[\spec \big( H_{[-\oN,\oN]}(x, \omega)\big), E^{(N')}_{j_0}(x, \omega) \big]\leq r''_0 \]
for any $x \in [\ux_1,\ox_1]$. Assume that $x_1 \in [\ux_1,\ox_1]\setminus\cB_{\oN,\omega}$.
Due to Lemma~\ref{lem:8.maxintervals} there exists a segment $\{E_{j_1}^{(\oN)}(x, \omega),\ux_1,\ox_1\}$
 such that conditions $(1)$, $(2)$, $(3)$ hold. Condition $(4)$ follows from
 $(3)$ combined with $(iii)$ and the Cauchy estimates for the derivatives. Condition $(5)$
 follows just from the definition of the center of localization
 and the first part in Lemma~\ref{lem:9.uncondcount}. For arbitrary $a\in [-\frac {3}{4}\oN ,\frac
{3}{4}\oN]$ the argument is similar.
\end{proof}

We now turn to the investigation of conditional spectral values. For technical reasons
we formulate this property in the following way which does not require $E_0$ to be a zero
of the first entry of $M_N$; in fact, it will be convenient to also state this at scale~$\ell$ rather than~$N$.
(NS) here stands for {\em non-spectral}, which refers to the fact that $E_0$ will be separated from the spectrum
of $H_N(e(x),\omega)$ (up to shifting the edges) uniformly in~$x$, see Proposition~\ref{prop:14.condspecperish}.

\begin{defi}\label{def:II0}
Let $\Omega_\ell$, $\cE_{\omega,\ell}$ be as in
Lemma~\ref{lem:3.waschno}. Let $\omega \in \tor_{c,a} \setminus
\mybigcup_{m=\ell-2,\ell-1,\ell} \Omega_m$, $x_0\in [0,1]$, and
$E_0\in \IR\setminus \cE_{\omega,\ell}$. By condition $(NS)$ we
mean that at least one of the Dirichlet determinants
\[f_{[1, \ell]}(e(x_0),
\omega, \cdot),\; f_{[1, \ell-1]}(e(x_0), \omega, \cdot),\; f_{[2,
\ell]}(e(x_0), \omega, \cdot),\; f_{[2, \ell-1]}(e(x_0), \omega,
\cdot)
\] has no zeros in $\cD(E_0, r_0)$,
$r_0:= \exp(-\ell^{\delta})$, where $0< \delta\ll 1$ is a parameter.
\end{defi}

Similarly to Corollary ~\ref{cor:9.uncondmod} one now has the
following statement.

\begin{lemma}\label{lem:14.condmod}
Assume that condition $(NS)$ holds. Then
\begin{equation}\label{eq:9.normuncond}
\mathop{\ell L}(E_0,\omega) -\ell^{2\delta}\leq\log\lVert
M_{\ell}(e(x_0),\omega,E)\rVert
\end{equation}
for any $|E-E_0| < r_0/2$. Conversely, assume that
\eqref{eq:9.normuncond} holds for any $|E-E_0| < r_0/2$. Then
condition $(NS)$ holds with $r_0$ replaced by
$r'_0:=\exp(-\ell^{4\delta})$.
\end{lemma}
\begin{proof}
If \eqref{eq:9.normuncond} fails for some
$|E-E_0| < r_0/2$,  then
\[
\log |f_{[a,\ell-b]}(e(x_0), \omega, E)| \leq \mathop{\ell
L}(E_0,\omega) -\ell^{2\delta}
\]
for any $a=1,2$, and $b=0,1$. Due to
Corollary~\ref{cor:2.uniexcepzero}, one has with $\eta :=r_0$
\[ (E_0-\eta,E_0+\eta)\cap\spec \big(H_{[a,\ell-b]}(x_0,\omega) \big)\neq\emptyset
\]
for any $a=1,2$, and $b=0,1$, contrary to the assumptions of the
lemma.
For the converse, let us assume that some determinant
$f_{[a,\ell-b]}(e(x_0), \omega, \cdot)$ has a zero at $E_1$, where
$|E_1-E_0|\le r'_0$. Just as in the proof of
Lemma~\ref{lem:9.iuncond} it follows that
\[\log|f_{[a,\ell-b]}(e(x_0), \omega, \cdot)| < \mathop{\ell L}(E_0,\omega) -\ell^{3\delta}\]
Consequently, if each entry of $M_\ell(e(x_0),\omega,\cdot)$ were to
exhibit such a zero, then
\[
\mathop{\ell L}(E_0,\omega) -\ell^{2\delta} > \log\lVert
M_{\ell}(e(x_0),\omega, E)\rVert
\]
which is a contradiction.
\end{proof}

In the following corollary we show that $(NS)$ has a natural
stability property.

\begin{corollary}\label{cor:14.condedges}
Assume condition $(NS)$. Then for any $-\ell^{\delta} < \ell' \leq
1 < \ell \leq\ell'' < \ell+\ell^{\delta}$,
 at least one of the
Dirichlet determinants
\[f_{[\ell'+1, \ell'']}(e(x_0),
\omega, \cdot),\; f_{[\ell'+1, \ell''-1]}(e(x_0), \omega, \cdot),\;
f_{[\ell'+2, \ell'']}(e(x_0), \omega, \cdot),\; f_{[\ell'+2,
\ell''-1]}(e(x_0), \omega, \cdot)
\]
has no zeros in $\cD(E_0, r'_0)$, where
$r'_0:=\exp(-\ell^{4\delta})$.
\end{corollary}
\begin{proof}
 For any $x$, $E$, and any integers
$N'_i,N''_i$, $i=1,2$ such that $N_1'<N_2'<N_2''<N_1''$ one has
\begin{equation}\label{eq:9.monddrel}
\big |\log\lVert M_{[N'_1,N''_1]}(e(x),E)\rVert-\log\lVert
M_{[N'_2,N''_2]}(e(x),E)\rVert\big | \leq
C(V,E)\max(N'_2-N'_1,N''_1-N''_2)
\end{equation}
Hence, the estimate \eqref{eq:9.normuncond} is stable under such
changes to the length of the monodromy matrix as long as the change
is much smaller than $\ell^{2\delta}$. In particular, this holds
with a change of size $\ell^\delta$ as stated in the corollary.
\end{proof}

The following proposition proves that conditional spectral values die out when we
pass from scale~$\ell$ to a larger scale~$N$.  Clearly, this is of great importance
as it ensures that our inductive procedure can be carried out with unconditional spectral
values~$E_0$ rather than conditional ones. The relevance of the unconditionality property here stems from
Proposition~\ref{prop:9.uncondreg} which shows that at larger scales eigenfunctions with eigenvalues close to~$E_0$
are localized away from the edges of the underlying interval. If this were not so, then
these eigenfunctions would not be stable when passing the next larger scale.
We note one disadvantage of the following proposition: for each~$x$ it gives a choice of four determinants.
We will subsequently see that {\em periodic} boundary conditions can be used uniformly for all~$x$.
This is the reason whey periodic boundary conditions appear in this paper at all.

\begin{prop} \label{prop:14.condspecperish} Let $E_0$ be arbitrary. Assume there is an
interval $(x'_0,x''_0)$ with $x''_0-x'_0 \ge \ell^{-\delta/3}$ such
that for any $x_0\in (x'_0,x''_0)$ condition $(NS)$ holds. Then
for any $|E-E_0|\le r_0/4$ and any $x\in \tor$, and $\ N \ge \ell^2$
at least one of the Dirichlet determinants
\[f_{[1, N]}(e(x),
\omega, \cdot),\; f_{[1, N-1]}(e(x), \omega, \cdot),\; f_{[2,
N]}(e(x), \omega, \cdot),\; f_{[2, N-1]}(e(x), \omega, \cdot)
\]
has no zeros in $\cD(E, r''_0)$, where $r''_0:=\exp(-N^{8\delta})$.
\end{prop}
\begin{proof} Let $|E-E_0|\le r_0/4$ be arbitrary. Note that $(NS)$ holds on the disk $\cD(E,r_0/4)$
for any $x_0\in (x'_0,x''_0)$.
 Let $x\in \tor$ be arbitrary. We invoke once again the avalanche principle
expansion~\eqref{eq:2.29}. Let $\ell^2\le N \le \ell^{10}$ be
arbitrary. By Corollary~\ref{cor:14.condedges}, one can pick $1\le a\le
\ell^{\delta}$, $N-\ell^{\delta}\le b \le N $, and arrange the
avalanche principle expansion so that the following conditions hold:
\begin{itemize}
\label{eq:10.partcond}
\item[(1)]
$ [a, b] = \mybigcup^t_{k=1}
\Lambda_k$, $\; \Lambda_k = [n_k+1, n_{k+1}]$
\item[(2)]
 for any $k=1,\ldots,t-1$ there exists $|n'_k-n_k|\le \ell^{\delta}$,
 such that $ x+n'_k\omega\in (x'_0,x''_0) $
\item[(3)] $|n_{k+1} - n_k - \ell| \le 2\ell^{\delta}$
\item[(4)]$n_1, \dots, n_t$ are adjusted to $\cD(e(x), r_0/C(V)) \times \cD(E,
r_0/2)$, $r_0=\exp(-\ell^{\delta})$ at scale $\ell^{\delta/6}$
\end{itemize}
In item $(4)$ we used Corollary~\ref{cor:18.2}, whereas for (2) we
invoked the dynamics: since $\omega\in \tor_{c,a}$, for any
$x\in\tor$ and any $s'\in \IZ , s' > 0$ there exists $s'\le s\le
s'^2$ such that
\begin{equation}\label{eq:14.diophappr}
\|x-s\omega\| \le 1/s
\end{equation}
In view of (1)--(4) we can apply the zero count of
Proposition~\ref{prop:add_zerosE}; the crucial observation here is
that due to Corollary~\ref{cor:14.condedges} one can pick the edge
points $n_k$ so that $ f_{[n_k+1, n_{k+1}]}(e(x),\omega,\cdot)$ has
no zeros in $\cD(E, r'_0)$, where $r'_0:=\exp(-\ell^{4\delta})$. In
conclusion,
\begin{equation}
\nn \nu_{f_{[a,b]}(e(x),\omega, \cdot)} (E, r'_0/2)=\sum^{t-1}_{k=0}
\nu_{f_{[n_k+1, n_{k+1}]}(e(x),\omega, \cdot)} (E_0,r'_0/2 ) =0
\end{equation}
This means that $(NS)$ holds for the entries of
$M_{[a,b]}(e(x),\omega,\cdot)$ on the disk $\cD(E,r_0'/2)$ for all
$x\in(x'_0,x''_0)$. By Corollary~\ref{cor:14.condedges} one can now
replace $[a,b]$ by $[1,N]$ whence the proposition for the range
$\ell^2\le N \le \ell^{10}$. For arbitrary $N\ge \ell^2$ one can use
an induction argument. Indeed, by the preceding any $x_0\in
(x'_0,x''_0)$ has the property that $(NS)$ holds for the entries
of $M_{[1,N^{(1)}]}(e(x),\omega,\cdot)$ on the disk
$\cD(E,r^{(1)})$, with $N^{(1)}:=\ell^{8}$ and
$r^{(1)}:=\exp\big(-(N^{(1)})^{\delta}\big)$. Therefore, one can apply the
very same argument with $N^{(1)}$ in the role of~$\ell$.
\end{proof}

We now discuss the relation between the unconditional Dirichlet
spectral values and the periodic spectrum. For this purpose, recall the
following relation from Proposition~\ref{th:2.avatrace} and~\ref{prop:LDEinEtrace} which follows from properties of the trace:
\begin{equation}\label{eq:13.hill}
\log |g_N(e(x),\omega,E)| = \log\|M_{2N}(e(x ),\omega,E)\| -
 \log \|M_N(e(x), \omega, E)\|] +O\bigl(\exp(-(\log N)^B)\bigr)
\end{equation}
 provided $\|N\omega\|\le \kappa_0(V,c,a,\gamma)$ and $E\in \IC
\setminus \cE_x$, where $\cE_x \subset \IC$, $\mes(\cE_x) \le
\exp(-(\log N)^B)$.
Here $g_N(e(x),\omega,E):=\det(H^{(P)}_N(x,\omega)-E)$ with
$H^{(P)}_N(x,\omega)$ the Schr\"odinger operator on $[1,N]$ with
periodic boundary conditions, and
\[
g_N(e(x),\omega,E)=\rtr( M_N(e(x),\omega,E))-2
\]
The reason for considering periodic boundary conditions is as follows: note that
Proposition~\ref{prop:14.condspecperish} shows that under the (NS) condition at scale $\ell$,
for each $x\in\tor$ it is possible to ``wiggle'' the boundary of $[1,N]$ slightly so as to ensure that the corresponding entry
has no zeros in a disk of energies.  The technically unpleasant feature here is that the ``wiggling''
or in precise terms, the choice of boundary conditions, depends on~$x$. However, we will now
see that~\eqref{eq:13.hill} implies that periodic boundary conditions
achieve the desired absence of zeros in~$E$ {\em uniformly} for all $x\in\tor$.

\begin{lemma}\label{lem:13.hillNclosetodir2N} Assume that for some
$x_0$, $E_0$ and with $2N$ in the role of $\ell$ condition $(NS)$
holds. Then
\item[(1)]
$ \nn  \cJ \bigl(\log |g_N(e(x_0), \omega, \cdot)|, E_0, r\bigr)\le
\exp(-(\log N)^C) $ for any $\exp(- N^{1/2})\le r\le
r_1:=\exp(-N^{10\delta})$
\item[(2)]
$ \nu_{g_N(e(x_0),\omega,\cdot)}(E,r_1/2)=0 $
\end{lemma}
\begin{proof} Since each logarithm involved in
\eqref{eq:13.hill} is subharmonic in $E$, one has
\[
  \cJ \bigl(\log |g_N(e(x_0), \omega, \cdot)|, E_0, r\bigr)\le
   \cJ \bigl(\log \|M_{2N}(e(x_0),\omega,\cdot)\| , E_0, r\bigr) +
C\exp(-(\log N)^B)
\]
 Therefore, the estimate in~$(1)$ is due to Proposition~\ref{prop:14.10}. Part~$(2)$ follows
from $(1)$ due to Lemma~\ref{lem:2.jensencor}.
\end{proof}

Property (2) in Lemma~\ref{lem:13.hillNclosetodir2N} simply says that the periodic
spectrum does not intersect the interval $(E'_0,E''_0)$. We record this fact as a separate
statement and definition.

\begin{corollary}\label{cor:14.1inlarge''}
 Using the notations of
Proposition~\ref{prop:14.condspecperish} one has
 with $(E'_0,E''_0):=(E_0-r_0/8,E_0+r_0/8)$
\begin{equation}\label{eq:13.spectrumfree}
 \spec\big ( H_N^{(P)}(x, \omega) \big)\cap (E'_0, E''_0)=\emptyset
\end{equation}
 for any $x\in \tor$ and any $N\ge \ell^2$ provided $\|N\omega\|\le
 \kappa_0(V,c,a,\gamma)$. An interval $(E'_0,
E''_0$) is called {\em spectrum free} if there exists $N_0$
depending on the usual parameters $a,c,V,\rho_0$ as well as on
$(E'_0, E''_0)$ such that~\eqref{eq:13.spectrumfree} holds for any
$N\ge N_0$ with $\|N\omega\| \le \kappa_0$ and any $x\in \tor$.
\end{corollary}

The following lemma establishes a crucial dichotomy between an interval
being spectrum free and the existence of a regular segment.

\begin{lemma}\label{lem:13.eitheruncondornonspectral}
Given a scale $\ell$, a parameter $0 < \delta \ll 1$, and intervals
$(E'_0, E''_0)$, $(x'_0,x''_0)$ with \[ E''_0-E'_0\ge \exp(-(\log
\ell)^A),\quad x''_0-x'_0 \ge \ell^{-\delta/3} \] either the interval \[ \big( E'_0 + \frac
{1}{4}\exp(-(\log \ell)^A),\, E''_0 -\frac {1}{4}\exp(-(\log
\ell)^A) \big) \] is spectrum free, or for any scale $\ell^2 \le N \le
\ell^{10}$ and any $\frac {1}{4} N \le a \le \frac {3}{4} N$ there
exists a regular $I_1$-segment \[ \bigl\{E_{j_1}^{(N)}(x, \omega),
\underline{x}_1, \overline{x}_1 \bigr\},\quad I_1\subset (E'_0, E''_0)
\] with \[ (\underline{x}_1, \overline{x}_1)\subset (x'_0,x''_0),\quad  a-2\ell \le \nu^{(N)}_j(x, \omega) \le a + 2\ell \]
\end{lemma}
\begin{proof}
Assume that there exist $x_0\in (x'_0,x''_0)$ and \[ E_0\in
\big(E'_0+\frac {1}{4}\exp(-(\log  \ell)^A),E''_0 -\frac
{1}{4}\exp(-(\log  \ell)^A)\big)\] such that condition $(NS)$ fails.
Then there exists $j_1$ such that \[ E_1:=E_{j_1}^{(\ell)}(x_0,
\omega) \in \big(E'_0+\frac {1}{8}\exp(-(\log  \ell)^A),\,E''_0 -\frac
{1}{8}\exp(-(\log  \ell)^A)\big)\] is an $r_0$-unconditional spectral
value, where $r_0:=\frac {1}{4}\exp(-\ell^{\delta})$. Due to
Proposition~\ref{prop:9.uncondreg} for any $\ell^2 \le N \le
\ell^{10}$ and any $\frac {1}{4} N \le a \le \frac {3}{4} N$ there
exists a regular $I_1$-segment \[ \bigl\{E_{j_1}^{(N)}(x, \omega),
\underline{x}_1, \overline{x}_1 \bigr\},\quad I_1\subset (E'_0, E''_0)
\] with \[ (\underline{x}_1, \overline{x}_1)\subset (x'_0,x''_0), \quad a-2\ell \le \nu^{(N)}_j(x, \omega) \le a + 2\ell \]
 Fix arbitrary
\begin{equation}\label{eq:E0choice} E_0\in \big(E'_0+\frac {1}{4}\exp(-(\log \ell)^A),E''_0 -\frac
{1}{4}\exp(-(\log  \ell)^A)\big)\end{equation} and assume that for any $x_0\in
(x'_0,x''_0)$ condition $(NS)$ holds. There exists $N\asymp\ell^C$
such that $\|N\omega\|\le \kappa_0$ with $\kappa_0$ as in Proposition~\ref{th:2.avatrace}. Then due to
Corollary~\ref{cor:14.1inlarge''} there exists a spectrum free
subinterval $(E', E'')\subset (E'_0, E''_0)$ containing $E_0$.
Since $E_0$ as in~\eqref{eq:E0choice} is arbitrary, the statement holds.
\end{proof}

To proceed, we need a version of Lemma~\ref{lem:8.maxintervals} for
the parametrization of the eigenvalues of the Schr\"odinger operator
with periodic boundary conditions. Let
\begin{equation}
\nn
 E_1^{(N,P)}(x,\omega)\le E_2^{(N,P)}(x,\omega)\le \dots \le
E_N^{(N,P)}(x, \omega)
\end{equation}
be the eigenvalues of $H_{[1,N]}^{(P)}(x,\omega)$. We define
segments of the graphs of $ E_1^{(N,P)}(\cdot,\omega)$ via approximation by the segments $\bigl\{E_{j_1}^{(2N)}(x,
\omega), \underline{x}_1, \overline{x}_1 \bigr\}$.

\begin{lemma}\label{lem:13.hillsegemnts}
Assume that $\|N\omega\|\le \kappa_0$, with $\kappa_0$
 as in Proposition~\ref{th:2.avatrace}. Then the following properties hold:
\begin{itemize}
\item[(1)]  Let
${\cE}^{(1)}_{2N, \omega}$ be as in
Proposition~\ref{prop:11.4rellich} with $2N$ in the role of $N$.
Assume that for some $x_0$, $j$ one has \[
\dist(E_j^{(N,P)}(x_0,\omega), {\cE}^{(1)}_{2N, \omega})\ge
2\exp(-N^{\delta}).\] Then there exists an unconditional segment
$\bigl\{E_{j_1}^{(2N)}(x, \omega), \underline{x}_1,
\overline{x}_1 \bigr\}$ with $x_0\in (\underline{x}_1,
\overline{x}_1)$ such that
\begin{equation}\label{eq:hillclosetodir2Na}
|E_j^{(N,P)}(x,\omega)-E_{j_1}^{(2N)}(x,\omega)|\le
\exp(-N^{2\delta})
\end{equation}
for any $x\in [\underline{x}_1, \overline{x}_1] $
\item[(2)] Let  $\bigl\{E_{j_1}^{(2N)}(x, \omega),
\underline{x}_1, \overline{x}_1 \bigr\}$ be a regular segment such
that $\frac {1}{4} N\le \nu^{(2N)}_{j_1}(\cdot, \omega) \le \frac {1}{2}
N$ on that segment. Then there exists $E_j^{(N,P)}(\cdot,\omega)$ such that \eqref{eq:hillclosetodir2Na}
holds for all $x\in (\underline{x}_1, \overline{x}_1)$.
\end{itemize}
\end{lemma}
\begin{proof} Due to Lemma~\ref{lem:13.hillNclosetodir2N} there
exists $j_1$ such that \eqref{eq:hillclosetodir2Na} holds for
$x=x_0$. Note that $E_{j_1}^{(2N)}(x,\omega)\notin {\cE}^{(1)}_{2N,
\omega}$. Hence there exists a segment $\bigl\{E_{j_1}^{(2N)}(x, \omega), \underline{x}_1, \overline{x}_1 \bigr\}$, with $x_0\in
(\underline{x}_1, \overline{x}_1)$. Once again due to
Lemma~\ref{lem:13.hillNclosetodir2N} for any $x$ there exists $j(x)$
such that
\begin{equation}\label{eq:hillclosetodir2Nb}
|E_j^{(N,P)}(x,\omega)-E_{j(x)}^{(2N)}(x,\omega)|\le
\exp(-N^{2\delta})
\end{equation}
Recall that for any $j\neq j_1$ and any $x\in [\underline{x}_1,
\overline{x}_1] $ one has
\begin{equation}\label{eq:2Nseparation}
|E_j^{(2N)}(x,\omega)-E_{j_1}^{(2N)}(x,\omega)|\ge \exp(-(\log N)^A)
\end{equation}
Combining \eqref{eq:2Nseparation} with \eqref{eq:hillclosetodir2Na}
for $x=x_0$ one concludes that $j(x)=j_1$ for any $x\in
[\underline{x}_1, \overline{x}_1] $ in \eqref{eq:hillclosetodir2Nb}.
Thus, \eqref{eq:hillclosetodir2Na} holds. It follows from
Lemma~\ref{lem:13.hillNclosetodir2N} and
\eqref{eq:hillclosetodir2Na} that each value
$E_{j_1}^{(2N)}(x,\omega)$ is unconditional. That proves the
first part. To establish the second part, let us note that since $\frac
{1}{4} N\le \nu^{(2N)}_{j_1}(\cdot, \omega) \le \frac {1}{2} N$, one has
\[
\|(H^{(P)}_{[1,N]}(x_1,\omega)- E_{j_1}^{(2N)}(x_1,\omega)) \psi_{j_1}^{(2N)}(x_1,\omega,\cdot)\| \le
\exp(-\gamma N/8)
\]
for any $x_1\in (\underline{x}_1, \overline{x}_1)$. Hence, for any
$x_1$ there exists $j(x_1)$ such that
\begin{equation}
\nn
 |E_{j(x_1)}^{(N,P)}(x_1,\omega)-E_{j_1}^{(2N)}(x_1,\omega)|\le
\exp(-\gamma N/8)
\end{equation}
Let $x_1=(\underline{x}_1 + \overline{x}_1)/2$. Then due to the lower bound
on the Dirichlet graphs, the last
estimate implies in particular that
\[ \dist(E_{j(x_1)}^{(N,P)}(x_1,\omega), {\cE}^{(1)}_{2N, \omega})\ge 2\exp(-N^{\delta})\]
 Due to first part there exists $j_2$ such
that
\begin{equation}
\nn
 |E_{j(x_1)}^{(N,P)}(x_1,\omega)-E_{j_2}^{(2N)}(x_1,\omega)|\le
\exp(-N^{2\delta})
\end{equation}
for any $x_1\in (\underline{x}_1, \overline{x}_1)$. Due to the separation property
of the segments one concludes that $j_2=j_1$ and we are done.
\end{proof}

The following proposition establishes a dichotomy which will be
 an essential ingredient in our proof of Theorem~\ref{th:1.2}.
It shows that either an interval does not intersect the infinite volume
spectrum, or it has to contain the graphs of both a positive slope
as well as of a negative slope {\em regular} segment at all sufficiently large scales. In view of Figure~2
this is of course important for the creation of resonances and thus also, gaps.
As already mentioned before, the regularity of the segments is essential for the rigorous
implementation of Figure~2.

\begin{prop}\label{prop:13.doubleresonancegraphs}
Given  $\ell$ large, a parameter $0 < \delta \ll 1$, and intervals
$(E'_0, E''_0)$, $(x'_0,x''_0)$ with $E''_0-E'_0\ge \exp(-(\log
\ell)^A)$,
 $x''_0-x'_0 \ge \ell^{-\delta/3}$,  there is the following dichotomy:
\begin{itemize}
 \item either the
interval $(E'_0 +$   $\frac
{1}{4}\exp(-(\log \ell)^A), E''_0 -\frac {1}{4}\exp(-(\log
\ell)^A))$ is spectrum free
\item  or  at some scale $\ell^2 \le N\le
\ell^{10}$ there exist a regular positive slope $I$-segment
$\bigl\{E_{j_1}^{(N)}(x, \omega), \underline{x}_1,
\overline{x}_1 \bigr\}$, as well as  a
regular negative slope $I$-segment $\bigl\{E_{j_2}^{(N)}$ $(x,
\omega), \underline{x}_2, \overline{x}_2 \bigr\}$
with  $I\subset (E'_0, E''_0)$ and
$(\underline{x}_1, \overline{x}_1)\subset (x'_0,x''_0)$, with the same
$I$.
\end{itemize}
\end{prop}
\begin{proof}
Assume that the first alternative does not hold. Then, due to Lemma~\ref{lem:13.eitheruncondornonspectral}
one can assume that for any $\ell^2 \le N_1 \le \ell^{4}$ there exists a
regular $I_1$-segment $\bigl\{E_{j_1}^{(2N_1)}(x, \omega),
\underline{x}_1, \overline{x}_1 \bigr\}$ with
\[ I_1\subset (E'_0, E''_0),\quad (\underline{x}_1, \overline{x}_1)\subset (x'_0,x''_0), \quad
 \frac
{1}{4} N_1\le \nu^{(2N_1)}_{j_1}(x, \omega) \le \frac {1}{2} N_1 \] Since
$\omega\in \tor_{c,a}$, one can choose $N_1$ so that $\|N_1\omega\|\le
\kappa_0(V,c,a,\gamma)$. Then due to the part~$(2)$ of
Lemma~\ref{lem:13.hillsegemnts} there exists $j$ such that
\begin{equation}\label{eq:hillclosetodir2Nc}
|E_j^{(N_1,P)}(x,\omega)-E_{j_1}^{(2N_1)}(x,\omega)|\le
\exp(-N_1^{2\delta})
\end{equation}
for any $x\in (\underline{x}_1, \overline{x}_1)$. Assume for
instance that $\bigl\{E_{j_1}^{(2N_1)}(x, \omega) \bigr\}$ is a
positive-slope segment. Let \[x_1:=(\underline{x}_1 +
\overline{x}_1)/2,\qquad E_1:= E_{j_1}^{(2N_1)}(x_1, \omega) \]
Without loss of generality (because of the periodicity) we may assume
that $x_1$ is close to the middle of $\tor\simeq[0,1]$.
 Due to the $1$-periodicity
$E_{j_1}^{(N_1,P)}(0,\omega)=E_{j_1}^{(N_1,P)}(1,\omega)$. Hence,
either $ E_1 \le E_{j_1}^{(N_1,P)}(0,\omega)$, or $E_1 \ge
E_{j_1}^{(N_1,P)}(1,\omega)$. Assume for instance that the latter case
takes place. Let \[x_2:=(x_1 + \overline{x}_1)/2,\quad E_2:=
E_{j_1}^{(2N_1)}(x_2, \omega),\quad \overline{E}_1:=
E_{j_1}^{(2N_1)}(\overline{x}_1,\omega) \] Then $E_1+ \exp((-\log
N_1)^{2A}) < E_2,\; E_2 + \exp((-\log N_1)^{2A}) < \overline{E}_1$.
Hence,
\begin{align*}
 E_{j_1}^{(N_1,P)}(1,\omega)+ \exp((-\log N_1)^{3A}) &<
E_{j_1}^{(N_1,P)}(x_2,\omega) \\
 E_{j_1}^{(N_1,P)}(x_2,\omega)+
\exp((-\log N_1)^{3A}) &<
 E_{j_1}^{(N_1,P)}(\overline{x}_1,\omega)
\end{align*}
 Since $E_{j_1}^{(N_1,P)}(\cdot,\omega)$
  is continuous, there exists $x_0\in (\overline{x}_1,1)$ such that
\begin{itemize}
\item[(i)]  $E_{j_1}^{(N_1,P)}(x_0,\omega)=E_{j_1}^{(N_1,P)}(x_2,\omega)$
\item[(ii)] for any $x\in (\overline{x}_1,x_0)$ one has
$E_{j_1}^{(N_1,P)}(x,\omega) >$ $E_{j_1}^{(N_1,P)}(x_2,\omega)$
\end{itemize}
Note that $x_0$ was chosen to be the first point to the right of~$\overline{x}_1$ where the
graph of $E_{j_1}^{(N_1,P)}(\cdot,\omega)$ hits the level~$E_{j_1}^{(N_1,P)}(x_2,\omega)$.
Due to part $(1)$ of Lemma~\ref{lem:13.hillsegemnts} there exists
 an unconditional segment $\bigl\{E_{j_2}^{(2N_1)}(x, \omega),
\underline{x}_2, \overline{x}_2 \bigr\}$ with $x_0\in
(\underline{x}_2, \overline{x}_2)$ such that
\begin{equation}\label{eq:hillclosetodir2Nd}
|E_j^{(N_1,P)}(x,\omega)-E_{j_2}^{(2N_1)}(x,\omega)|\le
\exp(-N_1^{2\delta})
\end{equation}
for any $x\in (\underline{x}_2, \overline{x}_2)$. Because of (ii) above this must be a negative-slope segment.
To see this, assume that it is a positive
slope segment. Note that since $|E_{j(x_1)}^{(2N_1)}(x_0,\omega)-E_2|$
$\le \exp(-N_1^{\delta})$ one has
\[ \dist(E_{j(x_1)}^{(2N_1)}(x_0,\omega), {\cE}^{(1)}_{2N_1, \omega})
\ge \exp(-(\log N_1)^{4A}) \] Therefore,
\[
|\partial_x E_{j_2}^{(2N_1)}(x,\omega)|\ge \exp(-(\log N)^B)
\]
for any $x\in [x^{-}_0,x^{+}_0]$, where $x^{\pm}_0:=x_0 \pm
\exp(-(\log N_1)^{5A})$. Due to part $(1)$ of
Lemma~\ref{lem:13.hillsegemnts} one has
\begin{equation}\label{eq:hillclosetodir2Ne}
|E_j^{(N_1,P)}(x,\omega)-E_{j_2}^{(2N_1)}(x,\omega)|\le
\exp(-N_1^{2\delta})
\end{equation}
for any $x\in [x^{-}_0,x^{+}_0]$. In particular,
$E_j^{(N_1,P)}(x^{-}_0,\omega) <$ $E_{j_1}^{(N_1,P)}(x_0,\omega)$
$=E_{j_1}^{(N_1,P)}(x_2,\omega)$, which contradicts~(ii) above.
Thus,  $\bigl\{E_{j_2}^{(2N_1)}(\cdot , \omega),
\underline{x}_2, \overline{x}_2 \bigr\}$ is indeed an unconditional,  negative slope,
segment. Applying Proposition~\ref{prop:9.uncondreg}
to $\bigl\{E_{j_1}^{(2N_1)}(\cdot, \omega), \underline{x}_1,
\overline{x}_1 \bigr\}$ and  $\bigl\{E_{j_2}^{(2N_1)}(\cdot, \omega), \underline{x}_2, \overline{x}_2 \bigr\}$ concludes
the argument. The case $ E_1 \le E_{j_1}^{(N_1,P)}(0,\omega)$ is treated analogously.
\end{proof}

\section{The proof of Theorem~\ref{th:1.sinai}}\label{sec:sinai}

In this section we prove Theorem~\ref{th:1.sinai} on Sinai's
parametrization of eigenfunctions. Sinai needed to assume
cosine-like potentials and his argument was perturbative. Our construction 
applies to general analytic potentials under the
condition that $L(\omega,E)>0$. We derive Theorem~\ref{th:1.sinai}
from the following detailed finite volume version. 

\begin{prop}\label{prop:14.basic} Assume that $L(\omega_0, E) \geq \gamma > 0$
 for some $\omega_0 \in \tor_{c, a}$
and any $E\in(E', E'')$. Then there exist $\rho^{(0)} =
\rho^{(0)}(V, c, a, \gamma) >0$ and $N_0 = N_0(V, c, a, \gamma)$
such that for any $N \geq N_0$ there exists a subset $\Omega_N \subset
\tor$ so that for all $ \omega \in \tor_{c, a}
\cap  (\omega_0 - \rho^{(0)}, \omega_0 + \rho^{(0)})
\backslash \Omega_N$ there exist $\cE_{N, \omega} \subset \IR$, $\cB_{N, \omega}\subset \tor$ such that
the following statements hold:
\begin{enumerate}
\item[(1)]  $ \mes(\Omega_N) \leq \exp(-(\log N)^{A_0} ), \,\, \compl(
 \Omega_N) \leq N^{C_0}\\ \mes(\cE_{N, \omega}) \leq \exp(-(\log N)^{A_0} ), \,\, \compl(
  \cE_{N, \omega}) \leq N^{C_0}\\ \mes(\cB_{N, \omega}) \leq \exp(-(\log N)^{A_0} ), \,\,
  \compl( \cB_{N, \omega}) \leq N^{C_0}$

\item[(2)]  For any $x \in \tor \backslash \cB_{N, \omega}$ one has
\[
\spec \big( H_{[-N, N]}(x, \omega) \big) \cap \bigl( (E', E'') \cap \cE_{N,
\omega} \bigr) = \emptyset,
\]

\item[(3)]  For any $\oN \ge N$ one has
\[
 \cS_{\oN, \omega} \backslash \cS_{N, \omega} \subset \cE_{N,
\omega}
\]
where
\[
\cS_{N, \omega}:=\bigcup_{x\in \TT} \spec \big( H_{[-N, N]}(x, \omega) \big)
\]
 Furthermore, assume that for
 some $N\ge N_0$,
  $$\omega \in \tor_{c, a} \cap (\omega_0 - \rho^{(0)}, \omega^{(0)} + \rho^{(0)})
   \backslash \mybigcup_{N' \ge N} \Omega_{N'}$$
   Then the following further properties hold:
\bigskip
\item[(4)] Let $x \in \tor \backslash \cB_{N, \omega}$. If some eigenvalue $E_j^{(N)}(x, \omega)$
 falls into the interval $(E', E'')$, then there exists $\nu_j^{(N)}(x,\omega)\in [-N,N]$
 such that
\[|\psi_{j}^{(N)}(x,\omega,n)|\le\exp\big(-\frac{\gamma}{2}|n-\nu_j^{(N)}(x,\omega)|\big)
\]
for all $|n - \nu_j^{(N)}(x,\omega)|\ge N^{1/2}$
\item[(5)]  Set $\hat {\Omega}_N:=\bigcup_{N' \ge N} \Omega_{N'}$,
 $\hat{\cB}_{N, \omega} := \mybigcup_{N' \ge N} \cB_{N', \omega}
$. Let $\omega \in \tor_{c,a} \backslash \hat{\Omega}_{N}$ and let
$x \in \tor \backslash \hat{\cB}_{N, \omega}$, $E_j^{(N)}(x,
\omega)\in (E',E'')$. Assume that
\[\nu_j^{(N)}(x,\omega)\in [-N + N^{1/2},N-N^{1/2}]\] Then for each
$N'\geq N$ there exists $j_{N'}$ such that
\begin{equation}\label{1.evstabilitywithscale}
\begin{split}
\bigl| E_{j_{N'}}^{(N')} (x, \omega) - E_j^{(N)}(x, \omega) \bigr| &<
\exp(-\frac{\gamma}{2}N^{\frac{1}{2}})\\
 \bigl|\partial_x
E_{j_{N'}}^{(N')} (x, \omega) -
\partial_xE_j^{(N)}(x, \omega) \bigr| &<
\exp(-\frac{\gamma}{2}N^{\frac{1}{2}}) \\
\bigl|\psi_{j_{N'}}^{(N')} (x, \omega, n) \bigr| &\leq \exp
\bigl(-\frac{\gamma}{2} |n - \nu_j^{(N)}(x,\omega)|\bigr)
\end{split}
\end{equation}
for any $ |n|\leq N'$, and any $ N' \geq N$;
\begin{align*} |E_{j_{N''}}^{(N'')} (x, \omega) -E_{j_{N'}}^{(N')} (x, \omega) |
 &\leq \exp(- \frac{\gamma}{2}(N')^{\frac{1}{2}}) \\
| \psi_{j_{N''}}^{(N'')} (x, \omega, n) - \psi_{j_{N'}}^{(N')}(x, \omega, n) |
&\leq \exp(-\frac{\gamma}{2}(N')^{\frac{1}{2}})
\end{align*}
for any $ |n| \leq N'$, and any $N \le N' \le N''$
In particular, the limits
\begin{align*}
E(x, \omega) &:= \lim_{N' \to \infty} E_{j_{N'}}^{(N')} (x, \omega)
\\
\psi (x, \omega, n) &:= \lim_{N' \to \infty} \psi_{j_{N'}}^{(N')} (x,
\omega, n),\;\; n \in \IZ
\end{align*}
exist,
\begin{equation}
\label{eq:1.eigenvaladd} | E(x, \omega) - E_j^{(N)} (x, \omega) |
\leq 2\exp(-\frac{\gamma}{2}N^{\frac{1}{2}}), \quad | \psi (x,
\omega, n) - \psi_{j_{N}}^{(N)}(x, \omega, n) | \leq
\exp(-\frac{\gamma}{2}(N)^{\frac{1}{2}})
\end{equation}
 and $\sum\limits_n |\psi(x, \omega, n) |^2 = 1$. Furthermore,
\begin{align}\label{eq:1.eigvalform}
H (x, \omega) \psi (x, \omega, \cdot) = E(x, \omega) \psi(x, \omega, \cdot)
\\\label{eq:1.expdecform}
\lim_{n \to \infty}\frac {1}{2n}\log(|\psi(x, \omega, n) |^2
+|\psi(x, \omega, n+1)|^2)= -L(\omega,E(x, \omega))
\end{align}
\item[(6)]  Let $\omega \in \tor_{c,a} \backslash \hat{\Omega}_{N}$,
$x \in \tor \backslash \hat{\cB}_{N, \omega}$,
 and $E_{j_m}^{(N)}(x, \omega) \in (E', E'')$, $m=1, 2$, $j_{1} \neq j_{2}$. Then
\[
|E_{j_1}^{(N)}(x, \omega) - E_{j_2}^{(N)}(x,\omega) | \geq \exp(-
N^{\delta})
\]
Let $E_m(x, \omega)$ be the eigenvalue of $H(x, \omega)$, defined in
(5) above, which obeys
\[
|E_m(x, \omega) - E_{j_m}^{(N)}(x, \omega) | <
2\exp(-\frac{\gamma}{2}N^{\frac{1}{2}})
\] for
$m = 1,2$. Then
\[
|E_1(x, \omega) - E_2(x, \omega) | > \frac{1}{2} \exp(-N^{\delta}).
\]
If $(E', E'') = (-\infty, +\infty)$, then for each $|j| \leq N/2$ one
has an ``almost Parseval identity''
$$0 \leq 1 - \displaystyle\sum\limits_m | \la \delta_j (\cdot), \psi_m(x,
\omega, \cdot) \ra |^2 \leq \exp(-\frac{\gamma}{8}N)$$ where
$\delta_j(\cdot)$ stands for the $\delta$--function at $n=j$. The
collection of all eigenfunctions $\psi_m(x, \omega, \cdot)$
obtained this way for $N, N+1, \dots$ form a complete orthonormal
system in $\ell^2$.
\item[(7)]   Let $\omega \in \tor_{c,a} \backslash \hat{\Omega}_{N}$,
$x \in \tor \backslash \hat{\cB}_{N, \omega}$.  Let
$\{\psi_m(x,\omega,\cdot)\} $ be all possible eigenfunctions of
$H(x, \omega)$ defined in (5) above for all scales $\ge N$. Let
$J_{x, \omega}$ be be the closure of the set of the corresponding
eigenvalues. If $E \in (E', E'') \backslash J_{x, \omega}$, then
$(H(x, \omega) - E)$ is invertible. In other words,
\[
\spec \big( H(x, \omega) \big) \cap (E', E'') = J_{x, \omega},
\]
and the functions $\psi_m(x, \omega, \cdot)$ form a complete orthonormal
system in the spectral subspace of $H(x, \omega)$ corresponding to
$(E', E'')$.

\item[(8)] Assume that $\nu_j^{(N)}(x,\omega)\in [-N +
N^{1/2},N-N^{1/2}]$. Then for any $k$ such that \[-N + N^{1/2}/2 <
\nu_j^{(N)}\xo + k < N - N^{1/2}/2\] there exists a unique
\[E^{(N)}_{j_k}(x + k\omega, \omega)\in \spec \big( H_{[-N, N]}(x + k\omega,
\omega) \big)\] such that
\begin{align}
& \bigl|E_j^{(N)}\xo - E^{(N)}_{j_k}(x + k\omega, \omega)\bigr|<
\exp\bigl(-\gamma N^{1/2}/4\bigr)\
,\label{eq:1.14}\\[5pt]
& E^{(N)}_{j_k} (x + k\omega, \omega) \notin \cE_{N,\omega}\ ,\label{eq:1.15}\\[5pt]
& \bigl|\nu^{(N)}_{j_k} (x + k\omega, \omega) - \bigl(\nu_j^{(N)}
\xo + k\bigr)\bigr| \le N^{1/2}/4\
 ,\label{eq:1.16}\\[5pt]
& -N + N^{1/2}/4 < \nu^{(N)}_{j_k}(x + k\omega, \omega) < N - N^{1/2}/4\ ,\label{eq:1.17}\\[5pt]
&\sum_{|m+k - \nu^{(N)}_j\xo| \le N^{1/2}/4} \bigl|\psi^{(N)}_{j_k}
(x + k\omega, m) - \psi_j^{(N)}(x, m + k) \bigr|^2 <
\exp\bigl(-\gamma N^{1/2}/8\bigr)\label{eq:1.18}
\end{align}
\item[(9)]   Let $\omega \in \tor_{c,a} \backslash \hat{\Omega}_{N}$,
$x_0 \in \tor \backslash \displaystyle\hat{\cB}_{N, \omega}$. Then  for any $x \in \tor$ one has
\[
\spec \big( H(x, \omega) \big) \cap (E', E'') = J_{{x_0}, \omega}
\]
In other words, $ \spec\big( H(x, \omega) \big) \cap (E', E'')$ is the same for
all $x \in \tor$
\item[(10)]   For any $x \in \tor$ and any $N$ one has
\[
\bigl[ (\spec ( H(x, \omega)) \backslash \cS_{N, \omega}) \cup (\cS_{N,
\omega} \backslash \spec (H(x, \omega) )) \bigr] \cap (E', E'') \subset
\mybigcup_{N' \geq N} \cE_{N', \omega}
\]

\item[(11)]
 If $\spec ( H(x, \omega))\cap (E', E'')\neq \emptyset$ for some $x$ then
\[
\mes (\spec (H(x, \omega))) \cap (E', E'') > 0
\]
If $(E', E'') = (-\infty, \infty)$ then
\[
\mes( \spec (H(x, \omega))) \geq \exp(-(\log N_0)^{C_1})
\]
\end{enumerate}
\end{prop}

Most of the preparatory work needed for the proof of this
proposition has already been done in the previous sections.
We need only few more auxiliary statements.

\begin{lemma}\label{lem:10.specn}
There exists $N_0=N_0(V,c,a,\gamma)$ such that for any $N \geq
N_0$, $\omega \in \tor_{c,a}$, $E\in(\uE,\oE)$ the following property
holds: if $\dist(E,\cS_{N,\omega})\geq\exp(-N^{\frac{1}{2}})$, then
$\dist(E,\cS_{\oN,\omega})\geq\frac{1}{2}\exp(-N^{\frac{1}{2}})$ for
any $\oN \geq N$.
\end{lemma}
\begin{proof}
Due to Corollary \ref{cor:2.uniexcepzero}, $\dist\big[E, \spec (
H_N(x,\omega)) \big]\geq \frac{1}{2}\exp(-N^{\frac{1}{2}})$ implies
\begin{equation}\label{10.lyap}
\log|f_{[1,N]}(e(x),\omega,E)| > \gamma N - N^{\frac{3}{4}}
\end{equation}
provided $N\geq N_0(V, c, a, \gamma)$.  Due to Lemma~\ref{lem:3.Green}, \eqref{10.lyap} implies
\[ \lvert (H_{[1,N]}(x,\omega)-E)^{-1}(m,n)\rvert\leq\exp(-\frac{\gamma}{4}N)
\]
for any $m,n\in[1,N]$, $|m-n|\geq\frac{N}{2}$. Therefore,
$\dist(E,\cS_{N,\omega})\geq\frac{1}{2}\exp(-N^{\frac{1}{2}})$
implies
\begin{equation}\label{10.green}
\lvert(H_{[N'+1,N'+N]}(x,\omega)-E)^{-1}(m,n)\rvert\leq\exp(-\frac{\gamma}{4}N)
\end{equation}
for any $x\in\tor$, $N'\in\IZ$, $m,n\in [N'+1,N'+N]$,
$|m-n|\geq\frac{N}{2}$. Let $N \leq\oN \leq
\exp(N^{\frac{3}{4}})$. Assume that
\begin{equation}\label{10.eigen}
H_{[1,\oN]}(x,\omega)\psi(x,n)=E\psi(x,n)
\end{equation}
Let $\mu \mathrel{\mathop{:}=} \max_{1\le n \le \oN} |\psi(x,n)| $.
Then, due to \eqref{10.green} and the ``Poisson formula'' one
obtains
\[ |\psi(x,n)| \leq 2\mu \exp(-\frac{\gamma}{4}N)\text{, for any } n \in [1,\oN]
\]
Hence,
\[ \mu \leq 2\mu \exp(-\frac{\gamma}{4}N)
\]
This yields $\mu=0$, and thus $\psi(x,n)=0$ for any $n\in [1,\oN]$.
Thus
\begin{equation}\label{10.specaux}
\dist(E,\cS_{N,\omega})\geq\frac{1}{2}\exp(-N^{\frac{1}{2}}) \Longrightarrow
E\not\in\spec  \big(H_{[1,\oN]}(x,\omega) \big) \quad \forall\; x \in \TT
\end{equation}
The lemma  follows from \eqref{10.specaux}.
\end{proof}

\begin{lemma}\label{lem:14.greendecay}
There exists $N_0=N_0(V,c,a,\gamma,A)$ such that for any $N \geq
N_0$, $\omega \in \tor_{c,a}$, $E\in(\uE,\oE)$ the following
property holds: if $\dist \big[ E,\spec\big( H_{[1,N]}(x,\omega) \big) \big] \geq\exp(-(\log
N)^A)$, then
\begin{equation}\label{14.greendecay1}
\lvert(H_{[1,N]}(x,\omega)-E)^{-1}(m,n)\rvert\leq\exp(-L(\omega,E)|m-n|
+(\log N)^{2A})
\end{equation}
for any $m,n\in [1,N]$.
\end{lemma}
\begin{proof} Due to Corollary \ref{cor:2.uniexcepzero}, $\dist\big[E, \spec \big(
H_N(x,\omega) \big) \big] \geq \exp(-(\log N)^A)$ implies
\begin{equation}\label{14.detlyaplowerest}
\log|f_{[1,N]}(e(x),\omega,E)| > L(\omega,E)N - (\log N)^{A+C}
\end{equation}
provided $N\geq N_0(V, c, a, \gamma,A)$. As in the proof
of Lemma~\ref{lem:3.Green} one sees that  \eqref{14.detlyaplowerest}
implies \eqref{14.greendecay1}. 
\end{proof}

\begin{lemma}\label{lem:14.eigenfunctionlowerest}
Given $\eps>0$ sufficiently small, there exists $N_0=N_0(V,c,a,\gamma,\eps)$
 such that for any $N \geq
N_0$, $\omega \in \tor_{c,a}$, $x\in \tor$, $E\in(\uE,\oE)$ the
following assertion holds. Assume that
\begin{equation}\label{14.generaleigenf}
H (x, \omega) \psi (\cdot) = E \psi(\cdot)
\end{equation}
for some $E$ and some function $\psi(n)$, $n\in \ZZ$. Assume
also that
\begin{equation}\label{14.geneigenfnontriv}
\max_{1\le n \le N} |\psi (n)| = 1
\end{equation}
Then
\begin{equation}\label{14.geneigenfnontrivmin}
\min_{1\le n \le N} \frac {1}{2}\log (|\psi (n)|^2+ |\psi
(n-1)|^2)\ge -N(L(\omega,E)+\eps)
\end{equation}
\end{lemma}
\begin{proof} Recall that for any $a<b$
\begin{equation}\label{14.basicpropagation}
\begin{bmatrix} \psi(b+1)\\ \psi(b) \end{bmatrix} =M_{[a, b]}(e(x), \omega, E)
\begin{bmatrix} \psi(a)\\ \psi(a-1) \end{bmatrix}
\end{equation}
where $M_{[a, b]}(e(x), \omega, E)$. Recall
also that due to the uniform upper estimate of
Lemma~\ref{lem:2.lipsub} for any $E\in(\uE,\oE)$ 
\[
\sup_{0<b-a\le N, x\in \tor} \|M_{[a, b]}(e(x), \omega, E)\| \le
\exp(N(L(\omega,E)+\eps))
\]
provided  $N\ge N_0=N_0(V,c,a,\gamma,\eps)$. Since
\[
\|M^{-1}\|=\|M\|
\]
for any unimodular matrix $M$, one obtains
\begin{equation}\label{14.propagatestmates}
\exp(-N(L(\omega,E)+\eps)) \Big\|\begin{bmatrix} \psi(a)\\ \psi(a-1)
\end{bmatrix}\Big\| \le \Big\|\begin{bmatrix} \psi(b+1)\\ \psi(b) \end{bmatrix}\Big\|
\le \exp(N(L(\omega,E)+\eps)) \Big\|\begin{bmatrix} \psi(a)\\ \psi(a-1)
\end{bmatrix}\Big\|
\end{equation}
for any $0<b-a\le N$. Due to the assumptions of the lemma there
exist $a\in [1,N]$ such that
\[
1\le \Big\|\begin{bmatrix} \psi(a)\\ \psi(a-1)
\end{bmatrix}\Big\|\le 2
\]
Therefore the statement follows from \eqref{14.propagatestmates}. 
\end{proof}

We now turn to the proof of Proposition~\ref{prop:14.basic}.

\begin{proof}[Proof of Proposition~\ref{prop:14.basic}]
 \begin{itemize}
 \item[(1)]
 Let $\Omega_N,\cE_{N,\omega}$ be the subset defined in Lemma~\ref{lem:3.waschno}.
 Let $\omega\in\tor_{c,a}$.
 Set \[
\cB_{N,\omega}=\big\{x\in\tor:\spec \big(H_{[-N,N]}(x,\omega)\big) \cap(E',E'')\cap\cE_{N,\omega}\neq \emptyset \big \}
\]
The subsets $\Omega_N, \cE_{N,\omega},\cB_{N,\omega}$ 
satisfy properties  $(1)$, $(2)$ of Proposition~\ref{prop:14.basic}.
\item[(3)] Clearly,  $\cS_{N,\omega}$
is a union of $(2N+1)$ closed intervals. Set
\[
\cE'_{N,\omega}=\{E\in\IR:E\not\in \cS_{N,\omega},
\dist(E,\cS_{N,\omega})\le \exp(-N^{\frac{1}{2}})\}
\]
Note that $\mes (\cE'_{N,\omega}) \lesssim N\exp(-N^{\frac{1}{2}}),
\, \compl (\cE'_{N,\omega}) \lesssim N$.
 With some abuse of notation we denote the set $\cE_{N,\omega}\cup \cE_{N',\omega}$ as $\cE_{N,\omega}$.
  Clearly $\cE_{N,\omega}$ obeys property (1) of Proposition~\ref{prop:14.basic}.
   Due to Proposition~\ref{lem:10.specn}, property (3) of Proposition~\ref{prop:14.basic} holds.
 \item[(4)]
Let $E_j^{(N)}(x,\omega)\in (E',E'')$ for some $x\in\tor \setminus
\cB_{N,\omega}$. Then due to Lemma~\ref{lem:3.waschno} there exists
$\nu_j^{(N)}(x,\omega)\in[-N,N]$ such that
\[
|\psi_j^{(N)}(x,\omega,n)|\le
\exp\bigl(-\frac{\gamma}{2}|n-\nu_j^{(N)}(x,\omega)|\bigr)
\]
provided $|n-\nu_j^{(N)}(x,\omega)|>N^{\frac{1}{2}}$. So, part (4)
of Proposition~\ref{prop:14.basic} is valid.
\item[(5)]

 Let $\omega \in \tor_{c,a} \backslash \hat{\Omega}_{N}$,
 $x\in\tor\setminus$ ${\hat{\cB}}_{N,\omega}$.
 Assume $\nu_j^{(N)}(x,\omega)\in[-N+N^{\frac{1}{2}}, N-N^{\frac{1}{2}}]$.
 Then due to standard perturbation theory of Hermitian matrices, for each $N\le N'\le N$ there exists an eigenvalue
 $E_j^{(N')}(x,\omega)$ of $H_{[-N', N']}(x,\omega)$ such that the estimates
 \eqref{1.evstabilitywithscale}
  hold. Assume that $x\in\tor\setminus\hat{\cB}_{N,\omega}$. Then
\eqref{1.evstabilitywithscale}
 applies for any $N'\ge N$. Therefore the limits
\begin{align*}
E(x,\omega) &=\lim\limits_{N' \to \infty} E_{j_{N'}}^{(N')}(x,\omega)
\\
\psi( x,\omega,n) &=\lim\limits_{N' \to
\infty}\psi_{j_{N'}}^{(N')}(x,\omega,n),\quad n\in\IZ
\end{align*}
exist, relations \eqref{eq:1.eigenvaladd}, \eqref{eq:1.eigvalform}
hold and
\begin{equation}\label{eq:14.normalization}
\sum\limits_n |\psi(x, \omega, n) |^2 = 1
\end{equation}
Pick  $N_1\in \ZZ$ so that $N\asymp(\log N_1)^B$. Since
$\omega \in \tor_{c,a} \backslash \hat{\Omega}_{N}$,
 $x\in\tor\setminus$ ${\hat{\cB}}_{N,\omega}$ due to
Proposition~\ref{prop:4.elimination} (applied with $t=2N$) one has with A=4B
\begin{equation}
\nn
\begin{split}
\dist \big[E_{j_{N_1}}^{(N_1)}(x,\omega),\spec \big( H_{[3N,N_1]}(x,\omega) \big) \big] \geq
\exp(-(\log N_1)^A), \\
\dist \big[ E_{j_{N_1}}^{(N_1)}(x,\omega),\spec \big(
H_{[-N_1,-3N]}(x,\omega)\big) \big]\geq \exp(-(\log N_1)^A)
\end{split}
\end{equation}
That implies
\begin{equation}\label{eq:14.specdistaux}
\begin{split}
\dist \big[ E(x,\omega),\spec \big( H_{[3N,N_1]}(x,\omega)\big) \big]\geq
\exp(-(\log N_1)^A), \\
\dist \big[ E(x,\omega),\spec \big( H_{[-N_1,-3N]}(x,\omega)\big) \big]\geq \exp(-(\log
N_1)^A)
\end{split}
\end{equation}
Due to Lemma~\ref{lem:14.greendecay} relation
\eqref{eq:14.specdistaux} in its turn implies the following
estimates
\begin{equation}\label{14.greendecay3}
\begin{split}
\lvert(H_{[3N,N_1]}(x,\omega)-E(x,\omega))^{-1}(m,n)\rvert\leq\exp(-L(\omega,E(x,\omega))|m-n|
+(\log N_1)^{2A}) \\
\lvert(H_{[-N_1,-3N]}(x,\omega)-E(x,\omega))^{-1}(m,n)\rvert\leq\exp(-L(\omega,E(x,\omega))|m-n|
+(\log N_1)^{2A})
\end{split}
\end{equation}
for any $m,n$ in the corresponding interval. Let $n$ be arbitrary
such that \[(\log N_1)^{3A}\le |n| \le N_1/2\] Applying 
Poisson's formula to $\psi(x, \omega, n)$ one obtains
\begin{equation}
\nn
\begin{split}
|\psi(x,\omega,n)| &\le \exp\big[ -L(\omega,E(x,\omega)) (|n|-3N -(\log
N_1)^{2A})  \big]  \\
&+\exp \big[-L(\omega,E(x,\omega)) (N_1 - |n| -(\log
N_1)^{2A}) \big] \\
&\lesssim \exp\big[-L(\omega,E(x,\omega)) |n|(1-o(1)) \big]
\end{split}
\end{equation}
as $n\to\infty$. 
This estimate implies
\begin{equation}\label{eq:14.upexpdec}
\begin{split}
 \limsup_{n \to \infty}\frac {1}{2n}\log(|\psi(x, \omega, n) |^2
+|\psi(x, \omega, n+1)|^2)\le -L(\omega,E(x, \omega))
\end{split}
\end{equation}
Note that using the above notations one has
\[
\max_{|n| \le N}|\psi(x, \omega, n) |\ge 1/N
\]
since $|\psi(x, \omega, n) |$ is normalized. Thus, in view of
Lemma~\ref{lem:14.eigenfunctionlowerest}, 
\begin{equation}
\nn
 \min_{-N\le n \le N_1} \frac {1}{2}\log (|\psi(x,\omega,n)|^2+ |\psi
(n-1)|^2)\ge -N(L(\omega,E)+\eps)
\end{equation}
provided $N_1\ge N_0(V,c,a,\gamma,\eps)$. Hence
\[
 \liminf_{n \to \infty}\frac {1}{2n}\log(|\psi(x, \omega, n) |^2
+|\psi(x, \omega, n+1)|^2)\ge -L(\omega,E(x, \omega))-\eps
\]
Since $\eps>0$ is arbitrary
\[
 \liminf_{n \to \infty}\frac {1}{2n}\log(|\psi(x, \omega, n) |^2
+|\psi(x, \omega, n+1)|^2)\ge -L(\omega,E(x, \omega))
\]
Thus
\[
 \lim_{n \to \infty}\frac {1}{2n}\log(|\psi(x, \omega, n) |^2
+|\psi(x, \omega, n+1)|^2)= -L(\omega,E(x, \omega))
\]
Similarly
\[
 \lim_{n \to -\infty}\frac {1}{2n}\log(|\psi(x, \omega, n) |^2
+|\psi(x, \omega, n+1)|^2)= -L(\omega,E(x, \omega))
\]
as claimed.

 \item[(6)] Let $\omega \in \tor_{c,a} \backslash \hat{\Omega}_{N}$,
 $x\in\tor\setminus\hat{\cB}_{N,\omega}$.
Assume that $E_{j_m}^{(N)}(x,\omega)\in(E',E'')$, and
\[ \nu_{j_m}^{(N)}(x,\omega) \in[-N+N^{\frac{1}{2}},
N-N^{\frac{1}{2}}],\;\;m=1,2\] Then due to
Proposition~\ref{prop:4.Ej_sep} one has
\begin{equation}\label{eq:11.esplit}
\bigl|E_{j_1}^{(N)}(x,\omega)-E_{j_2}^{(N)}(x,\omega)\bigr|>\exp(-N^\delta)
\end{equation}
Let $E_m(x,\omega), \psi_m(x,\omega,\cdot)$ be defined as in (5) for
$j=j_m$, $m=1,2$. Then, from~\eqref{eq:1.eigenvaladd}
and~\eqref{eq:11.esplit}, 
\begin{equation}\label{eq:10.efront}
\bigl|E_1(x,\omega) - E_2(x,\omega)\bigr| >
\frac{1}{2}\exp(-N^\delta)
\end{equation}
If $(E',E'')=(-\infty,\infty)$, then part (5) is applicable to each
eigenvalue $E_j^{(N)}(x,\omega)$, provided
\begin{equation}\label{eq:10.regev}
\nu_j^{(N)}(x,\omega)\in[-N+N^{\frac{1}{2}}, N-N^{\frac{1}{2}}].
\end{equation}
 Furthermore, let
$\psi_{j_m}^{(N)}(x,\omega,\cdot)$ be all eigenfunctions of
$H_{[-N,N]}(x,\omega)$ with
\[-N+N^{\frac{1}{2}}<\nu_{j_m}^{(N)}(x,\omega)<N-N^{\frac{1}{2}}\] Let
$\psi_m(x,\omega,\cdot)$ be the eigenfunction defined in part (5)
for $j=j_m$ where $1\le j\le m_0\le 2N+1$. Let $\delta_k(\cdot)$ be the
delta-function at $n=k$, $k\in [-N+4N^{1/2}, N-4N^{1/2}]$. Clearly,
one has \[|\la\delta_k, \psi_{j}^{(N)}(x,\omega,\cdot)\ra|\le
\exp(-\gamma N^{1/2}/2) \] for any $j$ with $\nu_j^{(N)}(x,\omega)\notin$
$[-N+N^{1/2}, N-N^{1/2}]$. Since $\big\{\psi_{j}^{(N)}(x,\omega,\cdot)\}_{j=1}^{2N+1}$,
form an orthonormal basis in the space of all
functions on $[-N,N]$, one has
\[
0\le 1-\sum\limits_{m=1}^{m_0}|\bigl\la\delta_k(\cdot),
\psi_m(x,\omega,\cdot)\bigr\ra|^2\le
\exp(-\frac{\gamma}{2}N^{\frac{1}{2}})
\]
 \item[(7)] Let $J_{x,\omega}$ be the closure of the set
of all eigenvalues $E_m(x,\omega)$ of $H(x,\omega)$ defined in~(5).
Let  $E_0\in(E',E'')$. Assume $\sigma_0:=\dist(E_0,J_{x,\omega})>0$.
Then due to the definition of the eigenvalues $E_m(x,\omega)$ one
has for any $N\ge N_0$
\[\qquad\qquad
\min\Big\{|E_0-E_j^{(N)}(x,
\omega)|:E_j^{(N)}(x,\omega)\in(E',E''),
\nu_j^{(N)}(x,\omega)\in[-N+N^{\frac{1}{2}},
 N-N^{\frac{1}{2}}]\Big\}\ge \sigma_0/2
\]
Let $\varphi(n), \, n\in\IZ$ be an arbitrary normalized
$\ell^2(\IZ)$ function supported on some interval $[-T, T]$. Let
$N>2T$. Then $|\la\varphi, \psi_{j}^{(N)}(x,\omega,\cdot)\ra|$ $\le
\exp(-N^{2/3})$ for any $j$ with $\nu_j^{(N)}(x,\omega)\notin$
$[-N+N^{1/2}, N-N^{1/2}]$. Hence,
 \[
\sum_{\nu_j^{(N)}(x,\omega)\in [-N+N^{1/2}, N-N^{1/2}]}|\la\varphi,
\psi_{j}^{(N)}(x,\omega,\cdot)\ra|^2 \ge 1/2
\]
Expanding $\varphi(\cdot)$ in the basis
$\{\psi_{j}^{(N)}(x,\omega,\cdot)\}_{j=1}^{2N+1}$  one obtains
 \[
 \|(H_{[-N, N]}(x,\omega)-E_0)\varphi\|^2\ge
(1/2)(\sigma_0/2)^2
\] Hence
\[\
\|(H(x,\omega)-E_0)\varphi\|^2 \ge (1/2)(\sigma_0/2)^2
 \]
 Since
$\varphi$ here is arbitrary one infers that 
\[
\dist (H(x, \omega), E_0) \gtrsim \sigma_0
\]
Thus part $(7)$ holds.
 \item[(8)] Part $(8)$ is due to
Lemma~\ref{lem:4.transl}.
 \item[(9)] Given $x$,  let $\cM (x, N, s)$ be the
collection of all $m$ such that the eigenfunction
$\psi_{j_m}^{(N)}(x,\omega,\cdot)$ defined as in $(5)$ obeys
\[ -N+sN^{\frac{1}{2}}<\nu_j^{(N)}(x,\omega)<N-sN^{\frac{1}{2}},\quad s=1,2,3\]
 Let $J_{x,\omega, s}$ be the closure of all eigenvalues
$E_m(x, \omega)$ defined in $(5)$ for all $N$ and $m\in \cM(x, N,
s)$. Then just as above one has $\spec( H(x, \omega))\cap (E',E'')$
$=J_{x,\omega, s}$. On the other hand, fix $x_0\in\tor$
$\setminus\hat{\cB}_{N,\omega}$, and take arbitrary $k$ such that
$-N/2 < \nu_j^{(N)}(x, \omega) + k < N/2$. Let $x=x_0+k\omega\;(\mod 1)$. Then
$x_0,x\in\tor\setminus\hat{\cB}_{2N,\omega}$. Due to part $(8)$ and
part $(5)$ for each $m\in \cM (x_0, 2N, 3)$, there exists $m'\in \cM
(x, 2N, 2)$ such that
\[E_{j_m}^{(2N)}(x_0,\omega)-E_{j_{m'}}^{(2N)}(x,\omega)\le
\exp(-N^{1/3})\]
 Hence
$J_{x_0,\omega, 3}\subset\{E\;:\; \dist(E,J_{x,\omega, 2}) \le
\exp(-N^{1/3})\}$. Thus
$$J_{x_0,\omega}=J_{x_0,\omega, 3}\subset J_{x,\omega, 2} = J_{x,\omega}$$
Switching the roles of $x_0$ and $x$ in this argument implies 
$J_{x_0,\omega}=J_{x,\omega}$ provided \[x_0\in\tor
\setminus\hat{\cB}_{N,\omega},\quad x=x_0+k\omega\;(\mod 1)\] for some $|k|\le
N/2$. Since $\hat{\cB}_{N',\omega}$ $\subset \hat{\cB}_{N,\omega}$
for $N'\ge N$ the claim is valid for any $x=x_0+k\omega\;(\mod 1)$ with
arbitrary $k$. Given arbitrary $x'$ and $\rho
>0$,   one can find $k$ such that with $x=x_0 + k\omega\;(\mod 1)$ one has
 $|x'-x| < \rho$. Then
\[ \spec ( H(x', \omega))\subset \{E: \dist\big( E, \spec (H(x,
\omega)) \big)\le \rho\}\] and \[\spec (H(x, \omega))\subset \{E:
\dist\big(E, \spec H(x', \omega)\big)\le \rho\}\] Since $\spec (H(x,
\omega)) =J_{x_0,\omega}$, and $\rho$ is arbitrary, one has $\spec
(H(x', \omega)) =J_{x_0, \omega}$. Thus $(9)$ is valid.

\item[(10)] Part $(10)$ follows from $(7)$ combined with part $(3)$.

\item[(11)] If $\spec( H(x, \omega))\cap (E', E'')\neq \emptyset$ for some
$x$, then $J_{x_0,\omega}\cap (E', E'')\neq \emptyset$ for
some $x_0\in\tor\setminus\hat{\cB}_{N,\omega}$ with some $N$.
Therefore, using the notations from the proof of part $(9)$ one has
\[J_{x_0,\omega, 1}\cap (E', E'')\neq \emptyset\]
Hence
$E_{j_{m}}^{(2N)}(x_0,\omega)\in (E', E'')$. It follows from 
Proposition~\ref{prop:11.4rellich} that  $\cS_{N, \omega}\cap (E',
E'')$ contains an $I$-segment, $|I|\ge \exp(-(\log N)^{C})$. Due to
part $(10)$ this implies \[\mes (\spec (H(x, \omega))) \cap (E', E'') \ge
\exp(-(\log N)^{C})/2\] which proves the first claim of $(11)$.
The second claim 
$(11)$ is implicit in the argument leading to the first part and we are done. 
\end{itemize}
\end{proof}

Before we proceed with the proof of Theorem~\ref{th:1.sinai} we make
the following remark which we will use in the proof of
Theorem~\ref{th:1.mainth}. It follows from inspection of the proof of assertion $(7)$ of
Proposition~\ref{prop:14.basic}.

\begin{remark}\label{rem:14.spectrumfree}
Let $\omega \in \tor_{c,a} \backslash \hat{\Omega}_{N_1}$,
$x_0\in\tor\setminus\hat{\cB}_{N_1,\omega}$, $E_0\in \IR$, and let
$\sigma_0>0$ be a constant. Assume that for any $N_2\ge N_1$ there
exists $N\ge N_2$ such that
\begin{equation}
\label{eq:17.nonspec1}
\begin{split}
\min\Big\{|E_0-E_j^{(N)}(x_0, \omega)|\::\: E_j^{(N)}(x_0,\omega)\in(E',E''),
\nu_j^{(N)}(x_0,\omega)\in[-N+N^{\frac{1}{2}},
 N-N^{\frac{1}{2}}] \Big\}\ge \sigma_0
\end{split}
\end{equation}
Then for any $x\in \TT$ one has 
$E_0\notin $ $\spec (H(x,\omega))$.
\end{remark}

\noindent To proceed we need the following general statement.

\begin{lemma}
\label{lem:17.dioph} Let $\omega\in \tor_{c,a}$. Given $x\in \tor$
and $N\in \ZZ$, $N>0$ set
\[
\mu(x,\omega,N):=\min_{t\in \ZZ, |t|\le N} \|x-t\omega\|
\]
Let $0<\alpha < \beta < 1$ be arbitrary. There exists
$N_0=N_0(c,a,\alpha,\beta)$ such that the following statement holds:
Let $x\in \tor$, and $x\notin \omega \ZZ\;(\mod 1)$. Assume that
$\|x-t_1\omega\|\le \exp(-|t_1|^{\beta})$ for some $t_1\in \ZZ$,
 $|t_1|\ge N_0$. Then there exists $s > t_1$ such that for any
 $s/4\le s_1 \le s$ holds
\[
\exp(-s^{\beta})\lesssim \mu(x,\omega,s_1)\lesssim \exp(-s^{\alpha})
\]
\end{lemma}
\begin{proof} Given arbitrary $N_0$ set ($\cR$ here stands for ``recurrence'') 
\[
\cR(x,\omega,N_0):=\{t\in \ZZ: |t|\ge N_0, \|x-t\omega\|\le
\exp(-|t|^{\beta})\}
 \]
If $t'\in \cR(x,\omega,N_0)$ and  $t\in \ZZ$, $|t|\le |t'|$,
$t\neq t'$, then $\|x-t\omega\|\ >$ $\|x-t'\omega\|$ provided $N_0$
is large enough. Indeed, otherwise
\[
\|(t-t')\omega\|\le 2\exp(-|t'|^{\beta})
\]
and $|t'-t|\le 2|t'|$ which contradicts the condition $\omega\in
\tor_{c,a}$. Hence for each $t'\in \cR(x,\omega,N_0)$ 
\[
\mu(x,\omega,|t'|)=\|x-t'\omega\|
\]
Assume first that $\cR(x,\omega,N_0)$ consists of a single integer
$t_1$. Let $T:=|t_1|$. Note that $\mu(x,\omega,T)>0$ since $x\notin
\omega \ZZ\;(\mod 1)$. Set $s:=[(\log \mu(x,\omega,T)^{-1})^{\frac {1}
{\alpha}}]$. Then 
\[
s\ge [ |t_1|^{\frac {\beta}{\alpha}}] > 4|t_1|=4T
\]
because of  $\alpha < \beta$. Let $t\in \ZZ$ be such that $T < |t|\le s$.
Since $t_1$ is the only element
of $\cR(x,\omega,N)$ one has $t\notin \cR(x,\omega,N)$.  Hence,
\[\|x-t\omega\| > \exp(-|t|^{\beta}) \ge  \exp(-s^{\beta})
\]
If $|t|\le T$, then  due to the definition of $\mu(x,\omega,T)$
one has
\[
\|x-t\omega\|\ge \mu(x,\omega,T)\ge \exp(-(s+1)^{\alpha})\asymp 
\exp(-s^{\alpha})
\]
Since $\mu(x,\omega,T)=\|x-t_{1}\omega\|\asymp \exp(-s^{\alpha})
 $ one concludes that
\[
\exp(-s^{\beta})\lesssim \mu(x,\omega,s_1)\lesssim \exp(-s^{\alpha})
\]
for any $T<s_1\le s$. Since $T < s/4$, the lemma follows  in this
case. Assume now that there exist at least two points in
$\cR(x,\omega,N)$, viz.~$t_i\in$ $\cR(x,\omega,N_0)$,
$i=1,2$, $t_1\neq t_2$. We can assume that $|t_1|< |t_2|$ and also
that
\[
|t_2|=\min \{|t|: t\in \cR(x,\omega,N_0)\setminus \{t_1\}, |t|\ge
|t_1| \}
\]
Then
\[
\|(t_1-t_2)\omega\|\le 2\max_i \|x-t_i\omega\| = 2 \|x-t_1\omega\|
\]
since
\[ \|x-t_1\omega\|=\mu(x,\omega,|t_1|)\ge \mu(x,\omega,|t_2|)=\|x-t_2\omega\|
\]
Since $\omega\in \tor_{c,a}$, that implies in particular
\[
|t_2-t_1|\ge \|x-t_1\omega\|^{-1/3}\ge  \exp(-\frac
{1}{3}|t_1|^{\beta})
\]
provided $N_0$ is large enough. Hence
\[
|t_2|\ge \frac {1}{2}\|x-t_1\omega\|^{-1/3}\ge \frac {1}{2}
\mu(x,\omega,|t_1|)^{-1/3}
\]
Define $T$ and $s$ just as before. Then $4|t_1|<s<|t_2|$. Moreover, if
$|t_1|<|t|\le s$, then $t\notin \cR(x,\omega,N_0)$ due to our choice
of $t_2$. Therefore, one can just repeat the argument from the case
$\cR(x,\omega,N_0)=\{t_1\}$.
\end{proof}
\begin{proof}[Proof of
Theorem~\ref{th:1.sinai}] We will follow the notations of
Proposition~\ref{prop:14.basic}. Set
\[
\Omega:=\bigcap_{N\ge N_0}\Omega_N, \quad
\cB_{\omega}:=\bigcap_{N\ge N_0}\bigcup_{|k|\le N}(\tilde {\cB}_{N,
\omega}+k\omega) \;(\mod 1)
\]
where
\[
\tilde {\cB}_{N, \omega}:= \{x\in \tor: \dist(x,\hat {\cB}_{N,
\omega}) \le \exp(-N^{1/5})  \}
\]
These sets are decreasing and it follows from part $(1)$ of Proposition~\ref{prop:14.basic}
that they are of Hausdorff dimension zero. Moreover, The set
$\cB_{\omega}$ is invariant under the shifts $x\mapsto 
x+m\omega$, $m\in \ZZ$. Let
\[\omega\in
(\tor_{c,a}\cap (\omega_0-\rho^{(0)},\omega_0+\rho^{(0)}))\setminus \Omega,
\]
Due to part $(11)$ of Proposition~\ref{prop:14.basic} $\mes (
\Sigma_{\omega}) >0$. Due to part $(5)$ of
Proposition~\ref{prop:14.basic} for any any $x \in \tor \backslash
\cB_{\omega}$ eigenfunctions $\psi_{j}(x,\omega,\cdot)$, $j=1,...$
are defined,
\[
\lim_{|N|\rightarrow \infty} \frac {1}{2N} \log
(|\psi_j(x,\omega,N)|^2 + |\psi_j(x,\omega,N+1)|^2) =
-L(\omega,E_j(x,\omega))
\]
and functions $\psi_{j}(x,\omega,\cdot)$ form an orthonormal basis
in $\ell^2(\ZZ)$. Moreover, the eigenvalues $E_j(x,\omega)$ are
simple. Let $E_0\in \IR$ be arbitrary. Assume that there exist
$x(k)\in \tor \setminus \cB_{\omega}$, $1\le k\le k_0\in \ZZ^+\cup\{\infty\}$ such that
\begin{itemize}
\item[(i)] for each $k$ there exists $j(k)$ such that
$E_0=E_{j(k)}(x(k))$
\item[(ii)]the orbits $\Gamma(x(k))$ are all different, i.e., 
for any $k_1\neq k_2$ and any $t\in \ZZ$ one has $x(k_2)\neq
x(k_1)+t\omega\;(\mod 1)$
\end{itemize}
To finish the proof of the theorem we have to evaluate $k_0$. There
exists $N_0$ such that $x(k)\in \tor \setminus \tilde
{\cB}_{N,\omega}$ for each $1\le k\le k_0$ and all $N\ge N_0$. It follows from parts
$(5)$ and (6) of Proposition~\ref{prop:14.basic} that for each $k$ there
exists $j_k$ such that
$-N+N^{\frac{1}{2}}<\nu_{j_k}^{(N)}(x,\omega)<N-N^{\frac{1}{2}}$ and
\[
|E_0-E_{j_k}^{(N)}(x(k),\omega)|\le \exp(-N^{1/3})
\]
Fix arbitrary $k$. Due to part $(8)$ of
Proposition~\ref{prop:14.basic} for each $s\in \ZZ$ with
$$-N+N^{\frac{1}{2}}<\nu_{j_k}^{(N)}(x,\omega)+s<N-N^{\frac{1}{2}}$$
(we call these $s$ admissible in this proof) there exists $j_{k,s}$ such that
\begin{equation}\label{eq:14.closeev}
|E_0-E_{j_{k,s}}^{(N)}(x(k)+s\omega,\omega)|\les \exp(-N^{1/3})
\end{equation}
Set $z_{k,s}:=e(x(k)+s\omega)$. Due to \eqref{eq:14.closeev} and
Corollary~\ref{cor:2.uniexcepzero} there exists $\zeta_{k,s}\in
\cD(z_{k,s},\exp(-N^{1/3}))$ such that
$f_N(\zeta_{k,s},\omega,E_0)=0$. Recall the following estimate (see
Remark~\ref{rem:5.totnumzerosupperbound})
\begin{equation}\label{eq:14.avernzero}
\cM_N(\omega, E,1/2,2)  := (2N+1)^{-1} \# \bigl\{z\;:\; z\in\cA_{\rho_0} , \; f_N\zoe =
0\bigr\}\le C(V)<\infty 
\end{equation}
Assume that $k_0> C(V)$. Then there exist $k_1\neq k_2$ and admissible 
$s_1$,$s_2$ such that $\zeta_{k_1,s_1}=\zeta_{k_2,s_2}$,
$|z_{k_1,s_1}-z_{k_2,s_2}|\les \exp(-N^{1/3})$. Hence
\begin{equation}\label{eq:14.closepointsontraj1}
\|x(k_1)-x(k_2)-s(k_1,k_2)\omega\|\les \exp(-N^{1/3})
\end{equation}
with some $s(k_1,k_2)\in \ZZ$, $|s(k_1,k_2)|\le 2N$. Due to
Lemma~\ref{lem:17.dioph} either there exists $N_1$, depending on $(
x(k_1)-x(k_2))$ such that for any $t\ge N_1$
\[
\|x(k_1)-x(k_2)-t\omega\|\gtrsim \exp(-t^{1/3})
\]
or there exists $N\ge N_1$, $t\in \ZZ$, $|t|\le N/4$ such that
\[
\exp(-N^{1/3}) \lesssim \|x(k_1)-x(k_2)-t\omega\|\lesssim
\exp(-N^{1/4})
\]
Recall that  $N'\le N''$ implies that $\tilde {\cB}_{N'',\omega}\subset
\tilde {\cB}_{N',\omega}$. Therefore, due to the argument which leads
to~\eqref{eq:14.closepointsontraj1}, one can assume that the latter
case takes place. Hence, we can
replace relation~\eqref{eq:14.closepointsontraj1} by the following
one
\begin{equation}\label{eq:14.closepointsontraj2}
\exp(-N^{1/3})\lesssim \|x(k_1)-x(k_2)-s(k_1,k_2)\omega\| \lesssim
\exp(-N^{1/4})
\end{equation}
where $|s(k_1,k_2)|\le N/4$. Since $x(k_i)\in \tor \setminus \tilde
{\cB}_{N,\omega}$, due to part $(5)$ of
Proposition~\ref{prop:14.basic} there exists a  segment
$\bigl\{E_{j_i}^{(N)}(\cdot, \omega)$ $, \underline{x}_i,
\overline{x}_ i\bigr\}$ such that $x(k_i)\in$
$(\underline{x}_i+\exp(-N^{1/5}),$ $\overline{x}_i-\exp(-N^{1/5}))$
and
\[
|E_{j(k_i)}(x(k_i),\omega)-E^{(N)}_{j_i}(x(k_i),\omega)|\le
\exp(- \gamma N^{\frac12} )
\]
$i=1,2$. Note that due to part $(5)$ of
Proposition~\ref{prop:14.basic} we can also assume that
\begin{equation}\label{eq:14.localcenter}
|\nu^{(N)}_{j_i}(x(k_i),\omega)| \le N/4
\end{equation}
$i=1,2$ since otherwise we can just replace $N$ by $N'=2N$. Set
$$\hat {x}(k_1):=x(k_2)+s(k_1,k_2)\omega \quad (\mod 1)$$  Then $\hat {x}(k_1)\in$
$(\underline{x}_1,$ $\overline{x}_1)$. Since $\exp(-N^{1/3})\lesssim
|x(k_1)-\hat {x}(k_1)|$ one has
\[
 |E^{(N)}_{j_1}(x(k_1),\omega)-E^{(N)}_{j_1}(\hat
 {x}(k_1),\omega)|\gtrsim \exp(-2N^{1/3})
\]
On the other hand due to part $(8)$ of
Proposition~\ref{prop:14.basic} and since $|s(k_1,k_2)|\le N/4$ 
\[
|E^{(N)}_{j_1}(\hat
 {x}(k_1),\omega)-E^{(N)}_{j_2}(x(k_2),\omega)|\lesssim \exp(-N^{1/2})
\]
Since $E_{j(k_1)}(x(k_1),\omega)=E_{j(k_2)}(x(k_2),\omega)$ we
arrive at a contradiction. Thus $k_0\le C(V)$. If $V(e(x))$ is a
trigonometric polynomial then $C(V)\le 2\deg ( V(e(\cdot)))$ due to
Remark~\ref{rem:5.totnumzerosupperbound}.
\end{proof}

\section{Elimination of triple resonances}\label{sec:tripleelim}

\noindent The goal of this section is to eliminate $\omega$ with the
property that for some $x\in\tor$ and some $0<m<m'$
\begin{equation}
  \label{eq:triples} \begin{aligned}
  |E_{j_1}^{(n)}(x+m\omega,\omega) - E_{j_0}^{(n)}(x,\omega)| &< \eps \\
 |E_{j_2}^{(n)}(x+m'\omega,\omega) - E_{j_0}^{(n)}(x,\omega)| &< \eps
\end{aligned}
\end{equation}
with three distinct $j_0, j_1, j_2$. Of course it will be important
to specify the mutual sizes of $n,m,m'$ and~$\eps$. Unlike the
previous elimination machinery based on the resultant of two
polynomials from Section~\ref{sec:resultant}, the elimination here
will be based on the implicit function theorem; this in turn will
require the lower bound on the slopes of the
$E_j^{(N)}(\cdot,\omega)$ that was obtained above 
at the expense of eliminating a small set of energies. Note that
there is no hope of eliminating the situation described
by~\eqref{eq:triples} unless the $E_j^{(N)}(x,\omega)$ truly depend
on~$x$ (for example, if the potential is constant -- this is of
course excluded in our case since we are assuming positive Lyapunov
exponents). We have chosen to present the elimination process
``abstractly'' at first, i.e., without any reference to the
$E_j^{(N)}$. Later we will apply the abstract elimination theorem,
see Proposition~\ref{thm:8.elim}, to the system~\eqref{eq:triples}.
We now begin
with a number of standard calculus lemmas that develop the implicit
function theorem in a quantitative way. We use the well-known idea
of basing the implicit function theorem on monotonicity arguments
(this is of course restricted to {\em scalar} implicit functions).
The first lemma is nothing but a careful statement of a quantitative
implicit function theorem.

\begin{lemma}
\label{lem:8.implic}\begin{itemize}
\item[(1)] Let $f\in C^1(\cR)$
where $\cR:=(a,b)\times (c,d)$. Assume that
\begin{equation}
\label{eq:8.dirpos} \mu:=\inf_{(x,y)\in\cR} \partial_y f(x,y) > 0,
\quad K :=\sup_{(x,y)\in\cR}|\partial_x f(x,y)|<\infty
\end{equation}
If $f(x_0,y_0)=0$ for some $(x_0,y_0)\in\cR$, then for any \[x\in
J_0:=(x_0-\kappa_0,x_0+\kappa_0)\cap (a,b),\qquad
\kappa_0:=h_0\mu K^{-1}, \quad h_0:=\min(y_0-c,d-y_0)\] there exists
a unique $y=\phi_0(x)\in(c,d)$ such that
$ f(x,\phi_0(x))=0 $.
Moreover,  $\phi_0\in C^1(J_0)$ and
$\sup_{x\in J_0} |\phi_0'(x)|\le K\mu^{-1}$.

\item[(2)] Assume in addition that the function $f(x, y)$ admits an analytic
continuation to the domain
\begin{equation}\label{eq:12.fdom}
\cT = \{(z, w) \in \IC^2 : \dist \bigr( (z, w), \cR \bigr) < r_0 \}
\end{equation}
for some $r_0 > 0$, and obeys
\begin{equation}\label{eq:12.fsup}
\sup\limits_{\cT} | f(z, w)| \leq 1
\end{equation} 
Then the implicit function $\phi_0$ has an 
analytic continuation to the rectangle $U=J_0\times (-r,r)$ where 
\begin{equation}\label{eq:rr0mu}
r := \min(r_0^3\mu^2,r_0)/8 \end{equation}
Furthermore, $f(z,\phi_0(z))=0$ for all $z\in U$ and  $\sup\limits_{z\in U} |\phi_0(z)| \le \max(|c|,|d|)+r_0$. 
\end{itemize}
\end{lemma}
\begin{proof}
Note that for any $x\in (a,b)$ one has $|f(x,y_0)|\le K|x-x_0|$. In
particular, for any $|x-x_0|<\kappa_0$
\[ |f(x,y_0)| < h_0 \mu \]
Given such $x$ consider the case $0<f(x,y_0)< h_0\mu$. Since $c\le
y_0-h_0<y_0+h_0<d$, we infer that
\[ f(x,y_0-h_0) < h_0\mu- h_0\mu =0.\]
Hence, there exists a unique $y=\phi_0(x)\in (y_0-h_0,y_0)$ such
that $f(x,\phi_0(x))=0$. If instead $-h_0\mu < f(x,y_0)\le 0$ then
there exists a unique $y=\phi_0(x)\in (y_0,y_0+h_0)$ such that
$f(x,\phi_0(x))=0$. It follows from the chain rule that $\phi_0\in
C^1(J_0)$ and $| \phi_0'(x)|\le K\mu^{-1}$; in fact,
\[
\phi_0'(x) = - \frac{\partial_x f(x,\phi_0(x))}{\partial_y
f(x,\phi_0(x))}
\]
for all $x\in J_0$. That finishes the proof of part $(1)$. To prove
part $(2)$ fix an arbitrary $x_1 \in J_0$ and set $y_1=\phi_0(x_1)$. Due to Cauchy's
estimates one has
 \[
|\partial^n_w f(x_1,y_1+w)|_{w=0}\le n! r_0^{-n}\qquad\forall\;n\ge0
\]
Since $f(x_1, y_1) = 0$, $\partial_w f(x_1,y_1+w)|_{w=0} \geq \mu$
the Taylor series expansion for $f(x_1,y_1+w)$ in the disk
$\cD(0,r_1)$, with $r_1 := \min(r_0^2\mu/2,r_0)$ yields
\[
|f(x_1, y_1+w)| \geq \mu r_1/2 
\]
for any $|w| = r_1$.  Furthermore, $w\mapsto f(x_1,y_1+w)$ has a simple zero at $w=0$ in the disk $|w|\le r_1$. 
 Applying Cauchy's estimate in the $z$-variable now implies that
\[
 |f(x_1+z, y_1+w)| \geq r_0^2\mu^2/8,\quad\forall\; |z|\le r, \; |w| = r_1
\]
where $r$ is as in the statement of the lemma. 
We now claim that 
\[
 \frac{1}{2\pi i}\oint_{|w|=r_1} \frac{f_w(x_1+z,y_1+w)}{f(x_1+z,y_1+w)}\, dw = 1
\]
for all $|z|< r$. Indeed, the integral on the left counts the number of zeros $f(x_1+z,y_1+\cdot)=0$
 inside the disk $|w|<r_1$ and with $z$ fixed. Since there is a unique zero at $w=0$ 
in this disk when $z=0$ is fixed, and since the integral is analytic in $|z|<r$ the claim follows.
By the residue theorem, the sought after implicit function is given by
\[
 \frac{1}{2\pi i}\oint_{|w|=r_1} w\frac{f_w(z,y_1+w)}{f(z,y_1+w)}\, dw = \phi_0(z)
\]
for all $|z-x_1|<r$. It is analytic and has all the desired properties. Covering all of $J_0$ with such disks
we obtain $\phi_0$ on $U$ by the uniqueness of analytic continuations. 
\end{proof}

The next lemma is a slight variant which does not require vanishing
at a point but only smallness. We of course reduce the latter case
to the former.

\begin{lemma}
\label{lem:8.ineqloc} Let $f\in C^1(\cR)$ be as in the previous
lemma and suppose $\mu$ and $K$ are as in~\eqref{eq:8.dirpos}.
Assume that $|f(x_1,y_1)| < \ve$ for some $(x_1,y_1)\in \cR$ and
$0<\ve\le h_1\mu$ where $h_1:=\min(y_1-c,d-y_1)/2$. Then
\begin{enumerate}
\item for any $x\in J_1:=(x_1-\kappa_1,x_1+\kappa_1)\cap (a,b)$ with
$\kappa_1:=h_1\mu K^{-1}$ there exists a unique $y=\phi_1(x)\in
(c,d)$ such that $f(x,\phi_1(x))=0$. Moreover,   $\phi_1\in C^1(J_1)$ and
\[\sup_{x\in J_1} |\phi_1'(x)|\le K\mu^{-1}\]
\item for any $x\in J_1$ and any $y\in (c,d)\setminus (\phi_1(x)-\ve
\mu^{-1},\phi_1(x)+\eps \mu^{-1})$ one has $|f(x,y)|\ge \ve$
\end{enumerate}
Assume in addition that $f$ admits
an analytic continuation to the domain \eqref{eq:12.fdom}
 and obeys condition \eqref{eq:12.fsup}. Then 
 $\phi_1$ admits an analytic continuation to the rectangle $V:=J_1\times (-r,r)$ with $r$  as in~\eqref{eq:rr0mu}. One has $\sup\limits_{z\in V} |\phi_1(z)| \leq \max(|c|,|d|)+r_0$
and $f(z,\phi_1(z))=0$ for all $z\in V$. 
\end{lemma}
\begin{proof}
Assume for instance  that $0\le f(x_1,y_1)<\ve$. Since $c<
y_1-h_1<d$, we conclude that $f(x_1,y_1-h_1)<\ve-\mu h_1\le 0$.
Hence, there exists a unique $\tilde y_1\in (y_1-h_1,y_1]$ such that
$f(x_1,\tilde y_1)=0$. By Lemma~\ref{lem:8.implic} with
$(x_0,y_0):=(x_1,\tilde y_1)$ there exists a $C^1$-function
$\phi_1(x)$ defined on the interval
\[J_2:= (x_1- \kappa_0,x_1+\kappa_0)\cap (a,b)\] with
\[ \kappa_0:=  h_0\mu K^{-1},\quad  h_0=\min(\tilde
y_1-c,d-\tilde y_1)\] such that $f(x,\phi_1(x))=0$ for any $x\in
J_2$. Moreover, $\sup_{x\in J_1}| \phi_1'(x)|\le K\mu^{-1}$. Note
first that \[ h_0 \ge \min(y_1-c,d-y_1)-h_1= h_1\] by construction.
So, $\phi_1(x)$ is defined on the interval~$J_1$. Clearly,
$|f(x,y)|\ge\ve$ for any $x\in J_1$ and any $y\in (c,d)\setminus
(\phi_1(x)-\eps\mu^{-1},\phi_1(x)+\ve\mu^{-1})$. That proves the
first part of the statement. The statement about the analytic continuation of $\phi_1$
follows from  part $(2)$ of
Lemma~\ref{lem:8.implic}. 
\end{proof}

We can now combine these local lemmas with a covering procedure to
obtain a global result. It is a quantitative version of the
following qualitative statement: suppose $f=f(x,y)$ is smooth on
some rectangle $\cR$ and $\partial_y f\ne0$ on~$\cR$. Then the set
where $|f|<\eps$ in $\cR$ with $\eps>0$ small is covered by a union
of neighborhoods of the (local) graphs $(x,\phi(x))$ where
$f(x,\phi(x))=0$. Figure~3 shows a possible form of the set
$U_f(h,\eps)$ appearing in the following proposition (indicated by
the shaded areas). The big rectangle is $\cR$ and the two horizontal
dashed  lines are at heights $c+h_1$ and $d-h_1$, respectively,
defining a smaller rectangle. Note that due to the fact that
$f(x,y)$ is increasing in $y$, the dashed areas cannot have any
$x$-projection in common. Also, they cannot ``die'' inside of the
smaller rectangle due to the previous lemmas (implicit function
theorem). Hence, they can only end on the boundaries of the smaller
rectangle. Also note that while each of the two shaded areas is
defined by graphs over~$x$, they are not graphs over~$y$. However,
we can cut them up into finitely many graphs over~$y$, see below.

\begin{prop}
\label{prop:8.ineqcover} Let $f\in C^1(\cR)$ and $\mu, K$ be as
in~\eqref{eq:8.dirpos}. Given $0<h_1< (d-c)/4$, and $0<\ve\le
h_1\mu$, there is a sequence of pairwise disjoint intervals
$\{J_i\}_{i=1}^{m}$ in $(a,b)$, with (provided they are arranged in
increasing order)
\[ \min_{2\le i\le m-1}|J_i|\ge \kappa_1:= h_1\mu
K^{-1}
\] and
satisfying the following properties:
\begin{enumerate}
\item for each $1\le i\le m$ there exists $\phi_i\in C^1(J_i)$  such that
\[f(x,\phi_i(x))=0\quad \forall\; x\in J_i,\quad \sup_{x\in J_i} |\phi_i'(x)|\le
K\mu^{-1}\]
\item the  set
\[ \cU_f(h_1,\ve) := \{(x,y)\in (a,b)\times (c+h_1,d-h_1)\::\:
|f(x,y)|<\ve\} \] satisfies
\[
\cU_f(h_1,\ve) \subset \bigcup_{i=1}^m \cS_x(\phi_i,\ve\mu^{-1})
\]
where for any $\sigma>0$,
\begin{align*}  \cS_x(\phi_i,\sigma) :=
 \{(x,y)\::\: x\in J_i,\; y\in
(\phi_i(x)-\sigma, \phi_i(x)+\sigma) \} \end{align*}
\end{enumerate}
Furthermore, assume in addition that $f$ admits
an analytic continuation to the domain \eqref{eq:12.fdom}
 and obeys condition \eqref{eq:12.fsup}. Then for each $i$ the function
 $\phi_i$ admits an analytic continuation to the domain
 \[
\cT_i= \{z=x+iy \;:\; x\in J_i,\; |y|<
 r\}
 \] with $r$ as in \eqref{eq:rr0mu}, and $f(z,\phi_i(z))=0$ as well as $ |\phi_i(z)| \leq \max(|c|,|d|)+r_0$ 
for all $z\in \cT_i$. 
\end{prop}
\begin{proof}
  By Lemma~\ref{lem:8.ineqloc}, for each $(x_1,y_1)\in \cU_f(h_1,\ve)$ there
  exists an interval $J(x_1,y_1):=J_1$ and a function $y=\phi_1(x)$ as described in that lemma.
  This defines a  collection \[\cC:= \{J(x_1,y_1)\}_{(x_1,y_1)\in \cU_f(h_1,\ve)}
  \]
Suppose $J(x_1,y_1), J(\tilde x_1,\tilde y_1)\in\cC$ have a nonempty
intersection and let $\phi_1$ and $\tilde\phi_1$ denote the
functions associated with $J(x_1,y_1)$ and  $J(\tilde x_1,\tilde
y_1)$, respectively. Then
\[\phi_1(x)=\tilde \phi_1(x) \quad \forall\; x\in J(x_1,y_1)\cap
J(\tilde x_1,\tilde y_1) \] due to the monotonicity of $y\mapsto
f(x,y)$. Define an equivalence relation on the intervals in~$\cC$ as
follows: $J, \tilde J\in\cC$ are equivalent iff they can be
connected by a chain of pairwise intersecting intervals in~$\cC$.
Then we find the $J_i$ in the statement above simply by taking the
union over all intervals in an equivalence class. The $\phi_i$ are
well-defined by the aforementioned uniqueness property of the
graphs.
The analytic continuation statement follows from Lemma~\ref{lem:8.ineqloc}.
\end{proof}

\begin{figure}[ht]
\centerline{\hbox{\vbox{ \epsfxsize= 13.0 truecm \epsfysize=6.0
truecm \epsfbox{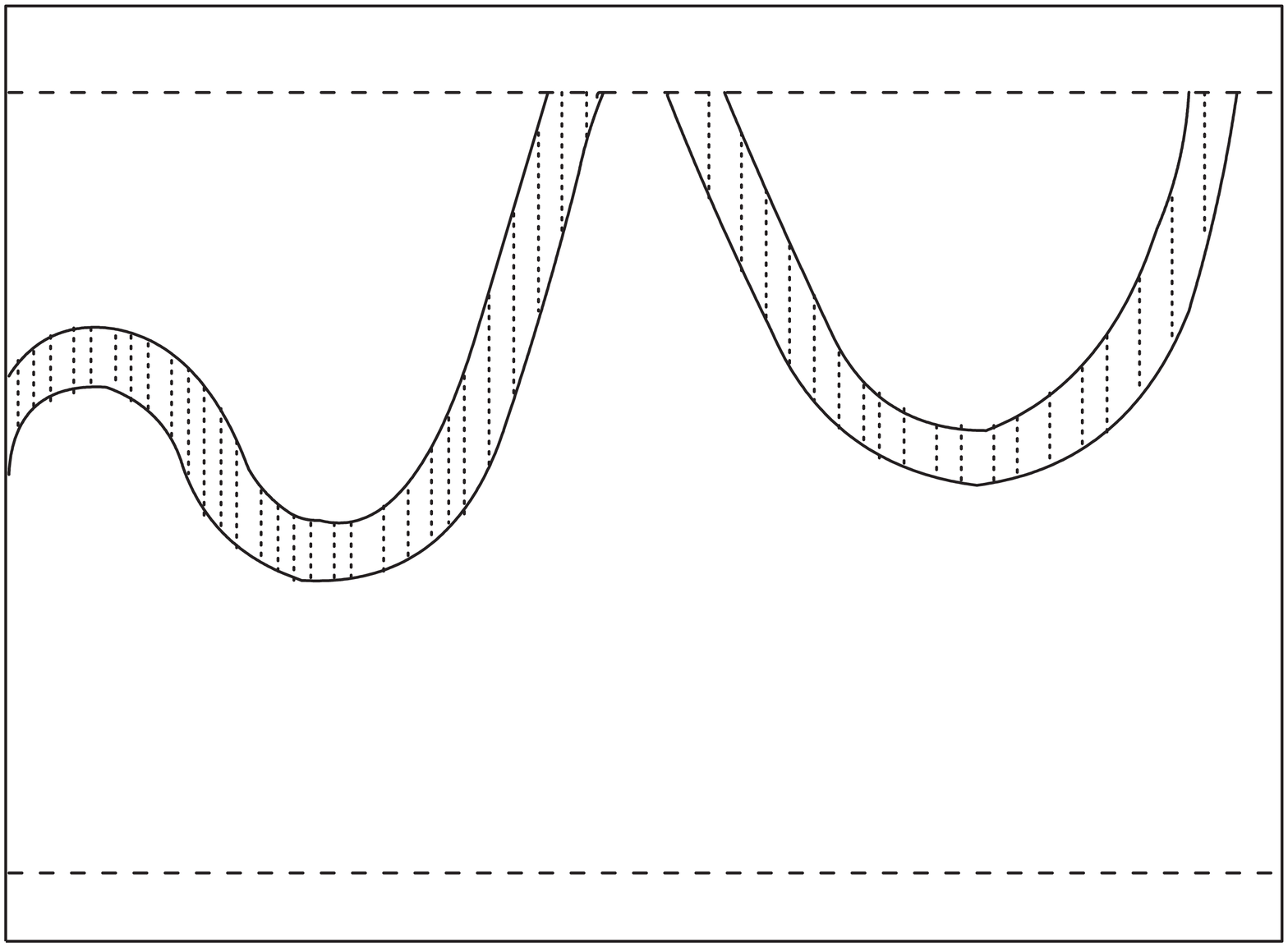}}}}
\caption{The set $\cU_f(h,\eps)$}
\end{figure}

In view of \eqref{eq:triples} we will need to apply the previous
proposition to a system $|f(x,y)|<\eps$, $|g(x,y)|<\eps$. More
precisely, we wish to eliminate a small set of~$y$ so that this
system fails for every~$x$ assuming suitable lower bounds on
$\partial_y f(x,y)$ and $\partial_y g(x,y)$ (one of these
derivatives will need to be much larger than the other). We shall
proceed by eliminating $x$ from the system $|f(x,y)|<\eps$,
$|g(x,y)|<\eps$. Note that this cannot be done on the basis of
Proposition~\ref{prop:8.ineqcover}  alone as we will need to invert
each function $y=\phi_i(x)$.  This will of course not always be
possible. However, by a Sard type argument we will be able to remove
a small set of~$y$ (in measure) so that the inversion can be carried
out for the remaining~$y$.

The precise formulation of this procedure is given by
Proposition~\ref{prop:8.ineqycover} below. We start with the
following lemma, which is a quantitative version of Sard's theorem.

\begin{lemma}
\label{lem:8.sardtype} Let $\phi$ be a real--valued function defined on
the interval $(\alpha,\beta)$ which admits an analytic continuation to the domain
\[
\cS = \{ z \in \IC \;:\; \dist (z, (\alpha, \beta) ) < r \}
\]
for some $r>0$ and satisfies $
\sup\limits_\cS |\phi(z) | \leq q
$ for some $q>0$. 
Then, given $0<\delta<\frac{q}{5r}$ there exist at most
$$n \le
12r^{-1}(\beta-\alpha+r) \log(qr^{-1} \delta^{-1})$$
 disjoint intervals
$I_j \subset (\alpha, \beta)$ such that
\begin{align}
\label{eq:int_delt}
|\phi' (x) | &\geq \delta, \qquad \forall x \in \mybigcup_jI_j\\
| \phi'(x)| &\leq 2\delta, \qquad \forall x \in (\alpha, \beta)
\backslash \mybigcup_jI_j
\end{align}
 In particular,
\begin{equation}\label{eq:total_var}
\int_{(\alpha, \beta) \backslash \mybigcup_j I_j} |\phi' (x) |\, dx
\leq 2\delta (\beta-\alpha)
\end{equation}
\end{lemma}
\begin{proof} Let $x_1 \in (\alpha, \beta)$ be arbitrary and define
 $(\alpha_1, \beta_1) := (x_1 - \frac{r}{8}, x_1 + \frac{r}{8})$. By a standard covering argument, it
suffices to consider $(\alpha_1, \beta_1)$ in the role of $(\alpha, \beta)$.
 The function $\phi'$ is  analytic on $\cD(x_1, r)$ and due to Cauchy's estimates
  satisfies 
\[
\max\limits_{\cD(x_1, 7r/8)} | \phi'(z) | \leq 8q r^{-1}\] Assume
first that there exists $\zeta_1\in \cD(x_1, r/8)$ such that
\begin{equation}\label{eq:12.aux}
  | \phi'(\zeta_1) | \geq 2\delta
\end{equation}
Then, due to the standard Jensen formula \eqref{eq:2.jensb} applied
to  $\phi'(z)\pm \delta$ in the disk $\cD(\zeta_1, 3r/4)$
one obtains
\begin{equation}
\nn
\begin{split}
\# \{z\in\cD(\zeta_1, r/4) \,:\, \phi'(z) \pm \delta =0\}& \leq
\int_0^1 \log |\phi'(\zeta_1+\frac {3}{4}r e(\theta))\pm \delta|\,
d\theta - \log|\phi'(\zeta_1)|\\
& \le \log (8qr^{-1}+\delta) - \log(2\delta)\le  \log(5q r^{-1} \delta^{-1})\\
&\le 2\log(q r^{-1} \delta^{-1})
\end{split}
\end{equation}
Hence there exist at most \[ n_1 \le  2\log(q r^{-1}\delta^{-1})+1 \le 3\log(q r^{-1}\delta^{-1})\] 
disjoint intervals $I_j \subset (\alpha_1, \beta_1)$ such that
\begin{align*}
|\phi' (x) | &\geq \delta, \,\,\, \forall x \in \mybigcup_j I_j\\
| \phi'(x)| &\leq \delta, \,\,\, \forall x \in (\alpha_1, \beta_1)
\backslash \mybigcup_jI_j
\end{align*}
If condition~\eqref{eq:12.aux} fails then
\[
|\phi'(x)|\leq 2\delta, \,\,\, \forall x\in(\alpha_1, \beta_1)\]
which finishes the proof.
\end{proof}

\begin{remark}
  \label{rem:scaling} In Lemma~\ref{lem:8.sardtype} we set up $\delta$ as a
  ``dimensionless'' quantity in the following sense: consider the scaling
  $\phi_\lambda(x):= \lambda^{-1}\phi(\lambda x)$ for any
  $\lambda>0$. Clearly, $\phi_\lambda'$ scales like $\lambda^0$ which
  is what we mean by~\eqref{eq:int_delt} and thus~$\delta$, being scaling invariant. Note, however, that $q$ scales like~$\lambda^{-1}$
which explains why the $q\delta^{-1} r^{-1}$ term inside the logarithm above is scaling invariant (which of course is necessary). 
  The somewhat strange scaling comes from the context of the implicit
  function theorem. Indeed, let $f_\lambda(x,y):= f(\lambda
  x,\lambda y)$ which simply means that we homothetically  scale the
  rectangle~$\cR$ to~$\lambda^{-1}\cR$. Then the implicit function $y=\phi(x)$ scales
  precisely as $\phi_\lambda(x)$ above. It is instructive to check
  our ``abstract'' results in this section against this scaling.
  For example, in Lemma~\ref{lem:8.implic} both $K$ and $\mu$ scale like $\lambda$,
  whereas $\phi_0'$ scales like $\lambda^0$ as required by the
  estimate $|\phi_0'(x)|\le K\mu^{-1}$. We will keep all ``abstract'' results in this section
scaling invariant. 
\end{remark}

 We are now able to
formulate the aforementioned ``$x=\psi(y)$'' version of
Proposition~\ref{prop:8.ineqcover}. In Figure~3 this corresponds to
removing those pieces from the two shaded areas where the graphs
defining the boundaries have horizontal tangents.

\begin{prop}
\label{prop:8.ineqycover} Let $f(x,y)$ be a real-valued function defined in
$ \cR = (a,b)\times (c,d)$. Assume that
$f$ admits an analytic continuation to the domain
\eqref{eq:12.fdom}
 and obeys condition \eqref{eq:12.fsup}. Let $\mu,K$ be as in~\eqref{eq:8.dirpos}.
 Apply
Proposition~\ref{prop:8.ineqcover} to~$f$ with parameters $h_1,\eps$
as specified there, and let $\phi_i$ be the resulting functions
defined on the intervals~$J_i$, $1\le i\le m$. Fix an arbitrary
\begin{equation}\label{eq:delta_cond} 
0<\delta< \frac{q}{5r},\qquad r=\min(r_0^3\mu^2,r_0)/8,\quad q:= \max(|c|,|d|)+r_0
\end{equation}
see~\eqref{eq:rr0mu}. 
Then for each $1\le i\le m$ there exist pairwise
disjoint intervals
\[ J_{i,j} \subset J_i,  \quad\forall\; 1\le j \le j(i) ,\]
so that for each of these intervals
\[
\phi_i\: :\: J_{i,j} \to I_{i,j}:=\phi_i(J_{i,j})
\]
is invertible. Denote the inverse by $\psi_{i,j}$. Then the
following properties hold:
\begin{enumerate}
\item The total number of intervals $J_{i,j}$ $($and thus of $I_{i,j}$ as well as of $J_{i,k}'$
which are defined to be the connected components of $J_i\setminus
\bigcup_{j=1}^{j(i)} J_{i,j}$ $)$ is at most
\[M:= 12\big((b-a)(r^{-1} + \kappa_1^{-1}) + 2\big)\, \log(q r^{-1}\delta^{-1})
\]
where $\kappa_1:=h_1\mu K^{-1}$. 
\item there exist  no more than $3M$ many intervals $I_\ell$ such
that $ \sum_{\ell} |I_\ell| \le 2(b-a)\delta + 3M\eps\mu^{-1}$ and
so that with  $\cU_f (h_1,\ve)$  as in
Proposition~\ref{prop:8.ineqcover},
\begin{equation}\label{eq:Ufminus}
\cU_f(h_1,\ve) \setminus \Big(\bigcup_{i=1}^m \bigcup_{j=1}^{j(i)}
J_{i,j}\times I_{i,j} \Big) \subset (a,b)\times  \bigcup_\ell I_\ell
\end{equation}
\item for  any $\sigma>0$, define
\[ \cS_y(\psi_{i,j},\sigma):=\{(x,y)\::\: y\in I_{i,j},\;|\psi_{i,j}(y)-x|<\sigma\}
\] Then
\begin{equation}\label{eq:Ufdurchschnitt} \cU_f(h_1,\ve) \bigcap \Big(\bigcup_{i=1}^m
\bigcup_{j=1}^{j(i)} J_{i,j}\times I_{i,j} \Big) \subset
\bigcup_{i=1}^m \bigcup_{j=1}^{j(i)}
\cS_y(\psi_{i,j},\ve\delta^{-1}\mu^{-1})
\end{equation}
\item $|\psi_{i,j}'(y)|\le \delta^{-1}$ for all $y\in I_{i,j} $ and all $i,j$
\end{enumerate}
\end{prop}
\begin{proof}
 By Proposition~\ref{prop:8.ineqcover} for each $i$ the function
 $\phi_i$ admits an analytic continuation in the domain
 \[
\cT_i= \{z=x+iy:x\in J_i, |y|\le
 r\}
 \] with $\sup\limits_{\cT_1} |\phi_0(z)| \leq q$.
  Applying Lemma~\ref{lem:8.sardtype} to
each $\phi_i$  produces pairwise disjoint intervals
\[ J_{i,j}\subset J_i, \quad 1\le j\le j(i)\le n_i\]
with, cf.~\eqref{eq:delta_cond}, 
\[
n_i\le 12 (|J_i| r^{-1}+1) \log(q r^{-1}\delta^{-1})
\]
such that
\begin{align}
 |\phi_i'(x)|&\ge\delta \quad \forall\; x \in \bigcup_{j=1}^{j(i)}\,
J_{i,j} ,\label{eq:lowerder}\\
 \sum_i \int\limits_{J_i\setminus \bigcup_{j=1}^{j(i)} J_{i,j}}
|\phi_i'(x)|\, dx & \le \sum_i 2|J_i|\delta 
\le 2(b-a)\delta 
\label{eq:mesini}
\end{align}
Since $J_i\ge \kappa_1$ for all but possibly two $i$, one has
\[
\sum_{i=1}^m n_i\le 12\big((b-a)(r^{-1} + \kappa_1^{-1}) + 2\big)\, \log(q r^{-1}\delta^{-1})
\]
as claimed in property~(1). To prove property~(2) of the
proposition, we define $I_\ell$  the be the collection of all
intervals arising as follows: $\eps\mu^{-1}$-neighborhoods of all
$\phi_i(J_{i,k}')$ (we refer to these as type~I), as well as
$\eps\mu^{-1}$-neighborhoods of each of the  two points in
$\phi_i(\partial J_{i,j})$ (these are type~II). By property~(1)
there are no more than $M$ type~I intervals, as well as at most $2M$
type~II intervals. Moreover, each type~II interval has measure not
exceeding $2\eps\mu^{-1}$, whereas all type~I intervals in total
have measure at most
\[
\sum_{i,k} \big( |\phi_i(J_{i,k}')| + 2\eps\mu^{-1}\big) \le
2(b-a)\delta + 2M\eps\mu^{-1},
\]
see \eqref{eq:mesini}. To prove~\eqref{eq:Ufminus}, observe from
Proposition~\ref{prop:8.ineqcover} that it suffices to prove the
inclusion
\[
 \bigcup_{i=1}^m \Big[\cS_x(\phi_i,\ve\mu^{-1}) \setminus \bigcup_{j=1}^{j(i)} J_{i,j}\times
I_{i,j} \Big] \subset (a,b)\times \Big[(c,d)\setminus \bigcup_\ell
I_\ell \Big]
\]
A point $(x,y)$ belongs to the set inside the brackets on the
left-hand side iff either of the following two scenarios occurs
\begin{itemize}
\item $x\in J_i\setminus \bigcup_{j=1}^{j(i)} J_{i,j}= \bigcup_{k}
  J_{i,k}'$ and $|y-\phi_i(x)|<\eps\mu^{-1}$
  \item $x\in J_{i,j}$ and $|y-\phi_i(x)|<\eps\mu^{-1}$ but $y\notin
  I_{i,j}=\phi_i(J_{i,j})$
\end{itemize}
In the first case, $y\in I_\ell$ where $I_\ell$ is a type-I
interval, whereas in the second case, $y\in I_\ell$ which is
type-II. In conclusion,
\[
\cS_x(\phi_i,\ve\mu^{-1}) \setminus \bigcup_{j=1}^{j(i)}
J_{i,j}\times I_{i,j}  \subset J_i\times \Big[(c,d)\setminus
\bigcup_\ell I_\ell \Big]
\]
which yields the desired property by taking unions over the pairwise
disjoint~$J_i$.
 On each interval
$I_{i,j}:=\phi_i(J_{i,j})$ an inverse function
$\psi_{i,j}=\phi_i^{-1}$ is defined, and moreover
\[ \sup\limits_{y\in I_{i,j}} |{\psi}_{i,j}'(y)|\le \delta^{-1} \]
from \eqref{eq:lowerder}. This establishes property~(4). Note that
if $y\in I_{i,j}$ and $|y-\phi_i(x)|<\ve \mu^{-1}$ for some $x\in
J_{i,j}$, then $|x-\psi_{i,j}(y)|\le\ve\delta^{-1}\mu^{-1}$. In
other words, in view of Proposition~\ref{prop:8.ineqcover} one has
\[ \cU_f(h_1,\ve) \cap
\Big(J_{i,j} \times I_{i,j} \Big) \subset
\cS_y(\psi_{i,j},\ve\mu^{-1}\delta^{-1})
\]
which implies property (3) above.
\end{proof}

We are now in a position to state and prove the main ``abstract''
elimination result of this section. It will be applied
to~\eqref{eq:triples}.

\begin{prop}
\label{thm:8.elim} Let $f(x,y)$ be a real-valued function defined in
$\cR=(a,b)\times (c,d)$. Assume that $f$ admits an
analytic continuation to the domain 
\[
\cT = \{(z, w) \in \IC^2 : \dist \bigr( (z, w), \cR \bigr) < r_0 \}
\]
for some $r_0 > 0$, and obeys $\sup\limits_{\cT} | f(z, w)| \leq 1$. 
Furthermore, let $g\in C^1(\wt\cR)$ be a real-valued function defined
on $\wt\cR:=(a,b)\times (c',d')$.  Assume that
\begin{align*} \mu &:= \inf_{(x,y)\in\cR}
\partial_y f(x,y)
> 0,\quad
\bar \mu := \inf_{(x,y)\in\wt\cR} \partial_y g(x,y) > 0 \\
K &:= \sup_{(x,y)\in\cR} |\partial_x f(x,y)| +
\sup_{(x,y)\in\wt\cR} |\partial_x g(x,y)| <\infty
\end{align*}
Set $r:=\min(r_0^3 \mu^2,r_0)/8$, $q:=\max(|c|,|d|)+r_0$ and pick $h_1\in (0,(d-c)/4)$, $\sigma_1\in(0, \frac{q}{5r})$, $\sigma_2\in(0, h_1\mu)$.
Define $\kappa_1:= h_1\mu K^{-1}$,  as well as 
\begin{equation}\label{eq:M_def}M:= 12\big((b-a)(r^{-1} + \kappa_1^{-1}) + 2\big)\, \log(q r^{-1}\sigma_1^{-1})
\end{equation}
and assume that
\[
 \bar\mu\ge \max\big[ 
 4(1+K\sigma_1^{-1}\mu^{-1})M, 2K\sigma_1^{-1} \big]
\]
Then  there exist subintervals $\{U_\ell\}_{\ell=1}^{\ell_0}$ of
$(c,d)$ and $\{V_k\}_{k=1}^{k_0}$ of $(c',d')$ so that
\begin{enumerate}
\item[(i)] $k_0\le \ell_0\le 3M$
\item[(ii)] $\sum_{\ell=1}^{\ell_0}|U_{\ell}|  \le 2
\sigma_1(b-a)+3M\sigma_2\mu^{-1}$, $\sum_{k=1}^{k_0} |V_k| \le
2\sigma_2(b-a) $
\item[(iii)] the intervals $U_\ell$ only depend
on the function $f$, the rectangle~$\cR$, and $\sigma_1,\sigma_2$
\item[(iv)] for any
\[ y\in (c+h_1,d-h_1)\cap (c',d')\setminus
\bigcup_{\ell,k} U_\ell\cup V_k
\]
and any $x\in(a,b)$ at least one of the following two inequalities
fails:
\begin{equation}
  \label{eq:dopp_eps} |f(x,y)|<\sigma_2,\quad |g(x,y)|<\sigma_2
\end{equation}
\end{enumerate}
\end{prop}
\begin{proof}
We apply Proposition~\ref{prop:8.ineqycover} to $f$ with
$\delta=\sigma_1,\eps:=\sigma_2$. This  produces intervals $J_{i,j}$
and $I_{i,j}$ as well as  functions~$\psi_{i,j}$ defined on
$I_{i,j}$ satisfying properties (1)--(4) in that proposition. First,
we define $\{U_\ell\}_{\ell=1}^{\ell_0}$ to be the same as the
$I_\ell$, see property~(2) of the proposition. Hence, properties
(i)--(iii) from above pertaining to~$\{U_\ell\}$ follow immediately
from Proposition~\ref{prop:8.ineqycover}. To define the $V_k$,
observe the following: by the chain rule, for any $y\in I_{i,j}\cap
(c',d')$
\begin{align*}
\frac{d}{dy} g(\psi_{i,j}(y),y) \ge \partial_y
g(\psi_{i,j}(y),y)-|\psi_{i,j}'(y)||\partial_x g(\psi_{i,j}(y),y)|
\ge \bar\mu - \sigma_1^{-1}K \ge \frac{\bar\mu}{2}
\end{align*}
since we are assuming that $\bar\mu\ge
2K\sigma_1^{-1}$. Hence, given
$\beta>0$ there exists $I_{i,j}'\subset I_{i,j}\cap (c',d')$ such
that
\begin{equation} |I_{i,j}'|\le 4\beta \bar\mu^{-1}
\label{eq:eta_est}
\end{equation}
and
\[ |g(\psi_{i,j}(y),y)|  >  \beta \qquad \forall \; y\in I_{i,j}\cap(c',d')\setminus I_{i,j}'  \]
Note that we allow for the possibility that
$I_{i,j}\cap(c',d')\setminus I_{i,j}'=\emptyset$ or
$I_{i,j}'=\emptyset$.  Now define the intervals
$\{V_k\}_{k=1}^{k_0}$ to be the entire collection
$\{I_{i,j}'\cap(c',d')\}_{i,j}$ (skipping empty $I_{i,j}'$), with
the choice of~$\beta:=\sigma_2(1+ K\sigma_1^{-1}\mu^{-1})$. Note
that $k_0\le j_0$ and property~(i) is done. Moreover, by~(1) of the
previous proposition,
\[
\sum_{k} |V_k| = \sum_{i,j} |I_{i,j}'| \le 4\beta\bar\mu^{-1} M \le
\sigma_2
\]
by our lower bound on~$\bar \mu$. This finishes property~(ii) above.
To prove property~(iv), let
\begin{equation}\label{eq:UlVk}(x,y)\in (a,b)\times \Big[ (c+h_1, d-h_1)\cap(c',d')
\setminus \bigcup_{\ell,k} U_\ell\cup V_k \Big]
\end{equation}
We can of course assume that $(x,y)\in \cU_f(h_1,\sigma_2)$ for
otherwise $|f(x,y)|\ge\sigma_2$. Because of~\eqref{eq:UlVk}
and~\eqref{eq:Ufminus}, we then have
\[
(x,y) \in \cU_f(h_1,\sigma_2)\bigcap\bigcup_{i=1}^m
\bigcup_{j=1}^{j(i)} J_{i,j}\times (I_{i,j}\setminus I_{i,j}')
\subset \bigcup_{i=1}^m \bigcup_{j=1}^{j(i)}
\cS_y(\psi_{i,j},\sigma_2\sigma_1^{-1}\mu^{-1}),
\]
where the inclusion is \eqref{eq:Ufdurchschnitt}. Now note the
following: fix $i,j$ so that $(x,y)\in J_{i,j}\times I_{i,j}$. Since
$ |f(x,y)|<\sigma_2$,
\[ |x-\psi_{i,j}(y)|\le\sigma_2\sigma_1^{-1}\mu^{-1}\] Hence,
\[ |g(x,y)|\ge |g(\psi_{i,j}(y),y)|-K| x-\psi_{i,j}(y)| >
\beta - K\sigma_2\sigma_1^{-1}\mu^{-1} >\sigma_2 \] provided we
choose $\beta=\sigma_2(1+ K\sigma_1^{-1}\mu^{-1})$ as before. In the
final step we used that \[(x,y)\in J_{i,j}\times (I_{i,j}\cap
(c',d')\setminus I_{i,j}').\]
 But this means that $|g(x,y)|\ge\sigma_2$
assuming~\eqref{eq:UlVk} and $|f(x,y)|<\sigma_2$, and the
proposition is proved.
\end{proof}

Next, we apply this result to establish the main elimination result
concerning system~\eqref{eq:triples}. We will use the notion of an
segment from Definition~\ref{def:8.1}. Let $\Big\{
E_j^{(m)}(\cdot,\omega), \ux',\ox'\Big\}$ be a segment. 
In the following corollary we will work with the following 
quantitative properties of segments:
\begin{itemize}
\item[(a)] $|\partial_x E_j^{(m)}(\cdot,\omega)|\ge
e^{-m^{\delta_1}}$, for any $x\in (\ux',\ox')$
\item[(b)] $\ox'-\ux'\ge e^{-m^{\delta_0}}$
\end{itemize}
where $0 < \delta_0$, $\delta_1\ll 1$ are some parameters. Note that
these are weaker than the ones implicit in~Definition~\ref{def:8.1}, cf.~\eqref{eq:EjtauN}. Furthermore, due
to Remark~\ref{rem:12.analytcontinu}, given $x$, $\omega$ the
function $E_j^{(m)}(\cdot,\cdot)$  admits an analytic continuation
to the polydisk
\[
\cP(x,\omega,m):=\{(z,w)\in \IC^2: |z-e(x)| < r(m), |w-\omega| <
r(m)\}
\]
where $r(m):=\exp(-m^{\delta_0})$. Moreover,
\[
\sup_{\cP(x,\omega,m)} |E_j^{(m)}(z,w)|\le C(V)
\]
Recall the sets $\Omega_{m}$, $\cE_{m,\omega}$, and
$\wt\cE_{m,\omega}$ from Section~\ref{sec:segments} (the latter obviously depends
on parameters $\delta_0,\delta_1$). In the following corollary,
\[
\wt\Omega_n:=\bigcup_{n\le m\le 100n} \Omega_{m}
\]

\begin{cor}
  \label{cor:triple_elim}
Choose $A\ge 10$, $d\ge 4$ arbitrary\footnote{In what follows, the parameter ``$d$'' is a large number that has nothing to do with
the parameter appearing earlier in this section in connection with the rectangle~$\cR$.}
 but fixed, as well as
$0<\eps_0<\frac{1}{20}$ and $0<2\delta_0\le \delta_1\le
\frac{\eps_0}{A}$.
  Let $N\ge N_0(V,\rho_0,a,c,\gamma,\delta_0,\delta_1, A,d,\eps_0)$ be large and set $n:=[(\log
  N)^A]$.
  Then there exist $\cB'_n, \cB''_n\subset \tor$ so that 
\begin{equation}
\nn
\begin{split}
  \mes(\cB'_n) &< N^{-\eps_0}, \quad \compl( \cB'_n) \le \exp((\log N)^{1/2}) \\
\mes(\cB''_n) &< N^{-d+3}, \quad \compl (\cB''_n) \le N^3
\end{split}
\end{equation}
   with the following property: for all  \[\omega\in\tor_{c,a}\setminus
  (\cB'_n\cup \cB''_n\cup\wt\Omega_n)\]
  and for each choice of $n\le n_1,n_2, n_3 \le 100 n$ and $1\le j_i\le 2n_i+1$, $i=1,2,3$, as well as every
   \[e^{n^{3\delta_1}}\le m_1\le \exp((\log N)^{1/4}),\quad N^{10\eps_0} \le m_2\le 2N \] the system
\begin{equation}
  \label{eq:triples2} \begin{aligned}
  |E_{j_2}^{(n_2)}(x+m_1\omega,\omega) - E_{j_1}^{(n_1)}(x,\omega)| &< N^{-d} \\
 |E_{j_3}^{(n_3)}(x+m_2\omega,\omega) - E_{j_1}^{(n_1)}(x,\omega)|
 &< N^{-d}
\end{aligned}
\end{equation}
has no solution with each $E_{j_1}^{(n_1)}(x,\omega)$,
$E_{j_2}^{(n_2)}(x+m_1\omega,\omega)$, and
$E_{j_3}^{(n_3)}(x+m_2\omega,\omega)$ being the evaluation of a
segment with the parameters $\delta_0,\delta_1$ as in $(a)$, $(b)$ above.
\end{cor}
\begin{proof} For each \[\omega\in \tor_{c,a}\setminus \wt\Omega_n, \qquad \wt\Omega_n=\bigcup_{n\le n'\le 100n}
\Omega_{n'}\] we enumerate all possible segments as follows:
the set
\[
[-C(V), C(V)] \setminus \bigcup_{n\le n'\le 100n } \cE_{n',\omega}
\]
can be written as the union of no more than $e^{2n^{\delta_1}}$
intervals $(\uE,\oE)$ of lengths~$e^{-n^{\delta_1}}$ (with $n$ large). Fixing such
an~$(\uE,\oE)$ one obtains no more than $e^{2n^{\delta_1}}$ many
 $I$-segments $\Big\{ E_j^{(n')}(\cdot,\omega), \ux,\ox\Big\}$ with $n\le n'\le 100n$ and $I=(\uE,\oE)$.
 In total, there are  no more than $e^{4n^{\delta_1}}$
many segments in this enumeration, each of which has slope
bounded below by~$s_0:= e^{-n^{2\delta_1}}$. Fix three $n'$ in the
specified range and denote them by $n_1, n_2, n_3$. In addition, fix
three segments from our list which we denote by
\[
\Big\{E_{j_1}^{(n_1)}(\cdot,\omega),\ux_1,\ox_1\Big\},\qquad
\Big\{E_{j_2}^{(n_2)}(\cdot,\omega),\ux_2,\ox_2\Big\},\qquad
\Big\{E_{j_3}^{(n_3)}(\cdot,\omega),\ux_3,\ox_3\Big\}
\]
Let $C(V)$ be a large enough. The functions
\begin{equation}
  \label{eq:fg_def9}
  \begin{aligned}
   f(x,\omega) &:= C(V)^{-1}(E_{j_2}^{(n_2)}(x+m_1\,\omega,\omega) -
E_{j_1}^{(n_1)}(x,\omega))\\
 g(x,\omega) &:=C(V)^{-1} (E_{j_3}^{(n_3)}(x+m_2\,\omega,\omega) -
 E_{j_1}^{(n_1)}(x,\omega))
\end{aligned}
\end{equation}
are defined on rectangles $\cR_f:=(x_1,x_2)\times (\omega_1,\omega_2)$  and
$\wt\cR_g:=(x_1,x_2)\times (\omega_1',\omega_2')$, respectively, where
\[
x_2-x_1\gtrsim \lambda_0 , \quad \omega_2-\omega_1\gtrsim \lambda_1/m_1,\quad
\omega_2'-\omega_1'\gtrsim \lambda_1/m_2
\]
with $\lambda_0:=e^{-n^{\delta_1}}$,
$\lambda_1:=e^{-n^{3\delta_1}}$. This is a consequence of the
stability of the segments under perturbations in $x,\omega$.

 Via an obvious covering argument, the total
number of such rectangles $\cR_f$ and $\wt\cR_g$ that we need to consider is
no larger than $ (\lambda_0\lambda_1)^{-1} m_1 $ and
$(\lambda_0\lambda_1)^{-1} m_2$, respectively, with $m_1$, $m_2$
fixed (and up to multiplicative constants). Similarly, the number of
choices of $f$ which we need to consider with $m_1$ fixed
 is no larger than
$e^{4n^{\delta_1}}(\lambda_0\lambda_1)^{-1} m_1 $ and that of all
possible $g$ is no larger than
$e^{4n^{\delta_1}}(\lambda_0\lambda_1)^{-1} m_2$. Finally, summing
over all admissible choices of $m_1, m_2$ as in the statement yields
\begin{equation}
  \label{eq:fgnums}
    F\les e^{4n^{\delta_1}}(\lambda_0\lambda_1)^{-1} \exp(2(\log N)^{1/4}),\quad G\les e^{4n^{\delta_1}}(\lambda_0\lambda_1)^{-1}
    N^2
\end{equation}
where $F,G$ denote the total number of $f$ and $g$, respectively,
that need to be considered in this enumeration. Note that by our
choice of $\delta_1$, one has
\[
e^{n^{10\delta_1}}\le \exp\big(({\log N})^{\frac{1}{2}} \big) \le N^\eps
\]
for any $\eps>0$ provided $N$ is large.
 We now verify the conditions of 
Proposition~\ref{thm:8.elim} for such  a fixed choice of $f,g$ living
on some pair $\cR=\cR_f,\;\wt\cR=\wt\cR_g$ as above. In the notation of that
proposition,
\[
 m_1\gtrsim  \mu \gtrsim m_1 s_0\gtrsim \sqrt{m_1},\quad m_2\gtrsim \bar\mu \gtrsim m_2 s_0\gtrsim \sqrt{m_2}, \quad K\lesssim 1,
 \quad r_0=e^{-n^{2\delta_1}}= s_0
\]
where we used that $m_1, m_2\gg s_0^{-1}$. 
Since $r_0\mu\gg1$, it follows that $r=\min(r_0^3\mu^2,r_0)/8\asymp r_0=s_0$. Moreover, due to  $r_0\le q\lesssim 1$, 
one has $\frac{q}{5r}\gtrsim 1$ so that the condition on $\sigma_1$ in Proposition~\ref{thm:8.elim}  turns into the  harmless $0<\sigma_1\lesssim 1$. 
We choose
$h_1:= (\omega_2-\omega_1)/8\gtrsim \lambda_1/m_1$ which implies that $\kappa_1=h_1\mu K^{-1}\gtrsim h_1 m_1 s_0\gtrsim \lambda_1 s_0$. 
Now set
$\sigma_1= N^{-2\eps_0}$ and $\sigma_2= N^{-d}$ which is admissible for  Proposition~\ref{thm:8.elim} by the preceding. Then the constant $M$ from \eqref{eq:M_def} satisfies
\[
M \lesssim  \lambda_1^{-1}s_0^{-1} \log(s_0^{-1} \sigma_1^{-1})  \lesssim \lambda_1^{-2} \log( \sigma_1^{-1})
\]
and our main condition
on~$\bar\mu$ reduces to
\[
\bar\mu\gg \sigma_1^{-1}\lambda_1^{-2} \log( \sigma_1^{-1})
\]
This holds because 
\[
 \bar\mu \gtrsim m_2 s_0 \gg \sigma_1^{-2}\lambda_1^{-2}
\]
Finally, since
\[
 \sigma_2 \ll \sigma_1 M^{-1}
\]
if follows that the first bound in part (ii) of Proposition~\ref{thm:8.elim} reduces to
\[
 \sum_{\ell} |U_\ell|\lesssim \sigma_1
\]
Applying  Proposition~\ref{thm:8.elim}, we define $\cB'_n$
and  $\cB''_n$ as the union over all possible choices of $f$ and $g$
as explained above of all intervals $U_\ell$ and $V_k$ respectively
as in Proposition~\ref{thm:8.elim}. Recall that due to property
(iii) in Proposition~\ref{thm:8.elim} the intervals $U_\ell$
depend on function $f$ only. Therefore $\cB'_n$ is a union of at
most $F$ intervals. The set $\cB''_n$ is a union of at most $FG$
intervals. Furthermore,
\begin{equation}
\begin{split}
\nn
 \mes (\cB'_n) &\les \sigma_1 F \le
N^{-2\eps_0}e^{4n^{\delta_1}} (\lambda_0\lambda_1)^{-1} \exp(2(\log N)^{1/4}) \les N^{-\eps_0}
\\ \mes(\cB''_n)&\les \sigma_2 FG \les  N^{-d} e^{8n^{\delta_1}}   (\lambda_0\lambda_1)^{-2}
   N^{2+2\eps_0} \les
N^{-d+3}
\end{split}
\end{equation}
as claimed.  The complexity bounds are as follows:
\begin{align*}
 \compl(\cB'_n) &\les MF\les \exp(\sqrt{N}) \\
 \compl(\cB''_n) &\les MFG \les N^3
\end{align*}
as desired. Finally, 
suppose~\eqref{eq:triples2} had a solution for some~$x$ satisfying
all the conditions stated above. Although the segments
from~\eqref{eq:triples2} do not necessarily belong to our list of
segments described in the beginning of the proof, locally
around~$x$ they would have to agree with some choice of segment
from our list. Therefore, for some choice of $f,g$ as
in~\eqref{eq:fg_def9} necessarily \[|f(x,\omega)|<\sigma_2, \quad
|g(x,\omega)|<\sigma_2\] contradicting that $\omega\not\in \cB_n$,
see Proposition~\ref{thm:8.elim}.
\end{proof}

\begin{remark}\label{rem:17.parameters}  The purpose of this remark is
to comment further on the set of exceptional~$\omega$ in Proposition~\ref{20.regpsitiveregnegative}. Recall that
we assume that
\begin{equation} \label{eq:17.sets1}
\omega\in \tor_{c,a}\setminus \cB(N), \quad \cB(N):=\cB'_{n}\cup
\cB''_n \cup \tilde {\Omega}_{n}
\end{equation}
with the effective bounds 
\begin{equation}
\nn
\begin{split}
  \mes(\cB'_n) & < N^{-\eps_0}, \quad \compl( \cB'_n) \le \exp((\log N)^{1/2}) \\
\mes(\cB''_n) &< N^{-d+3}, \quad \compl (\cB''_n) \le N^3
\end{split}
\end{equation}
Here $0 < \eps_0 \ll 1$, $d\ge 4$ are arbitrary, provided
\[
N\ge N_0(V,c,a,\gamma,\delta_0,\delta_1,A,d,\eps_0)
\]
and $\delta_0$, $\delta_1$, $A$ are  parameters as specified in Corollary~\ref{cor:triple_elim}. 
For the rest of this paper we fix all
the parameters except~$d$ in such a  way that the
corollary holds. Thus the statement of
Corollary~\ref{cor:triple_elim} holds as long as $N$ is large enough
depending on $d$. Recall also that
\[
\wt\Omega_n:=\bigcup_{n\le n'\le 100n} \Omega_{n'}
\]
where
\[
\mes (\Omega_{m})\le \exp(-(\log m)^A),\quad \compl (\Omega_{m})\le
m^C
\]
$n\asymp (\log N)^A$. It is convenient to use the outer
Hausdorff measures $\cH^\alpha_{r}(\cF)$, $\cF\subset \IR$
\[
\cH^\alpha_r (\cF):=\inf \Big\{ \sum_j |I_j|^\alpha\;:\;\cF\subset \bigcup_j I_j,\; \sup_j
|I_j|\le r \Big\}
\]
where $0<\alpha \le 1$, $r>0$ are arbitrary. For $d >7$ define
\begin{equation}\label{eq:15.hausparameters}
\alpha(d)=4/(d-3), \quad r(N)=\exp(-(\log \log N)^A)
\end{equation}
Then
\[
\cH^{\alpha(d)}_{r(N)}(\cB(N))\le \exp(-(\log \log N)^{B})
\]
where $B=A/2$. 
\end{remark}

\section{Resonances and the formation of pre-gaps}\label{sec:pregaps}

The main objective of this section is to establish the resonance
splitting picture for the Rellich parametrization of the eigenvalues
similar to the one described in Figure~2. For ease of notation, we mostly
drop $\omega$ from functions when it appears as an independent variable.

We begin with the following statement, which formalizes the idea that we can
make a positive slope $I$-segment intersect with a negative of the same scale by means of a shift
ot the form $m\omega$. Of crucial importance is the fact that the intersecting
segments can be chosen to be regular unless we are inside a spectrum free interval of energies
in the sense of Section~\ref{sec:formation}. 

\begin{lemma}
\label{prop:13.resexist}
\begin{itemize}
\item[(i)] Fix $\delta>0$ small and let
$\bigl\{E_{j_1}^{(\underline{N})} \, (x), \underline{x}_1,
\overline{x}_1 \bigr\}$ and $\{E_{j_2}^{(\underline{N})} \, (x),
\underline{x}_2, \overline{x}_2 \}$ be a positive-slope and a
negative-slope $I$-segment, respectively, where $I=[\underline{E},
\overline{E}]$, with $\overline{E} - \underline{E} >
\exp(-\underline{N}^\delta)$. Then for all $\underline{N}\ge N_0(\delta)$ there exists an integer $m \in
\bigl[\exp(\underline{N}^\delta), \exp(2\underline{N}^\delta)\bigr]$ and $x_0 \in
(\underline{x}_1, \overline{x}_1)$ such that
\begin{equation}
 \label{eq:13.1}
  E_{j_1}^{(\underline{N})} (x_0) = E_{j_2}^{(\underline{N})} (x_0 + m\omega).
\end{equation} Moreover\footnote{$\{\underline{x}_j,\overline{x}_j\}$ here is the set with these two points as elements.},
\begin{equation}
 \dist(x_0,\{\underline{x}_j,\overline{x}_j\})> C(V)^{-1}(\overline{E} - \underline{E})
\label{eq:x0_sep}
\end{equation}
for $j=1,2$. 
\item[(ii)] Given a scale $\ell$ and
interval $(E'_0, E''_0)$, $ E''_0-E'_0\ge \exp(-(\log \ell)^A)$,
either the interval \[ \Big(E'_0 +\frac {1}{4}\exp(-(\log \ell)^A),
E''_0 -\frac {1}{4}\exp(-(\log \ell)^A) \Big)
\] is spectrum free, or  at some scale $\ell^2 \le \underline{N}\le
\ell^{10}$ there exist a regular positive-slope I-segment $
\bigl\{E_{j_1}^{(\underline{N})} (x), \underline{x}_1, \overline{x}_
1\bigr\}$,
 a regular negative-slope I-segment
 $\bigl\{E_{j_2}^{(\underline{N})} (x), \underline{x}_2,
\overline{x}_2 \bigr\}$, $I = \bigl[\underline{E},
\overline{E}\bigr]\subset [E'_0,E''_0]$, $\overline{E} -
\underline{E}
> \exp(-{\underline{N}}^\delta)$, an integer $m \in [\exp(\underline{N}^\delta),
\exp{(2\underline{N}}^\delta)]$ and a point $x_0 \in
(\underline{x}_1, \overline{x}_ 1)$ such that conditions
\eqref{eq:13.1} and~\eqref{eq:x0_sep} hold.
\end{itemize}
\end{lemma}
\begin{proof}
\begin{itemize}
\item[(i)] Assume for instance, $\underline{x}_1 < \underline{x}_2$.
 Then necessarily also $\overline{x}_1 < \overline{x}_2$.
 Let $y_1 = \overline{x}_2 - \overline{x}_1$ and $y_2 = \underline{x}_2 - \underline{x}_1$.
  The function
$$
h(E)=\bigl(E_{j_2}^{(\underline{N})}\bigr)^{-1} (E) -
\bigl(E_{j_1}^{(\underline{N})}\bigr)^{-1} (E)
$$
is strictly decreasing and  satisfies
$
h(\overline{E}) = y_1, \,\,\,\, h(\underline{E}) = y_2
$. Let $\Delta E:=\overline{E}-\underline{E} $, $\overline{E}'=\overline{E}-\Delta E/4$ and $\underline{E}' = \underline{E}+ \Delta E/4$,  and 
define $y_1':= h(\overline{E}')$, $y_2':=h(\underline{E}')$.  Then
\[
 y_2'-y_1'> C(V)^{-1}(\overline{E}-\underline{E})  > C(V)^{-1} \exp(-\underline{N}^\delta)
\]
Hence, by the Diophantine nature of $\omega$, there exists $\exp(\underline{N}^\delta )\le m \leq
\exp(2\underline{N}^\delta )$ so that $\{ m\omega \} \in (y_1',
y_2')$. Consequently, there is a unique $E_0 \in (\underline{E},
\overline{E})$ so that $h(E_0) = \{m\omega\}$. Set $x_0 :=
\bigl(E_{j1}^{(\underline{N})}\bigr)^{-1} (E_0)$. By construction,
$\underline{x}_1 < x_0 < \overline{x}_1$ and
$$
E_{j_2}^{(\underline{N})} (x_0 + m\omega) = E_0 =
E_{j_1}^{(\underline{N})} (x_0)
$$
as desired. Moreover, \eqref{eq:x0_sep} follows from 
\[
 \dist(E_0,\{\overline{E},\underline{E}\})> \Delta E/4
\]

\item[(ii)]
This part follows from part $(i)$ due to
Proposition~\ref{prop:13.doubleresonancegraphs}.
\end{itemize}
\end{proof}


For the rest of this section we fix a regular positive-slope
$I$-segment $ \bigl\{E_{j_1}^{(\underline{N})} (x), \underline{x}_1,
\overline{x}_ 1\bigr\}$, a regular negative-slope $I$-segment
 $\bigl\{E_{j_2}^{(\underline{N})} (x), \underline{x}_2,
\overline{x}_2 \bigr\}$, $I = \bigl[\underline{E},
\overline{E}\bigr]$, $\overline{E} - \underline{E}
> \exp(-{\underline{N}}^\delta)$,
\begin{itemize}
\item[(i)]$|\partial_x E_{j_s}^{(\underline{N})}|\ge
\exp(-\uN^{\delta})$ for any $x\in (\underline{x}_s, \overline{x}_s)
$, $s=1,2$

\item[(ii)]
$\overline{x}_s -\underline{x}_s\ge \exp(-\uN^{\delta})$
\end{itemize}
where $\delta>0$ can be made arbitrarily small but fixed and $N$ large depending on~$\delta$. 
We also fix an integer $m \in [\exp{(\underline{N}}^\delta), \exp{(2\underline{N}}^\delta)]$
and a point $x_0 \in (\underline{x}_1, \overline{x}_ 1)$ such that
condition \eqref{eq:13.1} holds. As usual, we denote the eigenvalues of
$H_{[-N, N]} (x,\omega)$ by $E_j^{(N)} (x)$, $1\le j\le 2N+1$, and  a normalized
eigenfunction corresponding to $E_j^{(N)}(x)$ by $\psi_j^{(N)}(x, \cdot)$. As in the proof of
the previous lemma, 
$E_0:=E_{j_1}^{(\underline{N})}(x_0)$.  In the following proposition, we pass to a larger scale. 
The idea is as follows: since they are {\em regular} the segments $E_{j_1}^{(\underline N)}(\cdot)$ and $E_{j_1}^{(\underline N)}(\cdot +m\omega)$ correspond to eigenfunctions supported on $[-\underline N, \underline N]$, and $[m-\underline N, m+\underline N]$,
respectively which are exponentially small near the edges of these intervals. Hence, they each generate an approximate
eigenstate of the operator $H_{[-N,N]}(x,\omega)$ with eigenvalues close to~$E_0$.  Proposition~\ref{prop:13.resexistprop}
quantifies these qualitative properties. What remains to be done as far as Figure~2 is concerned is to show that there
is a true bottom arc $E^{-}$, i.e., an arc that achieves its maximum in the interior of the interval on which it is defined.  

\begin{prop}
 \label{prop:13.resexistprop}
Let $0<\eps_0\le\frac{1}{20}$ be arbitrary but fixed, and $0<\delta\ll \frac{\eps_0}{A}$ where $A=A(V,\gamma,a,c)$.  Set $N:=
[\exp(\underline{N}^{6\delta})]$, $n:=[(\log N)^A]$. 
We assume that \[\omega\in \tor_{c,a}\setminus (\cB'_{n}\cup \cB''_n
\cup \tilde {\Omega}_{n})\] as in Corollary~\ref{cor:triple_elim}. 
\begin{itemize}
 \item[(a)] For any interval $[N', N'']$, with $n\le N''-N'\le 100n$,
  $N^{10\eps_0} \le |N'| \le 2N $ and any $|x-x_0|<N^{-d-1}$ where $d$ is
  as in Corollary~\ref{cor:triple_elim}
 one has
   $$\spec (H_{[N',N'']} (x) ) \cap (E_0 - \kappa, E_0 + \kappa) = \emptyset$$
   with $\kappa := N^{-d}$

\item[(b)] For any $x_0 - N^{-d-1} \leq x \leq x_0 + N^{-d-1}$ there exists $j'_1$,
$j'_2$ (possibly depending on $x$) such that 
\[
|E_{j_1}^{(\uN)} (x)- E_{j'_1}^{(N)}(x)| \le \exp(-\uN^{1/3}) ,\quad \,\
|E_{j_2}^{(\uN)} (x+m\omega)- E_{j'_2}^{(N)}(x)|\le \exp(-\uN^{1/3})
\]
\item[(c)] For any $x_0 - N^{-d-1} \leq x_{-} \leq x_0-e^{-\underline{N}^{\frac13}}$  $($resp.~$x_0+e^{-\underline{N}^{\frac13}}\leq
x_+\leq x_0 + N^{-d-1})$ there exist $j'_-, j''_-$, respectively $j'_+,
j''_+$,  $($possibly depending on $x_{-},x_{+}$$)$ such that 
\begin{align*}
E_0 - C(x_0 -x_{-}) &<E_{j'_{-}}^{(N)}(x_{-})<
E_0-(x_0-x_{-})\exp(-\uN^{\delta})
\\E_0 - C(x_{+}
-x_0) &<E_{j''_{-}}^{(N)}(x_{+})<E_0-(x_0-x_{+})\exp(-\uN^{\delta})
\\
E_0+(x_0-x_{-})\exp(-\uN^{\delta})&<E_{j'_{+}}^{(N)}(x_{-})<E_0+C(x_0-x_{-})
\\
E_0+(x_{+}-x_0)\exp(-\uN^{\delta})&<E_{j''_{+}}^{(N)}(x_{+})<E_0+C(x_{+}-x_0)
\end{align*}
\item[(d)] If $E_j^{(N)} (x) \in (E_0 - N^{-d}/2, E_0 +N^{-d}/2)$ for some $x \in
[x_0 -N^{-d-1}, x_0 + N^{-d-1}]$ then
\[
\bigl | \psi_j^{(N)} (x, n') \bigr | \leq \exp({-\gamma}|n'|/4)\qquad\forall\; |n'| > N^{\frac{1}{2}}
\]

\item[(e)] If $E_j^{(N)} (x), E_{j'}^{(N)} (x) \in (E_0 - N^{-d}/2 , E_0 + N^{-d}/2)$, $j
\neq j'$, $x\in [x_-, x_+]$ then
\[
|E_j^{(N)} (x) - E_{j'}^{(N)}(x) | \geq \tau=\exp(-N^{11\eps_0})
\]
 \end{itemize}
\end{prop}
\begin{proof} To prove $(a)$ we use Corollary~\ref{cor:triple_elim}
Recall that we assume that $\omega\in \tor_{c,a}\setminus (\cB'_{n}
\cup \cB''_{n} \cup \tilde {\Omega}_{n})$. We set
$\delta_0:=\delta$, $\delta_1:=3\delta$, $\varepsilon_0:=1/50$,
$n_1:=\uN$, $n_2=n_1$, $n_3:=N''-N'+1$, $m_1:=m$, $m_2:=N'$. Then
all conditions of Corollary~\ref{cor:triple_elim} hold. Let
$|x-x_0|\le N^{-d-1}$. Then 
\[ |E_{j_1}^{(\underline{N})} (x) -
E_{j_2}^{(\underline{N})} (x + m_1\omega)|\le C(V)|x-x_0| \le
 N^{-d}
\]
 Due to Corollary~\ref{cor:triple_elim} one has
$|E_{j_1}^{(\underline{N})}(x)-E_{j}^{(n_3)}(x+m_2\omega)|$ $\ge
N^{-d}$ for any $1\le j\le 2n_3+1$ which proves $(a)$. Since
$\bigl\{E_{j_1}^{(\underline{N})} (x), \underline{x}_1,
\overline{x}_1 \bigr\}$ is a regular segment we have 
$\underline{N}^{1/2}\le$ $\nu_{j_1}^{(\underline{N})}(x)\le
\underline{N}-\underline{N}^{1/2}$. Therefore,
\[
\|(H_{[-N,N]}(x)-E_{j_1}^{(\underline{N})}
(x))\psi_{j_1}^{(\underline{N})} (x,\cdot)\|\le \exp(-\uN^{1/3})
\]
whence there exists $j'_1$ such that
\[
|E_{j_1}^{(\uN)} (x)- E_{j'_1}^{(N)}(x)| \le \exp(-\uN^{1/3})
\]
The proof of the statement  regarding $E_{j_2}^{(\uN)} (x)$
is similar. Thus $(b)$ holds. Because of~\eqref{eq:x0_sep} and
\begin{equation} \label{eq:13.slopes}
 \begin{aligned}
  \partial_xE_{j_1}^{(\uN)} (x) &> \exp(-\uN^\delta ), \,\, x \in [\ux_1, \ox_1]\\
  \partial_xE_{j_2}^{(\uN)} (x) & < \exp(-\uN^\delta ), \,\, x \in [\ux_2, \ox_2]
 \end{aligned}
\end{equation}
part $(c)$ follows from $(b)$. To validate $(d)$ and $(e)$ we invoke
Corollary~\ref{cor:sep}. Indeed,
 part $(a)$ of the current proposition implies that the
condition needed for Corollary~\ref{cor:sep} are met. Therefore,  $(d)$ and
$(e)$ hold.
\end{proof}

\begin{remark} \label{rem:18.rightwayelim} The previous proposition is based on the elimination
via Corollary~\ref{cor:triple_elim} and not via resultants as in Proposition~\ref{prop:2.zerosepar}. This is crucial, as the 
latter method of elimination requires the removal of some subset of
the energy~$E$.  Although the subset in
Proposition~\ref{prop:2.zerosepar} is small, its removal would
destroy the argument.
\end{remark}

The following lemma establishes the existence of the $E^-$ branch in Figure~2. 

\begin{lemma}\label{lem:13.localmax} For $N$ large, there exists $j_0$ with
\begin{equation}\label{eq:13.margins}
\bigl |E_{j_0}^{(N)}(x_0) - E_0 \bigr | < C(V)N^{-d-4}
\end{equation}
and an interval $[x_-, x_+]$,
\begin{equation}\label{eq:13.marginsx}
  0 < x_0-x_-\le N^{-d-4} ,\qquad 0 < x_+ - x_0 \leq
N^{-d-4}
\end{equation}
such that $E_{j_0}^{(N)}(x)$ assumes its maximum over
the interval $[x_-, x_+]$ at some point $x^{(0)}$ which is located
``properly inside the interval'', i.e., 
\begin{equation}\label{eq:13.maxininside}
x_- + C^{-1}N^{-d-13} < x^{(0)} < x_+ - C^{-1}N^{-d-13}
\end{equation}
\end{lemma}
\begin{proof}
To establish this lemma we use
Proposition~\ref{prop:13.resexistprop}
 and the property that the graphs of $E_j^{(N)}(x)$, $1\le j\le 2N+1$ never intersect. Set
\begin{align*}
x_- (k) &= x_0 - kN^{-d-8}, \,\,\, 1\le k\le N^4\\
E_-(k) &=E_{j_1}^{(\uN)}(x_- (k) ) \\
J_- (k) &= \bigl \{ j : \bigl | E_j^{(N)}(x_-(k)) - E_-(k)\bigr | <
\exp(-\uN^{\frac{1}{3}})\bigr \}
\end{align*}
where $k\ge1$. Due to part $(b)$ in
Proposition~\ref{prop:13.resexistprop} each set $J_-(k)$ is
non-empty. Furthermore, due to the fact $\partial_x E_{j_1}^{(\uN)}
(x)
> \exp(-\uN^\delta)$ one also has
\[E_- (k - 1) > E_-(k)+N^{-d-8}
\exp(-\uN^\delta)
\]
Let $j \in J_-(k)$ for some $k$. Assume that
\begin{equation}\label{eq:13.growdicht}
\max \bigl \{ E_j^{(N)}(x) \;:\; x \in [x_-(k), x_0] \bigr\} < E_-(k-N)-
\exp(-\uN^{1/3})
\end{equation}
Then clearly one has
\begin{equation}\label{eq:13.jparish}
j \not\in J_-(k') \text{\  for any\  } k' < k-N
\end{equation}
Since there are only $2N+1$ eigenvalues $E_j^{(N)} (x)$ in total
there exist $N^3\ge k_{-} \ge N^3/2$ and $j_- \in
J_-(k_{-})$ such that
$$
\max \bigl\{E_{j_-}^{(N)} (x) \;:\; x \in [x_-(k_-), x_0] \bigr\} \geq
E_-(k_- -N)-\exp(-\uN^{1/3})
$$
Analogously, we set 
\be
x_+(k) & := & x_0 + kN^{-d-11} , \,\,\,  1\le k\le  N^7,\nonumber \\
E_+(k) & := & E_{j_2}^{(\uN)} (x_+ (k)+m\omega ),\nonumber \\
J_+(k) & := & \bigl\{j \;:\; |E_j^{(N)}(x_+(k)) - E_+ (k) | <
\exp(-\uN^{\frac{1}{3}}) \bigr\} \nonumber \ee 
We want to show that
there exist $k_+$ and $j_+ \in J_+ (k_+)$
 such that the following conditions hold
\begin{enumerate}
\item[$(\alpha)$] $\max \bigl\{ E_{j_+}^{(N)} (x) : x \in [ x_0, x_+ (k_+) ] \bigr\}
\geq E_+ (k_+ - N) - \exp(-\uN^{\frac{1}{3}})$,

\item[$(\beta)$] $E_-(k_-) - N^{-d-11}   < E_+(k_+) < E_-(k_-) + CN^{-d-8}$.
\end{enumerate}
 To achieve this note that due to~\eqref{eq:13.slopes} and $k_- > N^3/2$, one has
\[ E_-(k_-) < E_0 - \frac{1}{2}N^{-d-5} \exp(- \uN^\delta )\]
On the other hand,
\[E_-(k_-)\ge E_-(N^3) > E_0 - CN^{-d-5}
\]
 Once again, due to~\eqref{eq:13.slopes}, \[ E_{j_2}^{(\uN)} (x_+(N^7)+m\omega) < E_0 - N^{-d-4} \exp(-\uN^\delta) \]
 Hence, there exists $k'_+ \le N^7$ such that $E_+(k'_+ +1) \leq E_-(k_-) < E_+(k'_+)$.
 Since
 \[\begin{split}
& E_0 - \frac{1}{2}N^{-d-5}\exp(- \uN^\delta )\ge E_-(k_-) \ge
E_{j_2}^{(\uN)} (x_+(k'_{+}+1)+m\omega)
> E_0 - Ck'_+ N^{-d-11}
\end{split} \]
one has $N^7\ge k'_+ \gtrsim N^6\exp(- \uN^\delta )$. Applying an argument
analogous to the one used to determine $k_-$, one
 finds $k_+ \in [k'_+ - N^3, k'_+]$ and $ j_+ \in J_+(k_+)$
  such that condition $(\alpha)$ holds. Furthermore,
\begin{equation}
\nn
\begin{split}
E_+(k_+) &> E_+ (k'_+) > E_-(k_-)\\
 E_+(k_+) &< E_+(k'_+ -
N^3) < E_+(k'_+) + CN^{-d-8}\\
&<E_+(k'_+ +1)+ CN^{-d-11} + CN^{-d-8}\\ 
&< E_-(k_-) + CN^{-d-11} + CN^{-d-8}
\end{split}
\end{equation}
Therefore, condition $(\beta)$ holds. To finish the proof of the
lemma assume for instance that
\begin{equation}\label{eq:13.left<right}
E_{j_-}^{(N)} (x_+ (k_+  ) ) \leq E_{j_+}^{(N)} (x_+ (k_+))
\end{equation}
Then
\begin{equation}
\nn
\begin{split}
&\max \bigl\{ E_{j_-}^{(N)} (x) \; :\; x \in [ x_+(k_+) - N^{-d-11},
x_+(k_+)] \bigr\} \\&\le E_{j_-}^{(N)} (x_+ (k_+  ) ) +CN^{-d-11} 
\leq E_{j_+}^{(N)} (x_+ (k_+)) +CN^{-d-11}\\
&\le
E_+(k_+)+\exp(-\uN^{\frac{1}{3}})+CN^{-d-11}\\
& \le E_{j_-}^{(N)} (x_- (k_-))
+ 3CN^{-d-8}
\end{split}
\end{equation}
Furthermore, 
\begin{equation}
\nn
\begin{split}
\max \bigl\{ E_{j_-}^{(N)} (x) \;:\; x \in [ x_-(k_-),x_-(k_-) + N^{-d-8}
] \bigr\}\le E_{j_-}^{(N)}(x_-(k_-))+CN^{-d-8}
\end{split}
\end{equation}
On the other hand,
\begin{equation}
\nn
\begin{split}
&\max \bigl\{ E_{j_-}^{(N)} (x) \;:\; x \in [ x_-(k_-), x_+ (k_+)] \bigr\}
\geq \max \bigl\{ E_{j_-}^{(N)} (x) \;:\; x \in [ x_-(k_-), x_0] \bigr\}\\ &
\geq
E_-(k_- -N)-\exp(-\uN^{1/3}) \ge E_{j_-}^{(N)}(x_-(k_-))+N^{-d-7}\exp(-\uN^{\delta})
-2\exp(-\uN^{1/3})\\
& \ge E_{j_-}^{(N)}(x_-(k_-))+\frac
{1}{2}N^{-d-7}\exp(-\uN^{\delta})
\end{split}
\end{equation}
Set $j_0 = j_-$, $[x_-, x_+]:= [ x_-(k_-), x_+ (k_+)] $. Then
\eqref{eq:13.maxininside}  holds. Note that
\[
0 <x_0-x_-(k_-)\le N^{-d-4}, \quad 0 < x_+ (k_+)- x_0 \le N^{-d-4}
\]
since $k_-\le N^4$, $k_+\le N^7$ as well as 
\[\begin{split}
|E_{j_-}^{(N)} (x_0)-E_0| &\le |E_{j_-}^{(N)}(x_0) - E_{j_-}^{(N)} (x_-(k_-))| +
|E_{j_-}^{(N)} (x_-(k_-))-E_-(k_-)| +|E_-(k_-)-E_0|\\& \le C(V)N^{-d-4}
\end{split}
\]
Hence
\eqref{eq:13.margins}, \eqref{eq:13.marginsx} also hold. Consider
now the case opposite to~\eqref{eq:13.left<right}. Since
the graphs of $E_{j_-}^{(N)} (x)$ and $E_{j_+}^{(N)} (x)$ do not intersect one
has in this case
\[E_{j_+}^{(N)} (x) < E_{j_-}^{(N)} (x)
\]
for all $x$. Hence,
\begin{equation}
\nn
\begin{split}
&\max \bigl\{ E_{j_+}^{(N)} (x) : x \in [ x_-(k_-),x_-(k_-) + N^{-d-13}]
\bigr\} \le \max \bigl\{ E_{j_-}^{(N)} (x) : x \in [ x_-(k_-),x_-(k_-) +
N^{-d-13} ] \bigr\} \\
&\le E_{j_-}^{(N)} (x_-(k_-)) +CN^{-d-12}\le E_-(k_-) +
\exp(-\uN^{1/3})+CN^{-d-12} \le E_+(k_+)+2N^{-d-11}\\ &\le E_{j_+}^{(N)} 
(x_+ (k_+)) +\exp(-\uN^{1/3}) + 2N^{-d-11}\le  E_{j_+}^{(N)} (x_+ (k_+))
+3N^{-d-11}
\end{split}
\end{equation}
Furthermore, 
\begin{equation}
\nn
\begin{split}
\max \bigl\{ E_{j_+}^{(N)} (x) : x \in [ x_+(k_+)-N^{-d-13},x_+(k_+)]
\bigr \} \le E_{j_+}^{(N)} (x_+ (k_+)) +CN^{-d-13}
\end{split}
\end{equation}
 On the other hand,
\begin{equation}
\nn
\begin{split}
&\max \bigl\{ E_{j_+}^{(N)} (x) : x \in [ x_-(k_-), x_+ (k_+)] \bigr\} \geq
\max \bigl\{ E_{j_+}^{(N)} (x) : x \in [x_0, x_+(k_+)] \bigr \} \ge
 \\& E_+(k_+
-N)-\exp(-\uN^{1/3})\ge E_+(k_+)+
N^{-d-10}\exp(-\uN^{\delta})-\exp(-\uN^{1/3})
 \ge \\ &E_{j_+}^{(N)} (x_+ (k_+)) + N^{-d-10}\exp(-\uN^{\delta}) -
2\exp(-\uN^{1/3})\ge  E_{j_+}^{(N)} (x_+ (k_+)) + \frac {1}{2}
N^{-d-10}\exp(-\uN^{\delta})
\end{split}
\end{equation}
Set $j_0 = j_+$. Thus~\eqref{eq:13.maxininside} and
\eqref{eq:13.marginsx} follow. Finally,
\begin{align*}|E_{j_+}^{(N)} (x_0)-E_0| &\le |E_{j_+}^{(N)} (x_0) - E_{j_+}^{(N)} (x_+(k_+))| +
|E_{j_+}^{(N)} (x_+(k_+))-E_+(k_+)| +|E_+(k_+)-E_0|\\
& \le C(V)N^{-d-4}
\end{align*}
and \eqref{eq:13.margins} holds as claimed.
\end{proof}

For convenience we summarize the conclusions  of
Lemma~\ref{lem:13.localmax} and
Proposition~\ref{prop:13.resexistprop} in the following corollary. 

\begin{corollary}
\label{cor:13.pregapped} Using the notations of
Proposition~\ref{prop:13.resexistprop} there exists $j_0$ and an
interval $[x_-, x_+] \subset [x_0 - N^{-d-4} , x_0 + N^{-d-4}]$ such
that

\begin{enumerate}
\item[(1)] The function $E^{(-)}(x) := E_{j_0}^{(N)} (x), x \in [ x_-, x_+]$
assumes its maximal value of the interval $[x_-, x_+]$ at some point
$x^{(0)} \in [x_- + CN^{-d-13}, x_+ - CN^{-d-13} ]$

\item[(2)] For any $x \in [x_-, x_+]$ the operator $H_{[-N, N]} (x)$
has no eigenvalues in the interval \[(E^{(-)}(x), E^{(-)} (x) +
\tau),\quad \tau = \exp(-N^{11\eps_0})\]

\item[(3)] The eigenfunction $\psi^{(-)} (x, n') := \psi_{j_0}^{(N)}  (x, n')$,
 obeys $| \psi^{(-)} (x, n') | \leq \exp(- \frac{\gamma}{4} |n'|)$, for $|n'| \ge N^{1/2}$.
\end{enumerate}
\end{corollary}

Next, we will establish that any energy  $E\in (E^{(-)}(x^{(0)}),E^{(-)}(x^{(0)})+\tau)$
generates two complex zeros of~$f_N(\cdot,E)$ with imaginary part bounded below. 
In order to do this, we apply the Weierstrass preparation theorem as in Corollary~\ref{cor:weierz_loc} to
$f_N (z, E)$ relative to the $z$-variable locally around the point
\[ (z^{(0)}, E^{(0)}), \quad z^{(0)} =
e(x^{(0)}), \quad E^{(0)} = E^{(-)}(x^{(0)})
\]
Thus, there exist a
polynomial $P(z, E) = z^k + a_{k-1} (E)z^{k-1} +\cdots+ a_0(E)$ with
$a_j(E)$ analytic in $\cD(E^{(0)}, r_1)$, where $r_1 \asymp \exp
(-N^{12\eps_0})$, and an analytic function $g(z, E)$, $(z, E) \in \cP
= \cD (z^{(0)}, r_1) \times \cD(E^{(0)}, r_1)$ such that:

\begin{enumerate}
\item[(a)] $f_N(z,E) = P(z,E)g(z,E)$,

\item[(b)] $g(z,E) \neq 0$ for any $(z,E) \in \cP$,

\item[(c)] For any $E \in \cD(E_0, r_1)$, the polynomial $P(.,E)$
has no zeros in $\IC \backslash \cD(z^{(0)}, r_1)$,

\item[(d)] $1\leq k = \deg P_N(\cdot, \omega, E) \leq (\log N)^{C_0}$
\end{enumerate}

Here the property $k \geq 1$ is due to the fact that $f_N(z^{(0)},
E^{(0)}) = 0$. We can now derive the following result:

\begin{lemma}
\label{lem:13.wmplzers} For any $E \in (E^{(0)} + r_1/4, E^{(0)} +
r_1/2)$, the Dirichlet determinant $f_N (\cdot, E)$ has at least two
complex zeros $\zeta^\pm = \zeta^\pm(E) = e(x(E) \pm iy(E)) \in
\cD(e(x^{(0)}), r_1)$, with $r_1/C_1 < |y(E)| < r_1$.
\end{lemma}

\begin{proof} For any $E \in \cD(E_0, r_1)$, the polynomial $P(\cdot, E)$ has
at least one zero $\zeta(E)$, with $\zeta(E) \in \cD(z^{(0)}, r_1)$. Let
$E_1 \in (E^{(0)} +r_1/4, E^{(0)} + r_1/2)$
 be arbitrary, and let $\zeta_1 = \zeta(E_1) = e(x_1 + iy_1)$.
 Note that $E^{(0)} - C r_1 \leq E^{(-)} (x_1) \leq E^{(0)}$.
 Therefore,
 \begin{equation}
\nn
\begin{split}
 E^{(-)}(x_1) + r_1/4 \le E^{(0)} +r_1/4 \le E_1 \le E^{(0)} +r_1/2 \le
 E^{(-)}(x_1) + Cr_1 \le E^{(-)}(x_1) +\tau/2
 \end{split}
 \end{equation}
Furthermore, recall that due to condition
 $(2)$ of Corollary~\ref{cor:13.pregapped}, the operator $H_N(x_1)$ has
 no eigenvalues in the interval $(E^{(-)} (x_1), E^{(-)} (x_1)+ \tau)$.
 In particular,
$$
\dist \big[ \spec \big( H_{[-N, N]} (x_1)\big), E_1 \big] \geq r_1/4.$$
Since $H_N(x_1)$ is self adjoint and \[ \|H_{[-N, N]} (x_1) - H_{[-N,
N]} (x_1 + iy_1) \| < C|y_1|, \] one has also
$$
0 = \dist \big[ \spec \big(H_{[-N, N]} (x_1 + iy_1)\big), E_1 \big] > r_1/4 - C|y_1|
$$
The determinant $f_N(e(x), E_1)$ assumes only real values for real
$x$. Therefore,  each complex
zero $\zeta_1$ produces a conjugate zero and we are done.
\end{proof}

The following is the main result of this section. Due to the fact that
an eigenfunction $\psi$ of $H_{[-N,N]}(x,\omega)$ which is very well localized
in~$[-N,N]$ remains close to an eigenfunction of this operator if it is translated
inside of the interval we obtain a whole sequence of complex zeros as in the previous lemma.
The parameters $\delta,\eps_0$
are as above. 

\begin{prop}\label{20.regpsitiveregnegative}
Using the notations of part $(ii)$
 in Lemma~\ref{prop:13.resexist}, assume that
 \[\big(E'_0+\frac{1}{4}\exp(-(\log \ell)^A),\;E''_0-\frac{1}{4}\exp(-(\log \ell)^A)\big)\]
is not spectrum-free. Then there exists $N=
[\exp(\underline{N}^{6\delta})]$  with $\ell^2 \le \underline{N}\le
\ell^{10}$ and an interval $(E',
E'')\subset (E'_0,E''_0)$, $E''-E'=$ $\exp(-N^{3\delta})$
 such that for any $E\in (E', E'')$ the Dirichlet determinant
$f_N(\cdot, \omega, E)$ has a sequence of zeros $\zeta^\pm_k = e(x_k
\pm iy_k)$, where $k$ runs in the interval $(-N+2N^{1/2}, N -
2N^{1/2})$, with $$\| x_k - x_0 - k\omega\| < \exp(-N^{1/8}), \quad
\exp(-N^{13\eps_0})< |y_k| \le \exp(-N^{10\eps_0})$$
\end{prop}
\begin{proof} Due to part $(a)$ of Proposition~\ref{prop:13.resexistprop} one has
$$\spec \big(H_{[N',N'']} (x,\omega) \big) \cap (E_0 - \kappa, E_0 + \kappa) = \emptyset$$
for any $|x-x_0|\le N^{-d-1}$ and any interval $[N', N'']$ with $$ 
n\le N''-N'\le 100n,\quad
  N^{10\eps_0} \le |N'| \le 2N $$
where $n$ is as above, $\kappa := N^{-d}$. Therefore, due to
Proposition~\ref{prop:add_zeros}, one obtains
\begin{equation}\label{eq:16.zerofree}
\nu_{f_{[N_1,N_2]}(\cdot, E)}(e(x), R_0)=0
\end{equation}
for any $$[N_1,N_2]\subset [-2N,2N]\setminus [-2N^{1/4},2N^{1/4}],\qquad
N_2-N_1\ge n,$$ and any $E\in (E_0 - \kappa/2, E_0 +
\kappa/2)$, where $R_0=:\exp(-(\log N)^C)$. Let $x^{(0)}$, $E \in
(E^{(0)} + r_1/2, E^{(0)} + r_1/2)$ be as in
Lemma~\ref{lem:13.wmplzers}. Then the Dirichlet determinant $f_N
(\cdot, E)$ has a complex zero $$\zeta^+(E) = e(x(E) + iy(E)) \in
\cD(e(x^{(0)}), r_1)$$ with $y(E)
> r_1/C_1$. Set $N_0:=[N^{1/3}]$. Then  \eqref{eq:16.zerofree} implies
in particular that $N_0$, $-N_0$ are adjusted to $(\cD(\zeta^+(E),
R_0),E^{(0)})$ at scale $\ell_0:=[N^{1/4}]$.  Due to
Proposition~\ref{prop:add_zeros}, with \eqref{eq:16.zerofree} taken
into account, one obtains
\begin{equation}\label{eq:16.zeroin}
\nu_{f_{[-N_0,N_0]}(\cdot, E)}(\zeta^+(E),r_0)\ge 1
\end{equation}
where $r_0=\exp(-\ell_0^{1/2})\asymp \exp(-N^{1/8})$. Let $k\in
[-N+2N^{1/2},N-2N^{1/2}]$ be arbitrary. Recall that
\[
f_{[a,b]}(e(x+k\omega+iy),E)=f_{[a+k,b+k]}(e(x+iy),E)
\]
 for any $x$, $y$, $E$ and any  $[a,b]$. Using once again
 Proposition~\ref{prop:add_zeros}  with \eqref{eq:16.zerofree}
 and \eqref{eq:16.zeroin} taken into account yields
\begin{equation}
\nn
 \nu_{f_{[-N,N]}(\cdot, E)}(\zeta^+(E)e(k\omega),r_0)\ge 1
\end{equation}
Hence there exists $\zeta^+_k(E)$ as claimed. The argument for
$\zeta^-_k(E)$ is similar.
\end{proof}

By a pre-gap we mean an interval of energies $(E',E'')$ as described in Proposition~\ref{20.regpsitiveregnegative}.
It remains to show that pre-gaps do not fill up again as we pass to larger scales. The following section makes this
precise. 

\section{Proofs of Theorems~\ref{th:1.mainth},~\ref{th:1.2}}\label{sec:proofs}

To prove Theorem~\ref{th:1.2} we use
Proposition~\ref{20.regpsitiveregnegative} in combination with the
analysis of the quantities
\begin{equation}\label{eq:20.avernzero}
\cM_N(E, R_1, R_2) = \tfrac{1}{N} \# \bigl\{z \in \cA_{R_1, R_2}:
f_N\zoe = 0\bigr\}
\end{equation}
from Section~\ref{sec:jensen}. Recall the following relations established for
these quantities, cf.~Lemma~\ref{lem:2.avernstab}:
\begin{equation}\label{eq:20.avernzel}
\begin{aligned}
\cM_N(\omega, E, R_1 + r_2, R_2 - r_2) & \le \cM_n(\omega, E, R_1 - r_2, R_2 + r_2) + n^{-1/4}\\[5pt]
\cM_n(\omega, E, R_1 + r_2, R_2 - r_2) & \le \cM_N(\omega, E, R_1 -
r_2, R_2 + r_2) + n^{-1/4}
\end{aligned}
\end{equation}
 for any $n > N_0(V,\gamma,a,c,\sigma)$, $N > \exp(\gamma\, n^\sigma)$,
$1- \rho^{(0)} < R_1 < R_2 < 1 + \rho^{(0)}$ where $N_0=N_0(V, c, a,
\gamma)$, $\rho^{(0)} = \rho^{(0)}(V, c, a, \gamma)
> 0$, $r_2 =
n^{-1/4}(R_2 - R_1)$, provided $r_2>\exp(-\gamma\, n^\sigma/100)$. Here $\sigma>0$ is small but fixed 
and $\gamma$ stands for the lower bound on the Lyapunov exponent
as usual. In what follows, $0<\eps_0\ll1$ is small and fixed as in the previous section and we will
set $\sigma:=20\eps_0$. 

\begin{lemma}\label{lem20:'densityofzerosfalls1}
Using the notations of Proposition~\ref{20.regpsitiveregnegative}
one has for any $E\in (E', E'')$, and any $N_1 > \exp(\gamma\,
N^{20\eps_0})$
\begin{equation}\label{densityofzerosfalls3}
\cM_{N_1}(E, R'_{-}(N), R'_{+}(N))\le \cM_N(E, R''_{-}(N),
R''_{+}(N)) - 2 + CN^{-1/4}
\end{equation}
where $R'_{\pm}(N)=1 \pm \exp(-N^{14\eps_0})$, $R''_{\pm}(N)=1 \pm
\exp(-N^{9\eps_0})$
\end{lemma}
\begin{proof} Due to Proposition~\ref{20.regpsitiveregnegative} for
any $E\in (E', E'')$, the Dirichlet determinant $f_N(\cdot, \omega,
E)$ has a sequence of zeros $\zeta^\pm_k = e(x_k \pm iy_k)$, where
$k$ runs in the interval $(-N+2N^{1/2}, N - 2N^{1/2})$,
\[ \| x_k - x_0
- k\omega\| < \exp(-N^{1/8}), \quad  \exp(-N^{13\eps_0})< |y_k| \le
\exp(-N^{10\eps_0})
\]
Since $\omega \in \tor_{c,a}$, all these zeros are different, i.e., 
$\zeta^\pm_k \neq \zeta^\pm_{k_1}$ if $k\neq k_1$. Hence, 
\[
\cM_N(E,R_{1,-},R_{1,+})\le
\cM_N(E,R''_{-}(N),R''_{+}(N))-2+2N^{-1/2}
\]
where $R_{1,\pm}=1 \pm \exp(-2N^{13\eps_0})$. Combining this estimate
and \eqref{eq:20.avernzel} one obtains \eqref{densityofzerosfalls3}.
\end{proof}

We can now prove Theorems~\ref{th:1.2} and~\ref{th:1.mainth}.

 \begin{proof} [Proof of Theorem~\ref{th:1.2}] Following the
 notations of Remark~\ref{rem:17.parameters} given $N$ we denote
 by $\cB(N)$ the set of exceptional values of $\omega$ defined in
 Corollary~\ref{cor:triple_elim}. Recall that due to Remark~\ref{rem:17.parameters}
 one has
\[
\cH^{\alpha(d)}_{r(N)}(\cB(N))\le \exp(-(\log\log N)^B)
\]
where $\cH^\alpha_{r}$ stands for the corresponding Hausdorff outer
measure, $r(N)=$ $\exp(- (\log \log N)^A)$, $\alpha(d)=$ $4/(d-3)$,
$A,B>1$ are constants,  and $d>7$ is arbitrary and provided $N\ge
N_0(V,c,a,\gamma,d)$. We need to iterate the result of
Lemma~\ref{lem20:'densityofzerosfalls1}. Set with some $0<\delta\ll \eps_0$ as in the previous section, 
$\oN(N,1):=N,\oN(N,t+1):=[\exp((\oN(N,t))^{\delta})]$ for all $t\ge1$, and
\[
\cB(N,T):=\bigcup_{1\le t \le T} \cB(\oN(N,t))
\]
Assume that $\omega \in \tor_{c,a}\setminus \cB(N,T)$. Given
arbitrary $1\le t\le T$ and an interval $(E_{t,1},E_{t,2})$ with \[
E_{t,2}-E_{t,1}\ge \exp(-(\log \oN(N,t))^C),
\] either
$(E_{t,1},E_{t,2})$ contains a spectrum-free subinterval $(E'_{t,1},E'_{t,2})$ with 
\[  E'_{t,2}-E'_{t,1}\ge \exp(-(\log
\oN(N,t))^C),\] or there exists a subinterval $(E'_t, E''_t)$ with  \[
E''_t-E'_t\ge \exp(-\oN(N,t+1)^{3\delta})\] such that
\begin{equation}\label{densityofzerosfalls2}
\cM_{\uN(N,t+1)}(E, R'_{-}(\oN(N,t)), R'_{+}(\oN(N,t)))\le \cM_N(E,
R''_{-}(\oN(N,t)), R''_{+}(\oN(N,t))) - 2 + C\oN(N,t)^{-1/4}
\end{equation}
for any $E\in (E'_t, E''_t)$. 
Hence, given an arbitrary interval $(E_{1},E_{2})$ with \[ E_{2}-E_{1}\ge
\exp(-(\log N)^C),\] either $(E_1,E_{2})$ contains a spectrum-free
subinterval $(E'_{T,1},E'_{T,2}) $ with \[ E'_{T,2}-E'_{T,1}\ge \exp(-(\log
\oN(N,T))^C),\] or there exists a subinterval $(E'_T, E''_T)$ with \[ 
E''_T-E'_T\ge \exp(-\oN(N,T+1)^{3\delta})\] such that
\begin{equation}\label{densityofzerosfalls4}
\cM_{\uN(N,T+1)}(E, R'_{-}(\oN(N,T)), R'_{+}(\oN(N,T)))\le \cM_N(E,
R''_{-}(N), R''_{+}(N)) - 2T + CN^{-1/4}
\end{equation}
for any $E\in (E'_T, E''_T)$. 
 Recall that due
to Lemma~\ref{lem:4.3} on translations of regular $I$-segments one
has in particular
\[
\cM_N(E, R_1, R_2)\ge 1- 2N^{-1/2}
\]
for any $E\in I$, and any $R_1$, $R_2$, provided there exists at least
one regular segment $\bigl\{E_{j}^{(N)} \, (x, \omega),
\underline{x}, \overline{x} \bigr\}$. On the other hand, $ \cM_N(E,
R_1, R_2)\le C(V)
$
for any $N$, $R_1$, $R_2$. In particular, using the notations of
\eqref{densityofzerosfalls3} with $T:=[C(V)/2]+1$ one has
\[
\cM_{\uN(N,T+1)}(E, R'_{-}(\oN(N,T)), R'_{+}(\oN(N,T))) \le \cM_N(E,
R''_{-}(N), R''_{+}(N)) - 2T + CN^{-1/4} \les N^{-1/4}
\]
for any $E\in (E'_T,E''_T)$. Combining this estimate with
\eqref{eq:20.avernzel} one obtains
\[
\cM_{N_1}(E, R'_{-}(\uN(N,T+1)), R'_{+}(\uN(N,T+1))) \le \cM_N(E,
R''_{-}(N), R''_{+}(N)) < 5N^{-1/4}
\]
for any $N_1\ge \uN(N,T+2)$. Consequently,  there is no regular
$I$-segment \[ \bigl\{E_{j}^{(N_1)} \, (x, \omega), \underline{x},
\overline{x} \bigr\}, \qquad I\subset (E'_T,E''_T)\] Due to
Proposition~\ref{prop:13.doubleresonancegraphs} the interval
$(E'_T,E''_T)$ is spectrum-free. Finally, set $\cB_{N,d}:=\cB(N,T)$.
Then
\[
\cH^{\alpha(d)}_{r(N)}(\cB_{N,d})\le T+1\le
C(V)
\]
and we are done. 
\end{proof}

\begin{proof}[Proof of Theorem~\ref{th:1.mainth}] It is convenient to split the proof
into a  few steps. We enumerate these steps as $a$, $b$, $c$, and $d$.
\begin{itemize}
\item[(a)] Using the
notations of Proposition~\ref{prop:14.basic} assume that $\omega \in
\tor_{c,a} \backslash \hat{\Omega}_{N_1}$,
$x\in\tor\setminus\hat{\cB}_{N_1,\omega}$ for some $N_1\ge N_0$.
Assume that $(E'_0,E''_0)$ is a spectrum free interval. Then for any
$N_2$ there exists $N\ge N_2$ such that
\[
\bigcup_{x\in \TT} \spec \big( H^{(P)}_{[-N+1,N]}(x,\omega) \big) \cap
(E'_0,E''_0)=\emptyset
\]
Let $E_0:=(E'_0+E''_0)/2$, $\sigma_0:=E''_0-E'_0$. Then
\begin{equation}
\label{eq:17.nonspec}
\begin{split}
\min\Big\{|E_0-E_j^{(N)}(x, \omega)| \::\: E_j^{(N)}(x,\omega)\in(E',E''),
\nu_j^{(N)}(x,\omega)\in[-N+N^{\frac{1}{2}},\;
 N-N^{\frac{1}{2}}] \Big\}\ge \sigma_0/8
\end{split}
\end{equation}
Indeed, assume that $|E_0-E_{j_0}^{(N)}(x, \omega)|\le \sigma_0/4$
for some $j_0$ with \[ E_{j_0}^{(N)}(x,\omega)\in(E',E''),\quad \nu_{j_0}^{(N)}(x,\omega)\in[-N+N^{\frac{1}{2}},N-N^{\frac{1}{2}}] \]
Then
\[
\|(H^{(P)}_{[-N+1,N]}(x,\omega)-E_{j_0}^{(N)}(x,\omega))\psi_{j_0}^{(N)}(x,\omega,\cdot)\|\le
\exp(-\gamma N^{1/2})\le \sigma_0/4
\]
Hence, \[\dist \big[ E_0,\spec \big( H^{(P)}_{[-N+1,N]}(x,\omega)\big) \big] \le
\sigma_0/2\] contrary to our assumption. Thus \eqref{eq:17.nonspec}
is valid. Due to Remark~\ref{rem:14.spectrumfree}  for any $x'\in
\TT$ and any $|E-E_0|\le \sigma_0/8$ one has  $E\notin \spec (
H(x',\omega))$.
\item[(b)] Assume that $\omega \in
\tor_{c,a} \backslash (\hat{\Omega}_{N_1}\cup \cB_{N_1,d})$, for
some $N_1\ge N_0$, where $\cB_{N_1,d}$ is defined as in
Theorem~\ref{th:1.2}. Let $(E',E'')$ be an arbitrary interval. Then
$(E',E'')$ contains a spectrum free interval $(E'_0,E''_0)$. Then
due to $(a)$, $(E'_0,E''_0)$ contains an interval $(E'_1,E''_1)$ such
that $\Sigma_{\omega}\cap (E'_1,E''_1)=\emptyset$
\item[(c)] Set $\Omega':=\bigcap_{N\ge N_0}\hat{\Omega}_{N_1}$.
Then $\Omega'$ has Hausdorff dimension zero. Furthermore, proceeding 
inductively,  for each $d\ge 7$ pick an arbitrary $N(d)$ large enough so
that $\cB_{N(d),d}$ is defined, and also $N(d) \ge \max
(\exp(\exp(d)),N(d-1))$. Set
\[
\Omega''_d:=\bigcup_{d'\ge d} \cB_{N(d'),d'}, \quad \Omega'':=
\bigcap_d \Omega''_d
\]
If $d'\ge d$ then $\alpha(d')\le \alpha(d)$ and $r(N(d'))\le r(N(d))$, due to 
$N(d)\le N(d')$. Hence
\begin{equation}
\nn \begin{split}
 \cH^{\alpha(d)}_{r(N(d))}(\Omega''_d)&\le \sum_{d'\ge
d}\cH^{\alpha(d)}_{r(d)}( \cB_{N,d'})\le \sum_{d'\ge d} \cH^{\alpha(d')}_{r(d')}(
\cB_{N,d'}) \\ &\le \sum_{d'\ge d}\exp(-(\log \log N(d'))^B)\le
\sum_{d'\ge d}\ (d')^{-B}\lesssim d^{-B+1}\le d^{-1}
\end{split}
\end{equation}
Therefore $\cH^{\alpha(d)}_{r(N(d))}(\Omega'')\lesssim d^{-1}$ for any $d$.
Since $\alpha(d),r(N(d))\to 0$ with $d\to
\infty$, the Hausdorff dimension of $\Omega''$ is zero.
\item[(d)] Set $\Omega:=\Omega'\cup \Omega''$. Then the Hausdorff
dimension of $\Omega$ is zero. Assume that $\omega \in
\tor_{c,a}\setminus \Omega$. Then there exist $N'_1$, $d_1$ such
that $\omega\notin \hat {\Omega}_{N'_1}$, and $\omega \notin
\cB_{N(d),d}$ for any $d\le d_1$. Pick arbitrary $d\ge d_1$ such
that $N(d)\ge N'_1$. Set $N_1:=N(d)$. Then $\hat
{\Omega}_{N_1}\subset \hat {\Omega}_{N'_1}$. Hence, $\omega \in
\tor_{c,a} \backslash (\hat{\Omega}_{N_1}\cup \cB_{N_1,d})$. Due to
$(b)$ any given interval $(E',E'')$ contains a subinterval
$(E'_1,E''_1)$ such that $\Sigma_{\omega}\cap
(E'_1,E''_1)=\emptyset$ as claimed. 
\end{itemize}
\end{proof}

\appendix

\section{Polynomials, Resultants, and Algebraic Functions}

We recall some basic facts on polynomials. They can be found in
\cite{Lan}, chapter 5.

\begin{itemize}
\item A polynomial $f(x), x=(x_j,\dots, x_\nu)$, of $\nu$ variables $x_j \in \IR, j-1,2\dots n$ is called irreducible if there is no factorization
\[
f(x) = g(x)h(x)
\]
with $g(x), h(x)$ being non-constant polynomials.

\item Each polynomials $f(x)$ can be factorized
\[
f(x) = \prod_{j=1}^m f_j(x)
\]
with $f_j(x)$  irreducible. This factorization is unique up to scalars.

\item Let
\[
f(x_1,\dots, x_\nu) = \sum\limits_{j=0}^{k} a_j(x_2, \dots,
x_\nu)x_1^j, \,\,\, \alpha_k \not\equiv 0
\]
\[
 g(x_1,\dots, x_\nu) = \sum\limits_{j=0}^{m} b_j(x_2, \dots, x_\nu)x_1^j, \,\,\, b_m \not\equiv 0
 \]
be two arbitrary polynomials of $\nu$ variables. The following determinant
 \begin{equation}\label{eq:B.res}
\Res(f, g) = \left|\begin{array}{ll} {\overbrace{\begin{array}{lll}
a_k & 0 & \cdots\\
a_{k-1}  & a_k & \cdots\\
a_{k-2} & a_{k-1} & \cdots\\
\cdots & \cdots & \cdots\\
a_0 & a_1 & \cdots\\
0 & a_0 & \cdots\\
0 & 0 & \cdots \\
\cdots & \cdots & \cdots \end{array}}^m} &
{\overbrace{\begin{array}{llll}
b_m & 0 & \cdots & 0\\
b_{m-1} & b_m & \cdots & \cdots\\
b_{m-2} & b_{m-1} &\cdots & \cdots\\
\cdots & \cdots & \cdots & \cdots\\
\cdots & \cdots & \cdots & \cdots\\
\cdots & \cdots & \cdots& \cdots\\
\cdots&\cdots&\cdots&\cdots\\
\cdots&\cdots&\cdots&\cdots\\
\end{array}}^k}
\end{array}\right|
\end{equation}
is called the resultant of $f$ and $g$ with respect to the first
variable. For arbitrary variable $x_j$ the resultant is defined
similarly.

\item Let $f(z) = z^k + a_{k-1} z^{k-1} + \cdots + a_0$, $g(z) =
z^m + b_{m-1} z^{m-1} + \cdots + b_0$, $a_i, b_j \in
\IC$. Let $\zeta_i$, $1 \le i \le k$ and $\eta_j$, $1 \le j \le m$
be the zeros of $f(z)$ and $g(z)$, respectively.  The resultant of
$f$ and $g$ satisfies
\begin{equation*}
\Res(f, g) = \prod_{i,j} (\zeta_i - \eta_j)
\end{equation*}
The discriminant of the polynomial $f$ is defined as
\begin{equation}
  \label{eq:discriminant}\disc f = \prod_{i \ne j}\bigl(\zeta_i - \zeta_j\bigr)\ .
\end{equation}
One has also
$$
\disc f = (-1)^{n(n-1)/2} \Res(f, f')\ .
$$

\item Let $R(x_2, \dots, x_\nu)$ be as in~\eqref{eq:B.res}. $R(x_2, \dots, x_\nu)$
is a polynomial $\deg(R) \le \deg(f)\deg(g)$. Recall that for
\begin{equation}\label{eq:B.deg}
f(x) = \sum\limits_\alpha c_\alpha x^\alpha
\end{equation}
The degree $\deg(f)$ is defined as $\max\{|\alpha|: c_\alpha \neq 0
\}$ Here
\[
x=(x_1, \dots, x_\nu), \;\alpha = (\alpha_1, \dots, \alpha_\nu),\;
x^{\alpha} = x_1^{\alpha_1} x_2^{\alpha_2} \dots x_\nu^{\alpha_\nu}.
\]

\item The main property of the resultant is as follows
\[
R(f, g) \equiv 0
\]
If and only if
\[
f=hf_1, \,\,\, g=hg_1
\]

\item In particular, if $f=\sum\limits_{j=0}^k a_j(x_1, \dots, x_\nu) x_1^j,\; a_k(x_2, \dots, x_\nu) \not\equiv 0$
is irreducible, then
\[
\Res(f, g) \equiv 0
\]
for any $g$.

\item Let $f(x,y), g(x,y)$ be polynomials of two variables, $x,y\in\IR$.
Bezout's Theorem says that if $f$ or $g$ is irreducible then the system
\[
f(x,y)=0,\quad
g(x,y)=0
\]
has at most $m\cdot n$ solutions, $m=\deg f, n=\deg g$.

\item Let $f(x,y)$ be a polynomial, and let  $\varphi(x)$
be a function defined on some interval $[a,b]$ such that
\begin{equation}\label{eq:B.alg}
f(x, \varphi(x)) = 0, \text{\ \ for any \ \ } x \in [a,b]
\end{equation}
has $m(m-1)$ solution at most. If for some $x_0 \in [a,b]$
\[
\partial_yf(x_0, y(x_0))\neq 0
\]
then $\varphi(x)$ is real analytic in some neighborhood of $x_0$,
\[
\partial_x\varphi = -\partial_xf/\partial_yf
\]
On the other hand if $\varphi(x)$ is real analytic on $[a,b]$ and
obeys~\eqref{eq:B.deg}, then there exists an irreducible polynomial
$f(x,y)$ such that
\[
f_1(x,\varphi(x)) = 0 \quad\text{\ \ \ for any \ \ \ } x \in [a,b]
\]
for any $\mu \in \IR$ the number of solutions of the equation
\[
\partial_x\varphi = \mu
\]
does not exceed $2 m(m-1)$, where $m=\deg f$.
\end{itemize}


\begin{thebibliography}{GolMolPasxxx}

\bibitem[AviJit]{AJ} Avila, A., Jitomirskaya, S. {\em  Solving the ten
martini problem}, preprint 2005.

\bibitem[AvrSim]{AvrSim} Avron, J., Simon, B. {\em Almost periodic Schr\"odinger operators. II. The integrated density of states.}
  Duke Math.\ J.~50  (1983),  no.~1, 369--391.

\bibitem[BelSim]{BelSim} B\'ellissard, J., Simon, B. {\em
Cantor spectrum for the almost Mathieu equation.} J.\ Funct.\
Anal.~48 (1982), no.~3, 408--419.

\bibitem[Bha]{Bhat} Bhatia, R. {\em Perturbation bounds for matrix eigenvalues.} Pitman research notes in mathematics series~162, Longman, 1987.

\bibitem[Bou1]{BouKAM} Bourgain, J. {\em H\"older regularity of integrated density of states for the
almost Mathieu operator in a perturbative regime.} Lett.\ Math.\ Phys.\ 51 (2000), no.~2, 83--118.


\bibitem[Bou2]{Boupos} Bourgain, J. {\em On the spectrum of lattice Schr\"odinger operators with deterministic potential. Dedicated to the memory of Thomas H.\ Wolff.}  J.\ Anal.\ Math.~87  (2002), 37--75. 


\bibitem[Bou3]{Bou} Bourgain, J. {\em Green's function estimates for lattice Schr\"odinger
operators and applications.} Annals of Mathematics Studies, 158.
Princeton University Press, Princeton, NJ, 2005.

\bibitem[BouGol]{BouGol} Bourgain, J., Goldstein, M.
{\em On nonperturbative localization with quasi-periodic potential.}
Ann.\ of Math.~(2) 152 (2000), no.~3, 835--879.

\bibitem[BouJit]{BourJit} Bourgain, J., Jitomirskaya, S. {\em Continuity of the Lyapunov exponent for quasiperiodic operators with analytic potential.} Dedicated to David Ruelle and Yasha Sinai on the occasion of their 65th birthdays.
J.\ Statist.\ Phys.\  108  (2002),  no.~5-6, 1203--1218.


\bibitem[Cha]{Chan} Chan, J. {\em Method of variations of potential of quasi-periodic Schr\"odinger equation}, preprint 2005, to appear
in Geom.\ Funct.\ Analysis. 

\bibitem[ChoEllYui]{CEY}  Choi, M.D., Elliott, G.A., Yui, N. {\em Gauss polynomials and the rotation algebra.} Invent.\
 Math.~99(2), 225--246 (1990)


\bibitem[Fal]{Falc} Falconer, d.\ J. {\em The geometry of fractal sets.} Cambridge Tracts in Mathematics, 85. Cambridge University Press, Cambridge, 1986.
 
\bibitem[Fed]{fed} Federer, H. {\em Geometric measure theory.}
Die Grundlehren der mathematischen Wissenschaften, Band 153
Springer-Verlag New York Inc., New York 1969.


\bibitem[FroSpe1]{FS1} Fr\"{o}hlich, J., Spencer, T.
{\em Absence of diffusion in the Anderson tight binding model for
large disorder or low energy.} Comm.\ Math.\ Phys.\ {\em 88} (1983),
151--189.

\bibitem[FroSpe2]{FS2} Fr\"{o}hlich, J., Spencer, T. {\em A rigorous approach to Anderson localization.}
Phys.\ Rep.\ {\em 103} (1984), no.\ 1--4, 9--25.

\bibitem[FroSpeWit]{FSW} Fr\"{o}hlich, J., Spencer, T., Wittwer, P.
{\em Localization for a class of one dimensional quasi-periodic
Schr\"{o}dinger operators.}  Commun.\ Math.\ Phys.\ {\em 132}
(1990), 5--25.


\bibitem[FurKes]{Furdes} F\"urstenberg, H., Kesten, H.
{\em Products of random matrices.}
 Ann.\ Math.\ Statist 31 (1960), 457--469.



\bibitem[GolSch1]{Gol Sch1} Goldstein, M., Schlag, W.
{\em H\"older continuity of the integrated density of states
 for quasiperiodic Schr\"odinger equations and averages of shifts of subharmonic functions.}
Ann.\ of Math.~(2) 154 (2001), no.~1, 155--203.


\bibitem[GolSch2]{Gol Sch2} Goldstein, M., Schlag, W.
{\em Fine properties of the integrated density of states and a
quantitative separation property of the Dirichlet eigenvalues},
to appear in Geom.\ Funct.\ Analysis.

\bibitem[GolSch3]{Gol Sch3} Goldstein, M., Schlag, W. {\em On the formation of gaps in the spectrum of Schr\"odinger operators with quasi-periodic potentials.}  Spectral theory and mathematical physics: a Festschrift in honor of Barry Simon's 60th birthday,  591--611, Proc.\ Sympos.\ Pure Math., 76, Part 2, Amer.\ Math.\ Soc., Providence, RI, 2007.

\bibitem[GolSch4]{Gol Sch4} Goldstein, M., Schlag, W. {\em On Schr\"odinger operators with dynamically defined potentials.}
  Mosc.\ Math.\ J.~5  (2005),  no.~3, 577--612.

\bibitem[Her]{herman} Herman, M.
{\em Une m\'{e}thode pour minorer les exposants de Lyapounov et
quelques exemples montrant le charact\`{e}re local d'un theoreme
d'Arnold et de Moser sur le tore de dimension $2$.} Comment.\ Math.\
Helv.\ 58 (1983), no.~3, 453--502.


\bibitem[Joh]{Joh} Johnson, R. {\em Cantor spectrum for the quasi-periodic Schr\"odinger equation.} J.\ Diff.\ Eq.~91,
88--110 (1991)

\bibitem[Lan]{Lan} Lang, S. {\em Algebra}, Third Edition, Addison Wesley, 1993.

\bibitem[Las]{Last2} Last, Y. {\em Zero measure spectrum for the almost Mathieu operator.} Comm.\ Math.\ Phys.~164 (1994), 421--432.

\bibitem[Lev]{levin} Levin, B.\ Ya. {\em Lectures on entire functions.} Transl.\ of Math.\ Monographs, vol.~150. AMS, Providence, RI, 1996.

\bibitem[Pui]{Puig}  Puig, J. {\em Cantor spectrum for the almost Mathieu operator.} Comm.\ Math.\ Phys.\ 244 (2004), no.~2, 297--309.

\bibitem[ReeSim4]{RS4}  Reed, M., Simon, B. {\em Methods of modern mathematical
physics. IV.} Academic Press [Harcourt Brace Jovanovich,
Publishers], New York-London, 1979.


\bibitem[Sin]{Sin1} Sinai, Y.\ G. {\em Anderson localization for one-dimensional difference Schr\"{o}dinger operator with quasi-periodic potential.}  J.\ Stat.\ Phys.\ {\em 46} (1987), 861--909.


\bibitem[Tod]{To}  Toda, M. {\em Theory of nonlinear lattices.} Second edition.
Springer Series in Solid-State Sciences, 20. Springer-Verlag,
Berlin, 1989.


\end{thebibliography}
\end{document}